\documentclass[12pt,oneside,a4paper]{article}

\usepackage[frenchb]{babel}
\usepackage{amsmath}
\usepackage{amsfonts}
\usepackage{amsthm}
\usepackage{amscd}
\usepackage[latin1]{inputenc}
\usepackage[T1]{fontenc}
\usepackage[autolanguage]{numprint}

\usepackage{graphicx}

\usepackage{dea}

\title{{\large \sc DEA de Mod\'elisation Stochastique et
Statistique \\ Universit\'e Paris-Sud XI}\\ \vspace{2.5cm} \'Etude
d'un mod\`ele de dynamique des populations\\ \vspace{4cm}}
\author{Sylvain Arlot\\  \textit{Stage encadr\'e par
Jean-Christophe Yoccoz} \vspace{1cm}}
\date{14 Septembre 2004}

\begin{document}
\pagenumbering{roman}

\maketitle \thispagestyle{empty}
\newpage

\tableofcontents
\newpage

\pagenumbering{arabic}
\section{Introduction}
La dynamique des populations est au c\oe ur de l'interface entre
syst\`emes dynamiques et biologie. Ainsi, l'un des mod\`eles
biologiques les plus simples et les plus importants --- le
mod\`ele logistique --- correspond \`a la dynamique des
polyn\^omes quadratiques, dont l'\'etude math\'ematique est des
plus int\'eressantes. Ce seul cas montre combien il est ais\'e
d'obtenir un comportement complexe sans introduire beaucoup de
complexit\'e dans le mod\`ele.

Mais le mod\`ele logistique est bien trop peu r\'ealiste pour que
sa complexit\'e dynamique puisse \^etre directement
interpr\'et\'ee dans le cadre d'une population biologique
r\'eellement observ\'ee. D'un point de vue math\'ematique, les
syst\`emes de dimension 1 pr\'esentent un nombre limit\'e de
dynamiques possibles. Il est donc int\'eressant de consid\'erer
des mod\`eles de dimension sup\'erieure tels que l'application de
H\'enon (qui est une petite perturbation de la dynamique d'un
polyn\^ome quadratique, en dimension 2), et qui sont encore mal
compris du point de vue th\'eorique.

Le pas suivant dans cette d\'emarche consiste en l'\'etude de
syst\`emes dynamiques de dimension infinie, qui seront de
<<petites perturbations>> des mod\`eles pr\'ec\'edents,
\latin{i.e.} le mod\`ele logistique. En effet, si l'on veut
int\'egrer le ph\'enom\`ene biologique de maturation des jeunes
individus, il est n\'ecessaire de consid\'erer la fonction
d'\'evolution de l'effectif en temps continu, et non seulement sa
valeur \`a un instant donn\'e, ce qui donne un syst\`eme dynamique
de dimension infinie (ou de grande dimension, si l'on discr\'etise
ce syst\`eme). Un autre ph\'enom\`ene int\'eressant \`a
consid\'erer est l'influence des rythmes saisonniers sur un tel
syst\`eme, lorsqu'il se combine avec cet effet de retard induit
par le temps de maturation des jeunes. Le mod\`ele que nous
consid\'erons combine ces deux effets avec une forme de
densit\'e-d\'ependance un peu diff\'erente de celle du mod\`ele
logistique.

\medskip

Nous commencerons par d\'efinir le mod\`ele \'etudi\'e, tel qu'il
a \'et\'e \'enonc\'e dans \cite{Yoccoz:toymodel}, puis sous une
forme l\'eg\`erement modifi\'ee, en motivant celle-ci aussi bien
par des raisons biologiques que des raisons de simplicit\'e
pratique. Nous verrons ensuite ce que l'on peut montrer simplement
par une \'etude th\'eorique \latin{a priori}, point de d\'epart
d'une \'etude future plus approfondie (mais surtout bien plus
difficile). La derni\`ere et plus importante partie de notre
\'etude sera consacr\'ee \`a l'analyse des r\'esultats de
simulations num\'eriques, en vue de comprendre l'influence des
param\`etres sur la dynamique du syst\`eme et d'analyser plus
finement un des attracteurs \'etranges que nous avons pu observer.
Ce travail est bien s\^ur loin d'\^etre complet, et se veut
surtout \^etre une introduction (et une motivation) pour de futurs
travaux, aussi bien math\'ematiques que biologiques.

\medskip

Je tiens \`a remercier particuli\`erement Jean-Christophe Yoccoz
pour le temps qu'il m'a consacr\'e, ses nombreux conseils et la
clart\'e des explications th\'eoriques qu'il m'a donn\'ees. Je
remercie \'egalement Gilles Yoccoz pour ses conseils
bibliographiques pour la partie biologique de ce m\'emoire.
\section{Description du mod\`ele}

\subsection{Mod\`ele initial continu}

Le mod\`ele suivant, d\'efini dans \cite{Yoccoz:toymodel},
d\'ecrit l'\'evolution temporelle d'une population de campagnols.

\begin{equation} \label{eq:modele}
N(t) = \int_{A_0}^{A_1} {S(a) m_{\rho}(t-a) N(t-a) m(N(t-a)) da}
\end{equation}

\begin{itemize}
\item $t$ est le temps (en ann\'ees),
\item $N$ la population active (\latin{i.e.} d'\^age sup\'erieur \`a $A_0$),
\item $A_0$ l'\^age de maturation,
\item $A_1$ l'\^age maximal,
\item $S$ le taux de survie,
\item $m_{\rho}$ le param\`etre de saison (d\'ecrit une probabilit\'e de
reproduction en fonction de la saison),
\item $m(N)$ le taux de f\'econdit\'e individuel annuel pour une population active de taille $N$.
\end{itemize}

On a choisi des formes simples pour les fonctions $S$, $m_{\rho}$
et $m$ :
\begin{equation}\label{eq:survie}
S(a) = 1-\frac{a}{A_1}
\end{equation}
\begin{equation}\label{eq:saison_discontinue}
 m_{\rho} (t) =
\begin{cases}
0 \text{ si } 0 \leq t \leq \rho \text{ mod. 1} \\
1 \text{ si } \rho \leq t \leq 1
\end{cases}
\end{equation}
\begin{equation}\label{eq:fecondite_C0}
m(N) = \begin{cases}
m_0  \text{ si } N \leq 1 \\
m_0 N^{-\gamma} \text{ sinon}
\end{cases}
\end{equation}

Cependant, pour \'eviter des art\'efacts d\^us \`a la
non-r\'egularit\'e du syst\`eme, il nous a sembl\'e pr\'ef\'erable
de remplacer $m$ et $m_{\rho}$ par des fonctions un peu plus
r\'eguli\`eres.

\subsection{Lissage de la f\'econdit\'e}

La fonction $m$ d\'efinie par \eqref{eq:fecondite_C0} est continue
mais pas $C^1$. Il y a une forte rupture de pente \`a la valeur
critique $N=1$. On pourrait ais\'ement <<recoller>> les deux
parties de la courbe pour obtenir une fonction $C^{\infty}$, mais
cela ne serait pas tr\`es pratique pour les simulations
num\'eriques. Nous utilisons ainsi une parabole interm\'ediaire qui
rend $m$ $C^1$.

\begin{equation}\label{eq:fecondite}
m(N) = \begin{cases}
m_0  \text{ si } N \leq N_{1,\gamma} \\
m_0 \times \left( A_{\gamma} + B_{\gamma} N + C_{\gamma} N^2
\right) \text{
si } N_{1,\gamma} < N \leq N_{2,\gamma} \\
 m_0 N^{-\gamma} \text{ si } N_{2,\gamma} < N
\end{cases}
\end{equation}

Pour garder un mod\`ele r\'ealiste et suffisamment proche du
mod\`ele initial, il faut conserver la d\'ecroissance de la
f\'econdit\'e $N \mapsto m(N)$, et s'assurer que la parabole
rejoint les valeurs extr\^emes $m_0$ et $m_0 N^{-\gamma}$
suffisamment pr\`es de $N=1$. Pour d\'efinir compl\`etement les
param\`etres $N_{1,\gamma}$, $A_{\gamma}$, $B_{\gamma}$,
$C_{\gamma}$ et $N_{2,\gamma}$, on impose \'egalement
$m(N_{2,\gamma}) = m_0 /2$. Cette derni\`ere contrainte permet de
limiter la zone interm\'ediaire, ce qui simplifiera notamment les
calculs explicites de la section \ref{sec:theorie}. On a ainsi les
conditions suivantes :
\begin{equation}
\left\{
\begin{aligned}
A_{\gamma} + B_{\gamma} N_{1,\gamma} + C_{\gamma} N_{1,\gamma}^2 &= 1 \\
B_{\gamma} + 2 C_{\gamma} N_{1,\gamma} &= 0 \\
N_{2,\gamma}^{-\gamma} &= \frac{1}{2} \\
A_{\gamma} + B_{\gamma} N_{2,\gamma} + C_{\gamma} N_{2,\gamma}^2 &= \frac{1}{2}\\
B_{\gamma} + 2 C_{\gamma} N_{2,\gamma} &= - \gamma
N_{2,\gamma}^{-\gamma-1} = - \frac{\gamma}{2 N_{2,\gamma}}
\end{aligned} \right. \end{equation}

\begin{figure}
\begin{center}
\includegraphics[height=7cm]{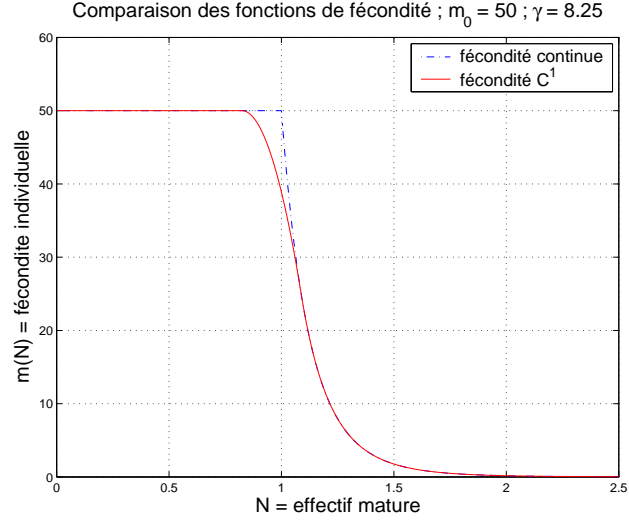}
\caption{\label{fig:fecondite}Comparaison des fonctions de
f\'econdit\'e.}
\end{center}
\end{figure}

De plus, on veut $C_{\gamma} \leq 0$ pour garantir la
d\'ecroissance de $m$. On choisit donc :
\begin{equation}
\left\{
\begin{aligned}
N_{2,\gamma} &= 2^{1/\gamma} \\
C_{\gamma} &= \frac{- \gamma^2}{8 \times 4^{1/\gamma}} \\
A_{\gamma} &= \frac{1}{2} \left( 1+\gamma - \frac{\gamma^2}{4} \right) \\
B_{\gamma} &= 2^{-1/\gamma} \times \left( \frac{\gamma^2}{4} - \frac{\gamma}{2} \right) \\
N_{1,\gamma} &= 2^{1/\gamma} \times \left( 1 - \frac{2}{\gamma}
\right)
\end{aligned} \right. \end{equation}

La figure~\ref{fig:fecondite} repr\'esente les deux fonctions $m$
--- continue et $C^1$ --- pour $\gamma = 8\virg 25$.

\subsection{Lissage du facteur saisonnier}

La fonction $m_{\rho}$ d\'efinie par \eqref{eq:saison_discontinue}
n'est pas continue, il est l\'egitime de vouloir consid\'erer un
facteur saisonnier un peu plus r\'egulier. On a choisi,
arbitrairement, de le rendre $C^1$ en effectuant le passage de 0 \`a
1 \`a l'aide d'un cosinus. Pour cela, on ajoute un param\`etre
$\epsilon$ qui est la dur\'ee du printemps et celle de l'automne. La
dur\'ee de l'\'et\'e est d\'esormais $1 - \rho - \epsilon$ et non plus $1
-  \rho$. Prendre $\epsilon = 0$ ram\`ene bien s\^ur au cas
pr\'ec\'edent.

\begin{figure}
\begin{center}
\includegraphics[height=7cm]{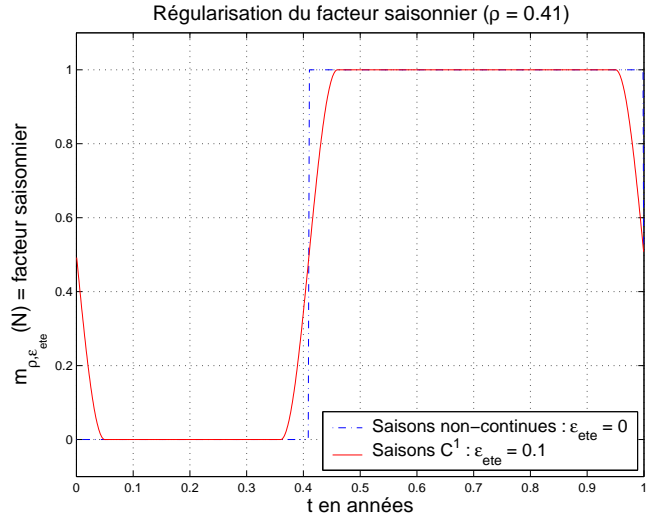}
\caption{\label{fig:saison}Comparaison des fonctions
$m_{\rho,\epsilon}$.}
\end{center}
\end{figure}

\begin{equation}\label{eq:saison}
 m_{\rho,\epsilon} (t\text{ mod. }1) =
\begin{cases}
\frac{1}{2}(1+ \cos(\pi\times(\frac{t}{\epsilon}+\frac{1}{2}))) &\text{ si } 0 \leq t < \epsilon/2 \\
0 &\text{ si } \epsilon/2 \leq t < \rho - \epsilon/2 \\
\frac{1}{2} \left( 1+ \cos\left (\pi\times(\frac{t-\rho}{\epsilon}-\frac{1}{2}) \right) \right) &\text{ si } \rho - \epsilon/2 \leq t < \rho + \epsilon/2 \\
1 &\text{ si } \rho + \epsilon/2 \leq t < 1 - \epsilon/2 \\
\frac{1}{2}(1+ \cos(\pi\times(\frac{t-1}{\epsilon}+\frac{1}{2})))
&\text{ si } 1-\epsilon/2 \leq t < 1
\end{cases}
\end{equation}

Cette d\'efinition n'\'etant valable que lorsque $\epsilon \leq \rho
\leq 1-\epsilon$, on posera $\epsilon = \min (\rho,1-\rho)$
lorsque ce n'est pas le cas \latin{a priori}.

La figure~\ref{fig:saison} repr\'esente $m_{\rho,\epsilon}$ pour
$\rho = 0\virg 41$ et deux valeurs de $\epsilon$. Remarquons enfin
que l'on aurait \'egalement pu rendre $m_{\rho,\epsilon}$
$C^{\infty}$ d\`es que $\epsilon >0$ en utilisant autre chose
qu'un cosinus. Le choix que nous avons fait tient compte de la
simplicit\'e des calculs num\'eriques futurs.

Pour les d\'etails concernant les simulations num\'eriques, voir
le paragraphe \ref{par:modele_discretisation} en annexe.

\subsection{Aspects biologiques}
\subsubsection{\'Etudes ant\'erieures}
\begin{figure}
\begin{center}
\includegraphics{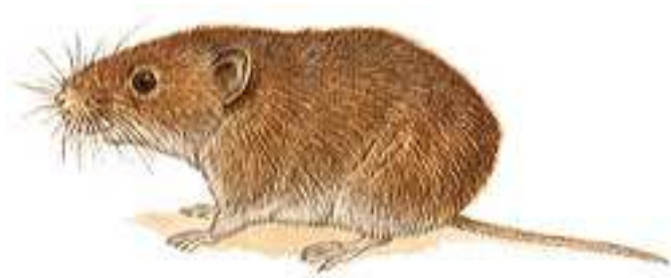}
\\ \copyright 2003
Missouri Botanical Garden\\
\includegraphics{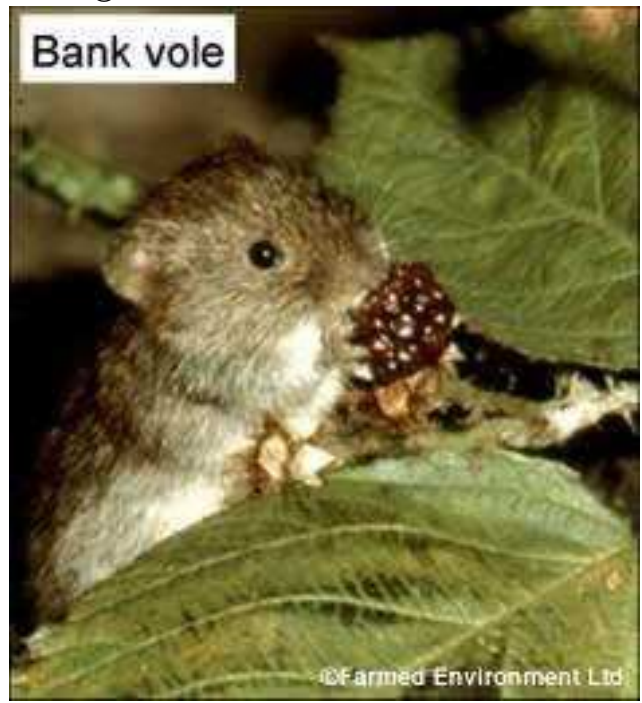}
\caption{\label{fig:dessin_photo} Campagnol rouss\^atre (bank
vole, \emph{Clethrionomys glareolus}). Ce petit rongeur, de la
famille des Microtin\'es, vit dans les for\^ets temp\'er\'ees et
se nourrit principalement de graines.}
\end{center}
\end{figure}

\paragraph{Esp\`eces concern\'ees}
Le mod\`ele que nous venons de d\'ecrire a \'et\'e \'elabor\'e en
relation avec l'\'etude de la dynamique de certaines populations
de petits rongeurs. Ceux-ci se caract\'erisent en effet par un
fort investissement dans la reproduction (port\'ees importantes et
maturit\'e sexuelle tr\`es rapide) et de grandes variations
annuelles de la taille de la population. Plusieurs esp\`eces de
campagnols ont ainsi \'et\'e \'etudi\'ees, notamment le campagnol
rouss\^atre \emph{Clethrionomys glaerolus}
(figure~\ref{fig:dessin_photo} ;
\cite{YoccozStenseth2000,CrespinYoccozetal2002,YoccozStensethetal2001})
et \emph{Microtus townsendii} \cite{LambinYoccoz2001}.

\begin{figure}
\begin{center}
\includegraphics[width=12cm]{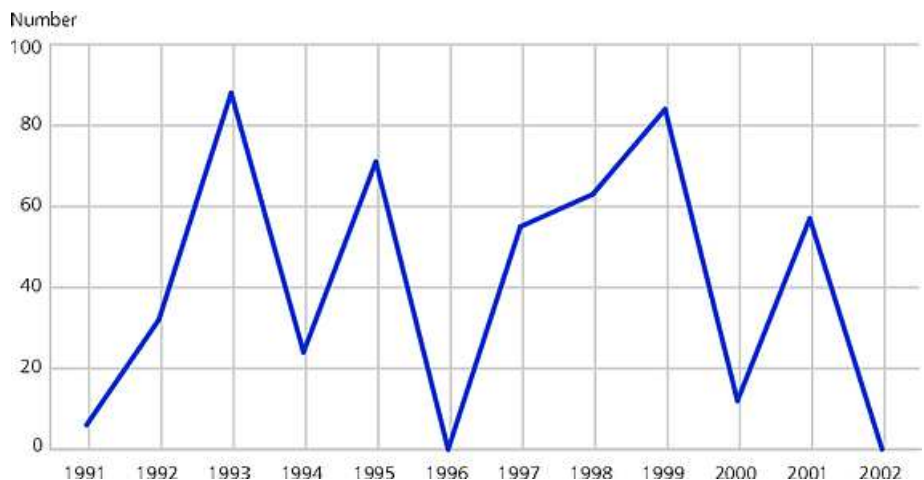}
\caption{\label{fig:popu_microtus_epiroticus} \emph{Microtus
epiroticus} \`a Svalbard : nombre d'individus captur\'es entre
1991 et 2002. (Source: NINA and University of Troms\o. 2003
\copyright\ Statistics Norway).}
\end{center}
\end{figure}

C'est plus particuli\`erement une population de \emph{Microtus
epiroticus}\footnote{sibling vole en anglais.}, introduite
accidentellement il y a une cinquantaine d'ann\'ees dans
l'archipel arctique de Svalbard (en Finlande, dans le Spitzberg),
qui est vis\'ee par ce mod\`ele. Ces campagnols poss\`edent en
effet une f\'econdit\'e est extr\^emement \'elev\'ee pour des
mammif\`eres. De plus, cette esp\`ece \'etant menac\'ee
d'extinction, l'\'etude de sa d\'emographie permettrait
\'egalement de mieux la prot\'eger. La
figure~\ref{fig:popu_microtus_epiroticus} montre ainsi de grandes
fluctuations de population, et des effectifs minimaux tr\`es
faibles, de l'ordre de quelques individus.

\paragraph{M\'ecanismes envisag\'es}
Plusieurs causes possibles \`a ces ph\'enom\`enes ont \'et\'e
\'etudi\'ees. Il a \'et\'e montr\'e exp\'erimentalement qu'une
augmentation de la \emph{quantit\'e de nourriture disponible}
augmente la densit\'e mais n'a pas d'influence sur la
densit\'e-d\'ependance (\emph{C. glaerolus}, Finlande
\cite{YoccozStensethetal2001}). Plus pr\'ecis\'ement, la
nourriture disponible est li\'ee aux variations inter-annuelles,
en jouant sur la survie de l'ann\'ee suivante ; on observe des
variations saisonni\`eres de taille comparable aux variations
inter-annuelles et qui se caract\'erisent par une survie plus
faible au printemps (\emph{C. glaerolus}, Belgique
\cite{CrespinYoccozetal2002}).

Le \emph{facteur climatique} semble lui aussi jouer un grand
r\^ole : une comparaison entre \emph{M. epiroticus} \`a Svalbard
et \emph{Chionomys nivalis}\footnote{snow vole en anglais.} dans
les Alpes fran\c{c}aises indique une relation entre un
environnement stable et un faible turnover\footnote{\latin{i.e.}
forte survie et faible f\'econdit\'e} \cite{YoccozIms1999}. En
effet, la population vivant dans l'Arctique, o\`u les conditions
hivernales sont tr\`es variables, a une reproduction extr\^emement
rapide, tandis que la population alpine, dont l'environnement est
stable, se reproduit peu et a une forte survie. On a ainsi
essay\'e d'inclure une stochasticit\'e environnementale et
d\'emographique dans les mod\`eles afin d'expliquer les
fluctuations de population observ\'ees (\emph{C. glaerolus}, Alpes
\cite{YoccozStenseth2000}).

Un ph\'enom\`ene important pourrait \^etre reli\'e \`a ces
facteurs environnementaux : la plasticit\'e de l'\^age \`a
maturit\'e. Il a ainsi \'et\'e montr\'e que les femelles n\'ees en
d\'ebut de saison de reproduction ont avantage \`a se reproduire
rapidement, malgr\'e le co\^ut \'elev\'e d'une reproduction
pr\'ecoce (\emph{M. townsendii}, Vancouver
\cite{LambinYoccoz2001}) : les femelles naissant plus t\^ot
peuvent se reproduire avant le fin de la saison de reproduction,
augmentant ainsi la contribution de leur m\`ere \`a la population
totale.

\subsubsection{Param\`etres du mod\`ele}
Les param\`etres de croissance et de reproduction des deux
populations\footnote{\emph{M. epiroticus} \`a Svalbard et \emph{M.
arvalis} en Finlande.} ont \'et\'e \'evalu\'es dans
\cite{YoccozImsSteen1993}. \`A partir de ces conclusions, nous
pouvons justifier le mod\`ele et le choix de param\`etre
<<typiques>>.

\paragraph{Densit\'e-d\'ependance}
La forme particuli\`ere de densit\'e-d\'ependance se justifie car
le facteur limitant est ici le nombre tr\`es restreint de sites de
reproduction. Ainsi, seule la population (femelle\footnote{Comme
souvent en dynamique des populations, seule la population femelle
est consid\'er\'ee dans la mesure o\`u elle est le facteur
limitant de la reproduction. Il n'est int\'eressant de
consid\'erer les m\^ales que si la femelle a peu de chances de
rencontrer un m\^ale (par exemple si la densit\'e de population
est faible), ou bien si le sex-ratio est loin de $1:1$.}) mature
doit entrer en ligne de compte, et en cas de surpopulation, les
quelques sites disponibles perdent beaucoup en qualit\'e. Ceci est
confirm\'e par les conclusions de \cite{YoccozStensethetal2001},
qui indiquent dans le cas de \emph{C. glaerolus} une
densit\'e-d\'ependance plus forte en \'et\'e, \`a cause de la maturation des
femelles. Nous choisissons de ne pas faire d\'ependre $\gamma$ de
la saison puisque la variation observ\'ee est li\'ee \`a la prise
en compte de la densit\'e totale, et non de la densit\'e de
femelles matures comme nous le faisons ici.

\paragraph{Saison de reproduction}
La p\'eriode de reproduction correspond \`a la saison de
croissance des plantes, c'est-\`a-dire du printemps \`a l'automne.
Nous consid\'erons avec ce mod\`ele un climat parfaitement
pr\'evisible, identique d'une ann\'ee sur l'autre. Il s'agit de
savoir si l'on peut observer un comportement chaotique dans des
conditions parfaitement stables.

Sa dur\'ee varie donc selon les lieux. \`A Svalbard comme dans les
Alpes, celle-ci dure de 3 \`a 4 mois (Juin \`a Septembre)
\cite{YoccozIms1999}, \latin{i.e.} $\rho \approx 0\virg 7$. Dans
la plupart des environnements temp\'er\'es, celle-ci est beaucoup
plus longue. Ainsi, en Belgique, \emph{C. glaerolus} se reproduit
en g\'en\'eral entre la deuxi\`eme semaine d'Avril et la fin du
mois d'Octobre \cite{CrespinYoccozetal2002}, mais peut varier de
telle sorte que l'on a $0 \virg 35 < \rho < 0 \virg 45$.

\paragraph{F\'econdit\'e}
Un \'el\'ement important du mod\`ele est la forte f\'econdit\'e
des campagnols. Pour ajuster les param\`etres du mod\`ele, nous
avons besoin d'\'evaluer le nombre de jeunes femelles par femelle
mature et par an, en l'absence de densit\'e-d\'ependance.

Pour \emph{M. epiroticus} \cite{YoccozImsSteen1993}, la p\'eriode
de gestation est de 20 jours, et la taille des port\'ees varie de
$4\virg 5$ \`a 7, celle-ci augmentant pour une m\^eme femelle au
fur et \`a mesure de ses reproductions. Le sex ratio est tr\`es
proche de $1:1$. Dans la mesure o\`u une femelle peut se
reproduire pendant la p\'eriode d'allaitement, on en d\'eduit une
valeur maximale $m_0 \approx 18 \times 6 \times 0 \virg 5 = 54$.

Dans le cas de \emph{M. townsendii} \cite{LambinYoccoz2001}, on
observe de 5 \`a 6 port\'ees par saison de reproduction (celle-ci
durant de Mars \`a Novembre, soit environ 8 mois : $\rho \approx 0
\virg 3$), chacune comportant 5 \`a 8 individus. On a donc dans ce
cas $m_0 \time 2/3 \approx 5 \times 6\virg 5 \times 0 \virg 5$
soit $m_0 \approx 24$. Il n'est pas \'etonnant de constater que
cette valeur est bien inf\'erieure \`a celle de \emph{M.
epiroticus}, qui repr\'esente un cas extr\^eme parmi les
mammif\`eres.

On peut donc prendre comme valeur $m_0 = 50$, mais une
f\'econdit\'e l\'eg\`erement inf\'erieure serait sans doute plus
r\'ealiste.

\paragraph{\^Age de premi\`ere reproduction}
L'\^age de maturit\'e des femelles est sup\'erieur \`a 17 jours,
auquel il faut rajouter la dur\'ee de gestation, soit 20 jours
suppl\'ementaires \cite{YoccozImsSteen1993}. L'\^age de premi\`ere
reproduction de \emph{M. epiroticus} est donc au minimum 37 jours,
\latin{i.e.} $A_0 \geq 0 \virg 10$. En moyenne, on observe
plut\^ot une premi\`ere port\'ee \`a un \^age d'environ 50 jours,
\latin{i.e.} $A_0 \approx 0\virg 14$.

Une telle pr\'ecocit\'e ne se retrouve pas chez les autres
Microtin\'es, \`a l'exception de \emph{M. arvalis} en Finlande. On
utilisera donc $0 \virg 1$ comme valeur minimale, tandis que $A_0
= 0 \virg 20$ (correspondant \`a 75 jours) est plus classique.

\paragraph{Survie}
Pour \emph{M. epiroticus} \`a Svalbard, la survie hivernale est de
l'ordre de $0 \virg 1$ (et tr\`es variable selon les ann\'ees), et
la survie estivale $0 \virg 85$ \cite{YoccozIms1999}. Le taux de
mortalit\'e par ann\'ee est donc de $0\virg 046$ en hiver et $0
\virg 52$ en \'et\'e. L'\^age maximal $A_1$ est toujours
inf\'erieur \`a deux ans.

Chez \emph{C. glaerolus}, en Belgique, les variations
saisonni\`eres de la survie ont \'et\'e \'etudi\'ees plus en
d\'etails \cite{CrespinYoccozetal2002}. Le taux de survie est plus
\'elev\'e en hiver ($0 \virg 95$ par semaine, soit $0 \virg 07$
par an) qu'en \'et\'e ($0 \virg 90$ par semaine, soit $0 \virg
004$ par an). Elle est \'egalement l\'eg\`erement diff\'erente
chez les femelles matures et immatures.

Le taux de mortalit\'e est suppos\'e dans le mod\`ele constant
avec l'\^age, et ind\'ependant des saisons, ce qui est loin
d'\^etre le cas en g\'en\'eral. L'\^age maximal de 2 ans est
\'egalement une l\'eg\`ere sur-estimation de ce qu'il est en
r\'ealit\'e. Pour plus de r\'ealisme, la fonction de survie est
sans doute l'un des premiers \'el\'ements du mod\`ele \`a
modifier.

\subsubsection{Probl\`emes pos\'es}

L'\'etude de ce mod\`ele n'a pas pour but de faire des
pr\'evisions pr\'ecises concernant l'avenir de la population de
\emph{Microtus epiroticus} \`a Svalbard. Nous nous efforcerons de
consid\'erer des param\`etres r\'ealistes pour de petits rongeurs,
pas n\'ecessairement \emph{M. epiroticus}. Il s'agit surtout
d'effectuer un travail th\'eorique dans un cadre assez simple,
afin de savoir si l'on peut observer une grande variabilit\'e
d'effectifs (voire une dynamique chaotique) dans un mod\`ele
compl\`etement d\'eterministe, dans un environnement r\'egulier.
Et si oui, quelles sont les facteurs biologiques d\'eterminants
(f\'econdit\'e, \^age de maturit\'e, dur\'ee de l'hiver) ? Un
autre objectif est de bien comprendre ce mod\`ele tr\`es simple
avant de le complexifier en introduisant d'autres m\'ecanismes
pouvant jouer un r\^ole dans la dynamique de cette population,
parmi ceux que nous avons \'evoqu\'es pr\'ec\'edemment.
\section{\'Etude th\'eorique} \label{sec:theorie}

Les r\'esultats de cette partie proviennent principalement de
\cite{Yoccoz:toymodel}, o\`u les fonctions $m_{\rho}$ et $m$
consid\'er\'ees \'etaient donn\'ees par \eqref{eq:saison_discontinue} et
\eqref{eq:fecondite_C0}. Nous avons consid\'er\'e ici un cadre un
peu plus g\'en\'eral, valable pour les fonctions $m_{\rho}$ et $m$
<<liss\'ees>>.

\subsection{D\'efinition du syst\`eme dynamique}

Pour $t_0 \in \R/\Z$, notons $Y_{t_0}$ l'ensemble des fonctions
continues $N$ sur $[-A_1 ; 0]$ \`a valeurs positives et v\'erifiant la
condition
\[ N(0) = \int_{A_0}^{A_1} S(a) N(-a)m(N(-a))m_{\rho}(t_0 - a) da.
\] Le syst\`eme dynamique d\'efini par \eqref{eq:modele} est donn\'e par
le semi-groupe $(T_s)_{s\geq 0}$ :
\begin{equation}\label{eq:semi_groupe} T^s (t,N) = ( t+s (\modulo
1) , N^s_t ) \end{equation}
\begin{equation}
N^s_t(-a) =
\begin{cases}
N(s-a) \text{ si } 0 \leq s \leq a \leq A_1& \\
\int_{A_0}^{A_1} S(b) N(s-a-b) m(N(s-a-b)) m_{\rho}(t+s-a-b) db
\text{ sinon}&
\end{cases}\end{equation} Ceci est bien d\'efini pour $0 \leq s
\leq A_0$, on l'\'etend \`a $s \geq 0$ par la propri\'et\'e de
semi-groupe. L'espace des phases est alors
\[ Y^{\sharp} = \{ (t,N) \telque \, t \in \R / \Z ,\, N \in Y_t
\}.
\] Dans la suite, on \'ecrira parfois $N^s$ au lieu de $N^s_t$
lorsque cela ne cr\'ee pas de confusion.

On s'int\'eresse \`a l'application $T^1 : Y_0 \rightarrow Y_0$
donnant l'\'evolution de la population d'une ann\'ee sur
l'autre\footnote{Le choix de $t_0=0$ est arbitraire, on pourrait
aussi bien consid\'erer $T^1 : Y_{t_0} \rightarrow Y_{t_0}$, qui
donnerait le m\^eme type de dynamique.}. \`{A} $N \in Y_{t_0}$
donn\'ee, on associe ainsi une unique fonction continue
$\overline{N}$ d\'efinie sur $[-A_1 ; +\infty [$ prolongeant $N$
et compatible avec $T$ (c'est-\`a-dire la solution de l'\'equation
\eqref{eq:modele}).

On munit $Y^{\sharp}$ de sa topologie naturelle, induite par la
topologie produit sur $\R / \Z \times \mathcal{C}([-A_1 ; 0])$,
l'espace $\mathcal{C}([-A_1 ; 0])$ des fonctions continues de
$[-A_1 ; 0]$ dans $\R$ \'etant muni de la topologie de la
convergence uniforme. Il d\'ecoule alors de la continuit\'e
uniforme de $N \rightarrow N \times m(N)$ la propri\'et\'e
suivante :
\begin{Pro}Pour tout $s\geq 0$, $T^s$ est un op\'erateur continu
$Y^{\sharp} \rightarrow Y^{\sharp}$.\end{Pro}

De plus, $N \rightarrow N \times m(N)$ \'etant
$K_f$-Lipschitzienne (la constante $K_f$ peut \^etre calcul\'ee
explicitement en fonction des param\`etres du mod\`ele),
l'application $T^s : Y_t \rightarrow Y_{t+s}$ est
$K$-Lipschitzienne, avec $K = \max(1,(A_1-A_0)\times K_f)$. Ceci
d\'ecoule directement de la d\'efinition de $N^s_t$. La constante
$K$ d\'epend uniquement des param\`etres du mod\`ele, et pas de $t
\in \R / \Z$.

\subsection{Existence d'un attracteur}

Nous allons montrer que pour des valeurs raisonnables des
param\`etres, un tel syst\`eme dynamique poss\`ede un attracteur,
ce qui n\'ecessite plusieurs lemmes techniques. Nous n'utiliserons
pas les formes explicites des fonctions $m_{\rho}$ et $m$ (pour
rester g\'en\'eraux, dans la mesure o\`u celles-ci pourraient
\^etre modifi\'ees ult\'erieurement), mais uniquement les
hypoth\`eses suivantes :
\begin{align}
m_0 \geq m(N) &\geq  \frac{m_0}{2}  &\text{   si } N \leq 1 \label{hyp_m0_1}\\
m_0 N^{- \gamma} \geq m(N) &\geq \left( \frac{1}{2} \wedge
N^{-\gamma} \right) m_0 &\text{   si } N \geq 1 \label{hyp_m0_2}
\\ 1 \geq m_{\rho}(t) &\geq 0 \, &\forall t \label{hyp_mrho_1}
\end{align}
\begin{equation}
m_{\rho}(t)=1 \text { sur un intervalle de longueur }
1-\rho-\epsilon \label{hyp_mrho_2}.
\end{equation}

On voit ais\'ement que les fonctions d\'efinies par
\eqref{eq:saison} et \eqref{eq:fecondite} (resp.
\eqref{eq:saison_discontinue} et \eqref{eq:fecondite_C0})
v\'erifient ces hypoth\`eses. Le symbole $\wedge$ est employ\'e
ici et dans la suite \`a la place de $\min$, de m\^eme que $\vee$
signifie $\max$.

Posons \begin{align} c_0 :&= \int_{A_0+\rho+\epsilon}^{A_0+1} S(a)
da \\ &= (1-\rho-\epsilon) \left( 1 - \frac{ 1 + \rho + \epsilon +
2A_0}{2A_1} \right).\end{align} Nous nous pla\c{c}ons d\'esormais
dans le cas o\`u les param\`etre v\'erifient les conditions
suivantes :
\begin{gather}
\gamma \geq 1 \\
A_1 \geq (2A_0) \vee (A_0+1)\\
c_0 \times m_0 >2 \label{cond_c0m0} \\
\rho + \epsilon < 1
\end{gather}

\begin{rem}
Ces conditions sont tr\`es raisonnables, et toujours v\'erifi\'ees
au cours des simulations que nous avons faites. En effet, si on
impose $\rho+\epsilon \leq \frac{6}{10}$, $A_1 \geq 2$, $A_0 \leq
1/2$, $\gamma\geq1$, alors $c_0 \geq \frac{14}{100}$ et donc $m_0
\geq 15$ suffit pour satisfaire \eqref{cond_c0m0}. Il n'y a donc
pas \`a s'inqui\'eter du manque de finesse de cette majoration.
\end{rem}

\begin{Le}\label{lemme:majoration}
Soit $t_0 \in \R / \Z$, $N \in Y_{t_0}$. On a alors, pour tout $0
\leq s \leq A_0$ : \[ N(s) \leq N_{\max} := m_0 \frac{A_1}{2}
\left( 1-\frac{A_0}{A_1} \right)^2\]
\end{Le}
\begin{proof}
On a toujours $N m(N) \leq m_0$, d'apr\`es \eqref{hyp_m0_1},
\eqref{hyp_m0_2}, et car $\gamma \geq 1$. Comme de plus $m_{\rho}
\leq 1$ \eqref{hyp_mrho_1}, on a \[ N(s) \leq m_0 \int_{A_0}^{A_1}
S(a)da = N_{\max}. \]
\end{proof}

\begin{Le} \label{lemme:minoration}
Soit $N \in Y_{t_0}$ telle que $N \leq N_{\max}$. On a $i(N)=
\min_{[-A_1 ; 0]} N >0$.
\begin{enumerate}\item Si $i(N) \leq N_{\max}^{1-\gamma}$, alors $N(s) \geq \frac{c_0
m_0}{2} i(N)$ pour $0 \leq s \leq A_0$. \item Si $i(N) \geq
N_{\max}^{1-\gamma}$, alors $N(s) \geq \frac{c_0 m_0}{2}
N_{\max}^{1-\gamma}$ pour $0 \leq s \leq A_0$. \end{enumerate}
\end{Le}
\begin{proof}
Commen\c{c}ons par montrer que $N m(N) \geq \frac{m_0}{2} \times
\left( i(N) \wedge N_{\max}^{1-\gamma} \right)$ sur $[-A_1 ; 0]$ :

si $N \leq 1$, \begin{equation}\begin{split} N m(N) &\geq N \times
\frac{m_0}{2} \\ &\geq i(N) \times
\frac{m_0}{2},\end{split}\end{equation}

et si $N \geq 1$, \begin{equation}\begin{split} N m(N) &\geq
\frac{m_0}{2} \times N^{1-\gamma} \\ &\geq \frac{m_0}{2} \times
N_{\max}^{1-\gamma}.\end{split}\end{equation}

Par cons\'equent, \[ N(s) \geq \frac{m_0}{2} \times (i(N) \wedge
N_{\max}^{1-\gamma}) \times \int_{A_0}^{A_1} S(a)
m_{\rho}(t_0+s-a) da. \] En fonction de la valeur de $t_0+s$, on
peut trouver un sous-intervalle de $[A_0 ; A_0+1] \subset [A_0 ;
A_1]$, de longueur $1-\rho-\epsilon$ sur lequel
$m_{\rho}(t_0+s-\cdot)$ vaut 1. La derni\`ere int\'egrale est donc
minor\'ee par la m\^eme int\'egrale restreinte \`a ce
sous-intervalle, qui est plus grande que $c_0$ car $S$ est
d\'ecroissante. Ceci ach\`eve la preuve du lemme.
\end{proof}
\begin{Cor}\label{corollaire:minoration}
Si $\gamma \geq 1$, $\frac{c_0 m_0}{2} >1$, $N \in Y_{t_0}$, alors
pour $s$ assez grand (d\'ependant de $N$), on a : \[ \frac{c_0
m_0}{2} N_{\max}^{1-\gamma} \leq N^s(a) \leq N_{max} ,\, \forall a
\in [-A_1 ; 0]. \]
\end{Cor}

Remarquons que l'on peut remplacer la condition \eqref{cond_c0m0}
par $c_0 m_0 >1$ dans le cas o\`u $m$ est d\'efinie par
\eqref{eq:fecondite_C0} (cf. \cite{Yoccoz:toymodel}). La constante
2 a \'et\'e choisie arbitrairement dans l'op\'eration de lissage
de $m$, celle-ci pourrait \^etre prise plus proche de 1 sans
difficult\'e suppl\'ementaire, mais toujours strictement
sup\'erieure \`a 1.

\begin{Le}\label{lemme:lipschitz}
Soit $N \in Y_{t_0}$. Posons $L=m_0 \left( 3-\frac{A_0}{A_1}
\right)$. Alors, si $0 \leq s_0 \leq s_1 \leq A_0$, on a \[
\absj{\overline{N}(s_1) - \overline{N}(s_0)} \leq L \absj{s_1 -
s_0}. \]
\end{Le}
\begin{proof}
Remarquons tout d'abord que $m_{\rho} \leq 1$, $0 \leq S \leq 1$
et \[ \absj{S(s_1-u) - S(s_0-u)} \leq A_1^{-1} \absj{s_1-s_0} .\]
Pla\c{c}ons-nous dans le cas o\`u $s_1-A_1 \leq s_0-A_0$ (c'est
vrai car on a suppos\'e $A_1 \geq 2 \times A_0$), et \'ecrivons la
d\'efinition de $\overline{N}(s_i)$ en fonction de $N$.

\begin{align}
\begin{split}
\absj{\overline{N}(s_1) - \overline{N}(s_0)} &= \Bigl| - \int_{s_0
-A_1}^{s_1 -A_1} S(s_0-u) N(u)m(N(u))m_{\rho}(t_0+u)du \\ &+
\int_{s_1 -A_1}^{s_0 -A_0} [S(s_1-u)-S(s_0-u)] N(u)m(N(u))
m_{\rho}(t_0+u)du \\ &+ \int_{s_0 -A_0}^{s_1 -A_0} S(s_1-u)
N(u)m(N(u))m_{\rho}(t_0+u)du \Bigr|
\end{split} \\
&\leq m_0 \absj{s_1-s_0} + m_0\left(1-\frac{A_0}{A_1}\right)
\absj{s_0-s_1} + m_0 \absj{s_1-s_0}
\\
&\leq L \absj{s_0-s_1}\end{align}
\end{proof}

Nous pouvons maintenant d\'efinir \begin{equation}\begin{split}
\mathcal{K}_{t_0} = & \Bigl\{ N \in Y_{t_0}; \, \forall s \in
[-A_1 ; 0] , \, c_0 m_0 N_{\max}^{1-\gamma} \leq N(s) \leq N_{\max} , \\
&\forall s_0 , s_1 \in [-A_1 ; 0], \, \absj{N(s_0) - N(s_1) } \leq
L \absj{ s_0 - s_1} \Bigr\} \end{split}\end{equation} qui est une
partie compacte de $Y_{t_0}$ pour la topologie de la convergence
uniforme, d'apr\`es le th\'eor\`eme d'Ascoli.

Les lemmes que nous venons de d\'emontrer peuvent se formuler de
la fa\c{c}on suivante :
\begin{Pro} \label{proposition:K0}
Soit $N \in Y_0$, $(T^s)_{s \geq 0}$ le semi-groupe d\'efini par
l'\'equation \eqref{eq:semi_groupe}. On se place dans les
conditions pr\'ec\'edemment \'enonc\'ees pour les diff\'erents
param\`etres du mod\`ele.
\begin{enumerate}
\item Si $N \in \mathcal{K}_0$, alors $N^s \in \mathcal{K}_s$ pour
tout $s \geq 0$. En particulier $T^1(\mathcal{K}_0) \subset
\mathcal{K}_0$. \item En g\'en\'eral, il existe $s_0 \geq 0$
(d\'ependant de $N$) tel que $N^s \in \mathcal{K}_s$ pour tout $s
\geq s_0$.\end{enumerate}
\end{Pro}
\begin{proof}
\begin{enumerate}\item $N^s \in Y_s$ par d\'efinition, $N^s \leq
N_{max}$ d'apr\`es le lemme \ref{lemme:majoration}, $N^s$ reste
$L$-lipschitzienne d'apr\`es le lemme \ref{lemme:lipschitz}. La
partie 2 du lemme \ref{lemme:minoration} donne la minoration, en
utilisant que $N \in Y_0$.
\item On utilise le corollaire \ref{corollaire:minoration} pour
montrer l'existence de $s_0$, le reste de la preuve \'etant
identique.
\end{enumerate}
\end{proof}

L'\emph{attracteur} du syst\`eme dynamique $(Y^{\sharp},(T^s)_{s
\geq 0})$ est d\'efini par: \begin{gather}
\label{eq:attracteur}\Lambda = \left\{
(t,N) \telque \, t \in \R / \Z , \, N \in \Lambda_t \right\} \\
\text{avec } \Lambda_t = \bigcap_{n \geq 0} T^n (\mathcal{K}_t)
\end{gather}

La propri\'et\'e suivante justifie l'appellation d'attracteur pour
$\Lambda$.
\begin{Pro}\label{proposition:lambda}
\begin{enumerate}\item
$\Lambda$ est une partie compacte de $Y^{\sharp}$. \item Pour tout
$s \geq 0$, $T^s(\Lambda) = \Lambda$. \item Pour tout voisinage
$U$ de $\Lambda$, et toute condition initiale $(0,N)$, $N \in
Y_0$, il existe $s_0$ (d\'ependant de $N$ et $U$) tel que
$T^s(0,N) \in U$ pour tout $s \geq s_0$.
\end{enumerate} \end{Pro}
\begin{proof}
\begin{enumerate}
\item D'apr\`es la continuit\'e de $T^1$ et la compacit\'e de $\mathcal{K}_t$,
$\Lambda_t$ est compact pour tout $t \geq 0$. De plus,
$T^{\epsilon} \xrightarrow[\epsilon \rightarrow 0]{} Id$
uniform\'ement sur $Y^{\sharp}$ (c'est une cons\'equence du lemme
\ref{lemme:lipschitz}, car on a alors
$\norm{T^{\epsilon}(N)-N}_{\infty} \leq L \epsilon$). Par
cons\'equent, $\Lambda$ est compact.
\item Par construction, $T^s(\{t\} \times \Lambda_t) = \{t+s\} \times
\Lambda_{t+s}$ pour tous $s,t \in \R/\Z$, d'o\`u $T^s(\Lambda) =
\Lambda$ pour tout $s \geq 0$.
\item On peut supposer $N \in \mathcal{K}_0$ d'apr\`es la proposition
\ref{proposition:K0}. La suite $T^n(0,N)$ est alors contenue dans
$\mathcal{K} = \{ (t,N), t \in \R / \Z , N \in \mathcal{K}_t \}$
qui est compact. Tout point d'accumulation de cette suite est
n\'ecessairement dans $\Lambda$, et donc $T^n(0,N) \in U$ pour $n$
entier assez grand.

De m\^eme, en consid\'erant la suite $\left( T^{\alpha n}(0,N)
\right)_{n \in \N}$ avec $\alpha>0$ r\'eel quelconque, on montre
que $T^{\alpha n}(0,N) \in U$ pour $n \geq n(\alpha,U)$ entier.
Munissons $Y^{\sharp}$ de la distance \[
d((s,N),(t,\widetilde{N})) = \absj{s-t} +
\norm{N-\widetilde{N}}_{\infty}, \] qui engendre bien la topologie
de $Y^{\sharp}$ pr\'ec\'edemment d\'efinie. Puisque les
\'el\'ements de $\mathcal{K}_0$ sont $L$-lipschitziens, on a pour
tous $s,t \geq 0$, \begin{align} d(T^s(0,N),T^t(0,N)) &=
\absj{s-t} + \norm{N^{s} - N^{t}}_{\infty} \\ &\leq \absj{s-t} +
\norm{\overline{N}(s+\cdot) - \overline{N}(t+\cdot)}_{\infty} \\
&\leq \absj{s-t}(1+L). \end{align}

De plus, comme $\Lambda$ est compact, il existe $\epsilon >0$ tel
que $\Lambda \subset \Lambda^{(\epsilon)} \subset U$ o\`u l'on a
not\'e $\Lambda^{(\epsilon)}$ l'\'epaississement de $\epsilon$ de
$\Lambda$ (c'est-\`a-dire l'ensemble des points situ\'es \`a
distance $< \epsilon$ d'un point de $\Lambda$)\footnote{on
recouvre $\Lambda$ par des boules contenues dans $U$, un nombre
fini suffit par compacit\'e, $\epsilon$ est alors le $\min$ des
rayons de ces boules.}.

Prenons $\alpha = \frac{\epsilon}{2(1+L)}$, alors pour tout $t
\geq n(\alpha,\Lambda^{(\epsilon/2)})$ r\'eel, $T^t(0,N) \in
\Lambda^{(\epsilon)} \subset U$ ce qui ach\`eve la preuve.
\end{enumerate}\end{proof}

Il est alors ais\'e de faire le lien avec les d\'efinitions
donn\'ees \ref{def:attracteur} et \ref{def:bassin} donn\'ees en
annexe, sous la forme du corollaire suivant.

\begin{Cor} Le compact $\Lambda$ d\'efini par l'\'equation
\eqref{eq:attracteur} est un attracteur pour le syst\`eme
dynamique $\left( \left(T^s\right)_{s\geq 0} , Y^{\sharp}
\right)$. Son bassin d'attraction est $Y^{\sharp}$. \end{Cor}

\begin{proof}
La seule difficult\'e suppl\'ementaire par rapport \`a la
proposition \ref{proposition:lambda} est qu'il faut montrer
l'existence d'un voisinage de $\Lambda$ revenant tout entier dans
lui-m\^eme en un temps fini $N$. Pour l'instant, nous savons
seulement que toute condition initiale arrive en temps fini dans
un voisinage donn\'e, mais ce temps peut \^etre arbitrairement
grand en fonction de la condition initiale d\`es que l'on est hors
de $\Lambda$.

D'apr\`es le $2.$ de la proposition \ref{proposition:lambda}, il
revient au m\^eme de consid\'erer le syst\`eme dynamique discret
$(T^1,Y^{\sharp})$. Nous nous placerons d\'esormais dans ce cas.
Soit $N_1 \geq A_1$ un entier et $\alpha>0$ tel que
$(1-\alpha)\frac{c_0 m_0}{2} > 1$. On d\'efinit alors l'ensemble
\[ V = \left(T^{N_1}\right) ^{-1} \left( \left\{
(t,N) \telque t \in \R/\Z,\, N \in Y_t,\, N>N_{\max}^{1-\gamma}
\times \frac{c_0 m_0}{2} \times (1-\alpha) \right\} \right)
\] Par continuit\'e de $T^{N_1}$, c'est un ouvert. Par
d\'efinition de $\mathcal{K}_t$, il contient $\{ (t,N) \telque t
\in \R/\Z,\, N \in \mathcal{K}_t \}$, et donc $\Lambda$.

Les lemmes \ref{lemme:majoration}, \ref{lemme:minoration} et
\ref{lemme:lipschitz} montrent que $T^{2N_1}(V) \subset \{ (t,N)
\telque N \in \mathcal{K}_t \}$, d'o\`u $T^{2N_1}(V) \subset V$.

On a \'egalement \[ \bigcap_{n\geq 0} T^n(V) \subset \bigcap_{n
\geq 0} T^n(\{ (t,N) \telque N \in \mathcal{K}_t \}) = \Lambda \]
Comme de plus $\Lambda \subset V$ et $T^n(\Lambda)=\Lambda$, on a
$\Lambda = \bigcap_{n \geq 0} T^n(V)$.

Le $3.$ de la proposition \ref{proposition:lambda} montre que si
$x=(t,N) \in Y^{\sharp}$, il existe un temps $t_0$ \`a partir
duquel $T^s(x) \in \Lambda$, et donc n\'ecessairement $\omega(x)
\subset \Lambda$.
\end{proof}

\subsection{Mod\`ele non-saisonnier}
Le cas $m_{\rho} \equiv 1$ ($\rho = 0$) peut \^etre trait\'e plus
en d\'etails, au moins pour de petites valeurs de $\gamma$.

Il existe une valeur d'\'equilibre (\latin{i.e.} une solution
constante en temps continu) \begin{equation} N_{eq} = m^{-1}\left(
\frac{2A_1}{(A_1-A_0)^2}\right) \end{equation} pourvu que la
quantit\'e \[ \frac{2A_1}{(A_1-A_0)^2} = \left(\int_{A_0}^{A_1}
S(a) da \right)^{-1}
\] soit plus petite que $m_0$. Si elle est plus petite que
$\frac{m_0}{2}$, lorsque $m$ est d\'efinie par
\eqref{eq:fecondite}, on peut r\'e\'ecrire \[ N_{eq} = \left[ m_0
\frac{(A_1-A_0)^2}{2A_1}\right]^{1/\gamma}. \] C'est toujours le
cas pour des valeurs raisonnables des param\`etres\footnote{e.g.
$\gamma\geq 1$, $A_1 \geq 2$, $A_1 \geq 2 A_0$, $m_0 \geq 8$.}.

Pour d\'eterminer la stabilit\'e de cet \'equilibre, d\'efinissons
\begin{align} F(\lambda) &= \int_{A_0}^{A_1} S(a) e^{-a \lambda}
da \\ \label{eq:Flambda} &= \left( \frac{1}{\lambda} \left( 1-
\frac{A_0}{A_1} \right) - \frac{1}{\lambda^2 A_1} \right) e^{-A_0
\lambda} + \frac{1}{\lambda^2 A_1} e^{- A_1 \lambda}, \end{align}
les valeurs propres de la diff\'erentielle\footnote{On consid\`ere
le syst\`eme sous la forme $\frac{dN}{dt}=f(N)$, au voisinage de
$N_{eq}$ (on peut expliciter $f$, au voisinage de $N_{eq}$, en
consid\'erant $T^{\epsilon}$ quand $\epsilon$ tend vers 0). Les
solutions de la forme $N_{eq}+h$ v\'erifient $\frac{dh}{dt} =
Df_{N_{eq}}h$. Si $\lambda$ est valeur propre de la
diff\'erentielle \`a l'equilibre, un vecteur propre associ\'e $h$
est n\'ecessairement sous la forme $h(t)=h_0 e^{\lambda t}$. En
introduisant cette relation dans \eqref{eq:modele}, on peut ainsi
d\'eterminer l'ensemble des valeurs propres.} \`a l'\'equilibre
sont les solutions de \begin{equation} \label{eq:val_propres}
F(\lambda) = \left[ \frac{2A_1}{(A_1-A_0)^2} (1-\gamma)
\right]^{-1} \egaldef c_{\gamma}.
\end{equation}

On peut alors d\'efinir des valeurs $\gamma_0(A_0) < \gamma_1(A_0)
< \cdots < \gamma_k(A_0) < \cdots$ telles que : \begin{itemize}
\item si $\gamma < \gamma_0(A_0) \egaldef 1 +
\frac{(A_1-A_0)^2}{2A_1}\absj{F(-i u_0)}^{-1}$, l'\'equilibre est
stable. \item si $\gamma_0(A_0) < \gamma < \gamma_1(A_0) \egaldef
1 + \frac{(A_1-A_0)^2}{2A_1}\absj{F(-i u_2)}^{-1}$, il y a
exactement deux valeurs propres (complexes conjugu\'ees)
instables. \item \ldots \item si $\gamma_{k-1}(A_0) < \gamma <
\gamma_k(A_0) \egaldef 1 + \frac{(A_1-A_0)^2}{2A_1}\absj{F(-i
u_{2k})}^{-1}$, il y a exactement $2k$ valeurs propres
instables\footnote{$k$ paires de valeurs propres complexes
conjugu\'ees.}. \end{itemize}
\begin{proof}[\'El\'ements de preuve:]
Les valeurs propres de la diff\'erentielle sont stables si et
seulement si leur partie r\'elle est n\'egative, puisque le
vecteur propre associ\'e est de la forme $t \rightarrow
\exp(\lambda t)$. Comme $c_{\gamma}$ est un r\'eel (n\'egatif si
$\gamma > 1$), ce sont des solutions de l'\'equation
$\Im(F(\lambda))=0$. Les valeurs des bifurcations correspondant
\`a $\lambda$ imaginaire pur, on s'int\'eresse \`a l'\'equation
$\Im(F(-iu))=0$, avec $u>0$.

On montre alors que l'ensemble de ces solutions peut s'\'ecrire $
u_0 < u_1 < u_2 < \cdots$ avec \begin{equation*} \left\{
\begin{aligned} \Re(F(-i u_{2k})) &<
\Re(F(-i u_{2k+2})) &< \cdots &< 0 \\ \Re(F(-i u_{2k+1})) &>
\Re(F(-i u_{2k+3})) &> \cdots &> 0 \end{aligned} \right.
\end{equation*} En effet, \begin{align*} \Im(F(-iu)) &= \frac{-1}{u^2 A_1}
\sin(A_1 u) + \frac{1}{u} (1-\frac{A_0}{A_1}) \cos(A_0 u) +
\frac{1}{u^2 A_1} \sin(A_0 u) \\ \Re(F(-iu)) &= \frac{-1}{u^2 A_1}
\cos(A_1 u) + \frac{1}{u} (1-\frac{A_0}{A_1}) \sin(A_0 u) +
\frac{1}{u^2 A_1} \cos(A_0 u) \end{align*} et donc l'ensemble des
solutions n'a pas de point d'accumulation (en 0, on le v\'erifie
par un d\'eveloppement limit\'e de $\Im(F(-iu))$ ; ailleurs, cela
d\'ecoule clairement de l'analycit\'e de la fonction qui n'est pas
identiquement nulle), ce qui permet d'\'enum\'erer les solutions.
Il faut v\'erifier par un calcul direct que les in\'egalit\'es
annonc\'ees sont vraies pour les premi\`eres valeurs de $k$. Pour
$k$ grand, le terme en $1/u$ est dominant, et donc $\Im(F(-iu))$
s'annule presque en m\^eme temps que $\cos(A_0 u)$, et en ce point
$\Re(F(-iu))$ se comporte comme $\frac{1}{u} (1 - \frac{A_0}{A_1})
\sin(A_0 u)$. On en d\'eduit l'alternance des signes et la
d\'ecroissance des valeurs absolues.

Nous ne consid\'erons que les $u_k$ tels que $\Re(F(-iu_k))<0$ car
$c_{\gamma}<0$. \`A chaque $u_{2k}$, on associe alors un
$\gamma_k$ tel que $c_{\gamma_k} = F(-i u_{2k})$, \latin{i.e.} \[
\gamma_k(A_0, A_1) \egaldef 1 + \frac{(A_1-A_0)^2}{2A_1}\absj{F(-i
u_{2k})}^{-1}. \] La d\'ecroissance des valeurs absolues des $F(-i
u_{2k})$ montre que les $\gamma_k$ sont ordonn\'es par ordre
croissant.

Nous venons de montrer que lorsque $\gamma$ varie, une paire de
valeurs propres traverse l'axe imaginaire en chaque $\gamma_k$, et
seulement en ces points-l\`a. Lorsque $\gamma$ tend vers 1 par
valeurs sup\'erieures, $c_{\gamma}$ d\'ecro\^it vers $-\infty$, et
donc les valeurs propres $\lambda$ doivent rendre $F(\lambda)$ de
plus en plus grand en module et n\'egatif. Or, le module de
$F(\lambda)$ est born\'e sur le demi-plan $\Re(\lambda)>0$
(d'apr\`es l'\'equation \eqref{eq:Flambda}), donc pour $\gamma$
assez proche de 1, toutes les valeurs propres ont une partie
r\'eelle n\'egative. C'est donc le cas pour tout $\gamma <
\gamma_0$.
\end{proof}

En r\'ealit\'e, $\gamma_0, \gamma_1, \ldots$ d\'ependent peu de
$A_0$, et leurs valeurs typiques sont $\gamma_0 \approx 6\virg 2$
et $\gamma_1 \approx 30$.

Lorsque $\gamma$ grandit et traverse $\gamma_0$ ($A_0$ \'etant
fix\'e), on s'attend \`a voir une bifurcation de Hopf (voir
annexe~\ref{annexe:Hopf}) : une orbite p\'eriodique attractive est
cr\'e\'ee au voisinage de l'\'equilibre pour $\gamma = \gamma_0$,
et attire toutes les solutions proches de l'\'equilibre (sauf
l'\'equilibre lui-m\^eme) pour $\gamma > \gamma_0$ proche de
$\gamma_0$.
\section{Simulations num\'eriques}
La mise en \oe uvre de simulations du mod\`ele \eqref{eq:modele} a
demand\'e un travail pr\'eliminaire de discr\'etisation et de mise
au point des param\`etres de simulation qui est d\'etaill\'e en
annexe~\ref{annexe:sim_preliminaire}. Les pr\'ecisions concernant
le traitement des donn\'ees sont donn\'ees en
annexe~\ref{annexe:traitement}, dans l'ordre de pr\'esentation des
r\'esultats. Elles sont cependant indispensables pour une r\'eelle
compr\'ehension de ceux-ci, car les nombreuses approximations qui
ont \'et\'e n\'ecessaires ont souvent une r\'eelle influence sur
les r\'esultats obtenus.

\subsection{Explorations de l'espace des param\`etres}
Pour commencer, on fait varier un param\`etre en gardant les
autres fixes, et l'on observe la fa\c{c}on dont la dynamique
stationnaire \'evolue. Trois param\`etres semblent d\'eterminants
: $A_0$ (qui introduit un effet de retard correspondant au temps
de maturation), $\rho$ (qui mesure l'importance du facteur
saisonnier) et $\gamma$ (qui traduit l'influence de la densit\'e
sur la f\'econdit\'e). Pour les autres param\`etres, on a fix\'e
les valeurs suivantes :
\begin{itemize}
\item $A_1 = 2$.
\item $m_0 = 50$ et la f\'econdit\'e est prise $C^1$.
\item $\epsilon_{ete} = 0\virg 1$ (ou 0 dans le premier cas, car le facteur
saisonnier n'avait pas encore \'et\'e r\'egularis\'e).
\end{itemize}
On pourra ainsi repr\'esenter chaque simulation par le triplet
$(A_0; \rho; \gamma)$ qui lui correspond.

\subsubsection{$A_0=0.18$, $\rho=0.41$, $\gamma$ variable}
\begin{figure}
\begin{center}
\includegraphics[width=\textwidth]{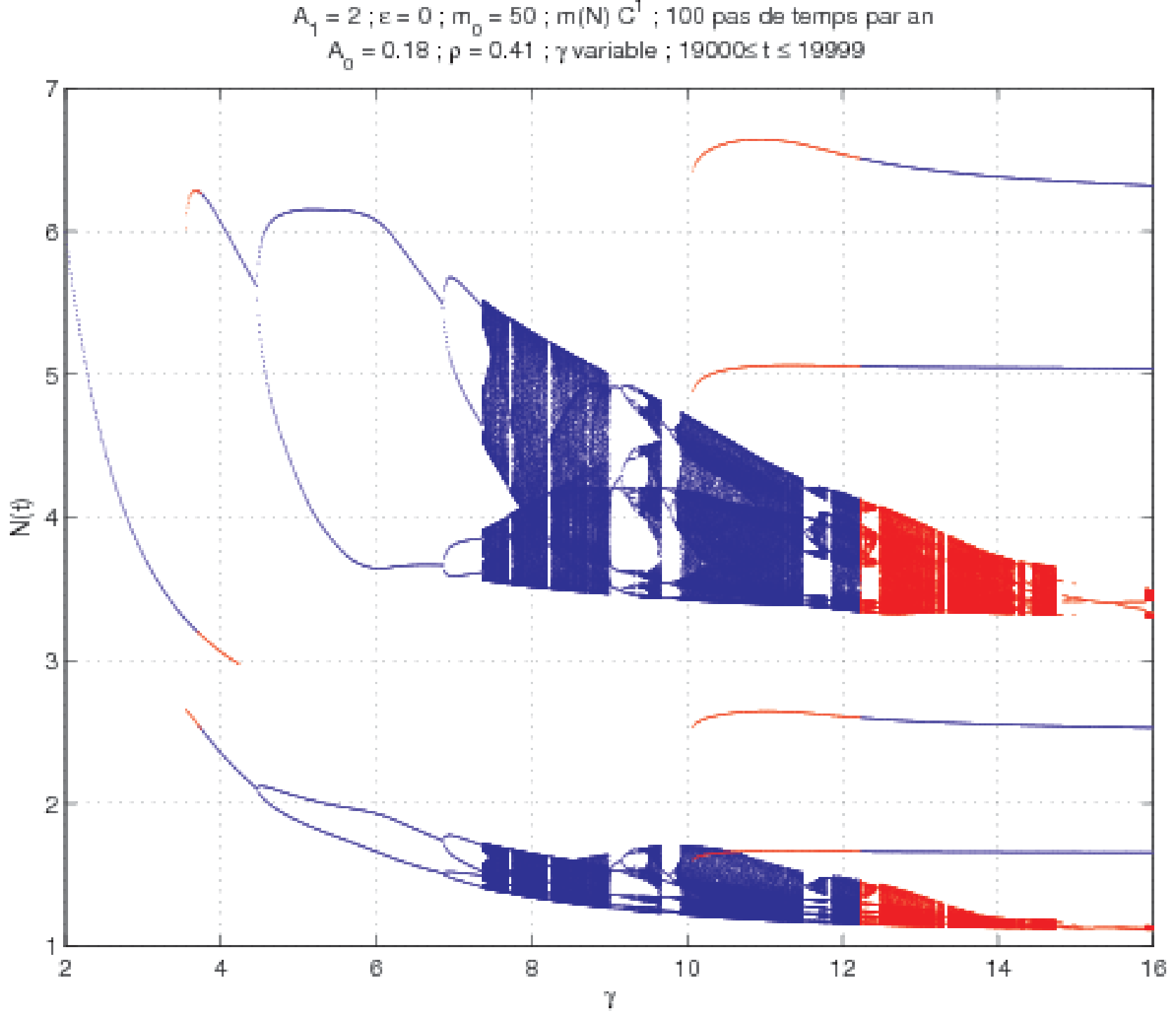}
\caption{\label{diag:0.18_0.41_gamma}Diagramme de bifurcation
$(0\virg 18; 0\virg 41; \gamma)_{2 \leq \gamma \leq 16}$. Pour
certaines valeurs de $\gamma$, plusieurs attracteurs coexistent.
On les a obtenus en utilisant la condition (I) (en bleu) et en
prolongeant aussi loin que possible (en rouge) certaines branches
interrompues dans le diagramme bleu. Noter que les saisons sont
discontinues ($\epsilon_{ete} = 0$).}
\end{center}
\end{figure}

Pour la premi\`ere exploration, nous sommes partis des valeurs
($0\virg 18$; $0\virg 41$; $8\virg 25$) et nous avons fait varier
$\gamma$. Contrairement aux simulations effectu\'ees
ult\'erieurement, le facteur saisonnier n'est pas r\'egularis\'e
(\latin{i.e.} $\epsilon=0$). Les r\'esultats sont repr\'esent\'es
sur un diagramme de bifurcation, figure~\ref{diag:0.18_0.41_gamma}
: pour chaque valeur de $\gamma$ sont repr\'esent\'ees les valeurs
de $N(t)$ aux temps entiers (\latin{i.e.} \`a la fin de
l'\'et\'e), en se limitant \`a $t$ assez grand (on a fix\'e
arbitrairement $19001 \leq t \leq 20000$). Les points bleus
correspondent \`a une m\^eme condition initiale (obtenue
al\'eatoirement), que nous notons (I) (voir
figure~\ref{fig:cond_I}). La partie rouge correspond \`a d'autres
simulations, d\'etaill\'ees ci-apr\`es.

\paragraph{Orbites p\'eriodiques attractives}

Pour $2 \leq \gamma \leq 7\virg 2$, le comportement observ\'e est
assez classique. On a d'abord un \'equilibre attractif, puis une
orbite de p\'eriode 2, et ainsi de suite avec des doublements de
p\'eriodes successifs (de plus en plus rapproch\'es) au fur-et-\`a
mesure que $\gamma$ se rapproche de la valeur limite
$\gamma_{\infty,1} \approx 7\virg 36$ . Il s'agit de bifurcations
par doublement de p\'eriode (voir la section
\ref{def:doublement_periode}) qui ont lieu pour $\gamma_0 <
\gamma_1 < \cdots < \gamma_n < \cdots < \gamma_{\infty,1}$, et qui
se traduisent par des doublements de p\'eriode successifs.

Pour $\gamma \geq \gamma_{\infty,1}$, le diagramme de bifurcations
permet de distinguer essentiellement deux comportements. D'une
part, il y a toujours des orbites p\'eriodiques attractives sur
certains intervalles de valeurs de $\gamma$, tout comme il y a des
<<fen\^etres de p\'eriodicit\'e>> dans le cas des polyn\^omes
quadratiques (voir section \ref{annexe:poly_quadratiques}). C'est
en particulier le cas pour $9\virg 00 \leq \gamma \leq 9\virg 39$,
$9\virg 67 \leq \gamma \leq 9\virg 86$ et $12\virg 23 \leq \gamma
\leq 16$.
On peut situer plus pr\'ecis\'ement ces fen\^etres \`a l'aide d'un
calcul de dimensions fractales.

\begin{figure}
\begin{center}
\includegraphics[height=7cm]{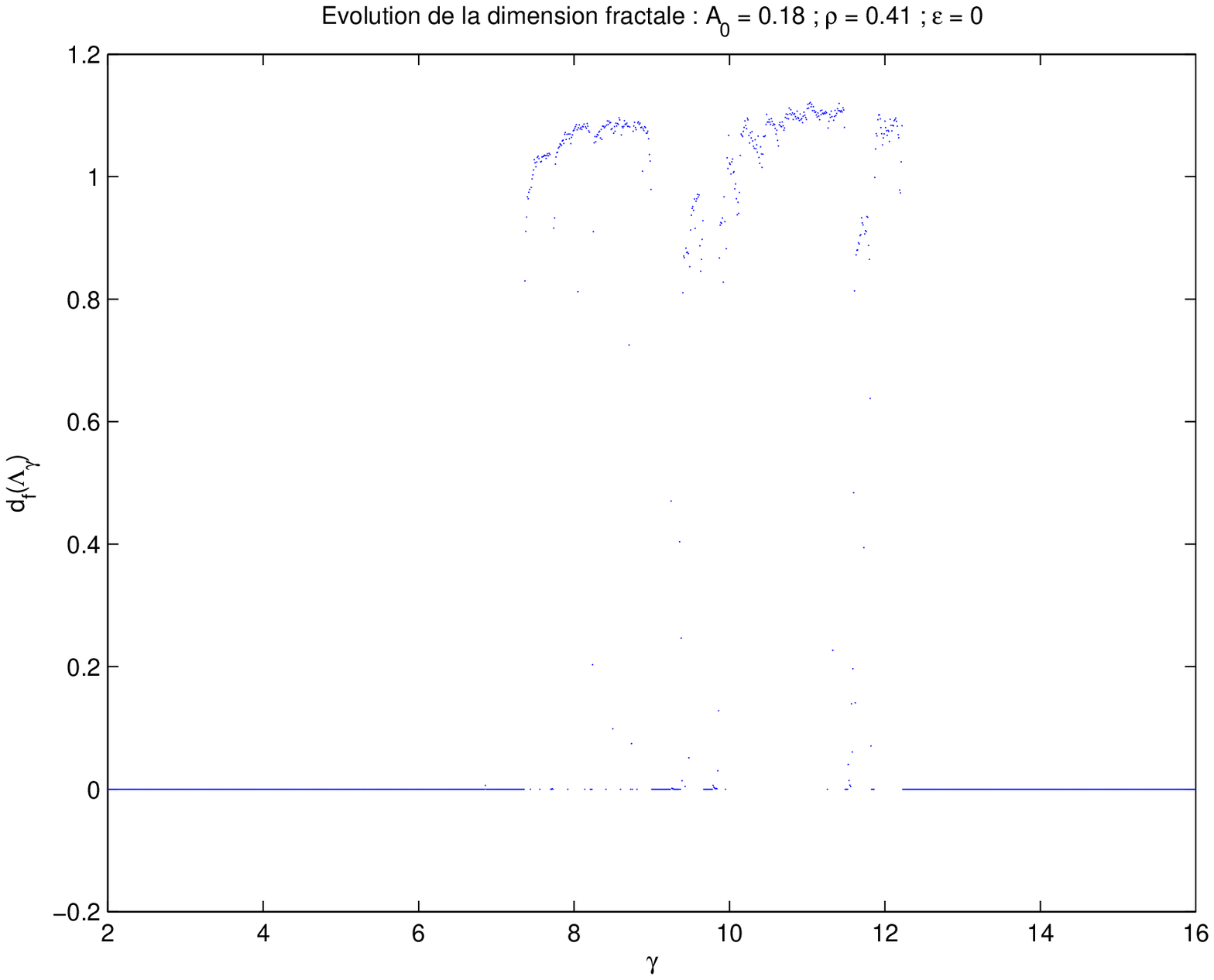}
\caption{\label{dim_f:0.18_0.41_gamma}Dimension fractale des
attracteurs $(0\virg 18; 0\virg 41; \gamma)_{2 \leq \gamma \leq
16}$. On s'est limit\'e \`a la condition initiale (I).}
\end{center}
\end{figure}

\paragraph{Dimensions fractales}
On peut ais\'ement calculer une valeur approch\'ee de la dimension
fractale des ensembles limites correspondant aux diff\'erentes
valeurs de $\gamma$. Le graphique obtenu est repr\'esent\'e
figure~\ref{dim_f:0.18_0.41_gamma}.

Les points o\`u la dimension fractale est nulle correspond aux
orbites p\'eriodiques attractives, ce qui nous permet de les
d\'etecter bien plus facilement qu'en observant le diagramme de
bifurcation. Dans le domaine interm\'ediaire ($7\virg 36 \leq
\gamma \leq 12\virg 23$), on constate qu'il y a alternativement
des attracteurs de dimension environ \'egale \`a 1 et des
fen\^etres de p\'eriodicit\'e. On peut raisonnablement penser
qu'il existe un ouvert dense dans l'espace des param\`etres pour
lequel il existe une orbite p\'eriodique attractive\footnote{qui
peut \'eventuellement coexister avec un autre attracteur.}. C'est
ouvert n'est en revanche certainement pas de mesure totale.

La pr\'ecision de ces calculs \'etant limit\'ee par le faible
nombre de points consid\'er\'es pour chaque ensemble, il n'est pas
ais\'e de d\'eterminer s'il y a ou non r\'eellement des
attracteurs de dimension non-enti\`ere. Si oui, elle n'est pas
tr\`es grande, certainement inf\'erieure \`a $1\virg 5$, et
probablement sup\'erieure \`a 1, dans la mesure o\`u le calcul
effectu\'e sous-estime l\'eg\`erement la dimension de l'attracteur
(en particulier, le petit nombre de points utilis\'es peut cr\'eer
de nombreux <<trous>> correspondant \`a la mesure physique sur
l'attracteur, et non \`a un trou r\'eel dans sa g\'eom\'etrie).

Il y a un autre argument th\'eorique en faveur d'une dimension
fractale sup\'erieure \`a 1 lorsqu'elle n'est pas nulle. En effet,
s'il y a une orbite p\'eriodique hyperbolique instable, sa
vari\'et\'e instable est contenue dans l'attracteur, qui doit donc
avoir une dimension au moins \'egale \`a 1. Il est donc
difficilement concevable que dans un grand domaine de l'espace des
param\`etres on puisse avoir un attracteur de dimension fractale
comprise strictement entre 0 et 1.

\paragraph{Attracteurs de type H\'enon}

\begin{figure}
\begin{center}
\includegraphics[height=7cm]{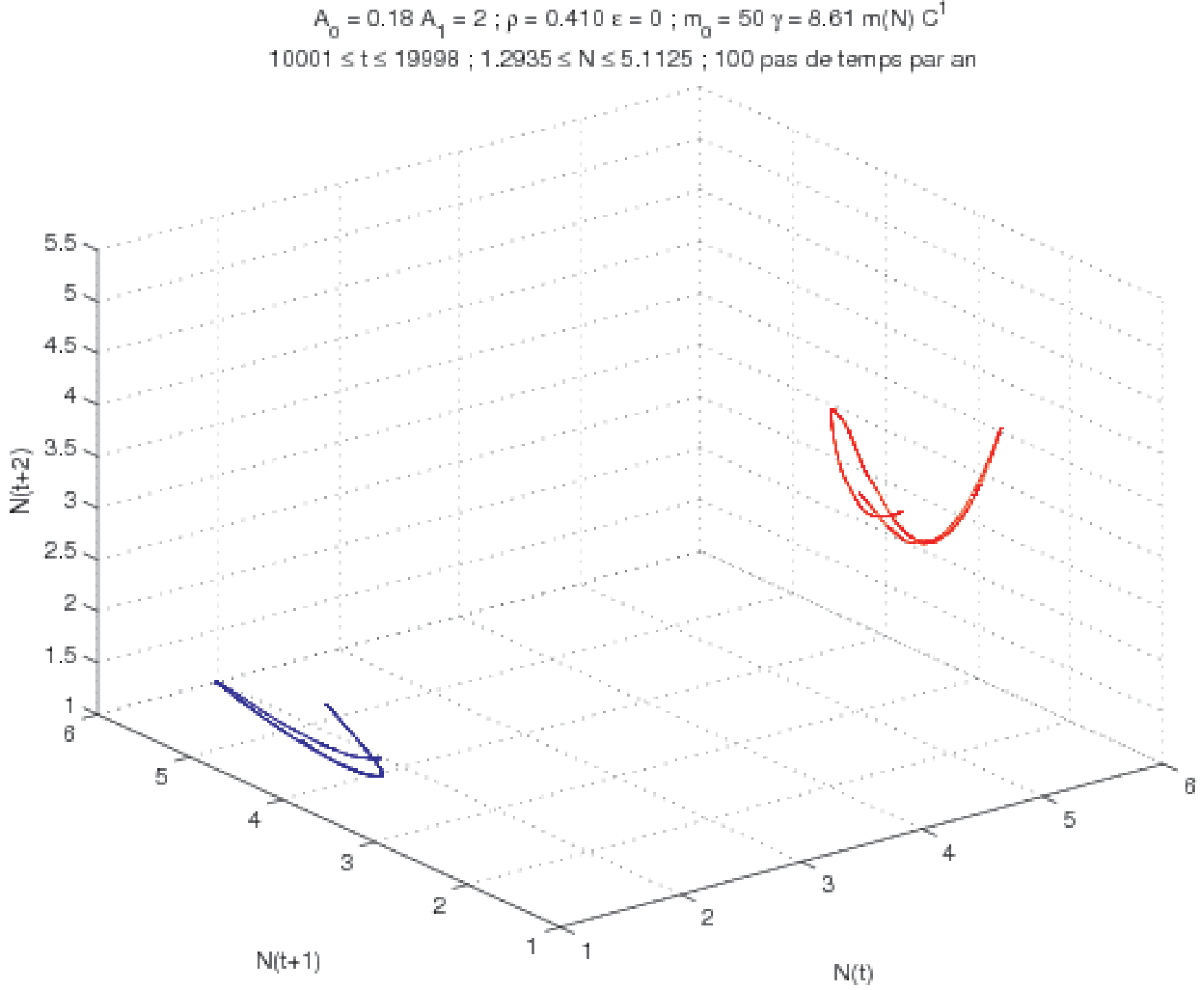}
\caption{\label{fig:0.18_0.41_8.61}$(0\virg 18; 0\virg 41; 8\virg
61)$. Condition initiale (I). Les deux composantes sont
invariantes par $T^2$, qui y semble topologiquement m\'elangeante.
Dimension fractale estim\'ee : $d_f \approx 1\virg 06$.}
\end{center}
\end{figure}

En-dehors des fen\^etres de p\'eriodicit\'e appara\^it un
comportement stationnaire non-p\'eriodique, le long d'un
attracteur qui semble constitu\'e de deux morceaux de courbes. Un
exemple est repr\'esent\'e avec la
figure~\ref{fig:0.18_0.41_8.61}, o\`u l'on a trac\'e les points
($N(t)$, $N(t+1)$, $N(t+2)$) pour chaque valeur enti\`ere de $t$
($10002 \leq t \leq 19999$). On s'int\'eresse alors \`a la
dynamique de l'application\footnote{Cette application n'est pas
parfaitement bien d\'efinie, le syst\`eme \'etudi\'e \'etant de
dimension infinie, alors que la visualisation consid\'er\'ee est
une projection de celui-ci en dimension 3.} $f$ : ($N(t)$,
$N(t+1)$, $N(t+2)$) $\mapsto$ ($N(t+1)$, $N(t+2)$, $N(t+3)$). Par
abus de notation, on \'ecrira $T^1$ au lieu de $f$, sans perdre de
vue que nous ne pouvons pas visualiser directement $T^1$. En
utilisant deux couleurs suivant la parit\'e de $t$, on constate
que chacune des deux parties de l'attracteur est envoy\'ee sur
l'autre. En revanche, il ne semble pas possible (pour cette valeur
de $\gamma$) de s\'eparer de la m\^eme fa\c{c}on l'attracteur en
un plus grand nombre de composantes. Il semble donc que
l'application $f^2$ restreinte \`a chacune des deux composantes de
l'attracteur soit topologiquement m\'elangeante (voir d\'efinition
\ref{def:melange}).

Un zoom sur l'attracteur fait appara\^itre des structures
semblables \`a celles de l'attracteur de H\'enon (voir
annexe~\ref{annexe:Henon}). Il est assez probable que la dynamique
soit du m\^eme type, mais cela n'\`a pas \'et\'e \'etudi\'e
pr\'ecis\'ement. La dimension fractale n'est pas clairement
diff\'erente de 1, mais elle est s\^urement sous-estim\'ee \`a
cause du petit nombre de points que nous avons calcul\'e. Si elle
s'av\'erait \^etre clairement diff\'erente de 1, cela renforcerait
l'hypoth\`ese d'une dynamique de type H\'enon.

\begin{figure}
\begin{center}
\includegraphics[width=\textwidth]{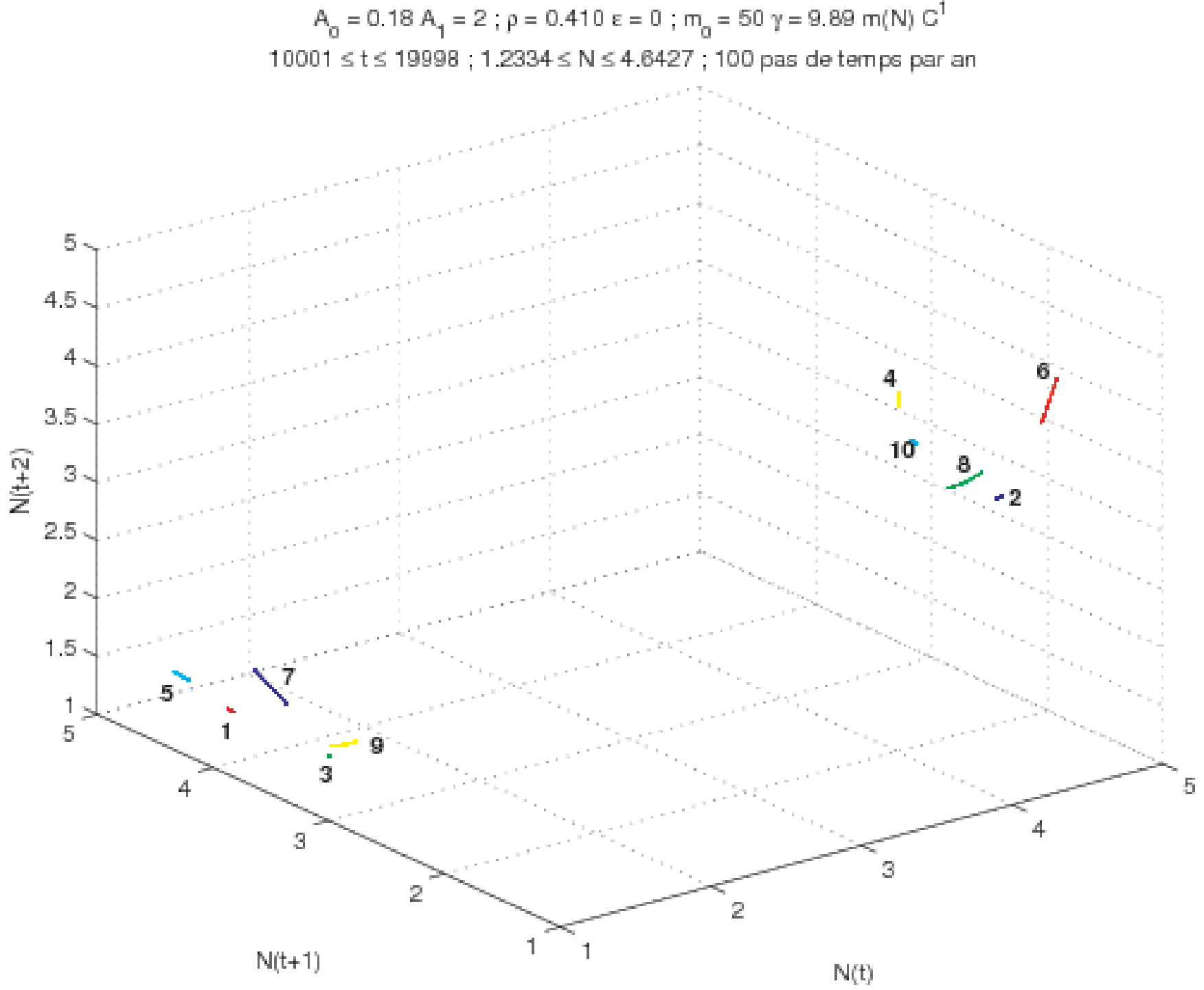}
\caption{\label{fig:0.18_0.41_9.89}$(0\virg 18; 0\virg 41; 9\virg
89)$. Condition initiale (I). Chacune des composantes est
invariante par $T^{10}$, qui y semble topologiquement
m\'elangeante. La num\'erotation des composantes correspond aux
valeurs de $t \modulo 10$ . Comme \`a la
figure~\ref{fig:0.18_0.41_8.61}, il y a deux groupes suivant la
parit\'e de $t$. Dimension fractale estim\'ee : $d_f \approx
0\virg 92$.}
\end{center}
\end{figure}

Il n'y a pas toujours deux composantes. On le voit sur le
diagramme de bifurcations notamment au voisinage des valeurs
<<limites>> de $\gamma$ (\latin{i.e.} juste apr\`es une zone o\`u
il y a une orbite p\'eriodique attractive). Ainsi, pour
$\gamma=8\virg 62$, on peut compter 10 composantes distinctes, et
$f^{10}$ semble topologiquement m\'elangeante sur chacune d'entre
elles (figure~\ref{fig:0.18_0.41_9.89}). Dans les deux cas
(figures \ref{fig:0.18_0.41_8.61} et \ref{fig:0.18_0.41_9.89}), il
semble donc qu'on ait une \emph{d\'ecomposition spectrale} (voir
th\'eor\`eme \ref{the:decompo_spectrale}) de l'attracteur
$\Lambda$ en un nombre fini de composantes (respectivement 2 et
10).

Lorsque $\gamma$ \'evolue\footnote{voir aussi l'animation
\fichier{film\_gamma\_100\_18\_200\_410\_00\_50\_\_1.avi}
(diagramme de bifurcation en quatre dimensions).}, on observe que
l'attracteur grandit petit-\`a-petit \`a l'int\'erieur d'un m\^eme
objet\footnote{\`A de l\'eg\`eres d\'eformations pr\`es, notamment
un changement de taille.}, que l'on visualise \`a peu pr\`es en
entier avec la figure~\ref{fig:0.18_0.41_8.61}. Tant que $\gamma
\leq 9\virg 96$, on remarque m\^eme que les orbites p\'eriodiques
attractives sont contenues dans le m\^eme objet. Lorsque
l'attracteur est continu mais en morceaux distincts, envoy\'es
p\'eriodiquement l'un dans l'autre, chaque morceau grandit et
ceux-ci fusionnent petit-\`a petit. Lorsque deux composantes
fusionnent, la stabilit\'e de chacune vis-\`a-vis de $f$
(compos\'ee le bon nombre de fois) semble instantan\'ement perdue,
et il y a alors m\'elange topologique \`a l'int\'erieur de chaque
composante. Cette fusion des composantes connexes s'apparente \`a
la cascade inverse qui suit la cascade harmonique directe, comme
c'est le cas pour la famille quadratique r\'eelle
(annexe~\ref{annexe:quadr_reelle}).

\paragraph{Discontinuit\'es du diagramme}

Les discontinuit\'es observ\'ees \`a $\gamma = 3\virg 70$ et
$\gamma = 12\virg 23$ sont plut\^ot surprenantes. L'hypoth\`ese la
plus plausible serait qu'il existe \`a ces valeurs de $\gamma$
deux attracteurs distincts, et la condition initiale (I) passe
brusquement du bassin de l'un au bassin de l'autre. Pour tester
cette hypoth\`ese, nous avons choisi comme nouvelles conditions
initiales les \'etats stationnaires obtenus de part et d'autre de
ces discontinuit\'es, et nous avons fait varier $\gamma$ pour
d\'eterminer s'il y a effectivement coexistence de deux
attracteurs pour certaines valeurs de $\gamma$.

\begin{figure}
\begin{center}
\includegraphics[height=7cm]{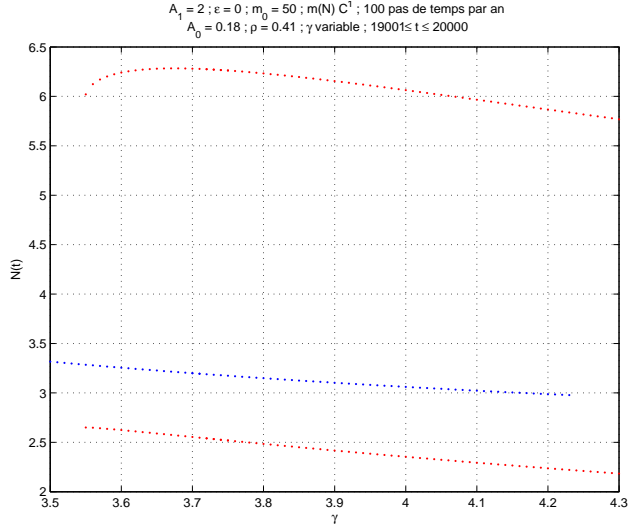}
\caption{\label{diag:0.18_0.41_gamma_3.70}Diagramme de bifurcation
$(0\virg 18; 0\virg 41; \gamma)_{3\virg 50 \leq \gamma \leq 4\virg
30}$. Deux attracteurs coexistent pour certaines valeurs de
$\gamma$ : un point fixe (au centre) et une orbite de p\'eriode 2
(en haut et en bas). Avec la condition (I), on passe de l'un \`a
l'autre pour $3\virg 70 < \gamma < 3\virg 75$.}
\end{center}
\end{figure}

En repartant de l'\'equilibre obtenu avec $\gamma = 3\virg 70$,
nous avons pu continuer la branche du diagramme jusqu'\`a $\gamma
= 4\virg 20$. Ensuite, on retrouve l'orbite de p\'eriode 2
d\'ej\`a trouv\'ee. Le m\^eme proc\'ed\'e nous permet de continuer
jusqu'\`a $\gamma = 4\virg 23$, mais nous n'avons plus pu
retrouver l'\'equilibre ensuite en utilisant le m\^eme
proc\'ed\'e. Dans l'autre sens, on prolonge le domaine o\`u se
trouve une orbite de p\'eriode 2 attractive jusqu'\`a $\gamma =
3\virg 60$. On peut ainsi tracer un nouveau diagramme de
bifurcations autour de ces valeurs, avec cette fois les deux
\'etats stationnaires (figure~\ref{diag:0.18_0.41_gamma_3.70}).

La seconde discontinuit\'e a lieu autour de $\gamma=12\virg 23$ :
on observe l'attracteur pour $\gamma = 12\virg 22$ et une orbite
attractive de p\'eriode 4 pour $\gamma = 12\virg 23$. En repartant
de l'\'etat final obtenu pour $\gamma = 12\virg 22$, on observe
pour les valeurs sup\'erieures de $\gamma$ (au moins jusqu'\`a 16)
un comportement similaire \`a ce qu'on constatait pour $\gamma
\leq 12\virg 22$, c'est-\`a-dire le m\^eme attracteur, avec
parfois des orbites p\'eriodiques attractives (mais dans un
domaine de valeurs de $N$ diff\'erent de l'orbite de p\'eriode 4).
Inversement, l'orbite 4-p\'eriodique attractive persiste jusqu'\`a
$\gamma \approx 10\virg 078$.

Les attracteurs ainsi d\'etect\'es sont repr\'esent\'es en rouge
sur la figure~\ref{diag:0.18_0.41_gamma}. Dans ces deux cas,
plusieurs attracteurs coexistent, mais il y a toujours des
discontinuit\'es dans le diagramme. Pour certaines valeurs de
$\gamma$ (envion $3\virg 60$ et $10\virg 078$), des orbites
p\'eriodiques deviennent attractives. A l'inverse, pour $\gamma
\approx 4\virg 23$, l'\'equilibre devient instable. Il pourrait
donc s'agir d'une bifurcation du type de $f_{(+1),\tau}$ (voir
annexe~\ref{annexe:bif1}, figure~\ref{fig:bif1_+1_diag}). Il est
\'egalement possible que le diagramme soit en r\'ealit\'e continu,
mais que les bassins d'attraction des orbites p\'eriodiques
attractives soient trop r\'eduits pour que l'on puisse les
atteindre par des simulations, avec la m\'ethode que nous avons
employ\'ee ici. Un petit travail th\'eorique serait n\'ecessaire
pour \'eclairer ce point.

\subsubsection{$A_0=0\virg 18$, $\rho=0\virg 30$, $\gamma$ variable}
\begin{figure}
\begin{center}
\includegraphics[width=\textwidth]{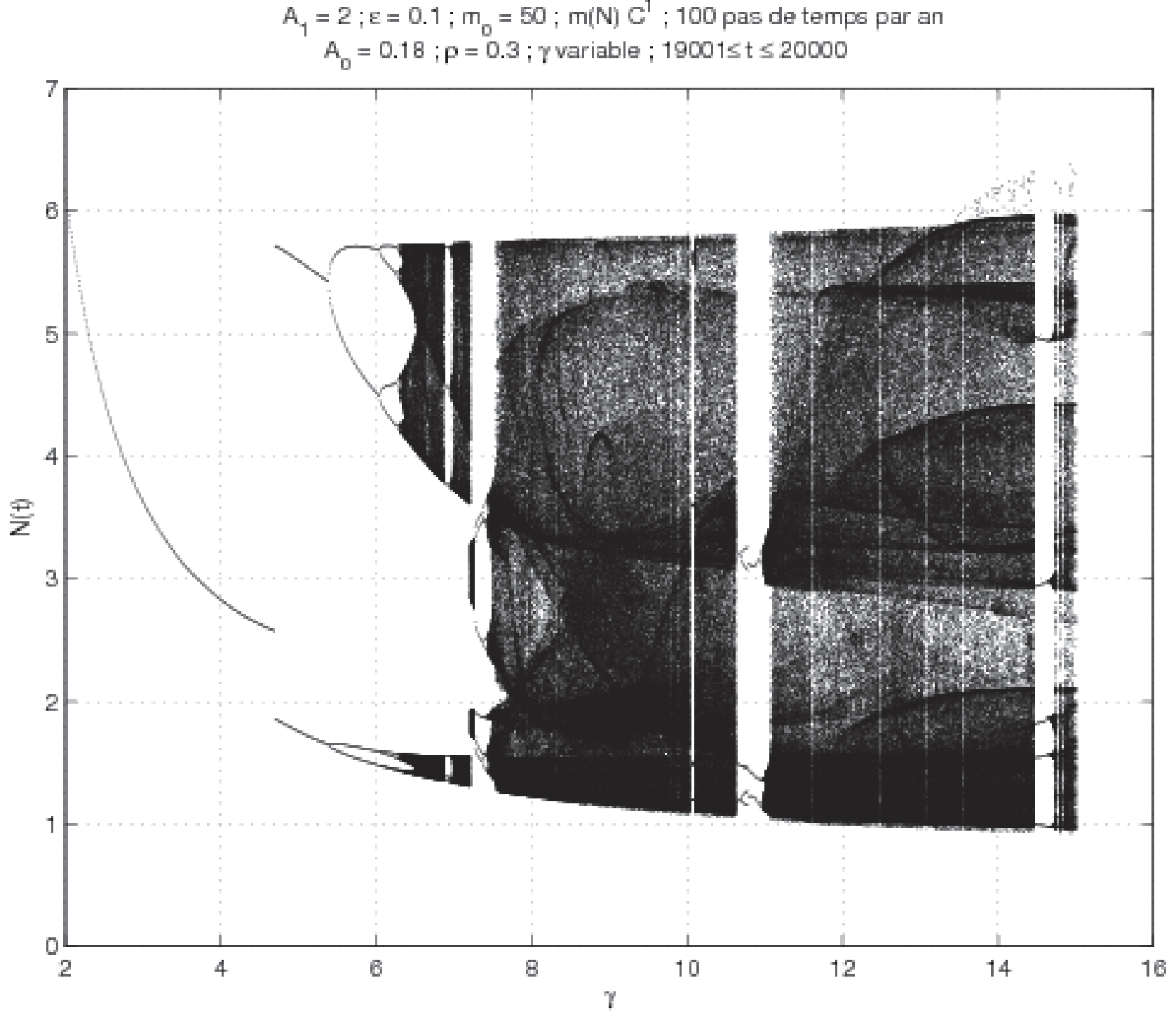}
\caption{\label{diag:0.18_0.30_gamma}Diagramme de bifurcation
$(0\virg 18; 0\virg 30; \gamma)_{2 \leq \gamma \leq 15}$. Les
saisons sont continues ($\epsilon_{ete} = 0\virg 1$). Condition
initiale (I).}
\end{center}
\end{figure}

Suite \`a une rapide exploration en faisant varier $\rho$, et au
vu de l'int\'er\^et des valeurs $(0\virg 18; 0\virg 30; 8\virg
25)$ (voir figure~\ref{fig:0.18_0.30_8.25}), nous avons effectu\'e
une deuxi\`eme exploration \`a $\gamma$ variable, autour de ces
nouvelles valeurs. Le diagramme de bifurcation ainsi obtenu est
repr\'esent\'e figure~\ref{diag:0.18_0.30_gamma}. Il est tr\`es
semblable au diagramme~\ref{diag:0.18_0.41_gamma}, mise \`a part
l'absence de deux composantes bien distinctes pour la plupart des
valeurs de $\gamma$, et le faible nombre de fen\^etres de
p\'eriodicit\'e.

\subsubsection{$A_0=0\virg 18$, $\rho$ variable, $\gamma=8\virg 25$}
\begin{figure}
\begin{center}
\includegraphics[width=\textwidth]{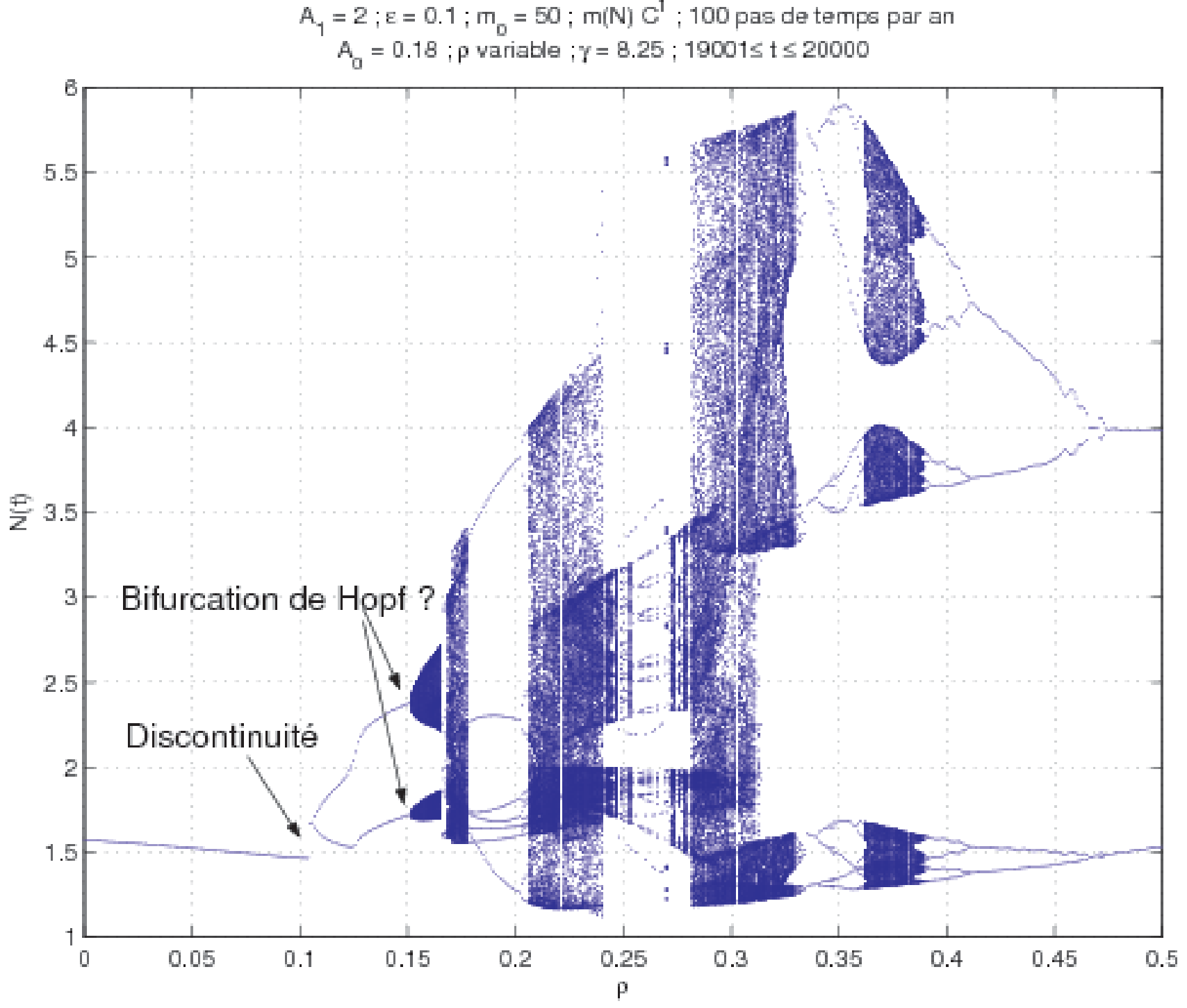}
\caption{\label{diag:0.18_rho_8.25}Diagramme de bifurcation
$(0\virg 18; \rho; 8\virg 25)_{0 \leq \rho \leq 0\virg 5}$. Les
saisons sont continues ($\epsilon_{ete} = 0\virg 1$). Condition
initiale (II).}
\end{center}
\end{figure}

\begin{figure}
\begin{center}
\includegraphics[width=\textwidth]{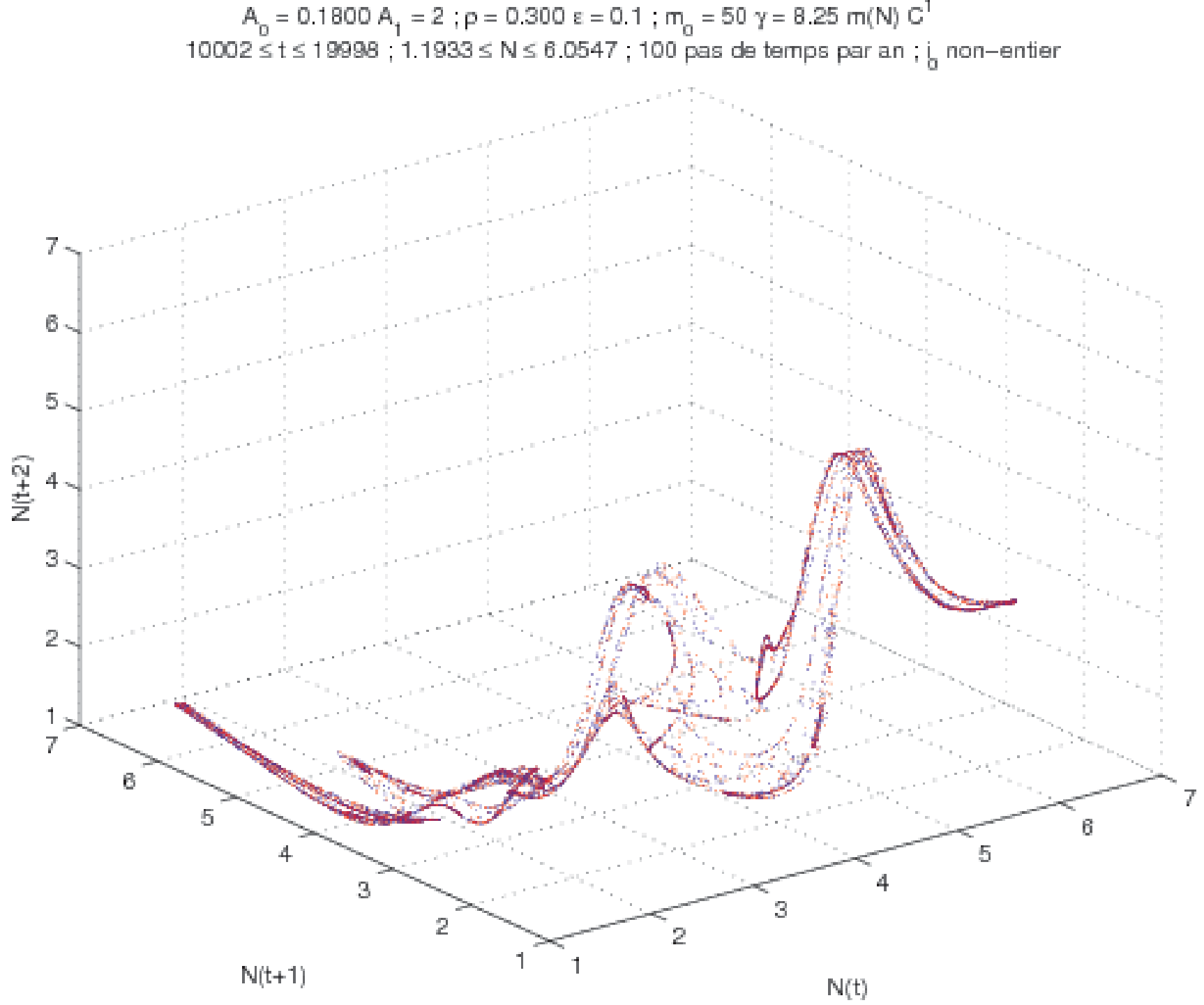}
\caption{\label{fig:0.18_0.30_8.25}Un comportement complexe :
$(0\virg 18; 0\virg 30; 8\virg 25)$. Les saisons sont continues
($\epsilon_{ete} = 0\virg 1$). Condition initiale (II). Il semble
que $T^1$ soit m\'elangeante sur cet attracteur. Dimension
fractale estim\'ee : $d_f \approx 1\virg 19$.}
\end{center}
\end{figure}

Le facteur saisonnier semble jouer un r\^ole d\'eterminant dans la
dynamique du syst\`eme (le mod\`ele non-saisonnier est
particuli\`erement simple, alors que pour des valeurs de $\rho$
plus proches de la r\'ealit\'e, on observe des comportements bien
plus complexes, \latin{e.g.} figure~\ref{fig:0.18_0.30_8.25}). Le
diagramme de bifurcation obtenu en faisant varier $\rho$ est
repr\'esent\'e figure~\ref{diag:0.18_rho_8.25}. Une autre
condition initiale, not\'ee (II), a \'et\'e utilis\'ee pour ces
simulations (voir figure~\ref{fig:cond_II}). On observe comme
pr\'ec\'edemment une discontinuit\'e dans le diagramme, pour
$\rho$ proche de $0\virg 1$, mais nous n'avons pas essay\'e de
prolonger les deux branches interrompues.

\paragraph{Bifurcation de Hopf}

\begin{figure}
\begin{center}
\includegraphics[width=\textwidth]{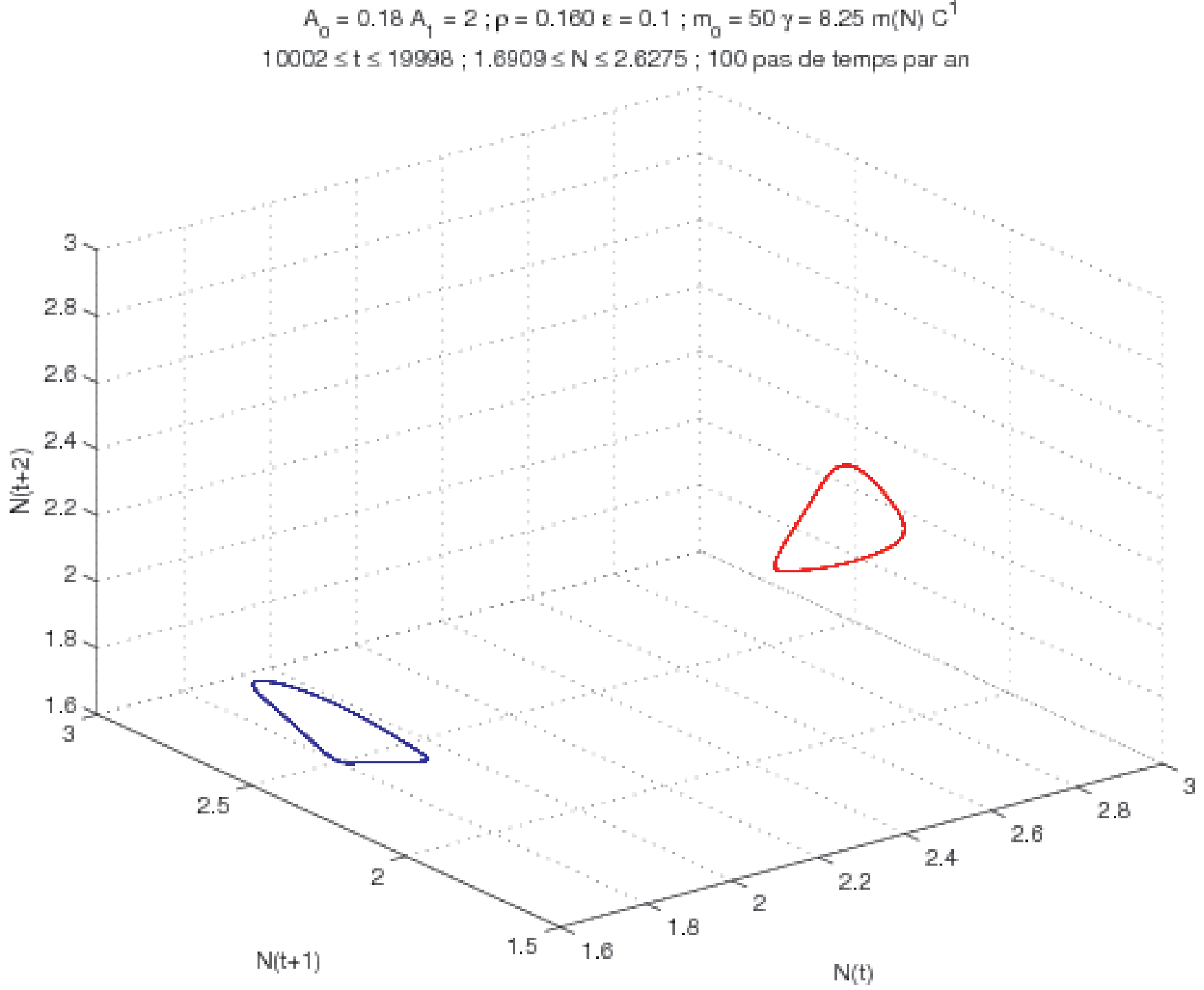}
\caption{\label{fig:0.18_0.16_8.25}Deux <<cycles>> : $(0\virg 18;
0\virg 16; 8\virg 25)$. Condition initiale (II). Chaque composante
est invariante par $T^2$. Dimension fractale estim\'ee : $d_f
\approx 0\virg 99$.}
\end{center}
\end{figure}

On constate sur ce diagramme un ph\'enom\`ene que nous n'avions
pas trouv\'e sur les diagrammes pr\'ec\'edents : une bifurcation
de Hopf pour $\rho \approx 0\virg 152$. On passe en effet d'une
orbite attractive de p\'eriode 2 \`a deux <<cycles>> attractifs
stables, qui persistent jusqu\`a $\rho \approx 0\virg 165$ (la
figure~\ref{fig:0.18_0.16_8.25} en repr\'esente un exemple).

Cet attracteur est de dimension fractale 1 et chaque lacet est
parfaitement connexe. Il n'est en revanche pas totalement certain
que $T^2$ soit bien topologiquement m\'elangeante sur chacun, bien
que l'on n'ait vu aucune p\'eriodicit\'e \'evidente. La dynamique
de $f$ sur ces cycles n'est pas forc\'ement simple (\latin{i.e.}
topologiquement conjugu\'ee \`a une rotation), et le cycle ne
co\"incide peut-\^etre pas exactement avec l'attracteur.
R\'epondre \`a ces questions demanderait une \'etude plus
pouss\'ee.

\subsubsection{$A_0$ variable, $\rho=0\virg 3$, $\gamma=8\virg 25$}
\label{sec:explor_A0_0.30_8.25}

\begin{figure}
\begin{center}
\includegraphics[width=\textwidth]{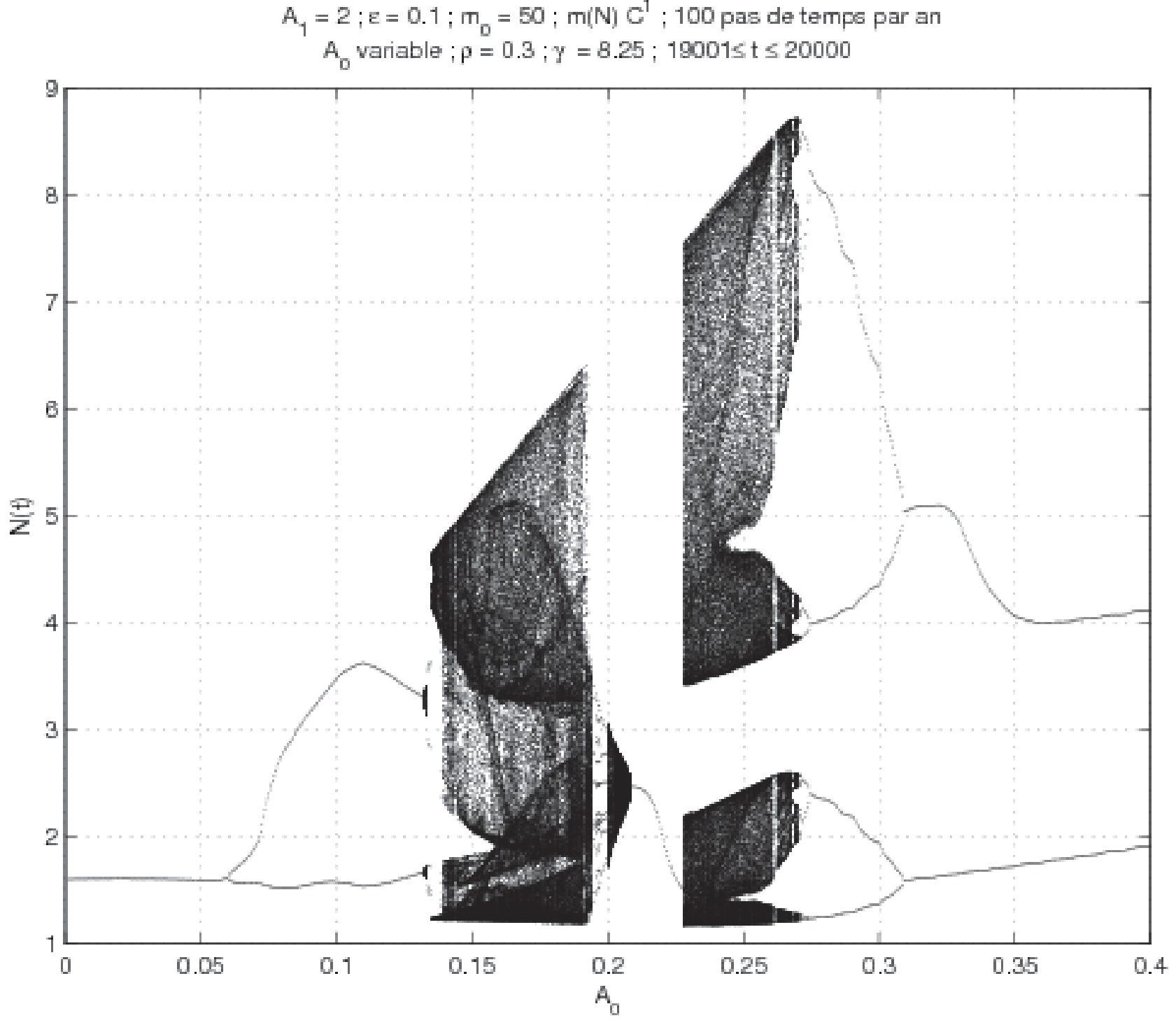}
\caption{\label{diag:A0_0.30_8.25}Diagramme de bifurcation $(A_0;
0\virg 30; 8\virg 25)_{0 \leq A_0 \leq 0\virg 4}$. Condition
initiale (II).}
\end{center}
\end{figure}

Pour cette derni\`ere exploration pr\'eliminaire, une difficult\'e
suppl\'ementaire a \'et\'e de trouver une fa\c{c}on de faire
varier la valeur de $A_0$ plus finement que le pas de
discr\'etisation. Celle-ci a conduit \`a utiliser la m\'ethode
\'evoqu\'ee \`a la fin du paragraphe
\ref{par:modele_discretisation}. Les r\'esultats obtenus sont
repr\'esent\'es dans le diagramme de bifurcations de la
figure~\ref{diag:A0_0.30_8.25}.

\paragraph{D\'ecomposition spectrale}

\begin{figure}
\begin{center}
\includegraphics[width=\textwidth]{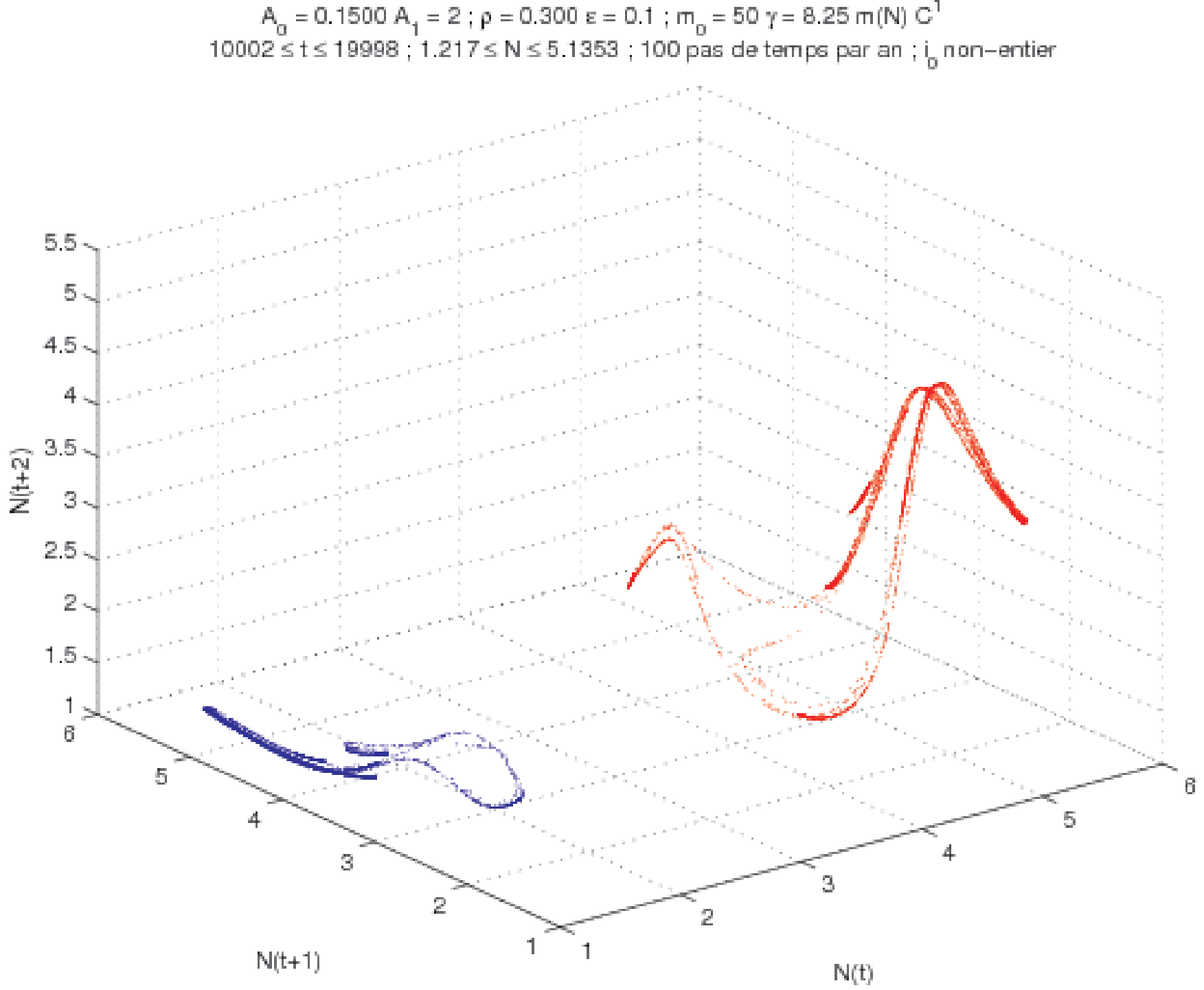}
\caption{\label{fig:0.15_0.30_8.25}Un comportement complexe :
$(0\virg 15; 0\virg 30; 8\virg 25)$. Condition initiale (II). Les
deux composantes sont invariantes par $T^2$, qui semble
topologiquement m\'elangeante sur chacune.}
\end{center}
\end{figure}

L'utilisation de 2 couleurs de visualisation montre que la
2-p\'eriodicit\'e de la figure~\ref{fig:0.15_0.30_8.25} ne se
retrouve pas \`a la figure~\ref{fig:0.18_0.30_8.25}. Il y ainsi
initialement deux composantes connexes bien distinctes, l'une
\'etant l'image de l'autre par l'application $T^1$. Sur chaque
composante, $T^2$ semble m\'elangeante. On a ainsi une
d\'ecomposition spectrale avec 2 composantes (voir th\'eor\`eme
\ref{the:decompo_spectrale}). Lorsque celles-ci fusionnent, on
perd cette 2-p\'eriodicit\'e et $T^1$ devient topologiquement
m\'elangeante. On retrouve le m\^eme comportement que
pr\'ec\'edemment avec le diagramme $(0\virg 18; 0\virg 41;
\gamma)$.

\subsection{\'Etude du cas $(0\virg 15; 0\virg 30; 8\virg 25)$}

\begin{figure}
\begin{center}
\includegraphics[width=\textwidth]{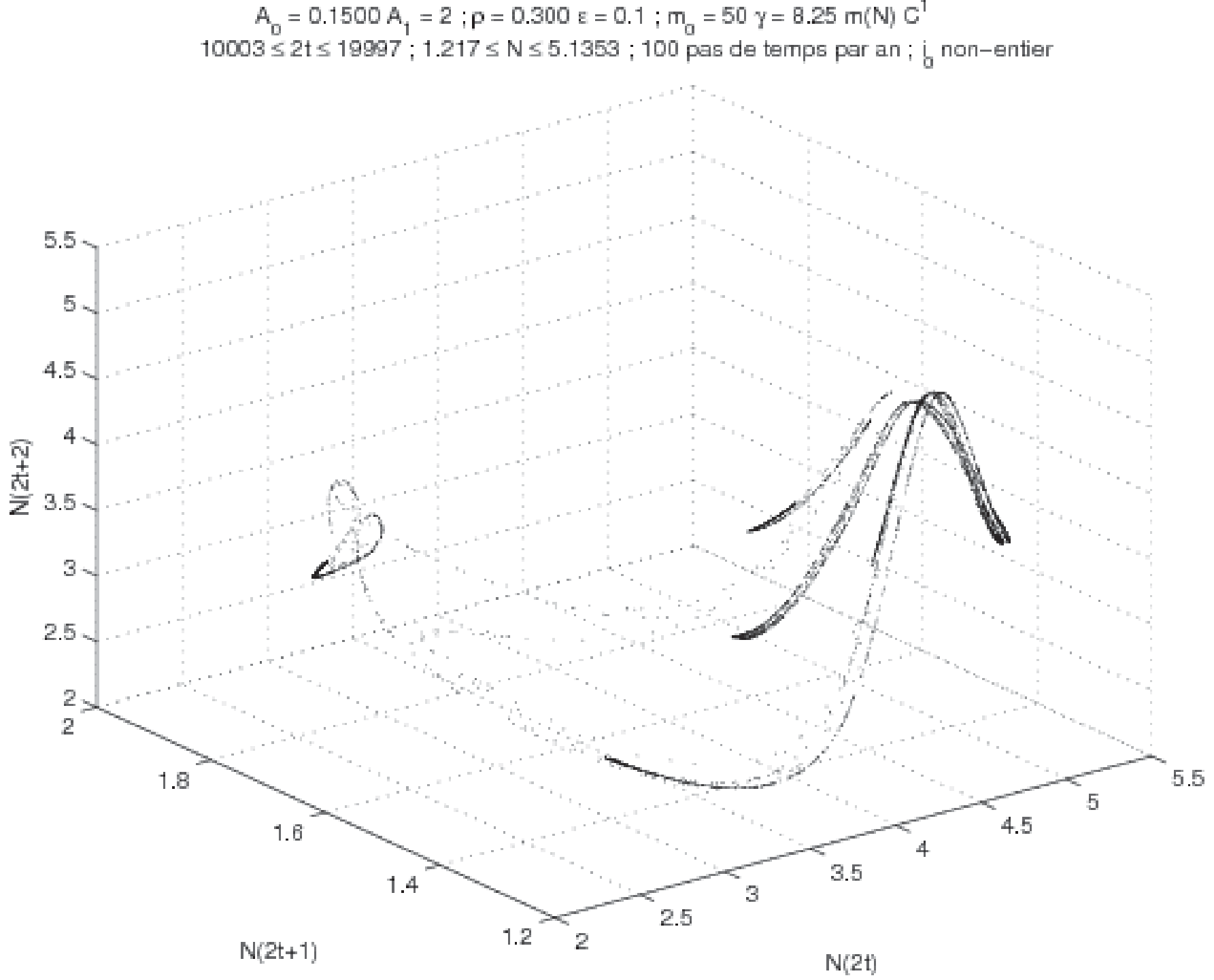}
\caption{\label{fig:0.15_0.30_8.25b} $(0\virg 15; 0\virg 30;
8\virg 25)$ Condition initiale (II). On n'a conserv\'e qu'une
composante connexe.}
\end{center}
\end{figure}

Essayons de comprendre la dynamique de la
figure~\ref{fig:0.15_0.30_8.25}. Nous avons vu qu'il y a deux
composantes connexes distinctes, il suffit donc de consid\'erer
l'une des deux pour comprendre la dynamique de $T^1$. Elle est
repr\'esent\'ee figure~\ref{fig:0.15_0.30_8.25b}.

\subsubsection{Visualisation en dimension 3}

La repr\'esentation que nous avons choisie (arbitrairement)
est-elle correcte ? Cette question est fondamentale : nous
projetons en effet un objet de dimension infinie\footnote{apr\`es
discr\'etisation, on se limite \`a une dimension finie tr\`es
grande, ici 201.} dans un espace de dimension 3.

\paragraph{Injectivit\'e de la projection}  \label{sec:injectivite} Pour tenter d'y
r\'epondre, nous pouvons \'evaluer la qualit\'e de la <<projection>>
$\pi$ : $\R^{201}$ $\rightarrow$ $\R^3$,
$x_{201}(t)=(N(t+k/100))_{k=0\ldots 200}$ $\mapsto$ ($N(t)$,
$N(t+1)$, $N(t+2)$)$=x_3(t)$ o\`u $N(t)$ d\'esigne la population
mature \`a l'instant $t$, en diff\'erents points de l'attracteur.
Nous voulons nous assurer que des points proches dans $\R^3$ sont
\'egalement proches dans $\R^{201}$, c'est-\`a-dire majorer
$\sup_{t \neq t^{\prime} \in \N}
\frac{\norm{x_{201}(t)-x_{201}(t^{\prime})}_{\R^{201}}}
{\norm{x_{3}(t)-x_{3}(t^{\prime})}_{\R^3}}$ pour diff\'erents
choix de normes ($L^1$, $L^2$ ou $L^{\infty}$).

\begin{figure}
\begin{center}
     \includegraphics[height=8cm]{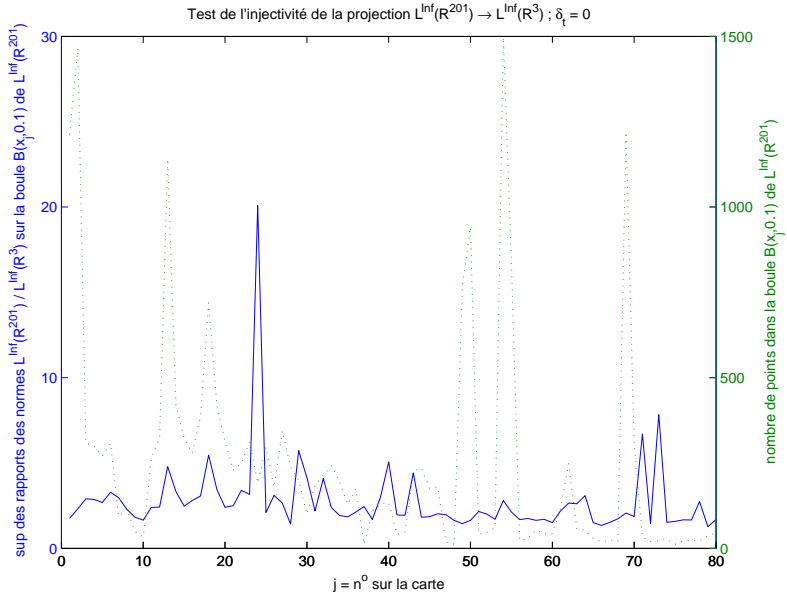}
\caption{\label{fig:injectivite_LInf_0} Injectivit\'e de la
projection : norme $L^{\infty}$, $\delta_t = 0$.}
\end{center}
\end{figure}

Le r\'esultat, repr\'esent\'e \`a la
figure~\ref{fig:injectivite_LInf_0}, montre que cette quantit\'e
est raisonnablement born\'ee. Avec les normes $L^1$ ou $L^2$, le
r\'esultat semble un petit peu meilleur, mais reste du m\^eme
ordre de grandeur. Une zone de l'attracteur semble en revanche
\^etre un peu moins bien repr\'esent\'ee par cette projection, il
s'agit du point 24 (et plus g\'en\'eralement des points 21 \`a
30). En se reportant \`a la figure~\ref{fig:80_points} o\`u sont
localis\'es ces points (voir section \ref{sec:dynamique80}), on
constate qu'il s'agit de la zone de pli. Une vue rapproch\'ee sur
cette zone de l'attracteur montre en effet des filaments
entrelac\'es, et certains rapprochements de filaments semblent
d\^us \`a la projection.

\begin{figure}
\begin{center}
     \includegraphics[height=7cm]{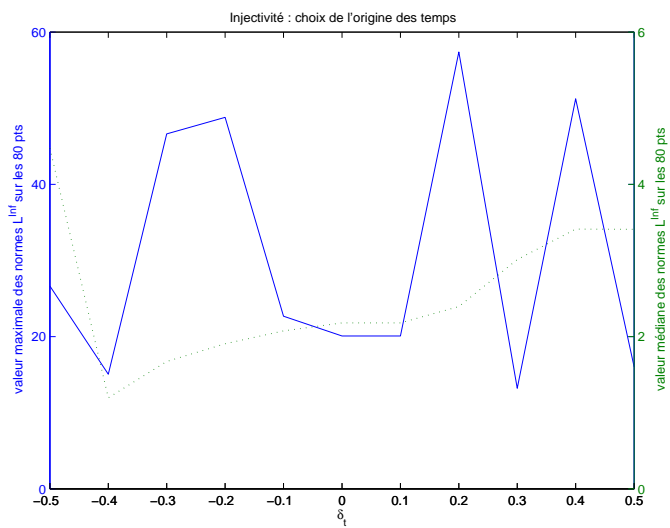}
\caption{\label{fig:injectivite_LInf_resume} Injectivit\'e de la
projection : norme $L^{\infty}$, $\delta_t$ variable.}
\end{center}
\end{figure}

\paragraph{Choix de l'origine des temps}\label{sec:origine_temps}
Nous l'avons arbitrairement fix\'ee \`a la fin de l'\'et\'e, mais
ce choix est-il judicieux ? Nous avons donc fait les m\^emes
calculs que pr\'ec\'edemment en d\'ecalant l'origine des temps. Il
semble que l'instant choisi initialement n'est pas mauvais. La
<<meilleure>> origine semble se situer autour de $\delta_t=0\virg
4$, mais la diff\'erence avec $\delta_t = 0$ n'est pas flagrante
(figure~\ref{fig:injectivite_LInf_resume}).

On peut expliquer ces r\'esultats en observant l'\'evolution en
temps continu de $N(t)$. En effet, la population mature atteint
tous les deux ans
--- un peu apr\`es le milieu de l'\'et\'e --- un maximum \'elev\'e,
suivi d'une chute brutale d'effectif. La valeur $\delta_t =
-0\virg 4$ correspond \`a l'instant du pic de population, qui est
suivi d'une simple diminution lin\'eaire de $N(t)$ (d\^ue \`a la
mortalit\'e naturelle, en l'absence de naissances), si bien que
les instants qui suivent sont encore des origines des temps de
bonne qualit\'e.

\begin{figure}
\begin{center}
     \includegraphics[width=\textwidth]{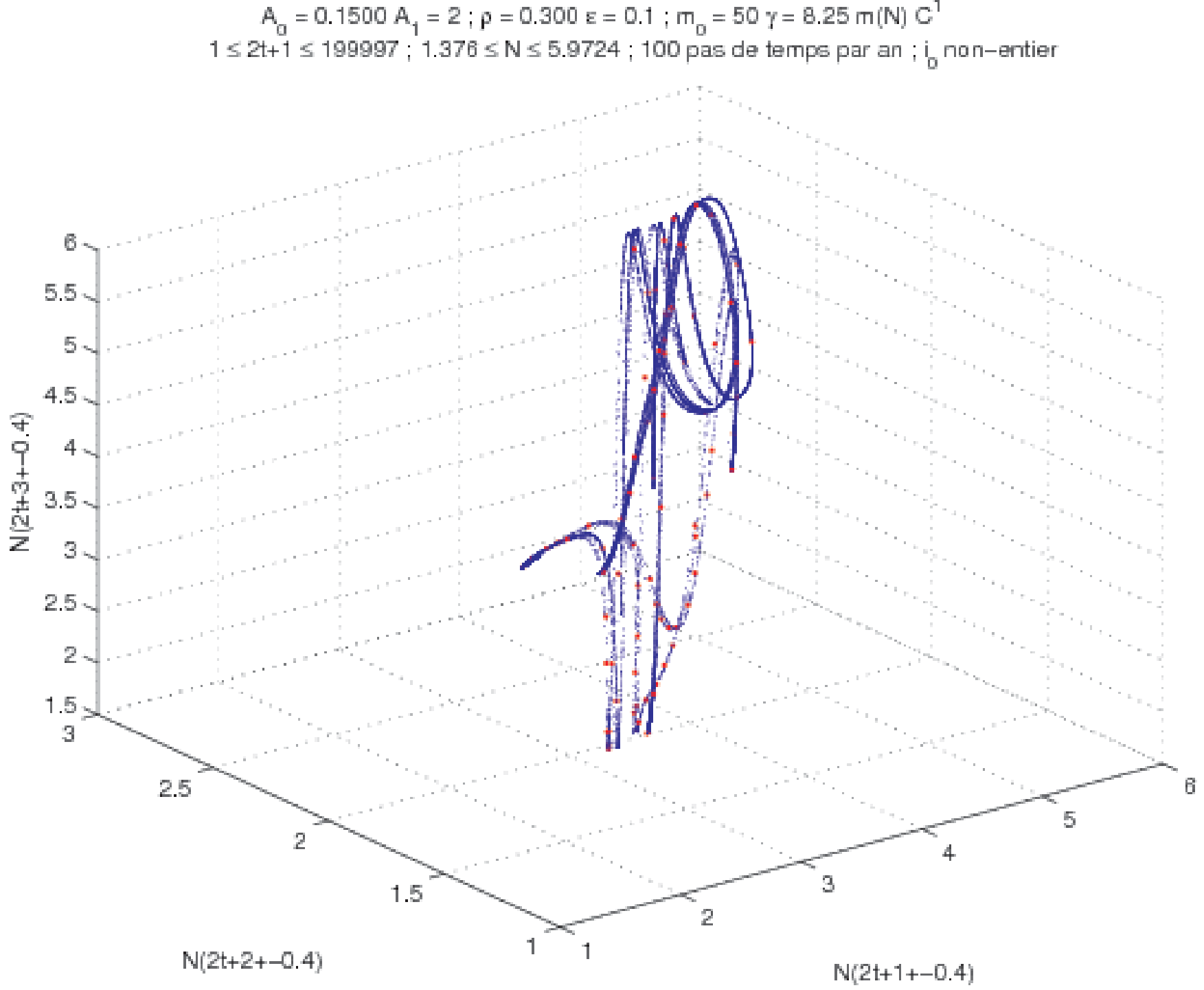}
\caption{\label{fig:delta_-0.4} Visualisation de l'attracteur
$(0\virg 15; 0\virg 30; 8\virg 25)$ avec une origine des temps
$\delta = -0\virg 4$.}
\end{center}
\end{figure}

Il est \'egalement int\'eressant, en vue de comprendre la
dynamique en temps continu du syst\`eme, de visualiser
l'\'evolution de l'attracteur tridimensionnel\footnote{l'animation
\fichier{film\_delta.avi} permet une bonne compr\'ehension de la
fa\c{c}on dont l'attracteur se d\'eforme, pour passer d'une
composante \`a l'autre quand $\delta_t$ varie de $-1$ \`a 0 ou de
0 \`a 1.} lorsque l'on fait varier l'origine des temps $\delta_t$
dans l'intervalle $[-1;1]$, la valeur $0$ correspondant \`a la fin
de l'\'et\'e. La figure~\ref{fig:delta_-0.4} en donne un exemple,
pour $\delta_t = -0\virg 4$.

\paragraph{\'Echelle logarithmique}
Une autre piste possible est de visualiser la projection
tridimensionnelle de l'attracteur suivant une \'echelle
logarithmique, c'est-\`a-dire de consid\'erer les points $(\log
N(t), \log N(t+1), \log N(t+2))$ pour $t$ entier grand.
\begin{figure}
\begin{center}
     \includegraphics[width=\textwidth]{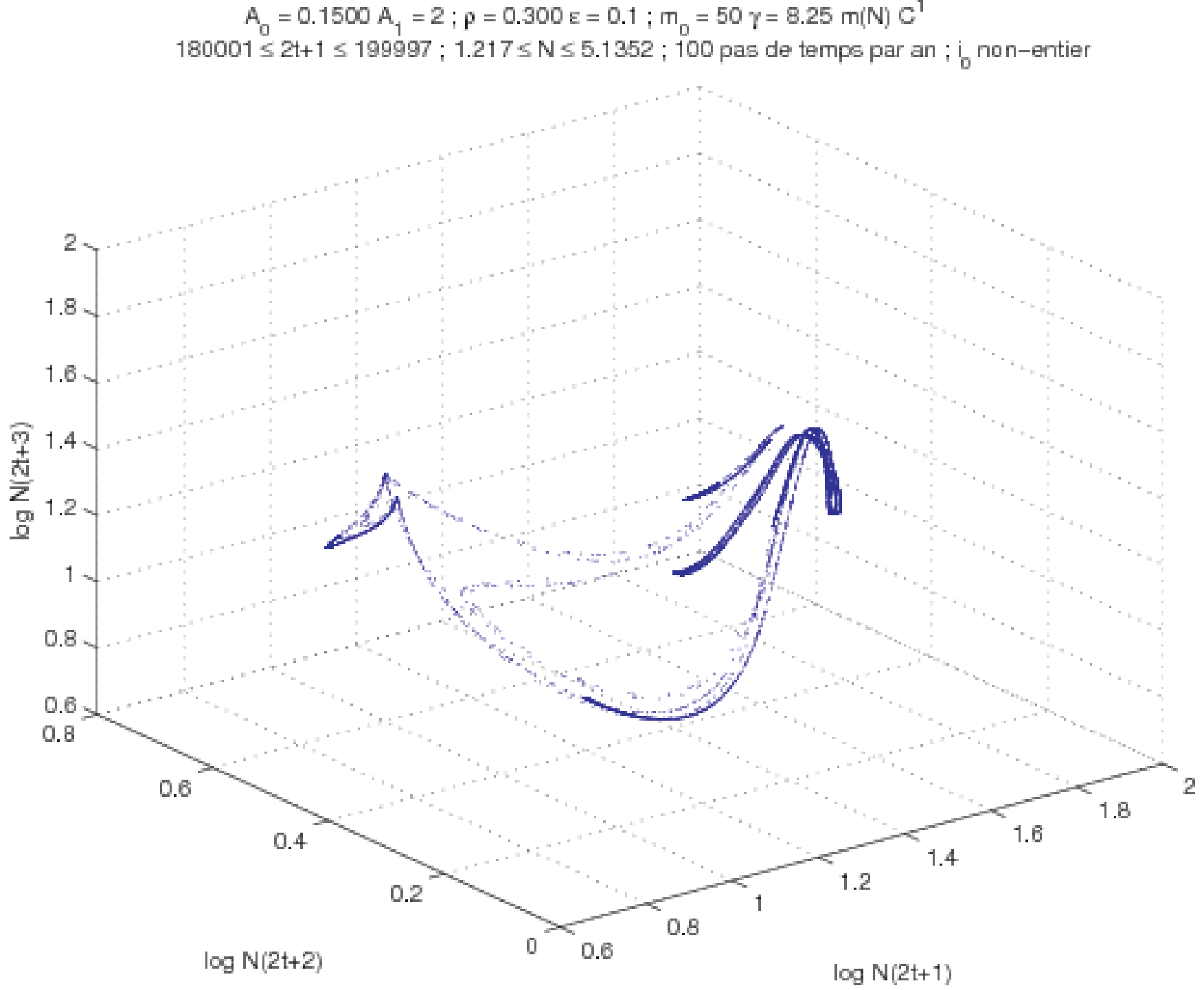}
\caption{\label{fig:log} Visualisation de l'attracteur $(0\virg
15; 0\virg 30; 8\virg 25)$ avec une \'echelle logarithmique.}
\end{center}
\end{figure}

\begin{figure}
\begin{center}
     \includegraphics[height=7cm]{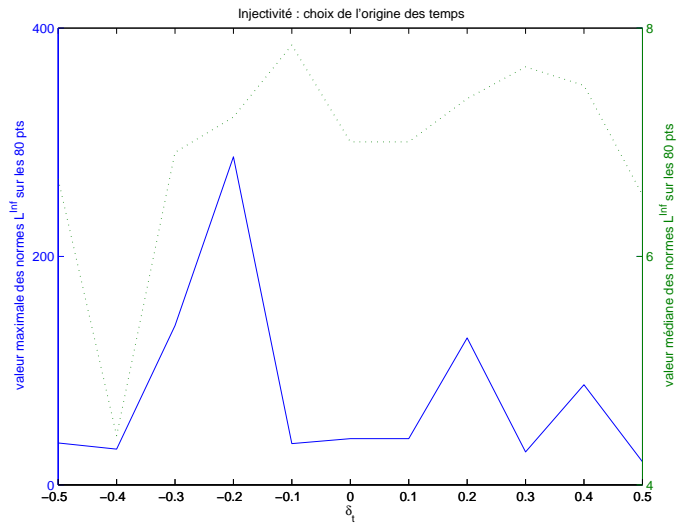}
\caption{\label{fig:injectivite_log} Qualit\'e de la projection en
\'echelle logarithmique.}
\end{center}
\end{figure}

Au vu de la figure~\ref{fig:log}, l'apport d'un tel changement
d'\'echelle n'est pas \'evident. L'attracteur est l\'eg\`erement
d\'eform\'e, mais garde le m\^eme aspect, et certaines zones
semblent toujours aussi <<emm\^el\'ees>>. Une \'evaluation
quantitative de la qualit\'e de cette nouvelle projection, comme
effectu\'e pr\'ec\'edemment, confirme l'aspect visuel : il n'y a
pas de gain significatif.

\subsubsection{G\'eom\'etrie de l'attracteur}
\begin{figure}
\begin{center}
     \includegraphics[width=\textwidth]{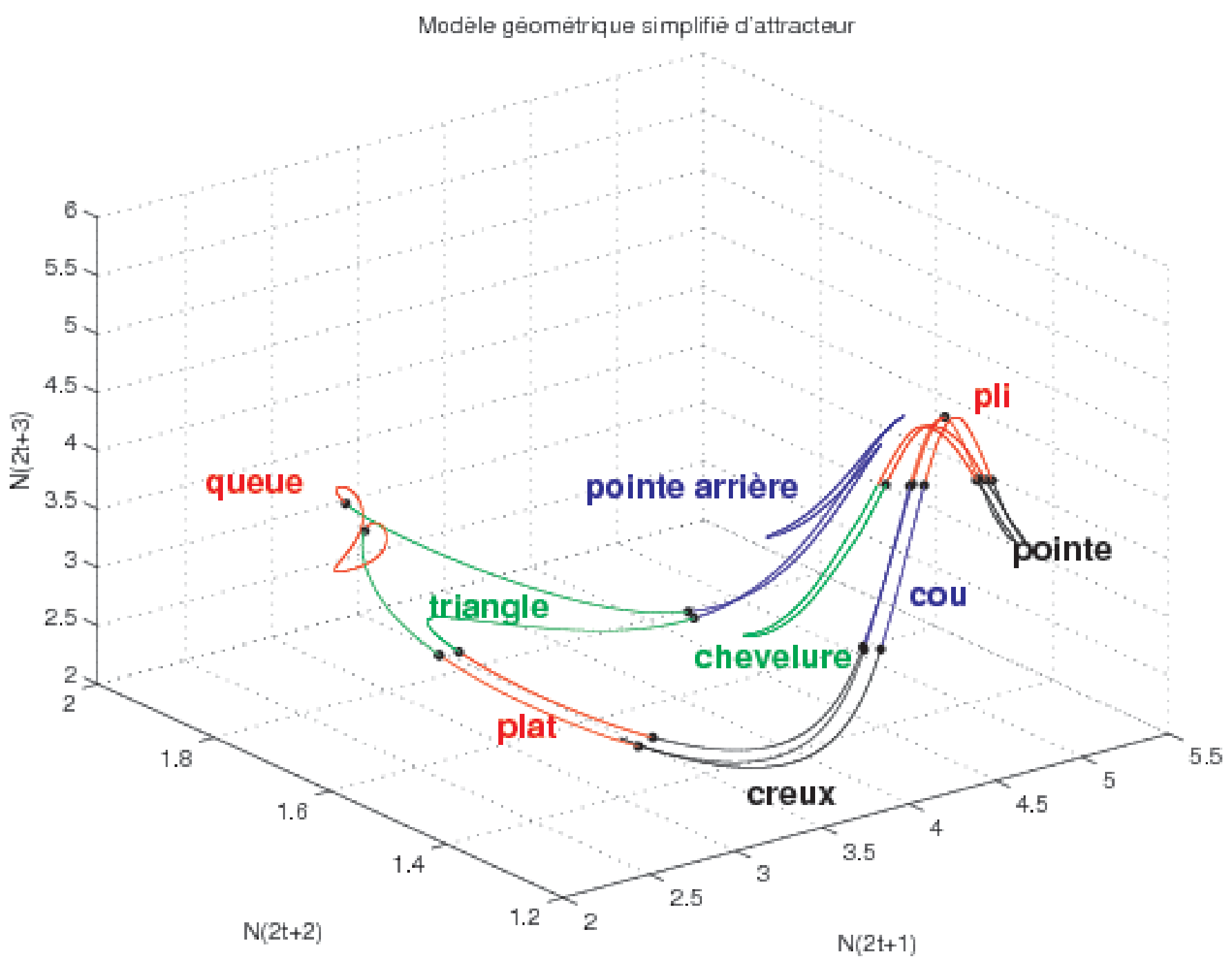}
\caption{\label{fig:geometrie} Repr\'esentation grossi\`ere de la
g\'eom\'etrie de l'attracteur $(0\virg 15; 0\virg 30; 8\virg
25)$.}
\end{center}
\end{figure}

La figure~\ref{fig:geometrie} repr\'esente de fa\c{c}on
simplifi\'ee la g\'eom\'etrie de l'attracteur de la
figure~\ref{fig:0.15_0.30_8.25b}, en distinguant neuf r\'egions
principales. Celles-ci sont nomm\'ees en fonction de leur forme et
de leur position dans l'attracteur, identifi\'e au corps d'un
animal dont la t\^ete serait situ\'ee \`a droite
(chevelure-cou-pli-pointe) et la queue \`a gauche.

Quatre ph\'enom\`enes principaux semblent pouvoir expliquer la
complexit\'e de l'attracteur que nous observons. Les deux premiers
\'etaient d\'ej\`a pr\'esents dans le sol\'eno\"ide : un tr\`es
fort pincement et un \'etirement. Le troisi\`eme est \'egalement
pr\'esent dans l'attracteur de H\'enon : un pli (il y en a
peut-\^etre plusieurs ici). Le quatri\`eme semble nouveau, et
ressemble \`a un ou plusieurs <<embranchements>>. Avec la
num\'erotation introduite dans la section pr\'ec\'edente, on peut
en situer trois : 67--80, 63--55, 58--74. Cela ne signifie pas
pour autant que ces embranchements sont distincts.

\begin{figure}
\begin{center}
     \includegraphics[width=8cm]{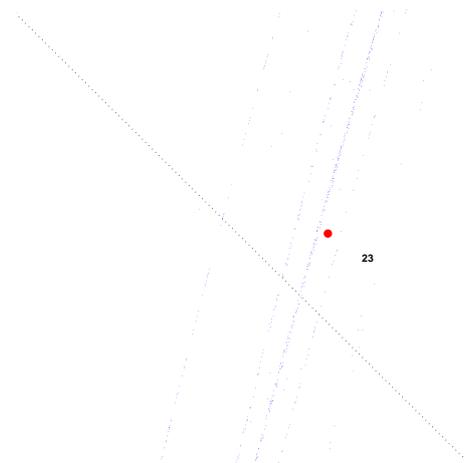}
\caption{\label{fig:zoom23} Zoom sur un filament, au voisinage du
point 23.}
\end{center}
\end{figure}

\begin{figure}
\begin{center}
\begin{tabular}{c@{}c}
     \includegraphics[width=7cm]{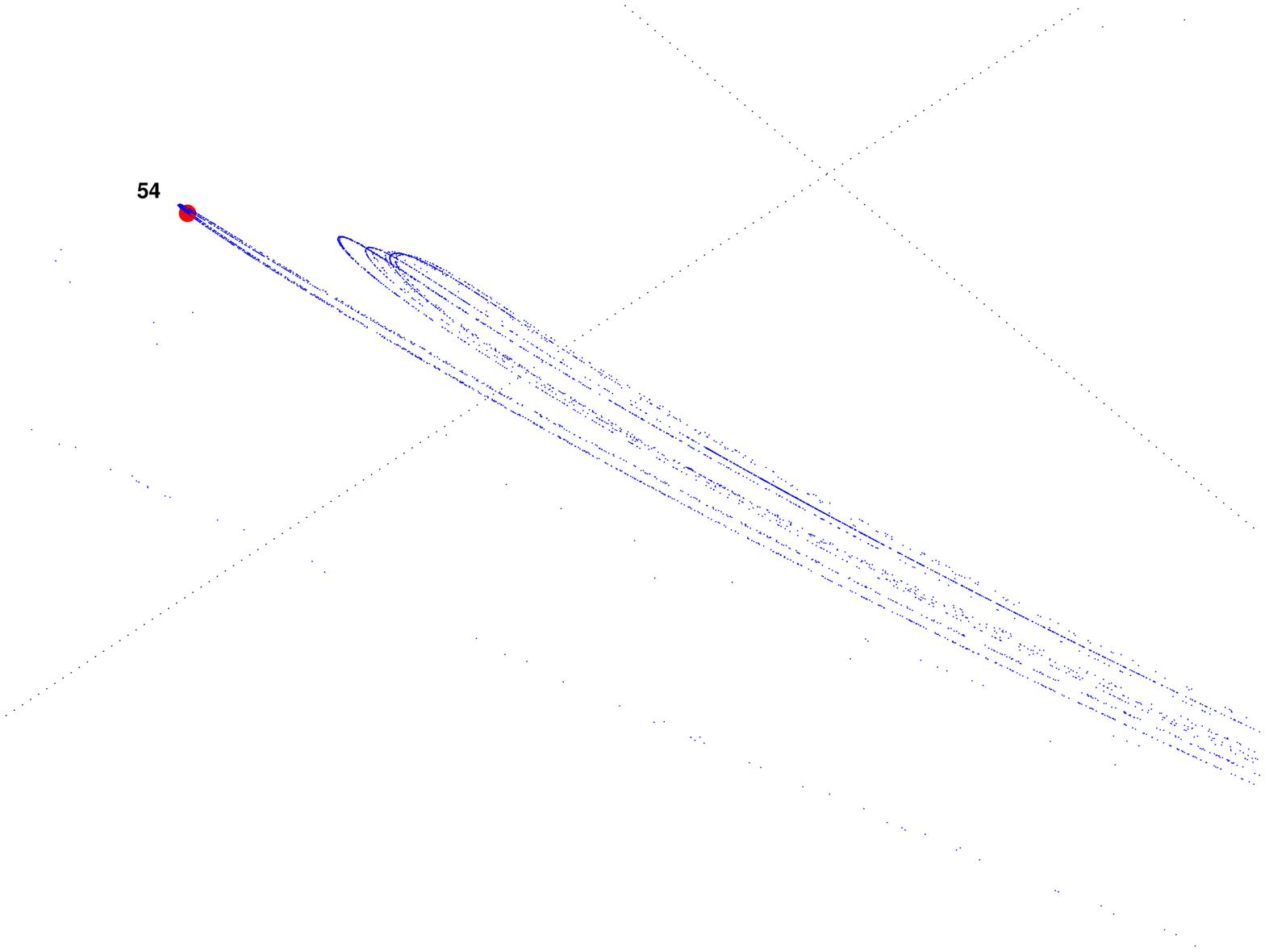}
     &
     \includegraphics[width=7cm]{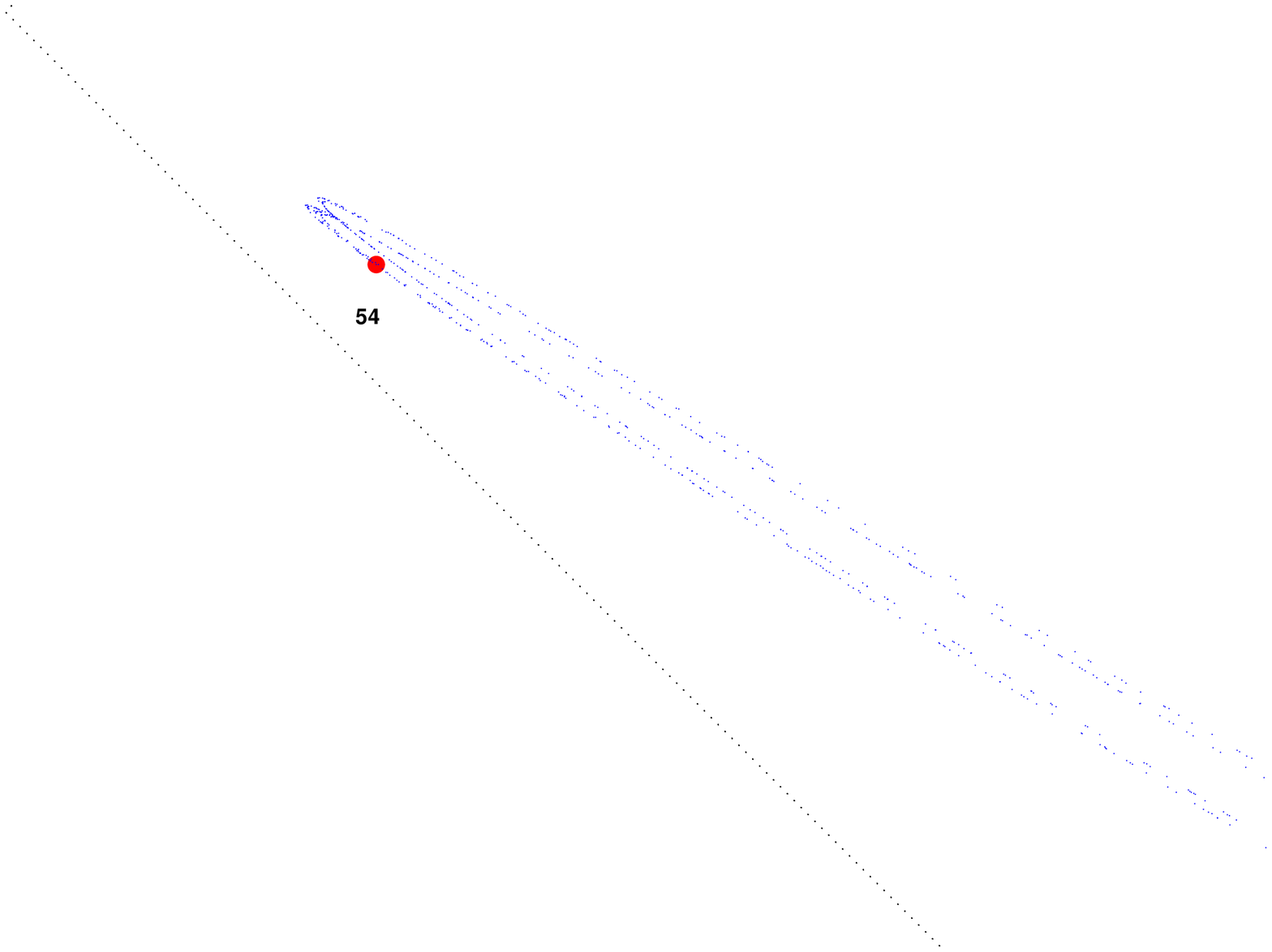}
\end{tabular}
\caption{\label{fig:zoom54} Zooms successifs sur une pointe, au
voisinage du point 54.}
\end{center}
\end{figure}

Localement, l'attracteur ressemble au produit d'une droite et d'un
ensemble de Cantor (figure~\ref{fig:zoom23}), sauf en certains
points o\`u l'on observe des <<pointes>>
(figure~\ref{fig:zoom54}). Par ces aspects, il ressemble beaucoup
\`a l'attracteur de H\'enon\footnote{voir
annexe~\ref{annexe:Henon}.}.

\subsubsection{Dynamique sur l'attracteur} \label{sec:dynamique80}
\begin{figure}
\begin{center}
\includegraphics[height=15cm,angle=90]{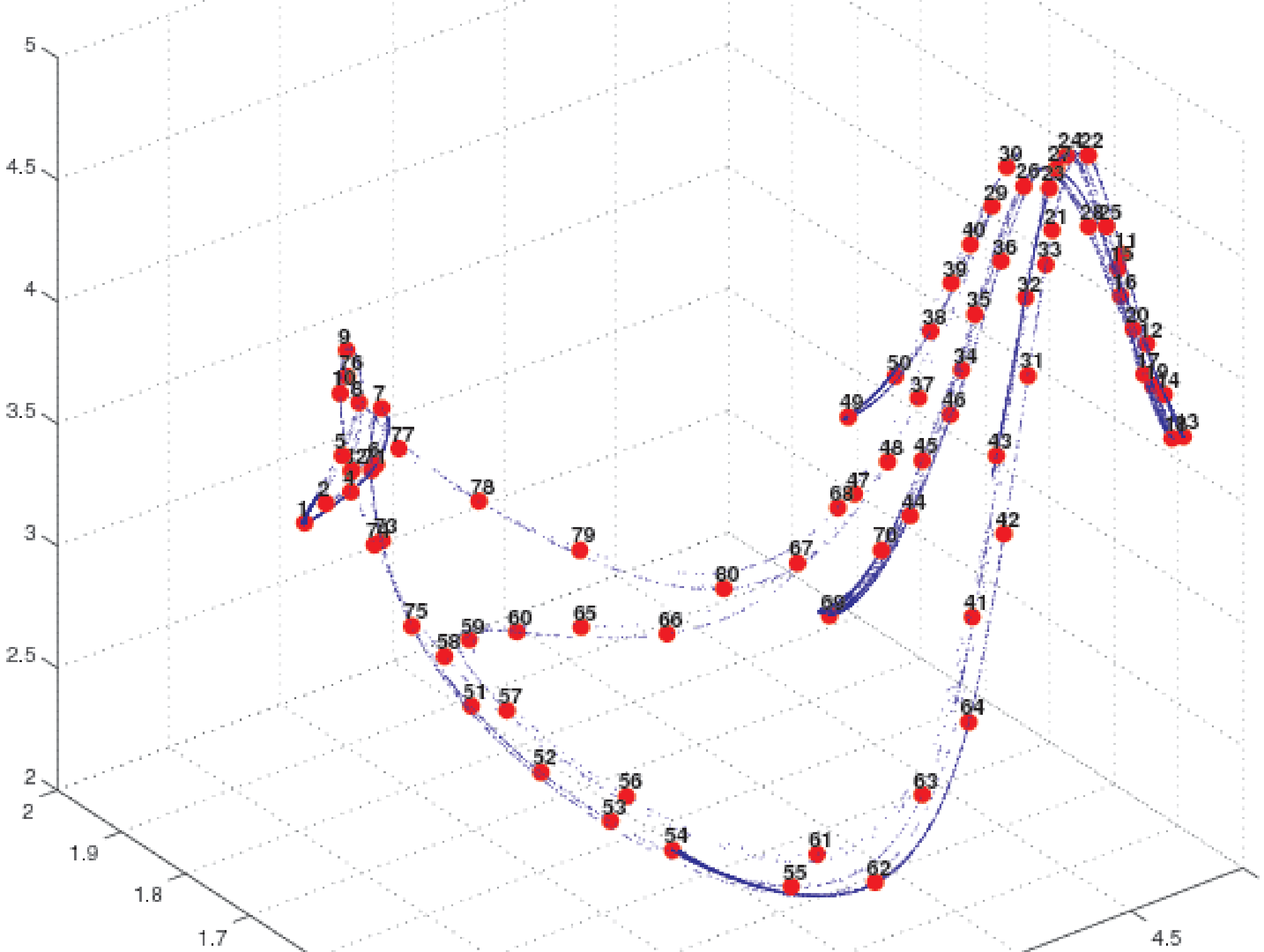}
\caption{\label{fig:80_points} Position des 80 points choisis sur
l'attracteur.}
\end{center}
\end{figure}

Consid\'erons l'application $T^2$. Comment agit-elle sur les
points de l'attracteur ? Pour essayer de le comprendre, 80 points
ont \'et\'e choisis\footnote{Ce choix a \'et\'e fait
arbitrairement, en essayant de r\'epartir ces points
uniform\'ement suivant la mesure de Hausdorff sur l'attracteur, et
non la mesure physique.} sur l'attracteur, num\'erot\'es de 1 \`a
80, comme repr\'esent\'e sur la figure~\ref{fig:80_points} (des
vues plus rapproch\'ees sont en annexe~\ref{annexe:details}). La
position des images directes et r\'eciproques de ces 80 points est
indiqu\'ee dans le tableau \ref{tab:dynamique80}.

\begin{table}
\begin{center}
\begin{tabular}{|r@{ $\rightarrow$ }c@{ $\rightarrow$ }l|r@{ $\rightarrow$ }c@{ $\rightarrow$ }l|r@{ $\rightarrow$ }c@{ $\rightarrow$ }l|r@{ $\rightarrow$ }c@{ $\rightarrow$ }l}
\hline
54--55 & 1 & 49 &35 & 11 & 16--20 \\
54--55 & 2 & 50 &36 & 12 & 26 \\
56--61 & 3 & 39 &26--27 ($9^{3+}$) & 13 & 69 \\
55--62 & 4 & 38 &27--15 & 14 & 46--34 \\
53 & 5 & 39--40 &39 & 15 & 28--16 \\
62 & 6 & 39--40 &32--23 & 16 & 28--27 \\
62--64 & 7 & 39--40 &21 & 17 & 35 \\
52 & 8 & 30 &24 & 18 & 69--69 \\
51--75 & 9 & 40 &29 & 19 & 46--34 \\
59--60 & 10 & 68 &40 & 20 & 26--27 \\
\hline
\end{tabular}
\end{center}

\begin{center}
\begin{tabular}{|r@{ $\rightarrow$ }c@{ $\rightarrow$ }l|r@{ $\rightarrow$ }c@{ $\rightarrow$ }l|r@{ $\rightarrow$ }c@{ $\rightarrow$ }l|r@{ $\rightarrow$ }c@{ $\rightarrow$ }l}
\hline
46 & 21 & 17 &45--46 & 31 & 24 &43--43 & 41 & 53--55 \\
34 & 22 & 18--13 &37 & 32 & 16--28 &45 & 42 & 62--55 \\
50 & 23 & 18--17 &46--45 & 33 & 15--17 &49 & 43 & 41--42 ($1^{3+}$) \\
31 & 24 & 18 &12--13 & 34 & 22 &13--14 & 44 & 62--64 \\
39--38 & 25 & 20--19 &17 & 35 & 11 &17--18 & 45 & 42 \\
12 & 26 & 13--12 &14--15 & 36 & 12 &19--18 & 46 & 21 \\
20--20 & 27 & 13--13 &07--71 & 37 & 32 &76--76 & 47 & 61--63 \\
15--12 & 28 & 14--15 &04 & 38 & 24--25 &07--06 & 48 & 63--41 \\
07--06 & 29 & 19 &03 & 39 & 15 &01 & 49 & 43 \\
08 & 30 & 18--19 &09 & 40 & 20 &02 & 50 & 23 \\
\hline
\end{tabular}
\end{center}

\begin{center}
\begin{tabular}{|r@{ $\rightarrow$ }c@{ $\rightarrow$ }l|r@{ $\rightarrow$ }c@{ $\rightarrow$ }l|r@{ $\rightarrow$ }c@{ $\rightarrow$ }l|r@{ $\rightarrow$ }c@{ $\rightarrow$ }l}
\hline
41--63 & 51 & 09--08 &68 & 61 & 03--07 &63--64 & 71 & 80 \\
41--63 & 52 & 08 &70--45 & 62 & 06 &65 & 72 & 79 \\
41--63 & 53 & 05 &43--41 & 63 & 07--03 &63--41 & 73 & 78 \\
69 & 54 & 01--05 &44--46 & 64 & 75 &66 & 74 & 77 \\
41--42 ($1^{3+}$) & 55 & 04--02 &77--78 & 65 & 72 &64 & 75 & 76 \\
67--80 & 56 & 08--03 &77--76 & 66 & 74 &75 & 76 & 47--67 \\
80 & 57 & 09--08 &71--07 & 67 & 56--61 &74 & 77 & 65--66 \\
66--67 & 58 & 09--76 &10 & 68 & 61 &73 & 78 & 60 \\
79 & 59 & 10--09 &13 & 69 & 54 &72 & 79 & 59 \\
78 & 60 & 10--72 &18--17 ($2^{3+}$) & 70 & 62--54 &71 & 80 & 57 \\
\hline

\end{tabular}
\end{center}
\caption{\label{tab:dynamique80} Dynamique sous $T^2$ des 80
points sur l'attracteur. On a not\'e $1^{3+}$ pour $T^6(1)$ ;
54--55 pour un point situ\'e entre 54 et 55, plus proche de 54 ;
76--76 pour un point proche de 76 mais diff\'erent de 76.}
\end{table}

On peut penser qu'il y a un point fixe dans la r\'egion du point
15, une orbite de p\'eriode 2 dans la r\'egion 42--31 et une
orbite de p\'eriode 3 dans le <<triangle>> 66--74--77. Nous
\'etudierons le point fixe plus en d\'etails dans une des sections
suivantes. Partant de ces trois orbites remarquables, on peut
essayer de comprendre sch\'ematiquement la dynamique interne \`a
l'attracteur.

\paragraph{Sortie de la <<t\^ete>>}
L'\'equilibre (entre pli et pointe) est r\'epulsif avec expansion
n\'egative (la valeur propre dominante est n\'egative), tandis que
l'orbite de p\'eriode 2 (dans le cou, entre les points 42 et 31)
est r\'epulsive positive. Dans cette r\'egion de la t\^ete et de
la pointe arri\`ere, on peut sch\'ematiser la dynamique de la
fa\c{c}on suivante :
\begin{itemize}
\item Pointe arri\`ere $\rightarrow$ creux-cou-pli-pointe (sans inversion haut-bas).
\item Chevelure $\rightarrow$ cou-pli-pointe (sans inversion).
\item Bas du cou (sous l'orbite de p\'eriode 2 = 42--31) $\rightarrow$ plat-creux-bas du cou (sans inversion).
\item Haut du cou (au-dessus de la p\'eriode 2) $\rightarrow$ haut du cou-pli (sans inversion).
\item Pli $\rightarrow$ bas de la pointe (avec inversion).
\item Haut de la pointe $\rightarrow$ bas du pli (avec inversion).
\end{itemize}

Ainsi, si l'on part d'assez haut, on arrive en bas de la pointe
puis en bas de la chevelure (69) et enfin dans le creux (54). Si
l'on part plus bas, on arrive directement dans le creux, sans
\^etre pass\'e par la pointe.

\paragraph{Sortie du <<corps>>}
Dans le <<corps>> (queue, triangle, plat, creux), l'orbite de
p\'eriode 3 (aux milieux des sommets du triangle) joue un r\^ole
majeur. Elle est r\'epulsive positive, et agit sur le triangle
comme une rotation d'angle $2\pi/3$ (75--72  $\rightarrow$ 76--79
$\rightarrow$ 67--59) combin\'ee avec un peu d'expansion. On sort
ainsi du triangle soit par le bas de la pointe arri\`ere, soit par
le plat ou le creux, soit par la queue. La dynamique se
sch\'ematise alors ainsi :
\begin{itemize}
\item Bas de la pointe arri\`ere $\rightarrow$ plat-creux.
\item Creux $\rightarrow$ droite de la queue (1--7)-c\^ot\'e du
triangle (75--73) (avec inversion droite-gauche).
\item Plat $\rightarrow$ gauche de la queue (sans inversion haut-bas).
\item Queue $\rightarrow$ pointe arri\`ere.
\end{itemize}

On sort ainsi du corps, pour y revenir rapidement (si l'on arrive
trop bas dans la pointe arri\`ere, soit pr\`es du triangle, soit
\`a la pointe 49), ou (le plus souvent) apr\`es un passage dans le
pli et \'eventuellement la pointe.

Ce bref aper\c{c}u de la dynamique nous permet de comprendre
comment s'instaure le m\'elange et le chaos de la dynamique sur
l'attracteur.

\subsubsection{Visualisation en temps continu}
Il est int\'eressant, du point de vue math\'ematique comme du
point de vue biologique, de mettre en relation les diff\'erentes
r\'egions de l'attracteur (en dimension 3) avec la dynamique en
temps continu dans ces r\'egions. La forme de $N(t)$ pour $t\in
[t_0 - 5 ; t_0 +5]$ peut ainsi \^etre mise en correspondance avec
les points de l'attracteur et leurs r\'egions de provenance et de
destination par $T^2$. On a repr\'esent\'e ces donn\'ees pour les
80 points de la <<carte>> en annexe~\ref{annexe:details}.

\paragraph{Aspect g\'en\'eral}
La dynamique en temps continu garde en g\'en\'eral un aspect
identique dans toutes les zones de l'attracteur : un pic chaque
ann\'ee pour $t\approx 0\virg 6 \modulo 1$ (croissance tr\`es
rapide suivie d'une d\'ecroissance lin\'eaire moins brusque),
\'elev\'e les ann\'ees impaires (entre $2\virg 5$ et $6\virg 5$),
plus faibles voire inexistant les ann\'ees paires (entre 1 et 2).
Juste avant les maxima se trouvent des minima locaux plus ou moins
bas : ceux-ci sont toujours l\'eg\`erement sup\'erieurs \`a 1
(mais inf\'erieurs \`a $1\virg 5$) avant un maximum faible, mais
pouvant atteindre $0 \virg 5$ avant un maximum \'elev\'e. Il y a
ainsi une pseudo-p\'eriodicit\'e de 2 ans, avec une tr\`es forte
variation d'amplitude (les maxima \'etant \`a $6\virg 5$, les
minima entre $0\virg 5$ et 1). On retrouve en partie l'aspect
g\'en\'eral de la figure~\ref{fig:popu_microtus_epiroticus}
(alternance de maxima et de minima, avec un facteur allant
jusqu'\`a 8 entre les deux).

Ce comportement s'explique par l'alternance entre une explosion de
la population d\^ue \`a la tr\`es forte f\'econdit\'e, qui est
ainsi suivie d'une chute lin\'eaire de la population mature (d\^ue
\`a la mortalit\'e naturelle, en l'absence quasi-totale de
naissances).

\paragraph{Diff\'erences entre r\'egions}
Les r\'egions de l'attracteur se diff\'erencient par l'amplitude
des pics (il y a un facteur 2 entre les amplitudes possibles des
maxima \'elev\'es) et des creux (inf\'erieurs ou sup\'erieurs \`a
1), ainsi que par l'amplitude relative des maxima secondaires
(inexistants ou bien valant jusqu'\`a 2).

Les pics faibles se trouvent au niveau de la queue (lorsqu'il se
produit \`a $t=0$) ou du plat (lorsqu'il se produit \`a $t=2$).

Les pics tr\`es \'elev\'es ($N\geq 6$) se trouvent dans le pli, la
pointe, le haut de la pointe arri\`ere et le cou, avec des maxima
secondaires quasi-inexistants et des minima inf\'erieurs \`a 1. Le
deuxi\`eme pic est l\'eg\`erement inf\'erieur au premier dans la
pointe, mais le rapport s'\'equilibre quand on se rapproche du
pli.

On peut encore affiner cette analyse en s'aidant des figures
plac\'ees en annexe~\ref{annexe:details}, qui permettent de faire
la diff\'erence entre filaments dans une r\'egion \`a l'aide de
l'\'evolution en temps continu dans le pass\'e ou dans le futur
proche.

\subsubsection{Dimension fractale de l'attracteur}
\begin{figure}
\begin{center}
\includegraphics[width=\textwidth]{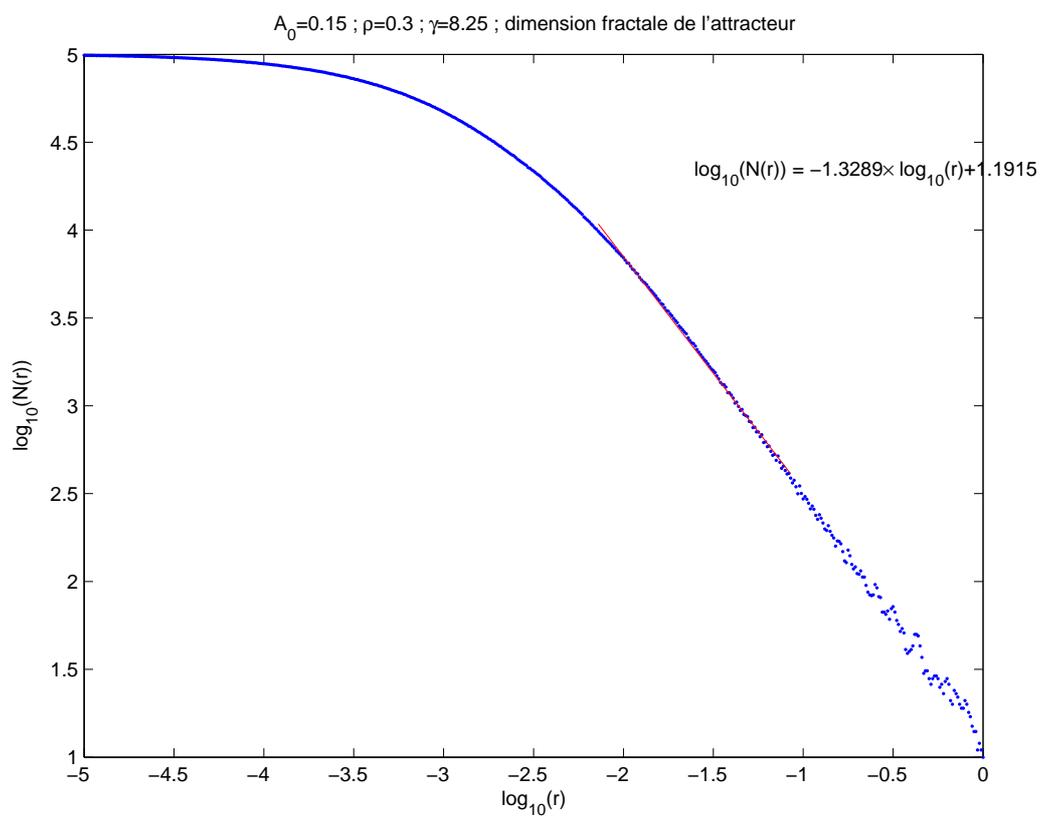}
\caption{\label{fig:dim_fract_0.15_0.30_8.25} \'Evaluation de la
dimension fractale de l'attracteur.}
\end{center}
\end{figure}

Pour l'\'evaluer, nous avons calcul\'e pour diff\'erentes valeurs
de $r$ le nombre de bo\^ites de c\^ot\'e $r$ (et dont les points
ont des coordonn\'ees qui sont des multiples entiers de $r$)
contenant des points de l'attracteur (avec %
\nombre{100000}
 points). En notant $N(r)$ ce nombre de points, on a alors trac\'e
$\log_{10}(N(r))$ en fonction de $\log_{10}(r)$. Tant que $N(r)$
est assez petit devant %
\nombre{100000}
 et assez grand devant 1, les points
obtenus sont presque align\'es, et la pente (en valeur absolue) de
la droite de r\'egression est une bonne estimation de la dimension
fractale de l'attracteur. Le r\'esultat est repr\'esent\'e \`a la
figure~\ref{fig:dim_fract_0.15_0.30_8.25}. On obtient donc une
dimension de l'ordre de $1\virg 33$. Cette valeur correspond bien
\`a l'impression visuelle que l'on a : localement, l'attracteur
semble \^etre le produit d'une droite et d'un ensemble de Cantor
de dimension proche de $1/3$, soit une dimension fractale
d'environ $4/3$ (voir figure \ref{fig:zoom23}). Il faut bien s\^ur
prendre ce r\'esultats avec beaucoup de pr\'ecautions, dans la
mesure o\`u cette r\'egression est faite dans la zone qui semble
--- visuellement --- pr\'esenter une <<bonne pente>> (voir annexe~\ref{annexe:dimfract}).

Une des cons\'equences de cette \'evaluation est la justification
\latin{a posteriori} de la possibilit\'e d'utiliser trois
dimensions seulement pour visualiser l'attracteur. En effet, le
th\'eor\`eme de Whitney \ref{the:whitney}\footnote{voir
annexe~\ref{annexe:whitney}.} affirme qu'un objet de dimension $d$
peut \^etre visualis\'e avec $N$ dimensions pourvu que $N>2d$.
L'attracteur ayant une dimension strictement comprise entre 1 et
$1\virg 5$, l'utilisation de 3 dimensions pour le visualiser
semble raisonnable\footnote{mais bien s\^ur, comme nous avons
calcul\'e la dimension fractale de la projection tridimensionnelle
de l'attracteur, nous n'avons pas la dimension fractale de
l'attracteur lui-m\^eme mais une l\'eg\`ere sous-estimation de
celle-ci.}.

\subsubsection{Sensibilit\'e aux conditions initiales}
Pour l'\'evaluer, il est int\'eressant de regarder la dynamique
future et pass\'ee d'une petite boule centr\'ee sur un point de
l'attracteur.

\begin{figure}
\includegraphics[width=\textwidth]{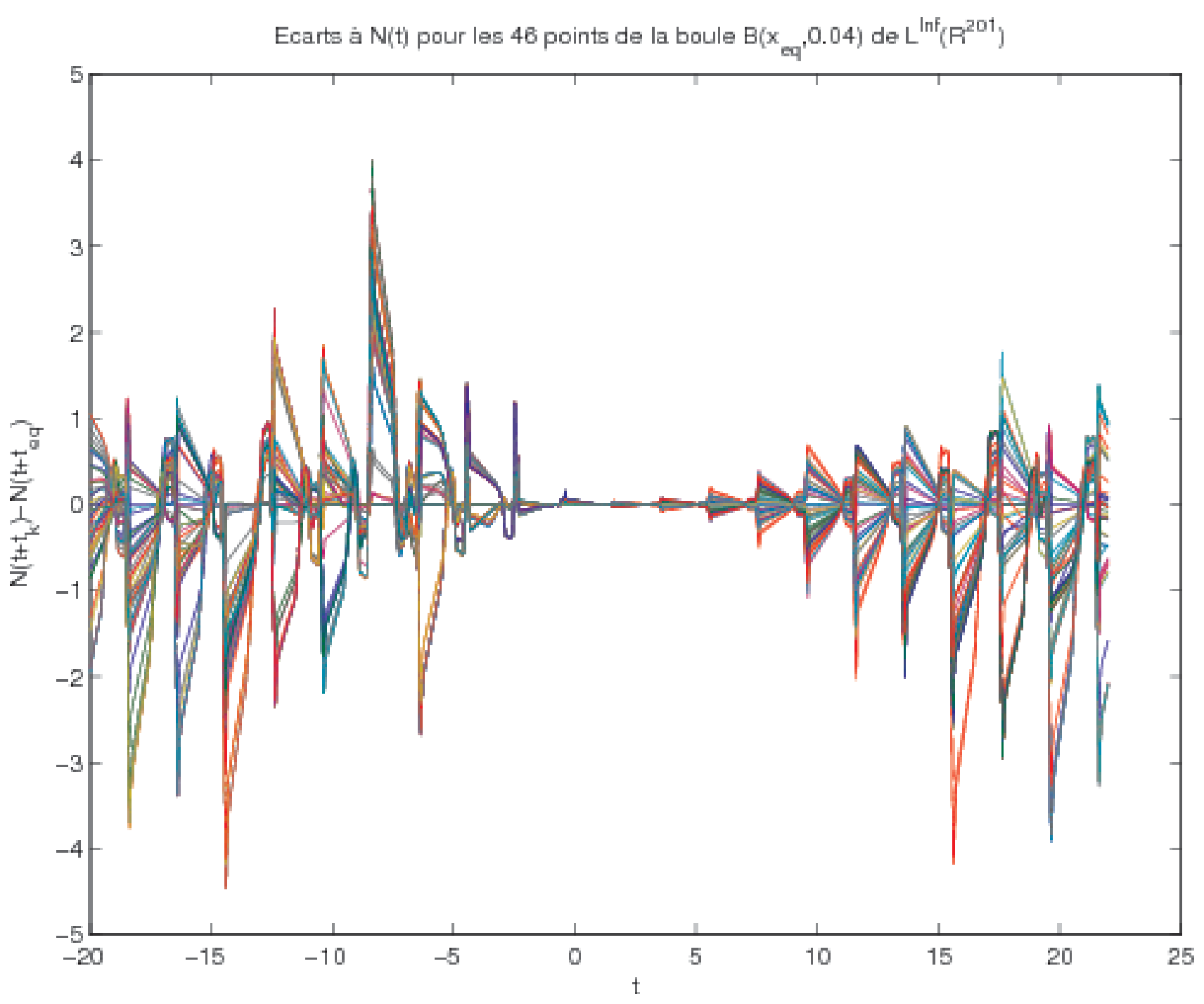}
\caption{\label{fig:passe_futur_eq_continu} Dynamique pass\'ee et
future de points au voisinage de l'\'equilibre. Les points ont
\'et\'e choisis pour leur proximit\'e de l'\'equilibre dans
l'intervalle de temps $[0;2]$.}
\end{figure}

\begin{figure}
\includegraphics[width=\textwidth]{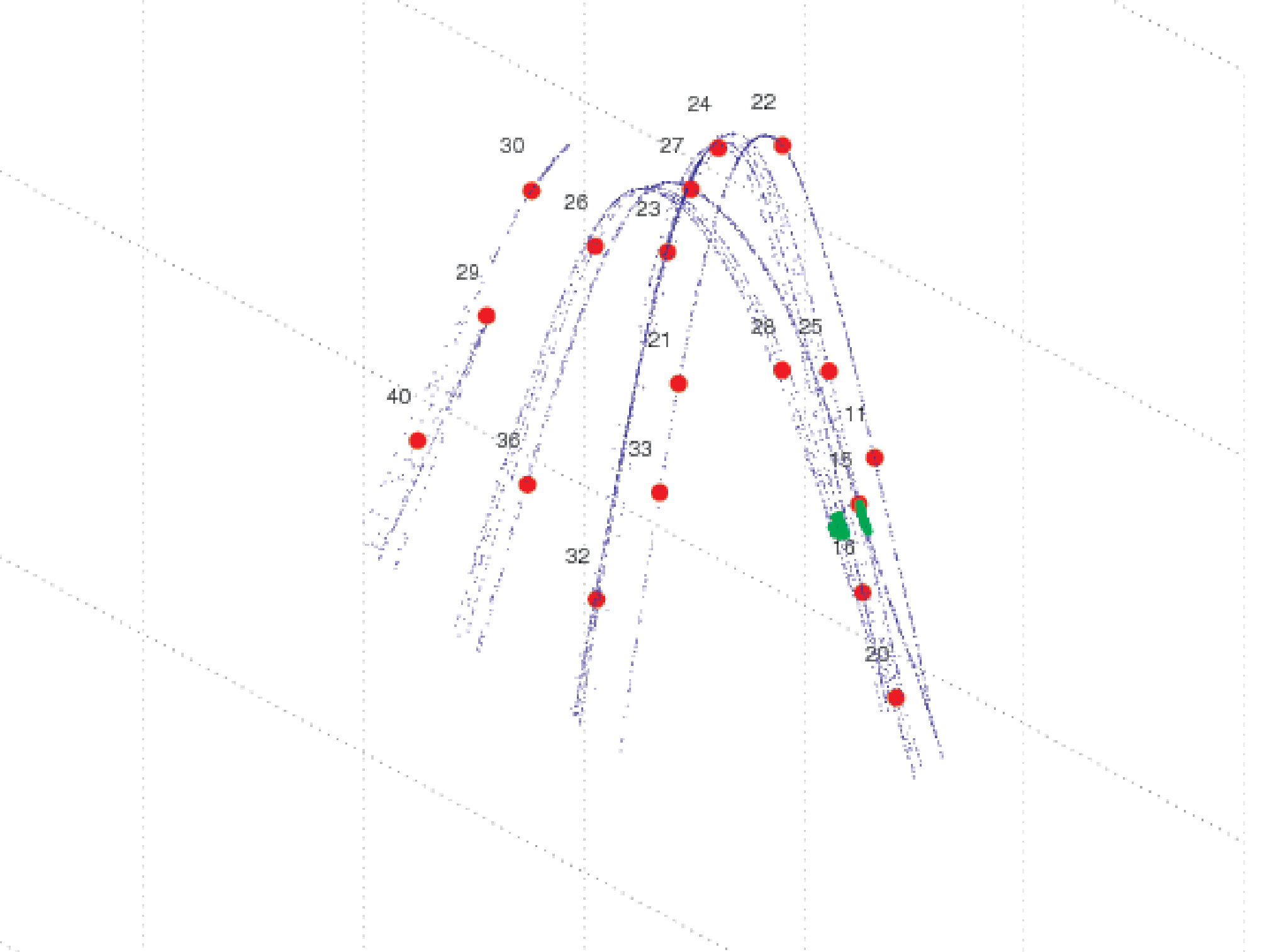}
\caption{\label{fig:passe_futur_eq_3d} Position des points au
voisinage de l'\'equilibre : l'ensemble des points de la boule
consid\'er\'ee est repr\'esent\'e en vert.}
\end{figure}

Les figures \ref{fig:passe_futur_eq_continu} et
\ref{fig:passe_futur_eq_3d} permettent de juger du r\'esultat au
voisinage de l'\'equilibre. En se r\'ef\'erant \`a l'exemple du
sol\'eno\"ide\footnote{voir annexe~\ref{annexe:solenoide}.}, on
peut tenter d'interpr\'eter la forme des courbes obtenues.

Tout d'abord, il y a clairement une tr\`es forte sensibilit\'e aux
conditions initiales, en tout point de l'attracteur, aussi bien
dans le pass\'e que dans le futur. Dans le cas de la
figure~\ref{fig:passe_futur_eq_continu}, on constate ainsi des
\'ecarts de l'ordre de 4 en moins de 15 ans, aussi bien dans le
pass\'e que dans le futur, alors que les courbes \'etaient
initialement s\'epar\'ees de moins de $0\virg 04$.

Les dynamiques futures des diff\'erents points se r\'epartissent
de fa\c{c}on \`a peu pr\`es homog\`ene, au moins au cours des 10
premi\`eres ann\'ees. Ceux-ci se s\'eparent en effet selon leur
r\'epartition initiale dans la direction instable, puisque la
direction stable est contract\'ee dans le futur.
L'homog\'en\'eit\'e dans le futur traduit une r\'epartition \`a
peu pr\`es uniforme des points d'une orbite dans la direction
instable. Cette propri\'et\'e est \`a rapprocher du cas du
sol\'eno\"ide, o\`u la mesure physique (qui donne la r\'epartition
des points d'une orbite sur l'attracteur) poss\`ede une densit\'e
par rapport \`a la mesure de Lebesgue dans la direction instable.

La dynamique pass\'ee semble bien diff\'erente, les diff\'erentes
courbes se s\'eparant <<par paquets>>, et non plus de fa\c{c}on
homog\`ene. Ainsi, en \`a peine 3 ans, on observe d\'ej\`a une
diff\'erence de 1 pour l'une des courbes, tandis que la plupart
des points ont une orbite toujours tr\`es proche de l'\'equilibre.
Ceci nous donne des informations sur la mesure physique dans la
direction stable, puisque la direction instable est contract\'ee
dans le pass\'e. Ainsi, comme dans le cas du sol\'eno\"ide, il
semble que la mesure physique poss\`ede une densit\'e par rapport
\`a la mesure de Hausdorff sur un ensemble de Cantor de dimension
fractale $0\virg 3$.

Les r\'esultats observ\'es ici sont cependant moins clairs que
dans le cas du sol\'eno\"ide, pour $\absj{t} \geq 10$. Ceci est
sans doute d\^u au passage des orbites dans un pli, ph\'enom\`ene
qui ne se produit pas dans le cas du sol\'eno\"ide.

On observe le m\^eme type de r\'esultat en de nombreux autres
points de l'attracteur, d'autant plus nettement que l'on n'est pas
au voisinage d'une <<pointe>>.

\subsubsection{Point fixe, vari\'et\'e instable}
\begin{figure}
\begin{center}
\includegraphics[height=10cm]{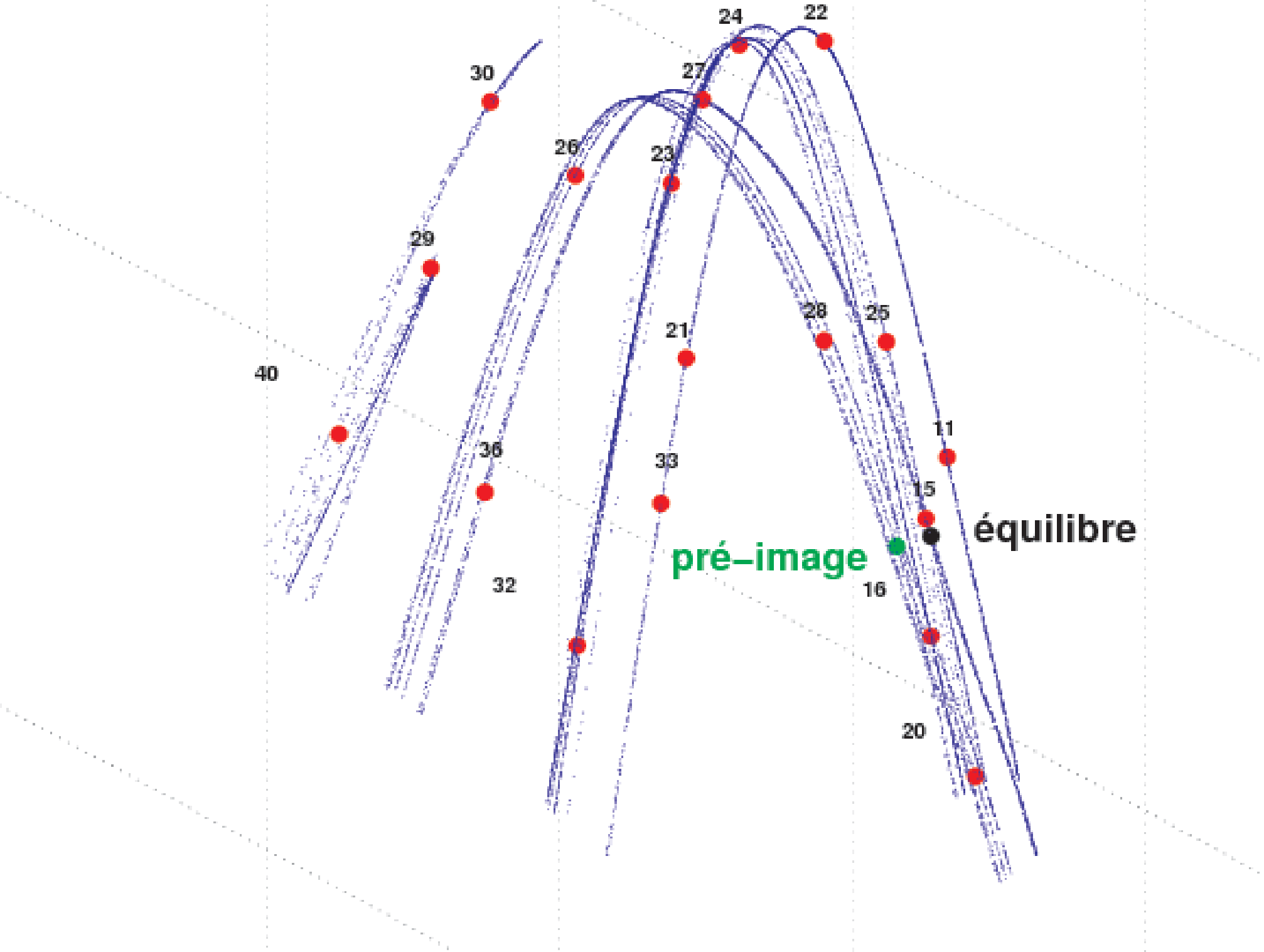}
\caption{\label{fig:equilibre_position} Position de l'\'equilibre,
et de sa <<pr\'eimage>>.}
\end{center}
\end{figure}

\begin{figure}
\begin{center}
\includegraphics[height=7cm]{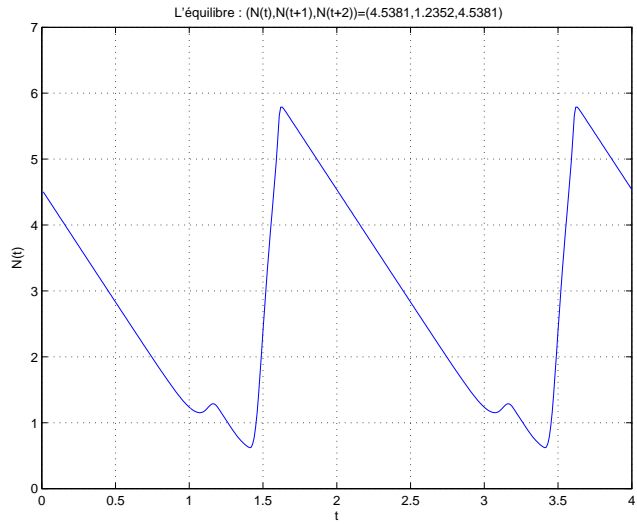}
\caption{\label{fig:equilibre_continu} L'\'equilibre : temps
continu.}
\end{center}
\end{figure}

Un point fixe (instable, bien s\^ur) a \'et\'e rep\'er\'e sur
l'attracteur. La figure~\ref{fig:equilibre_position} donne sa
localisation approximative en dimension 3, et la fonction continue
correspondante est repr\'esent\'ee
figure~\ref{fig:equilibre_continu}.

\begin{figure}
\begin{center}
\includegraphics[height=7cm]{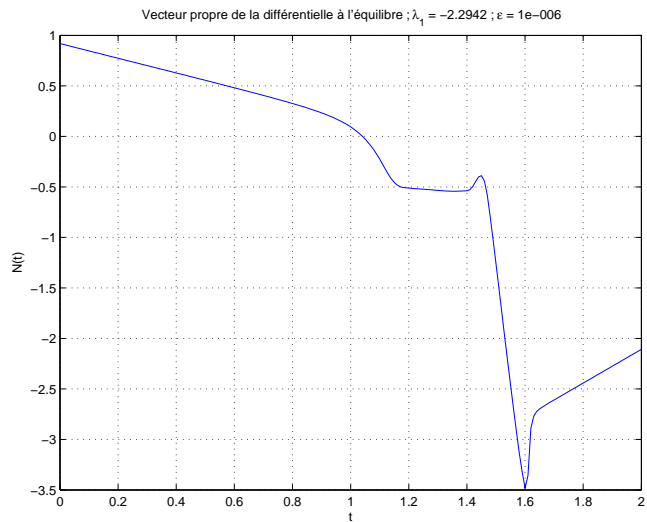}
\caption{\label{fig:diff_equilibre_vp1} Diff\'erentielle \`a
l'\'equilibre : premier vecteur propre.}
\end{center}
\end{figure}

\begin{figure}
\begin{center}
\includegraphics[height=7cm]{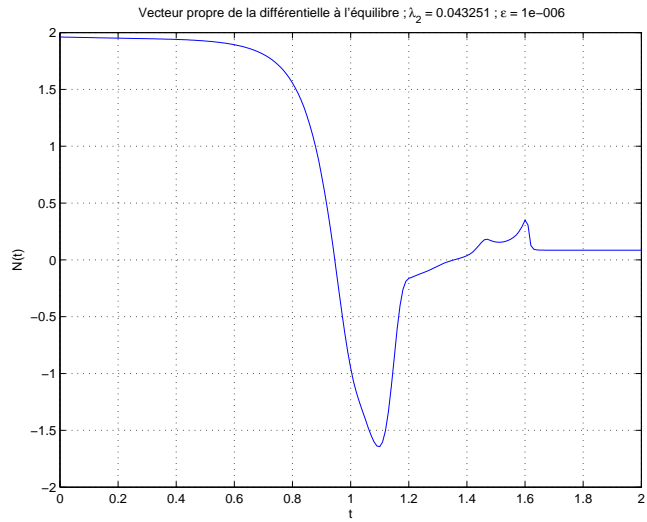}
\caption{\label{fig:diff_equilibre_vp2} Diff\'erentielle \`a
l'\'equilibre : second vecteur propre.}
\end{center}
\end{figure}

\begin{figure}
\begin{center}
\includegraphics[height=7cm]{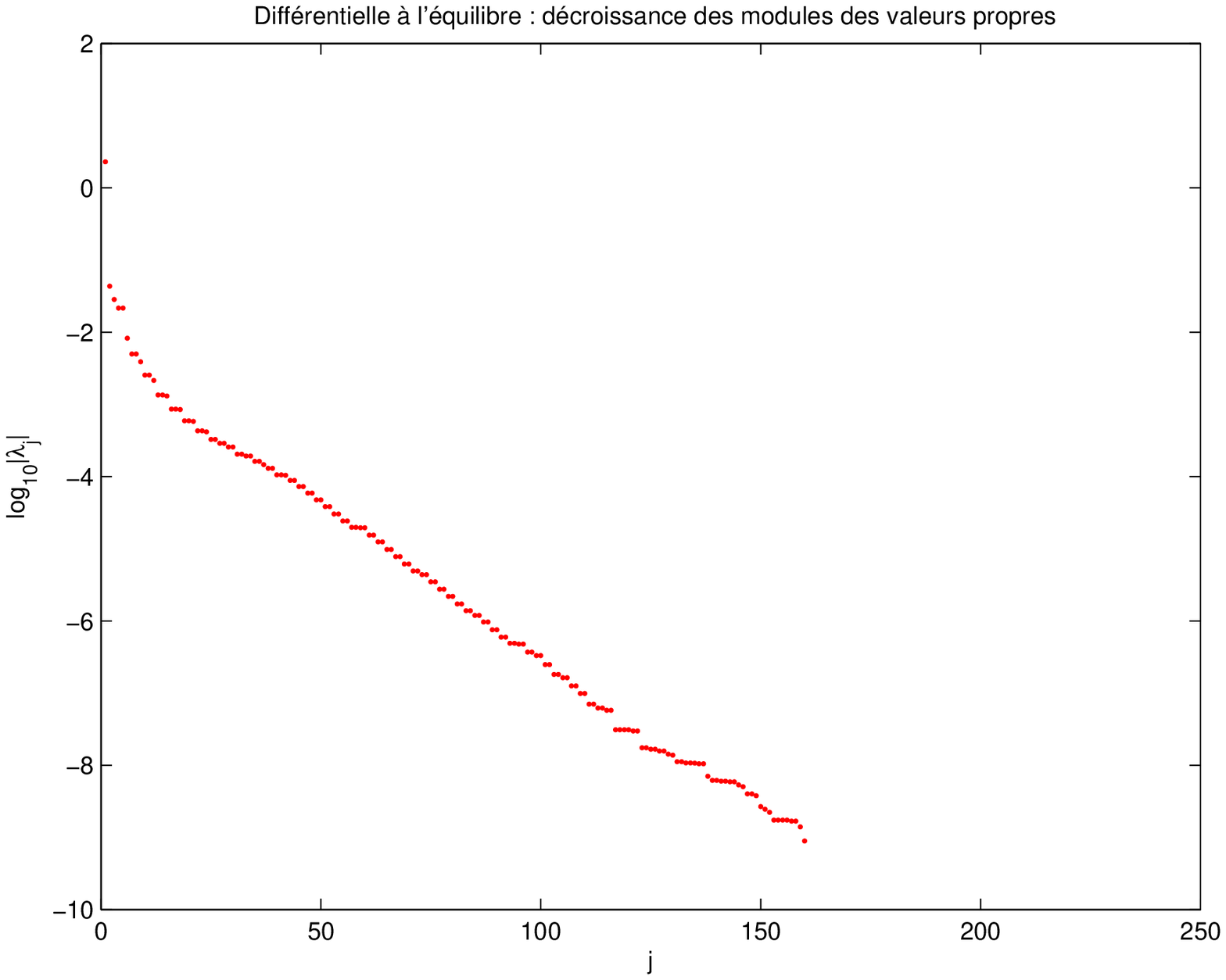}
\caption{\label{fig:diff_equilibre_val_propres} Diff\'erentielle
\`a l'\'equilibre : d\'ecroissance des valeurs propres.}
\end{center}
\end{figure}

\paragraph{Diff\'erentielle de $T^2$ \`a l'\'equilibre} Celle-ci nous fournit beaucoup de
renseignements sur le syst\`eme dynamique. On peut calculer
ais\'ement ses valeurs propres. Une seule est de module
strictement sup\'erieur \`a 1 et correspond \`a l'expansion dans
la direction instable : $\lambda_1 \approx -2\virg 29$. Le vecteur
propre associ\'e est trac\'e figure~\ref{fig:diff_equilibre_vp1}.
La seconde plus grande valeur propre (la pr\'ecision de ce calcul
est faible) est $\lambda_2 \approx 0\virg 043$, et le vecteur
propre associ\'e est repr\'esent\'e
figure~\ref{fig:diff_equilibre_vp2}. Les modules des valeurs
propres suivantes d\'ecroissent ensuite rapidement, comme le
montre la figure~\ref{fig:diff_equilibre_val_propres}.

Le point fixe est donc hyperbolique, puisqu'il n'a aucune valeur
propre de module proche de 1. De plus, on constate que l'expansion
est assez raisonnable (elle est sans doute un peu plus forte dans
certaines zones, un peu moins dans d'autres, mais reste de cet
ordre de grandeur), tandis que la contraction est beaucoup plus
forte. Cela nous donne un argument suppl\'ementaire pour penser
que 3 dimensions suffisent \`a repr\'esenter l'attracteur : les
valeurs propres suivantes ayant un module encore plus petit,
l'attracteur n'est vraiment \'etendu que dans 2 ou 3 dimensions,
les autres \'etant peu importantes.

\begin{figure}
\begin{center}
     \includegraphics[width=\textwidth]{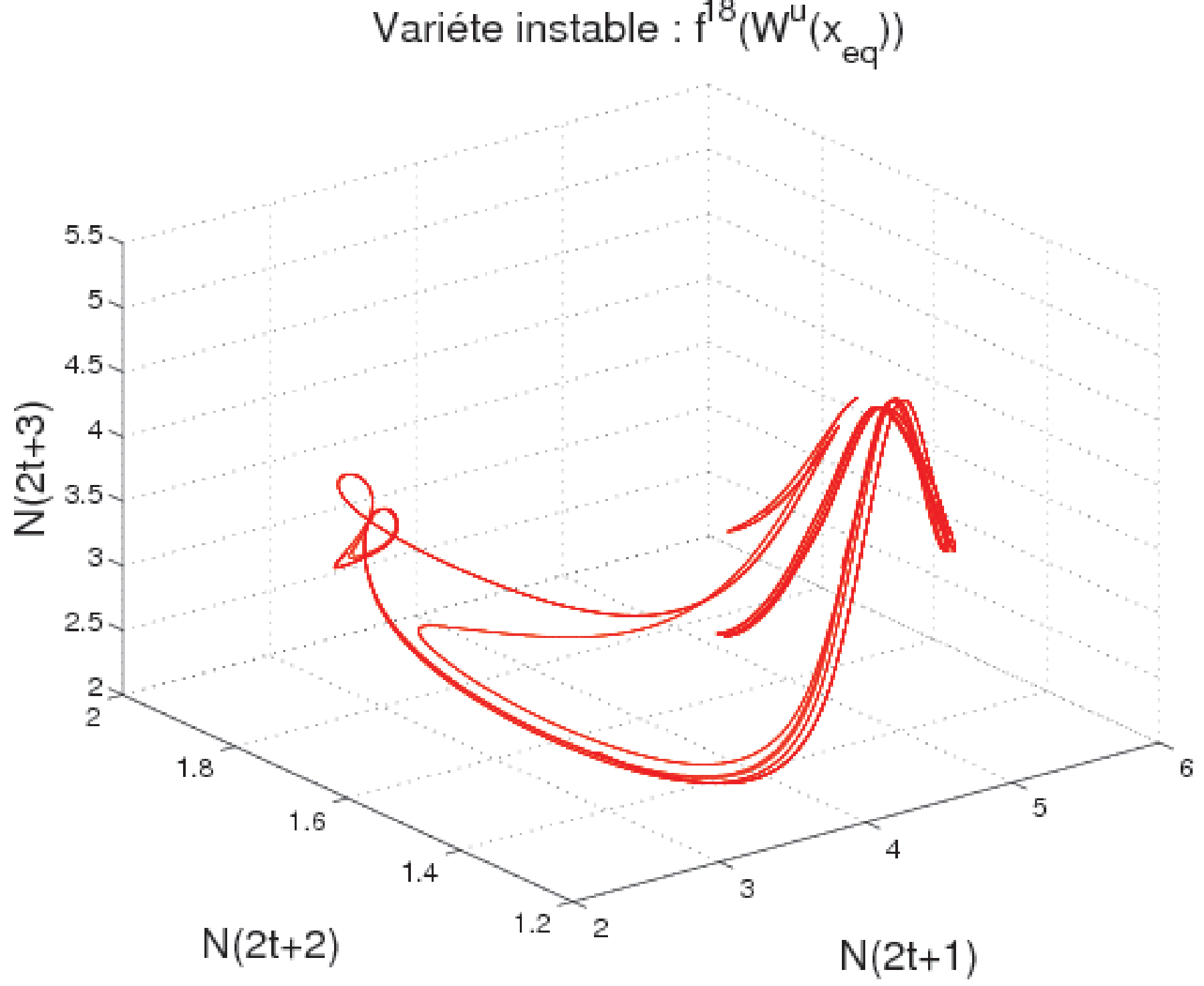}
\caption{\label{fig:variete_instable} Vari\'et\'e instable.}
\end{center}
\end{figure}

\paragraph{Vari\'et\'e instable} \label{sec:variete_instable}
On peut d\'eterminer la vari\'et\'e instable en regardant les
images par $f=T^2$ d'un segment situ\'e dans la direction instable
au voisinage de l'origine (figure~\ref{fig:variete_instable}).

\begin{figure}
\begin{center}
\begin{tabular}{c@{}c}
     \includegraphics[width=6.5cm]{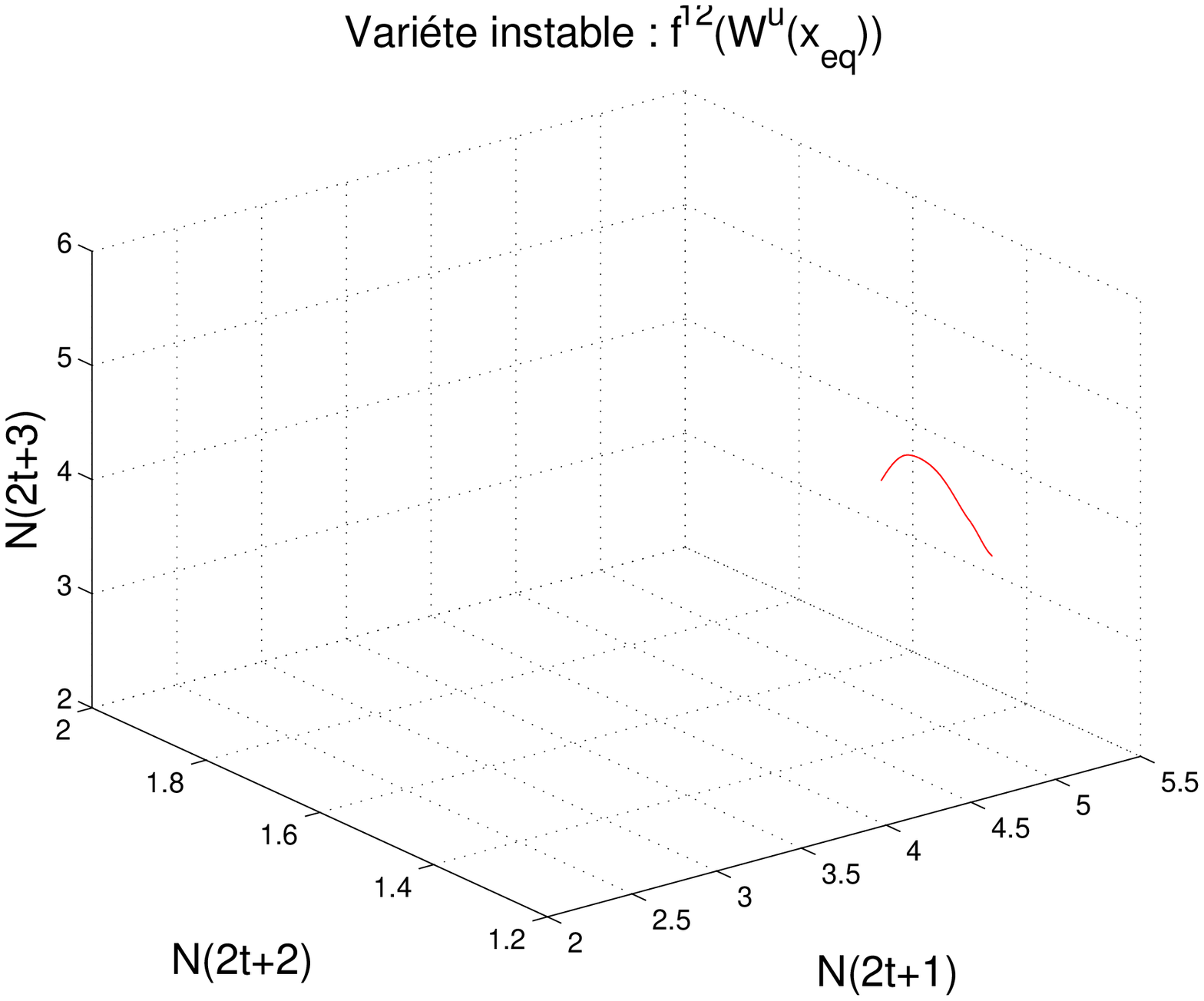}
     &
     \includegraphics[width=6.5cm]{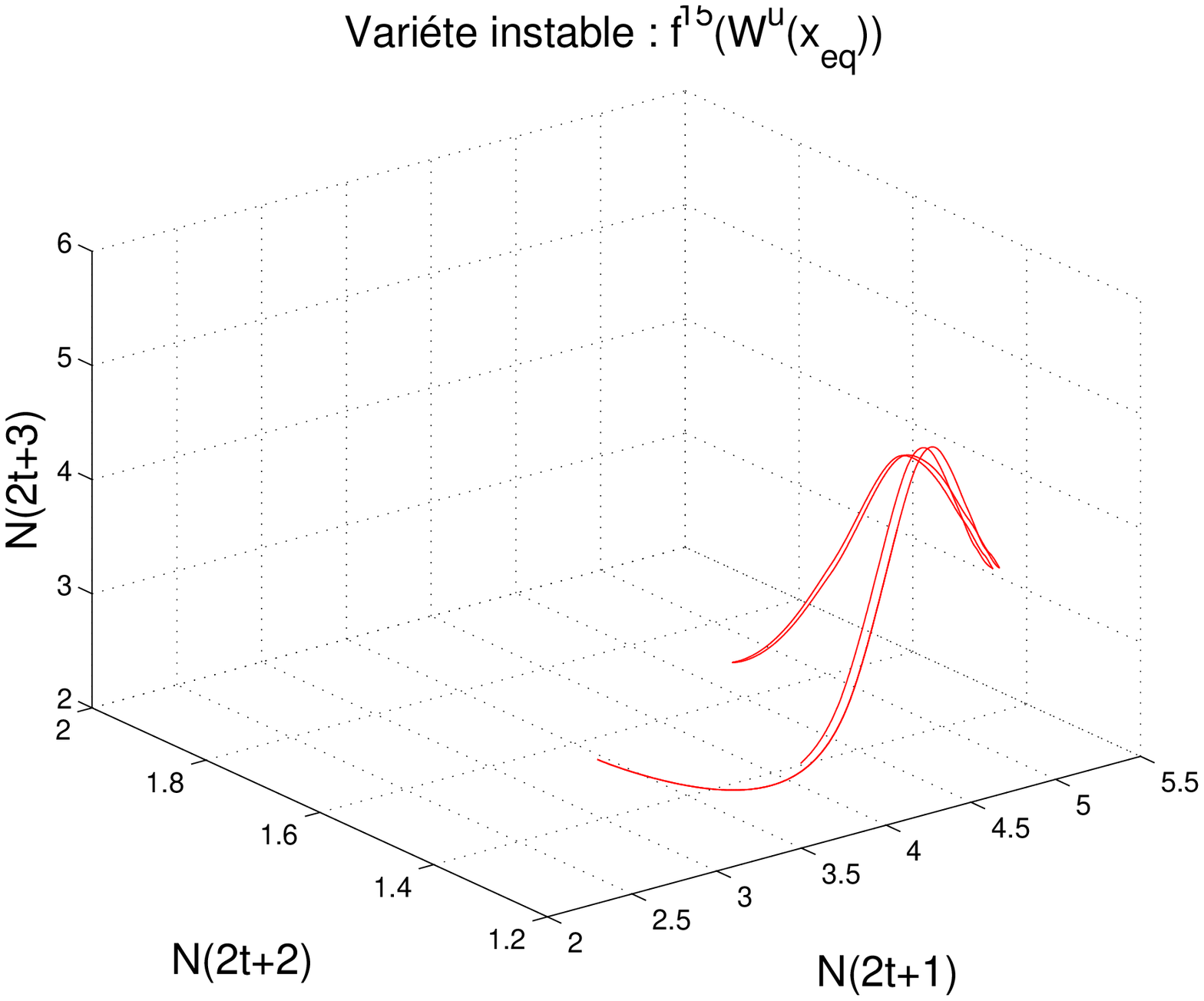}
          \\
     (a) $n=12$ & (b) $n=15$ \\
     \includegraphics[width=6.5cm]{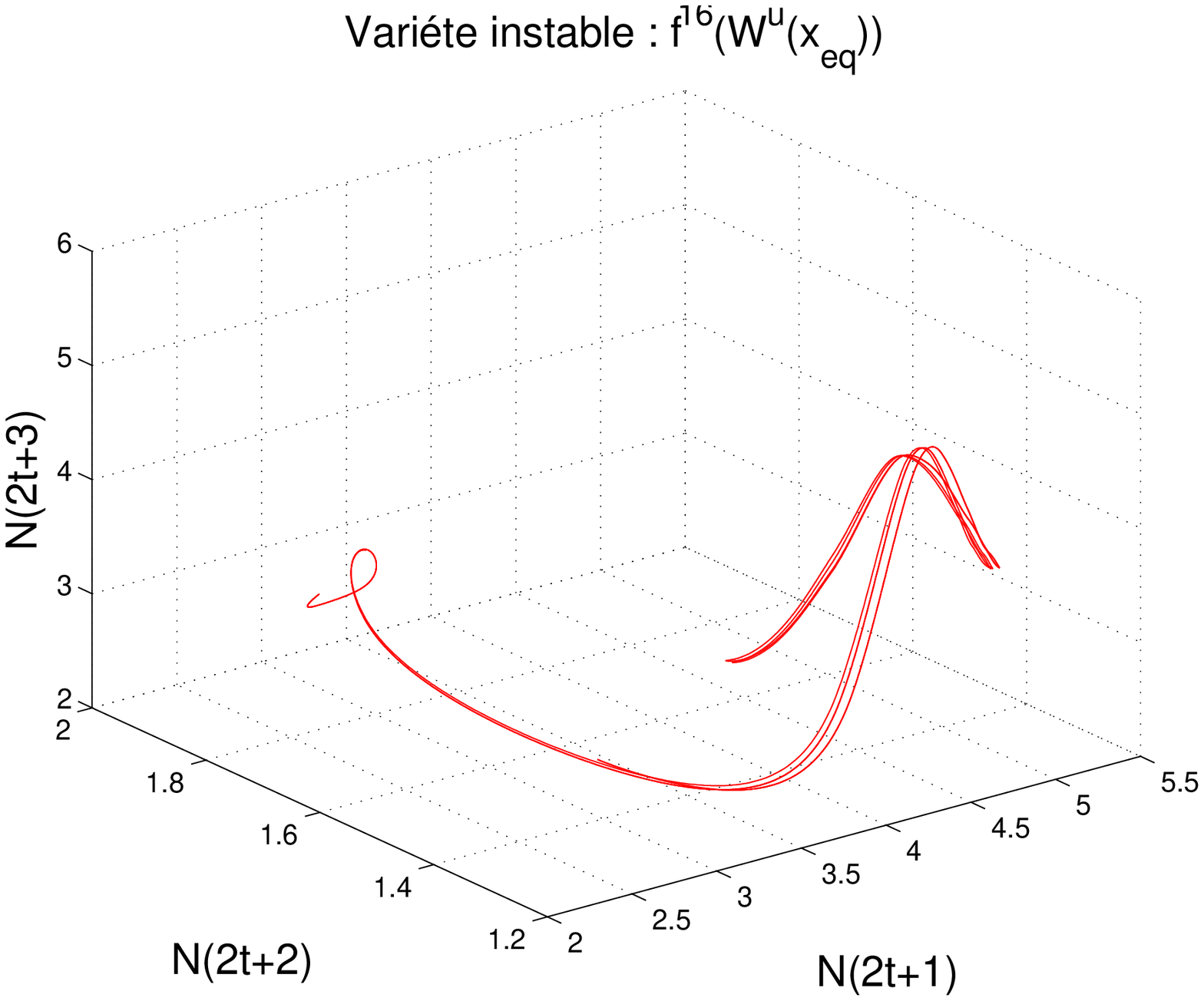}
     &
     \includegraphics[width=6.5cm]{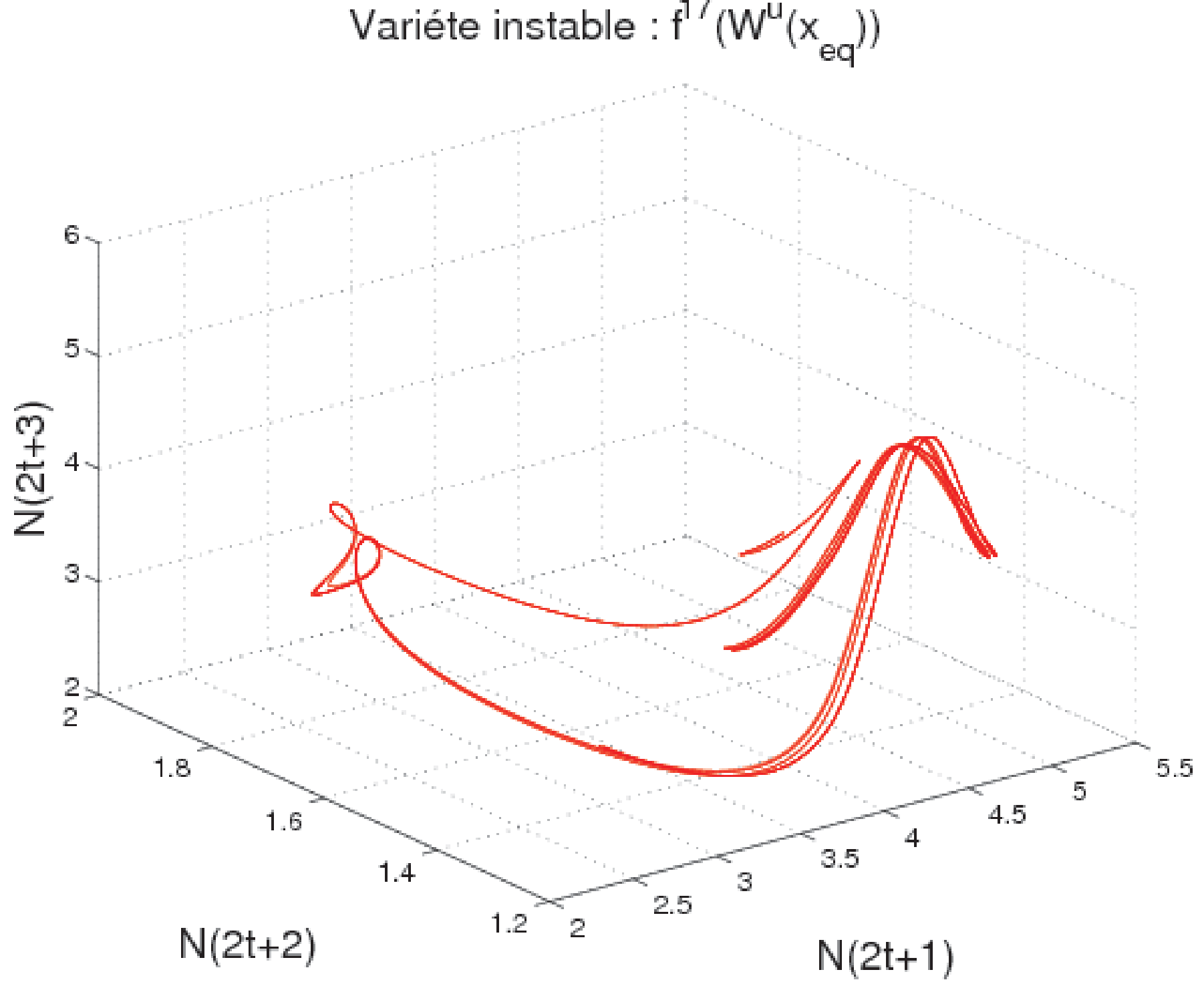}
          \\
     (c) $n=16$ & (d) $n=17$
\end{tabular}
\caption{\label{fig:var_u_deploiement} D\'eploiement de la
vari\'et\'e instable \`a l'\'equilibre $f^n(W^u(x_{eq}))$.}
\end{center}
\end{figure}

Il est int\'eressant, d'un point de vue dynamique, de visualiser
comment celle-ci se d\'eploie \`a l'int\'erieur de l'attracteur,
\`a la fois d'un point de vue dynamique et d'un point de vue
g\'eom\'etrique. En effet, la vari\'et\'e instable est une courbe
continue, ce qui nous donne une id\'ee plus pr\'ecise de la
g\'eom\'etrie de l'attracteur que lorsque nous ne disposons que
d'un nuage de points. Ceci est fait avec une
animation\footnote{consulter les fichiers
\fichier{var\_u\_0.1\_2\_1\_18.avi} et
\fichier{var\_u\_0.1\_2\_18.avi}. } dont la
figure~\ref{fig:var_u_deploiement} donne des extraits.

On peut d\'ecrire $f^n(W^u(x_{eq}))$ pour les valeurs successives
de $n$ de la fa\c{c}on suivante :
\begin{itemize}
\item $1 \leq n \leq 11$ : une petite portion rectiligne autour de
l'\'equilibre.
\item $n=12$ : un filament pli\'e, dans la direction de la pointe
interm\'ediaire situ\'ee dans le cou (43)
(figure~\ref{fig:var_u_deploiement}a).
\item $n=13$ : le filament est un peu \'etendu et atteint le bas de la pointe.
\item $n=14$ : le bas de la chevelure est atteint, ainsi que le cou.
\item $n=15$ : le creux (avec la pointe 54) est atteint. Les filaments sont
d\'edoubl\'es (figure~\ref{fig:var_u_deploiement}b).
\item $n=16$ : un c\^ot\'e du triangle (71--75) et la droite de la queue sont atteints. Les
filaments se d\'edoublent ailleurs
(figure~\ref{fig:var_u_deploiement}c).
\item $n=17$ : un autre c\^ot\'e du triangle (76--80) et la pointe arri\`ere sont
atteints (figure~\ref{fig:var_u_deploiement}d).
\item $n=18$ : le dernier c\^ot\'e du triangle (59--66) est
atteint. C'est la derni\`ere r\'egion importante \`a \^etre
touch\'ee. Notons tout de m\^eme que le filament 58--74 (c'est
l'un des embranchements) n'est pas encore atteint.
\end{itemize}

\subsubsection{Formation du pli}
\begin{figure}
\begin{center}
\begin{tabular}{c}
     \includegraphics[height=8.5cm]{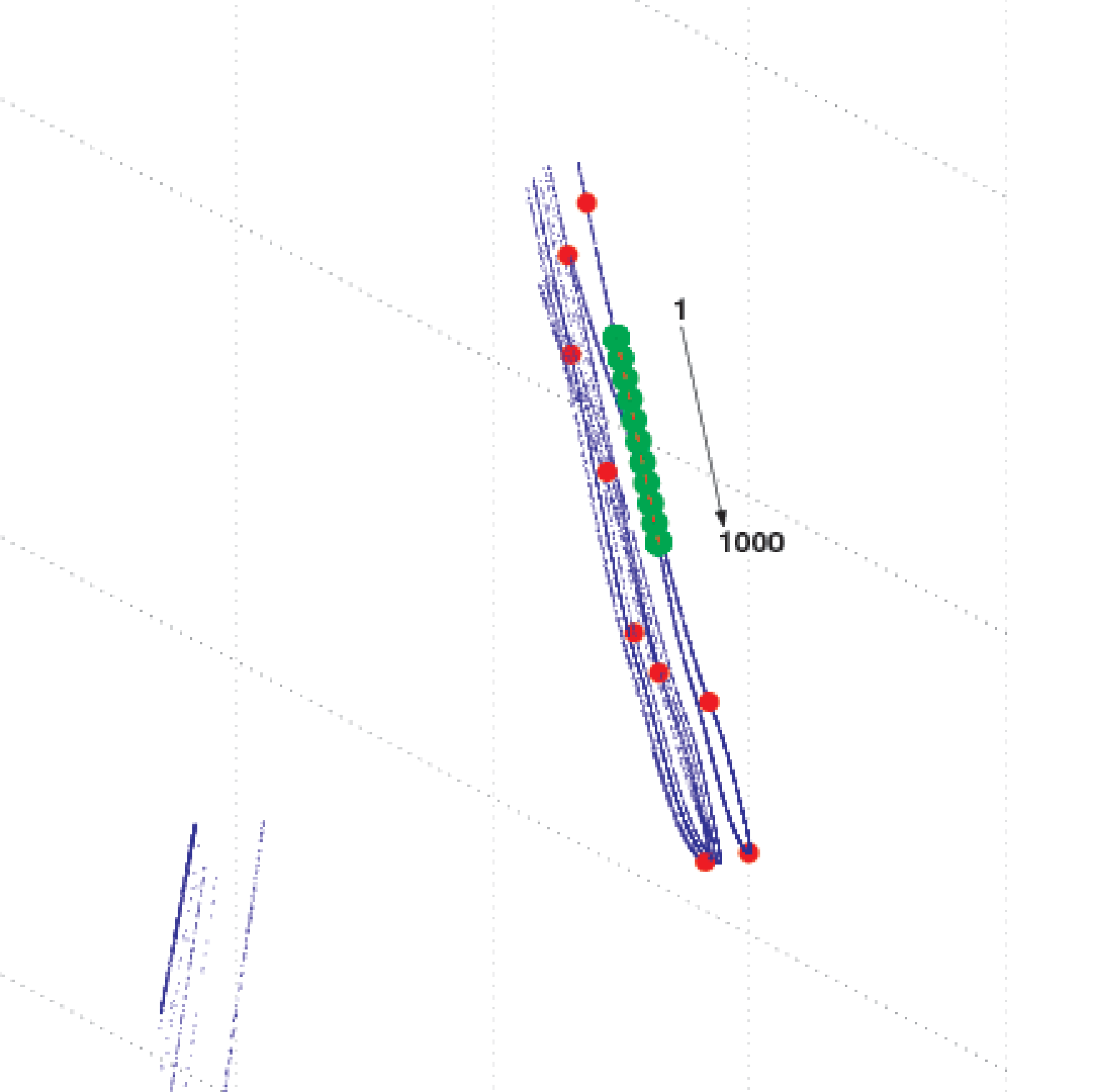} \\
    (a) $x_1 , \ldots , x_{1000}$ \\
     \includegraphics[height=8.5cm]{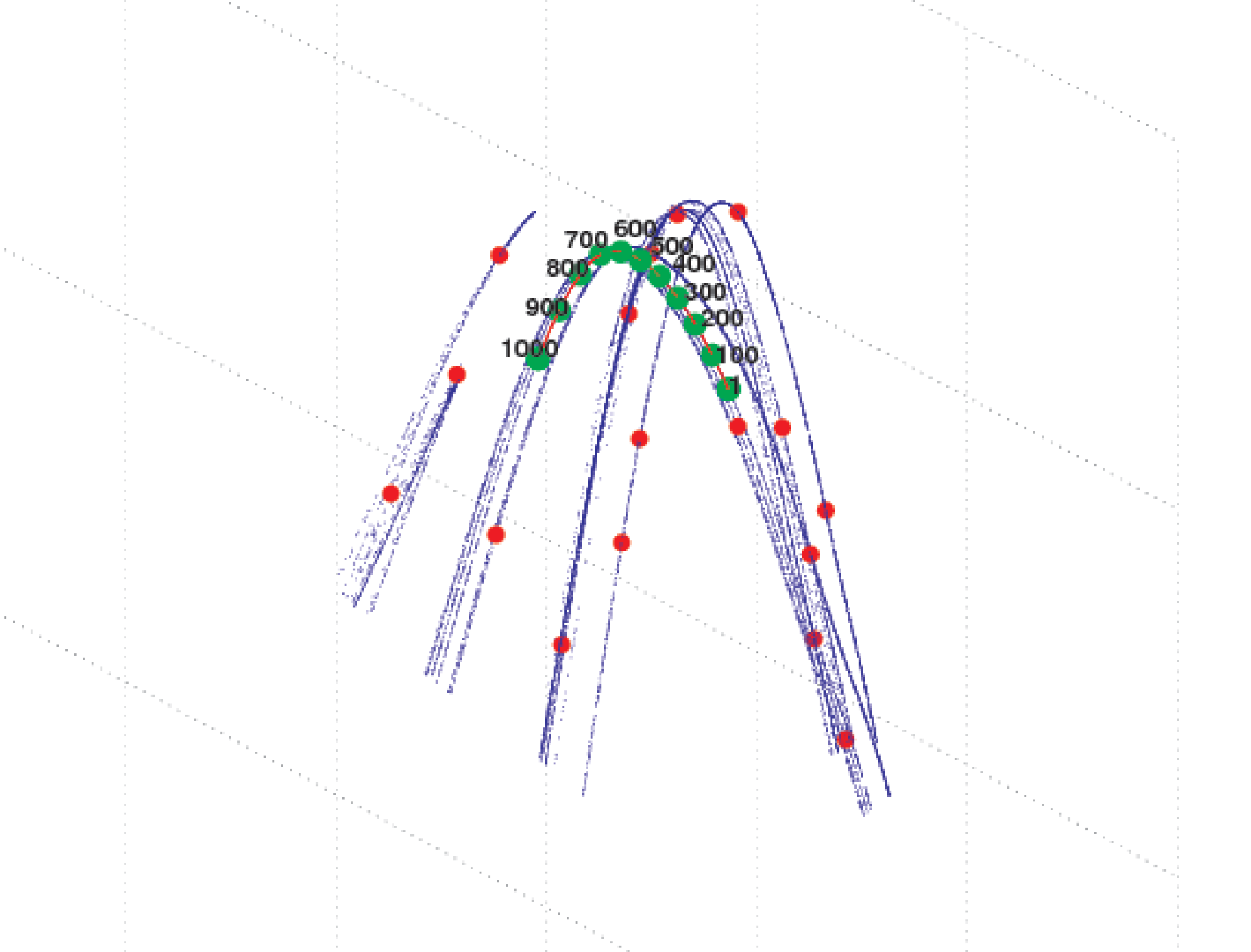} \\
    (b) $T^2(x_1), \ldots, T^2(x_{1000})$
\end{tabular}
\caption{\label{fig:pli} Localisation du pli et de sa pr\'eimage.}
\end{center}
\end{figure}

Une des caract\'eristiques du syst\`eme est l'existence d'un
pli\footnote{Il y en a peut-\^etre plusieurs, mais il n'est pas
\'evident de distinguer si deux r\'egions pli\'ees sont ou non
ind\'ependantes. Cela demanderait une \'etude plus approfondie.}.
Sous l'action de $T$, la r\'egion quasi-rectiligne de la
figure~\ref{fig:pli}a se courbe progressivement pour arriver dans
la r\'egion <<pli\'ee>> de la figure~\ref{fig:pli}b.

\begin{figure}
\begin{center}
\begin{tabular}{c@{}c}
     \includegraphics[width=7.5cm]{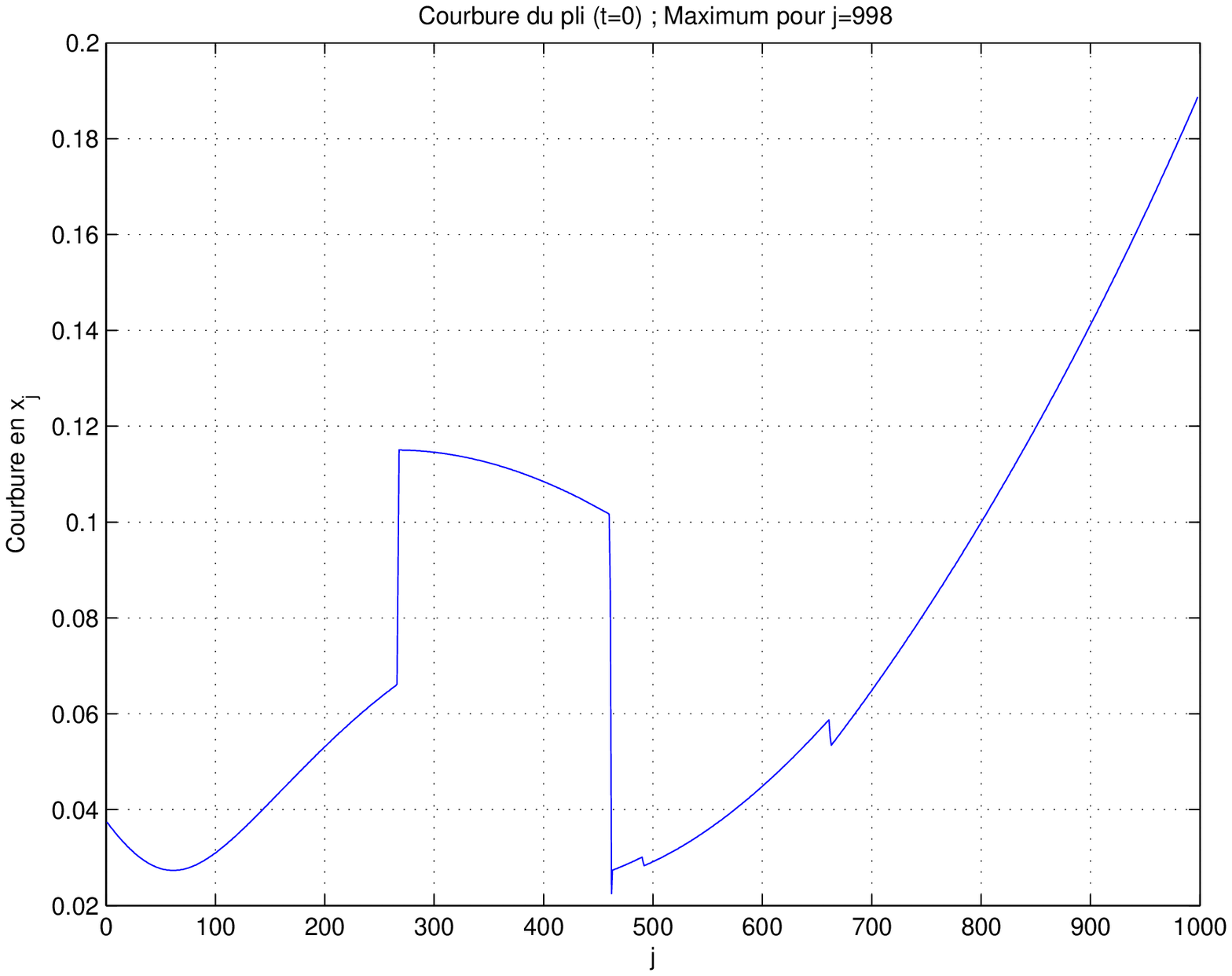}
     &
     \includegraphics[width=7.5cm]{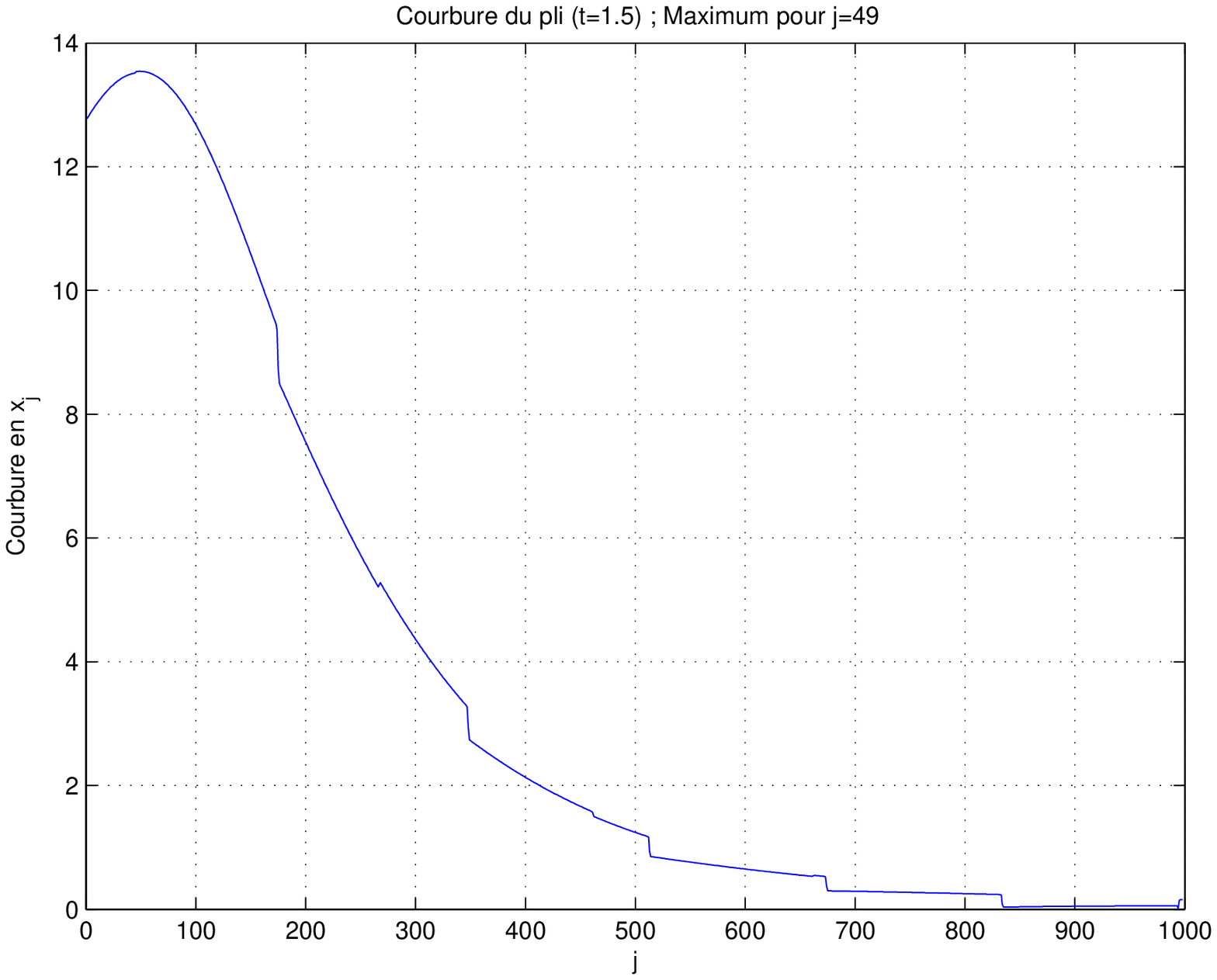}
     \\
     (a) $t=0$ & (b) $t=1\virg 5$
     \\
     \includegraphics[width=7.5cm]{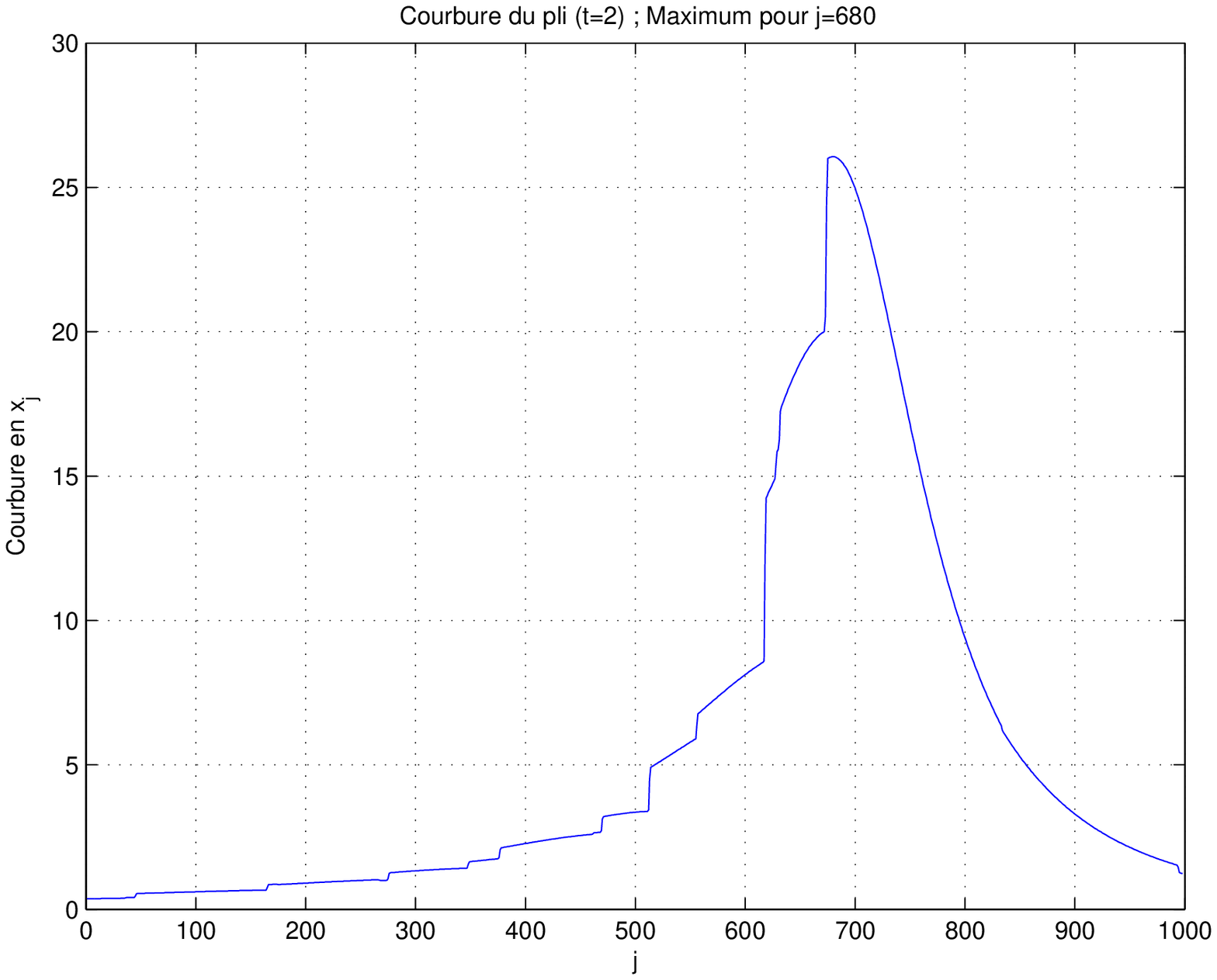}
     &
     \includegraphics[width=7.5cm]{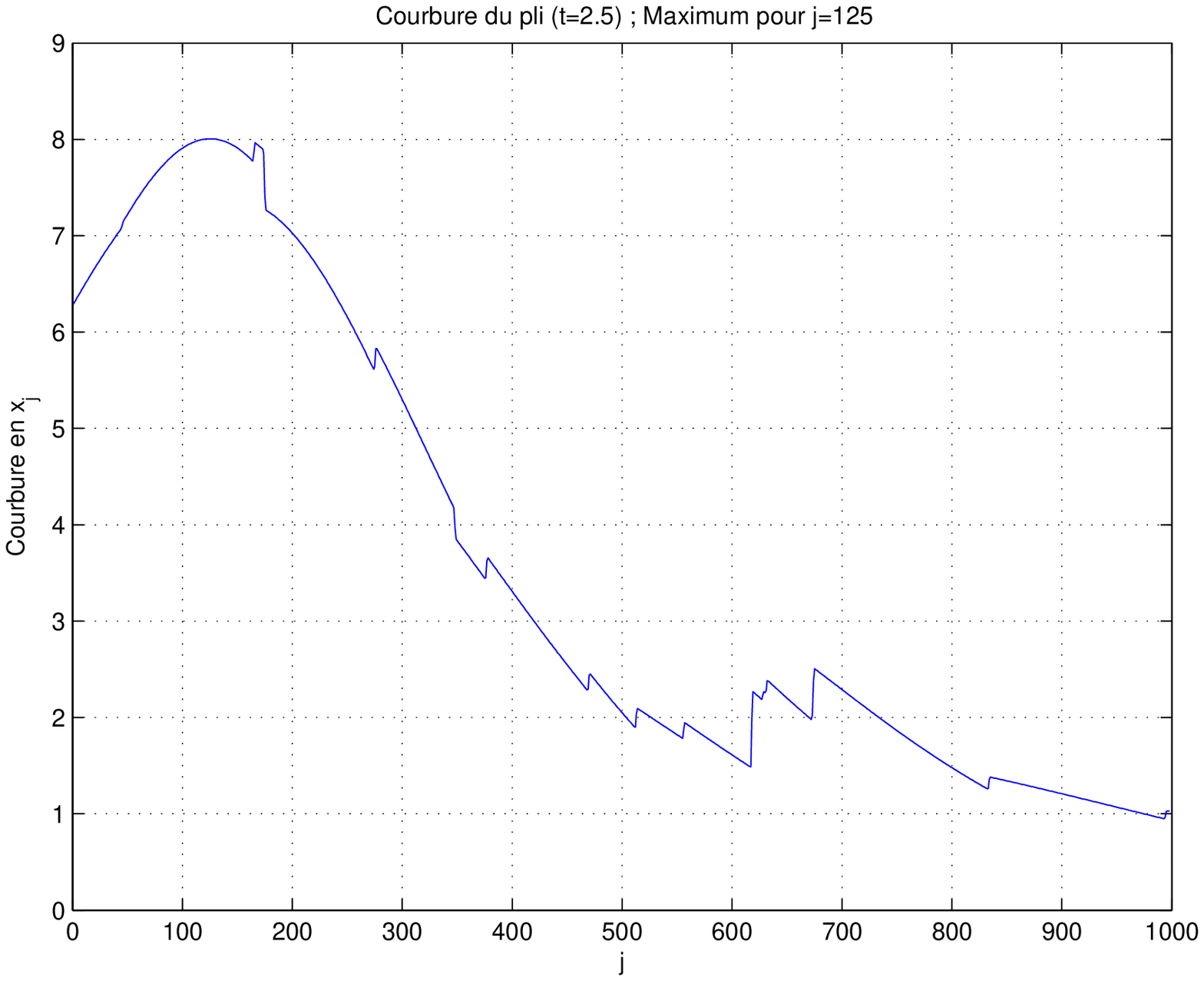}
     \\
     (c) $t=2$ & (d) $t=2\virg 5$
     \\
     \includegraphics[width=7.5cm]{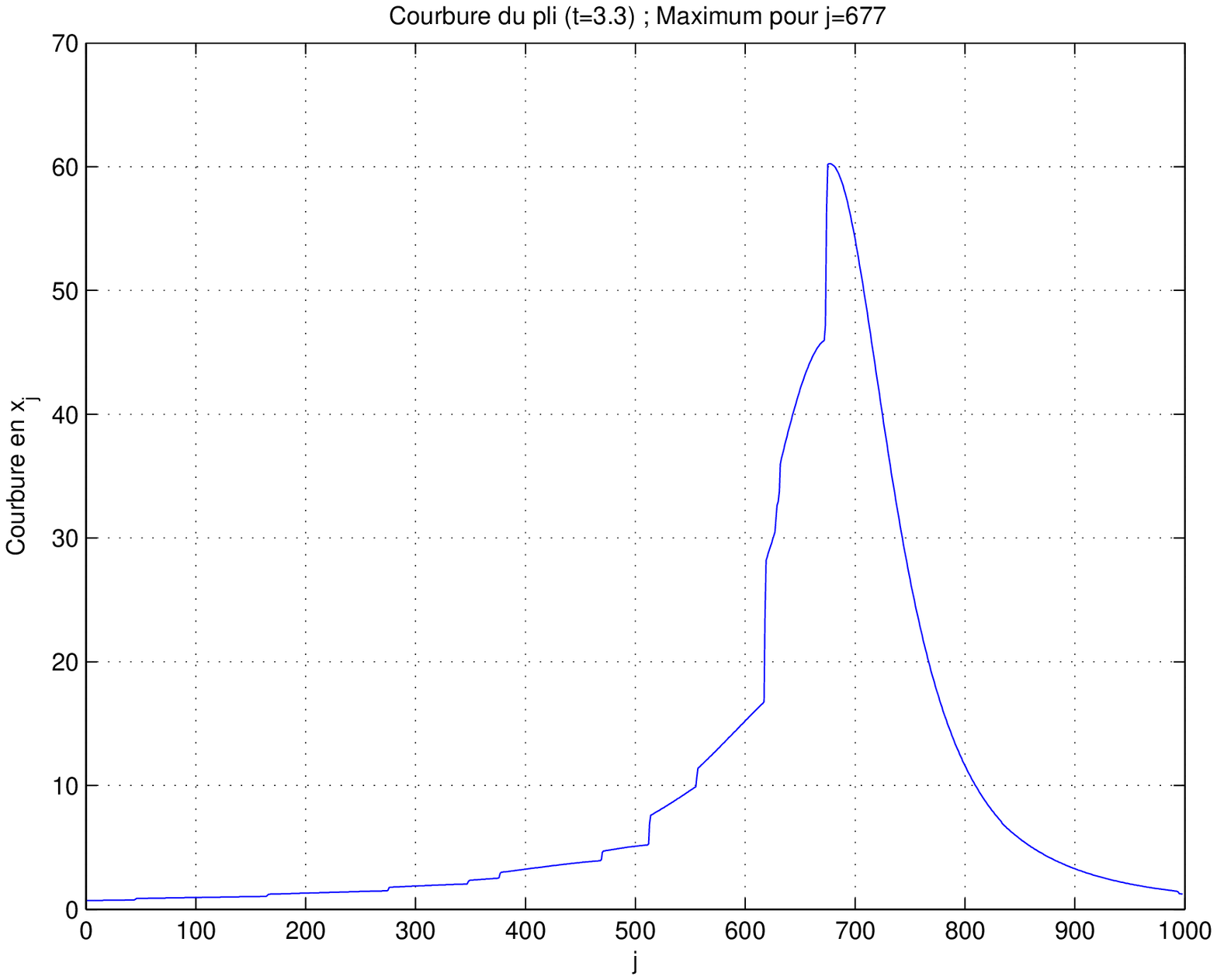}
     &
     \includegraphics[width=7.5cm]{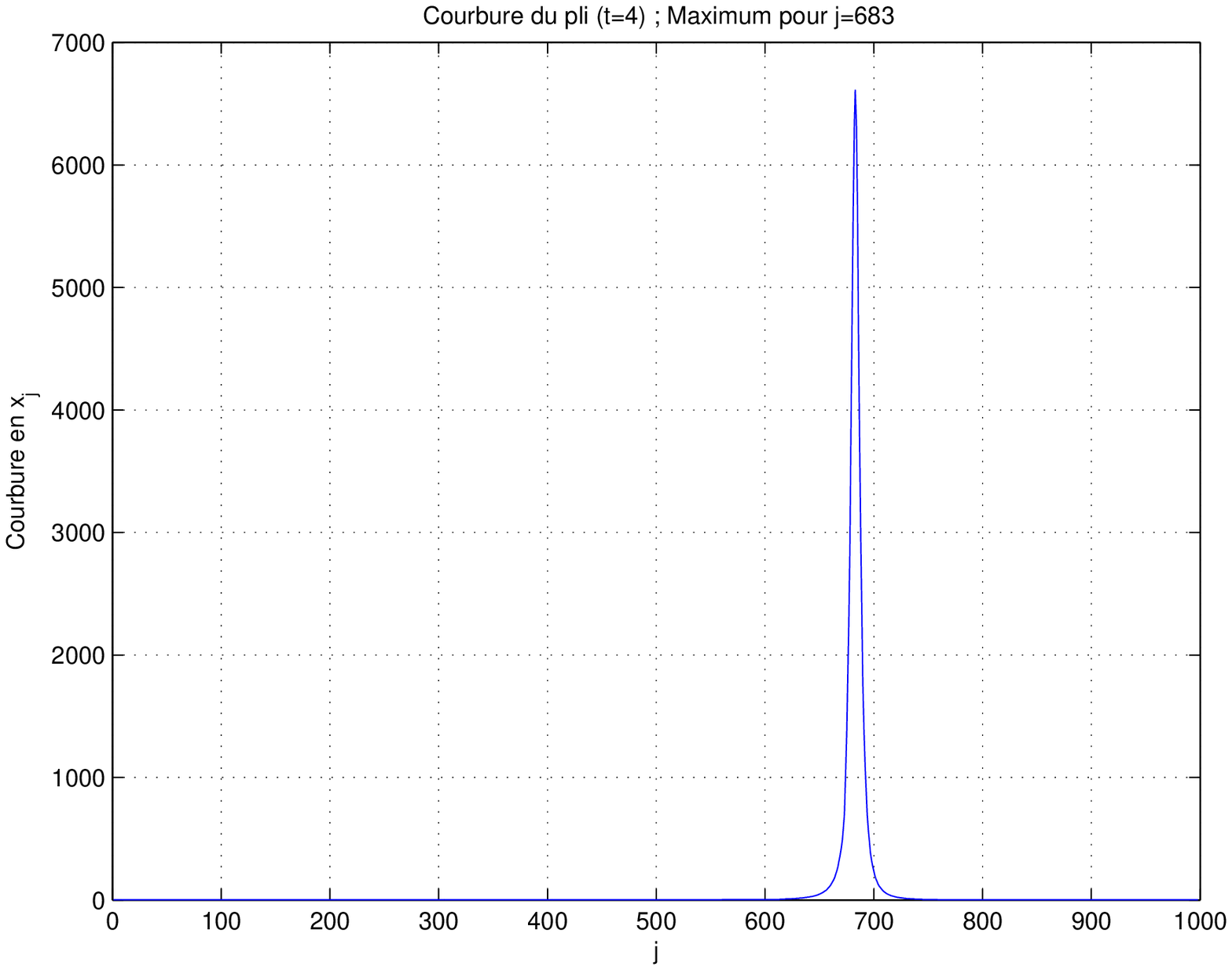}
     \\
     (e) $t=3\virg 3$ & (f) $t=4$
\end{tabular}
\caption{\label{fig:courbure_pli} Formation du pli : courbure.}
\end{center}
\end{figure}

\begin{figure}
\begin{center}
\begin{tabular}{c@{}c}
     \includegraphics[height=7cm]{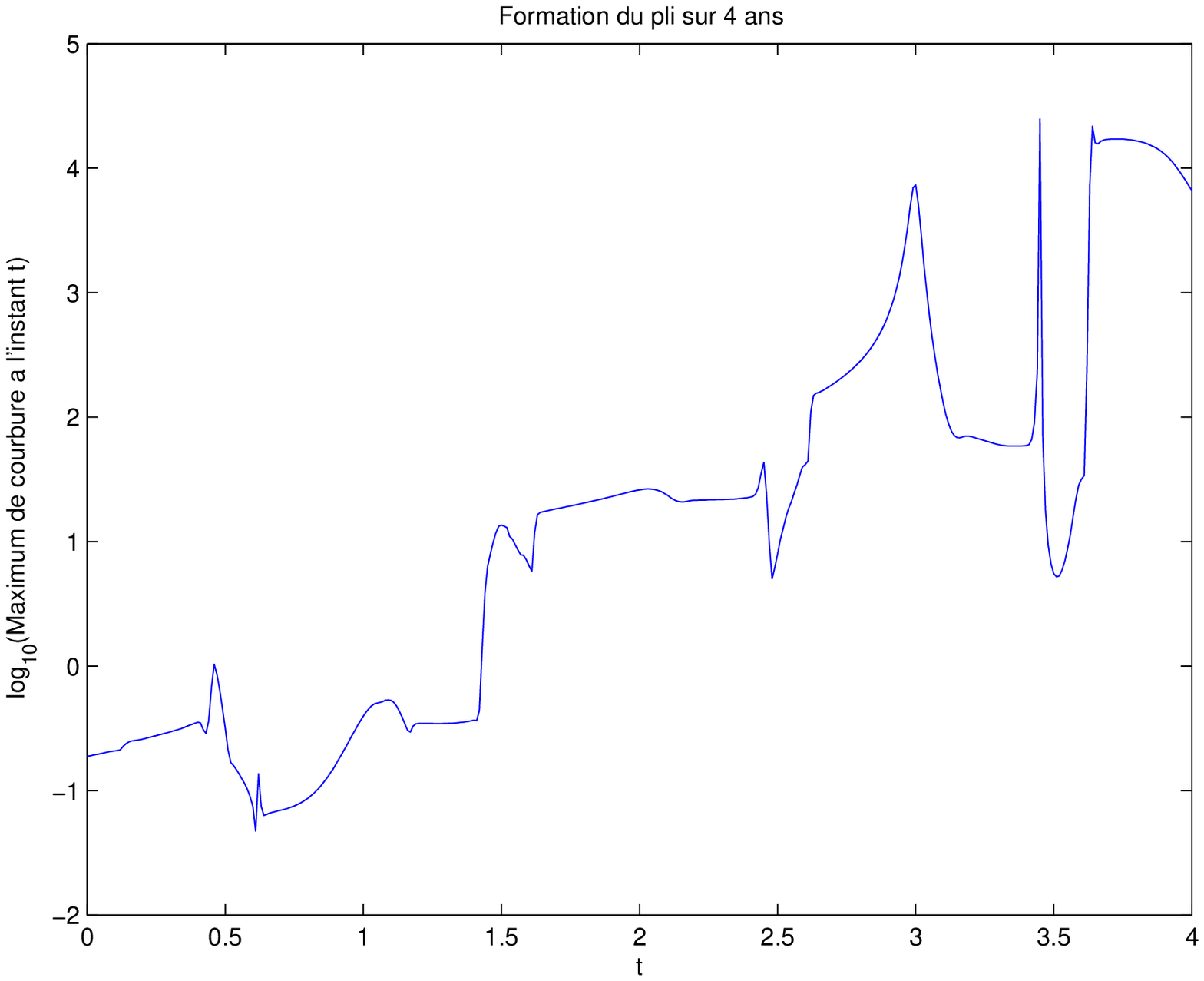}
     \\
     (a) Valeur maximale
     \\
     \includegraphics[height=7cm]{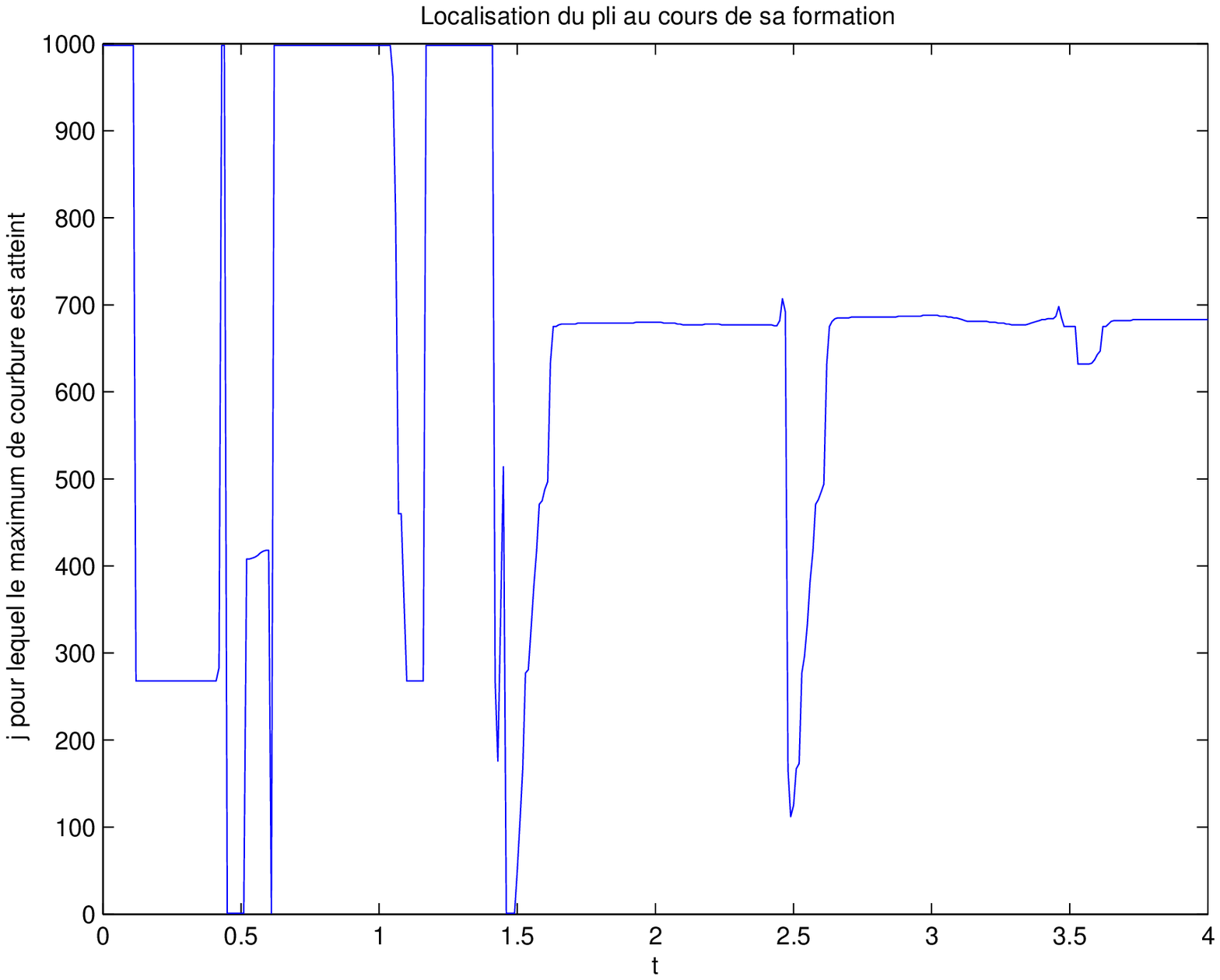}
     \\
     (b) Localisation du maximum
\end{tabular}
     \caption{\label{fig:courbure_maxi}Formation du pli : maximum de courbure.}
\end{center}
\end{figure}

\paragraph{Approche g\'eom\'etrique}
On peut visualiser la formation de ce pli en calculant la courbure
au niveau du pli (sur un m\^eme filament de l'attracteur) \`a
diff\'erents instants. Le r\'esultat\footnote{Cette figure est
extraite d'une animation donnant plus d'informations sur la
formation g\'eom\'etrique du pli.} est reproduit
figure~\ref{fig:courbure_pli}. Hormis quelques irr\'egularit\'es,
il se forme clairement un pli pour $1\virg 5 < t < 1 \virg 6$, et
celui-ci s'accentue fortement pour former un pli tr\`es marqu\'e
\`a $t=4$. L'\'evolution du maximum de courbure est report\'ee
figure~\ref{fig:courbure_maxi}a. Cette \'etude permet \'egalement
de localiser tr\`es pr\'ecis\'ement l'endroit pli\'e, en notant
\`a chaque instant la position du maximum de courbure sur le
segment. On constate avec la figure~\ref{fig:courbure_maxi}b que
celui-ci est situ\'e au point $690$ du segment initialement
choisi.

\paragraph{\'Etude des discontinuit\'es de la courbure}
Les diff\'erentes figures montrant la formation du pli
pr\'esentent cependant quelques irr\'egularit\'es. Ainsi, on peut
se demander \`a quoi sont d\^ues les discontinuit\'es de la
courbure figure~\ref{fig:courbure_pli}. Au vu de l'impr\'ecision
de la m\'ethode de localisation du pli, il est quasiment certain
que la ligne bris\'ee consid\'er\'ee saute d'un filament \`a
l'autre. Il est \'egalement possible que le manque de
r\'egularit\'e des fonctions du mod\`ele (qui sont $C^1$ et non
$C^2$) engendre des ruptures de courbure au niveau de l'attracteur
lui-m\^eme.

Pour tester cette seconde hypoth\`ese, nous avons tent\'e
d'utiliser la vari\'et\'e instable globale de l'\'equilibre, qui
devrait nous assurer que l'on consid\`ere un seul filament. Les
r\'esultats obtenus sont repr\'esent\'es en
annexe~\ref{annexe:details}, avec les figures \ref{fig:pli2} \`a
\ref{fig:courbure_maxi2}. Il y a toujours des discontinuit\'es
dans la courbure le long du filament. On peut sans aucun doute
attribuer le pic de courbure aux environs de $j=500$ \`a un saut
d'un filament \`a un autre (figure~\ref{fig:courbure_pli2}). Il
semble tr\`es difficile d'extraire avec suffisamment de
pr\'ecision un seul filament de l'attracteur dans la zone de pli.
Cette t\^ache est peut-\^etre simplement rendue impossible par
l'aspect fractal de celui-ci. En revanche, la similitude des
figures \ref{fig:courbure_pli} et \ref{fig:courbure_pli2} confirme
nos observations sur le processus de pliage. La
figure~\ref{fig:courbure_maxi2}b est d'ailleurs plus simple \`a
interpr\'eter, puisque l'on a uniquement le pic artificiel \`a
$j\approx 500$ et le pic r\'eel \`a $j \approx 875$.

\begin{figure}
\begin{center}
     \includegraphics[width=10cm]{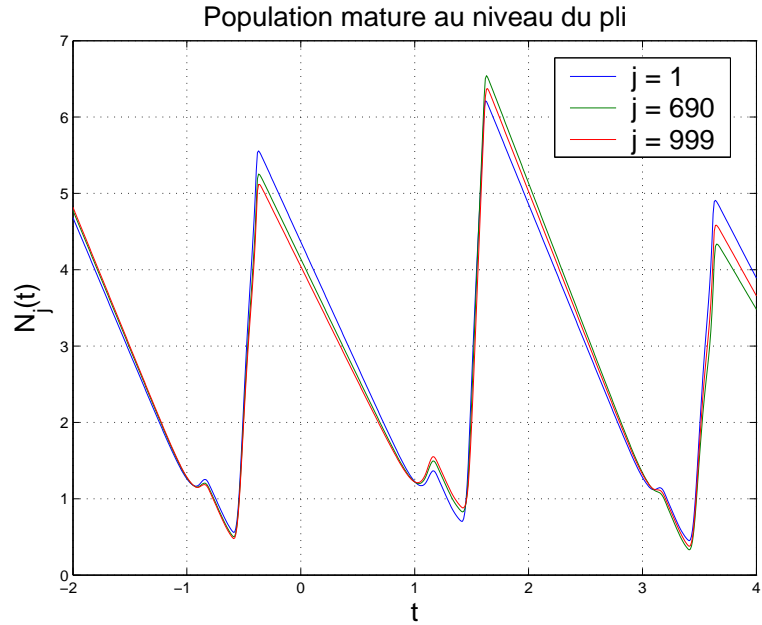}
\caption{\label{fig:pli_mature} \'Evolution en temps continu au
niveau du pli.}
\end{center}
\end{figure}

\begin{figure}
\begin{center}
     \includegraphics[height=7cm]{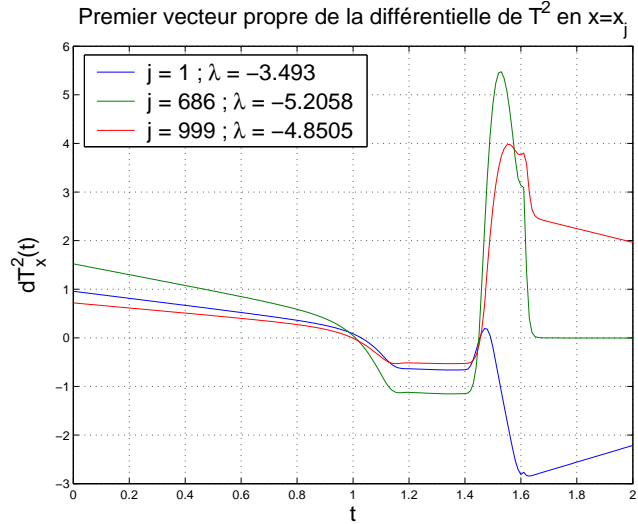}
\caption{\label{fig:pli_diff_avant_vp1} Diff\'erentielle de $T^2$
sur le filament non-pli\'e : premier vecteur propre.}
\end{center}
\end{figure}

La figure~\ref{fig:courbure_maxi} pr\'esente \'egalement des
irr\'egularit\'es facilement interpr\'etables. En effet, la
courbure ne cro\^it pas toujours avec le temps, et le maximum de
courbure se d\'eplace par moments, m\^eme apr\`es la formation du
pli \`a $t=1\virg 6$. Ces irr\'egularit\'es peuvent s'expliquer
par le non-uniformit\'e du pincement et de l'\'etirement dans la
dynamique. Ainsi, le pincement fait augmenter fortement la
courbure, et l'\'etirement entra\^ine au contraire sa diminution.
Sur deux ans, cela se traduit par une forte augmentation de la
courbure, mais en temps continu, il est des p\'eriodes o\`u la
tendance s'inverse. Un argument en faveur de cette hypoth\`ese est
que ces p\'eriodes co\"incident avec les instants des pics de
population, c'est-\`a-dire les instants d'expansion maximale (voir
figures \ref{fig:pli_mature} et \ref{fig:pli_diff_avant_vp1}).

\paragraph{Approche dynamique}
On peut aborder la question du pli d'un point de vue dynamique et
non plus purement g\'eom\'etrique en observant l'\'evolution de la
population en temps continu (figure~\ref{fig:pli_mature}). Les
points $x_j$ du filament non-pli\'e dans $\R^3$ correspondent \`a
l'intervalle $[-2;0]$.

\begin{figure}
\begin{center}
\begin{tabular}{c@{}c}
     \includegraphics[width=6.5cm]{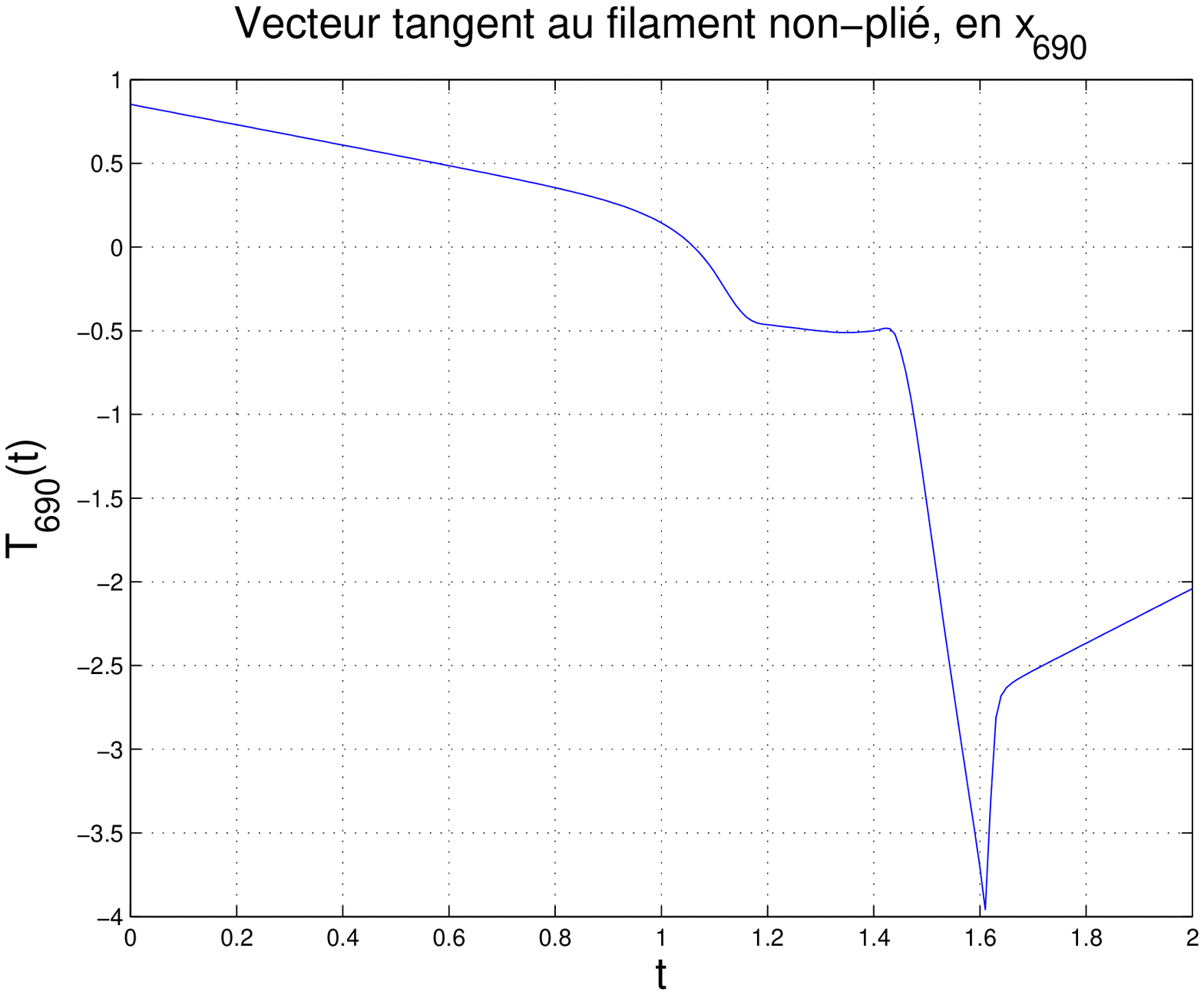}
     &
     \includegraphics[width=6.5cm]{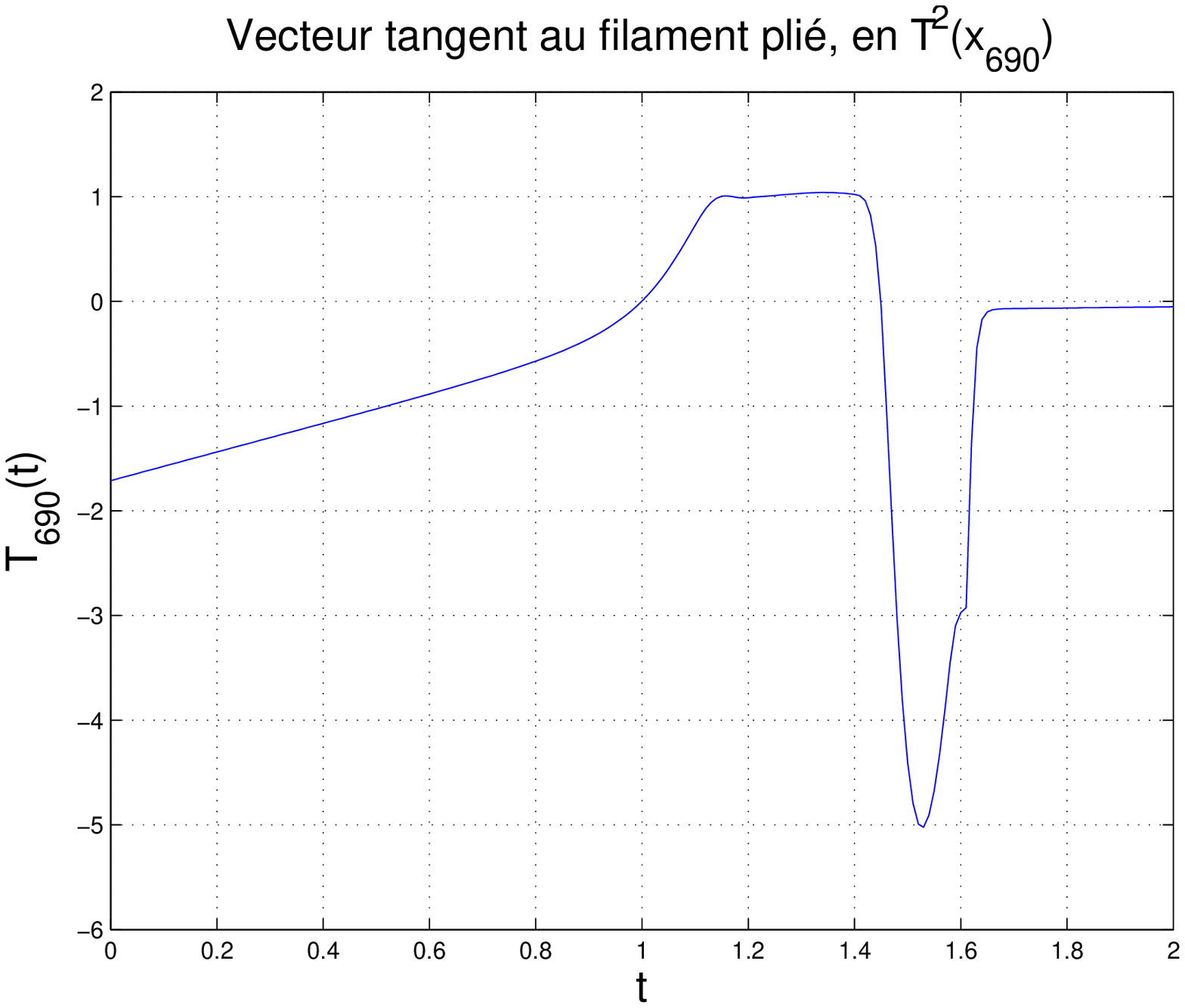}
          \\
     (a) en $x_{690}$ ($t=0$) & (b) en $T^2(x_{690})$ ($t=2$)
\end{tabular}
\caption{\label{fig:pli_tangent} Vecteurs tangents au filament.}
\end{center}
\end{figure}

Il est \'egalement int\'eressant de consid\'erer la
diff\'erentielle de $T^2$ le long du pli pour mieux comprendre ce
qui se passe dans cette r\'egion de l'attracteur. Les vecteurs
tangents (figure~\ref{fig:pli_tangent}) nous indiqueront alors
comment $T^2$ plie le filament lui-m\^eme.

\begin{figure}
\begin{center}
     \includegraphics[height=7cm]{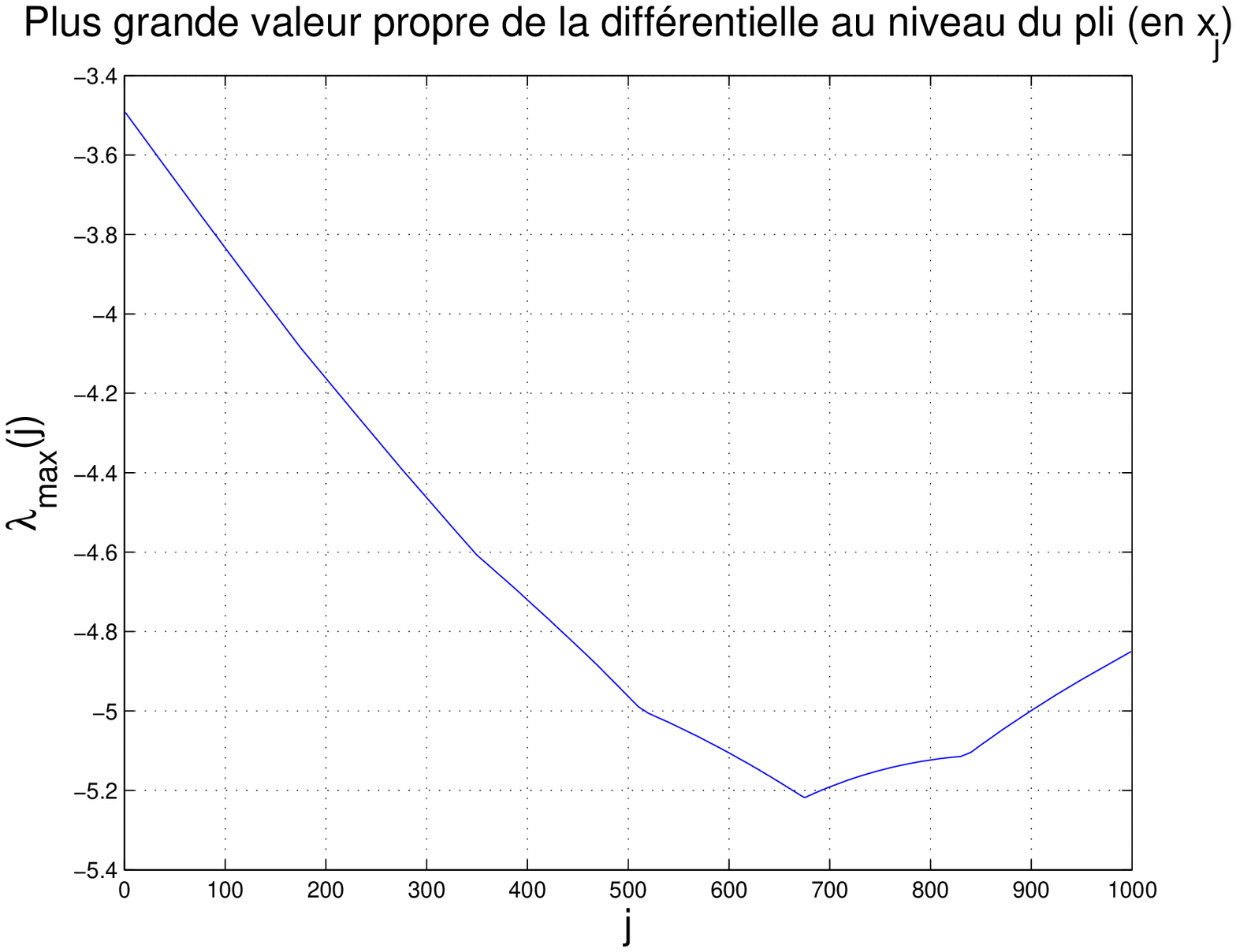}
\caption{\label{fig:pli_diff_avant_valp} Diff\'erentielle de $T^2$
sur le filament non-pli\'e : premi\`ere valeur propre.}
\end{center}
\end{figure}

\subparagraph{Filament non encore pli\'e} \`A $t=0$ (\latin{i.e.}
sur le filament non-pli\'e $x_1 , \ldots x_{1000}$), la plus
grande valeur propre de la diff\'erentielle\footnote{On a
calcul\'e ici les valeurs et vecteurs propres de la
diff\'erentielle $D$ de $T^2$, et non ceux de $\trans{D}D$, qui
permettent d'\'evaluer exactement les directions contract\'ees ou
dilat\'ees, dans la mesure o\`u
$\prodscal{x}{\trans{D}D(x)}=\norm{D(x)}^2$. Ces deux r\'esultats
ont cependant de bonnes chances d'\^etre tr\`es semblables.} de
$T^2$ est n\'egative et comprise entre $-5 \virg 2$ et $-3 \virg
5$ (figure~\ref{fig:pli_diff_avant_valp}). Elle atteint un minimum
en $x_{675}$, c'est-\`a-dire \`a proximit\'e du futur maximum de
courbure $x_{690}$.

La seconde valeur propre est tr\`es \'eloign\'ee de $1$, et reste
comprise entre $0 \virg 04$ et $0 \virg 07$. Le premier vecteur
propre correspond donc \`a la direction instable. Il est
repr\'esent\'e en diff\'erents points du filament sur la
figure~\ref{fig:pli_diff_avant_vp1}. On peut le d\'ecomposer en
trois parties : (1) $0 \leq t \leq 1 \virg 4$ ; (2) $1 \virg 4
\leq t \leq 1 \virg 6$ ; (3) $1 \virg 6 \leq t \leq 2$. Il faut
tenir compte de la normalisation du vecteur propre pour
interpr\'eter correctement son \'evolution le long du pli. La
premi\`ere partie varie tr\`es peu le long du filament, ce qui
confirme l'hypoth\`ese de formation du pli \`a $t \approx 1 \virg
4$ que l'on a formul\'ee d'apr\`es la
figure~\ref{fig:courbure_maxi}a. La troisi\`eme partie est
rectiligne et remonte en m\^eme temps que le pic de la deuxi\`eme
partie. Elle traverse 0 en $x_{686}$, \`a l'instant pr\'ecis o\`u
le pic atteint son maximum. Le pliage semble ainsi correspondre
\`a une forte expansion localis\'ee en $t=1 \virg 5$,
accompagn\'ee d'une absence d'expansion (ou plut\^ot une
contraction puisque les autres vecteurs propres ont de tr\`es
petites valeurs propres) sur l'intervalle $[1 \virg 6 ; 2]$. Ces
deux ph\'enom\`enes sont difficiles \`a distinguer dans la mesure
o\`u la normalisation $L^1$ du vecteur propre entra\^ine une
compensation entre la hauteur du pic et celle de la partie
rectiligne. On reconna\^it tout de m\^eme une forme de pliage en
nous limitant \`a la partie (3) du vecteur propre, l'expansion
ayant lieu dans des sens oppos\'es de part et d'autre du pli.

On remarque \'egalement que la direction instable \`a $t=0$
(figure~\ref{fig:pli_diff_avant_vp1}) est semblable au vecteur
tangent \`a $t=2$ (figure~\ref{fig:pli_tangent}b). Cela signifie
que le filament est suffisamment transverse \`a la vari\'et\'e
instable en $x_{690}$.

\begin{figure}
\begin{center}
     \includegraphics[height=7cm]{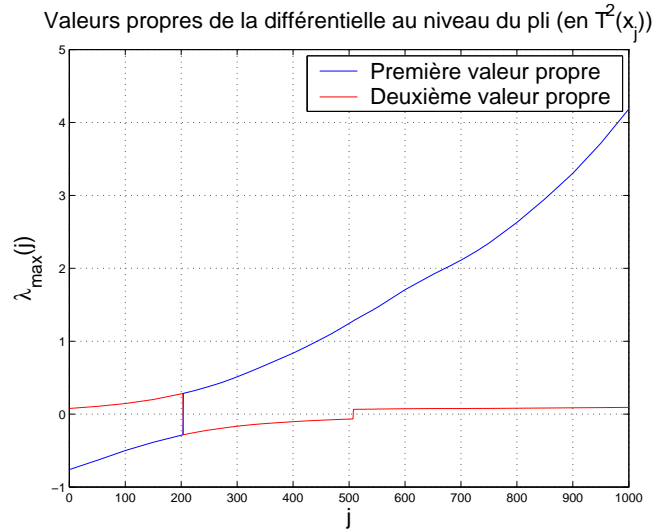}
\caption{\label{fig:pli_diff_valp2} Diff\'erentielle de $T^2$ sur
le filament pli\'e : premi\`eres valeurs propres.}
\end{center}
\end{figure}

Les quelques irr\'egularit\'es observ\'ees sur les courbes des
figure~\ref{fig:pli_diff_avant_valp} et \ref{fig:pli_diff_valp2}
semblent confirmer l'hypoth\`ese avanc\'ee lors de l'\'etude de la
courbure du pli : la ligne bris\'ee consid\'er\'ee doit <<sauter>>
d'un filament \`a l'autre, deux filaments c\^ote \`a c\^ote
n'ayant pas pr\'ecis\'ement la m\^eme courbure.

\subparagraph{Filament pli\'e} \`A $t=2$ (\latin{i.e.} sur le
filament pli\'e $T^2(x_1) , \ldots T^2(x_{1000})$), les deux plus
grandes valeurs propres de la diff\'erentielle de $T^2$ sont
repr\'esent\'ees figure~\ref{fig:pli_diff_valp2}. Elle est
l\'eg\`erement sup\'erieure \`a $2$ au voisinage de
$T^2(x_{690})$, mais il n'y pas d'expansion tout au long du
filament. Celle-ci est en effet de module inf\'erieur \`a 1 entre
$T^2(x_1)$ et $T^2(x_{444})$ o\`u elle traverse le cercle unit\'e,
puis cro\^it lorsqu'on se d\'eplace le long du filament pli\'e
vers $x_{1000}$ o\`u elle d\'epasse 4. Il n'y a donc pas uniforme
hyperbolicit\'e sur l'attracteur.

\begin{figure}
\begin{center}
     \includegraphics[height=7cm]{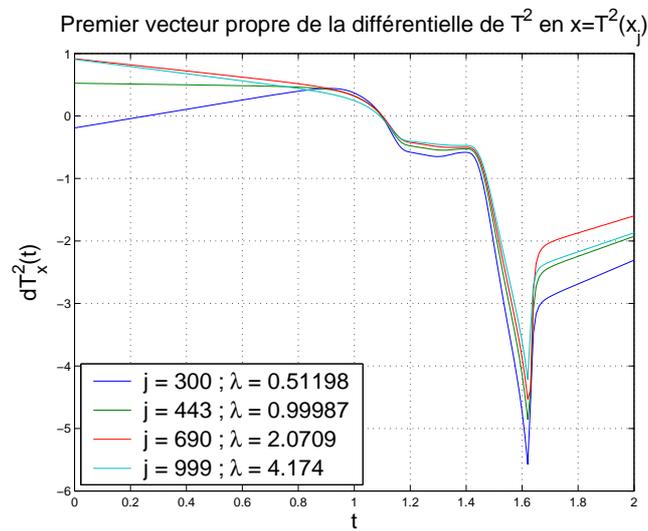}
\caption{\label{fig:pli_diff_vp1} Diff\'erentielle de $T^2$ sur le
filament pli\'e : premier vecteur propre.}
\end{center}
\end{figure}

Dans le domaine o\`u les valeurs propres sont toutes de module
strictement inf\'erieur \`a 1, la courbe $\lambda_{max}(j)$ n'est
pas continue. En tra\c{c}ant aussi la seconde valeur propre, on
constate qu'il s'agit d'un \'echange dans l'ordre des deux
premiers vecteurs propres. Tant que la premi\`ere valeur propre
n'est pas suffisamment grande devant la seconde, le premier
vecteur propre ne suffit pas \`a d\'ecrire la diff\'erentielle de
$T^2$. Il est cependant int\'eressant de le consid\'erer autour de
$T^2(x_{690})$ (figure~\ref{fig:pli_diff_vp1}). Il a alors le
m\^eme aspect qu'\`a la figure~\ref{fig:pli_diff_avant_vp1}, avec
un pic autour de $t=1 \virg 6$.


\begin{figure}
\begin{center}
\begin{tabular}{c@{}c}
     \includegraphics[width=6.5cm]{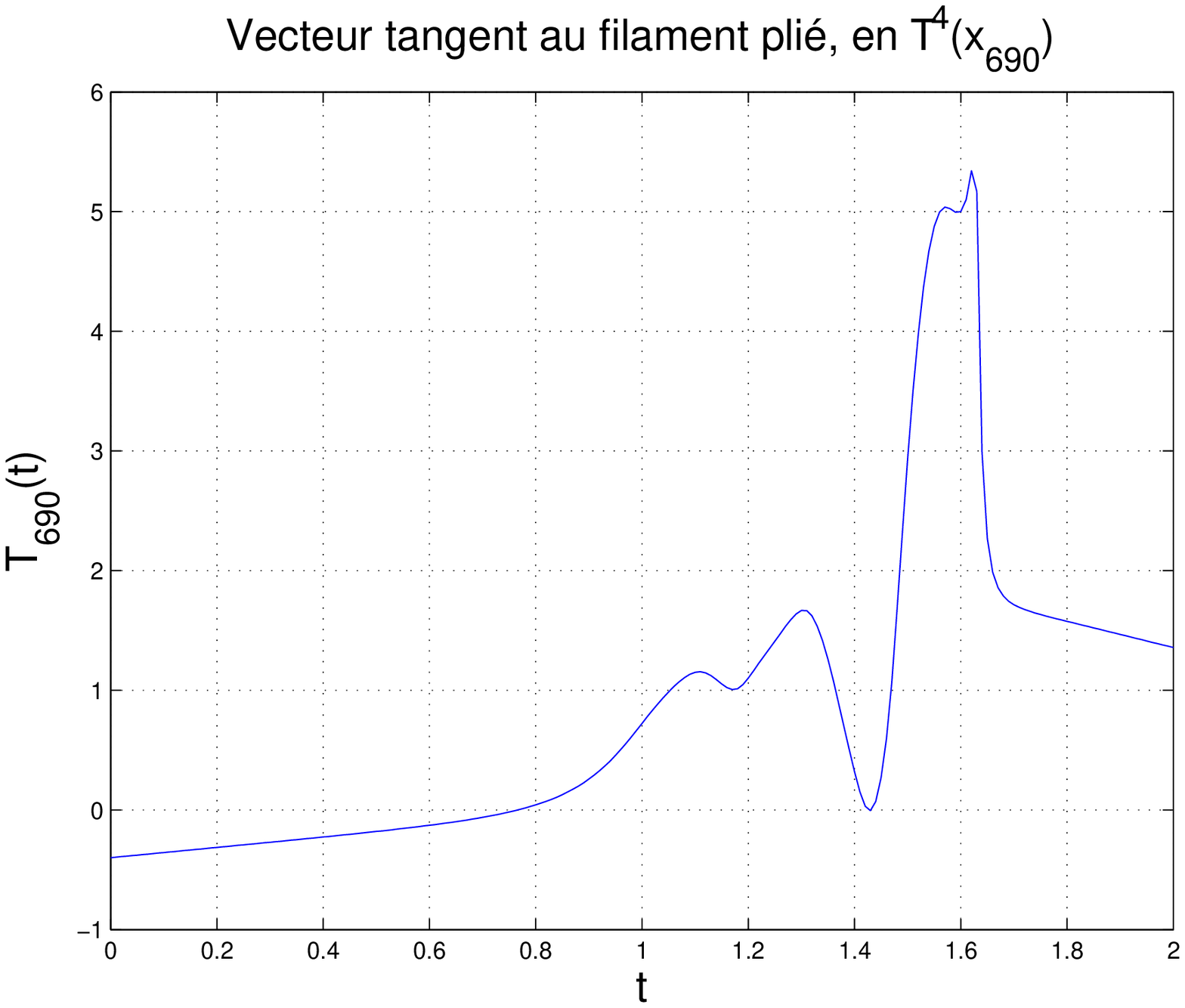}
     &
     \includegraphics[width=6.5cm]{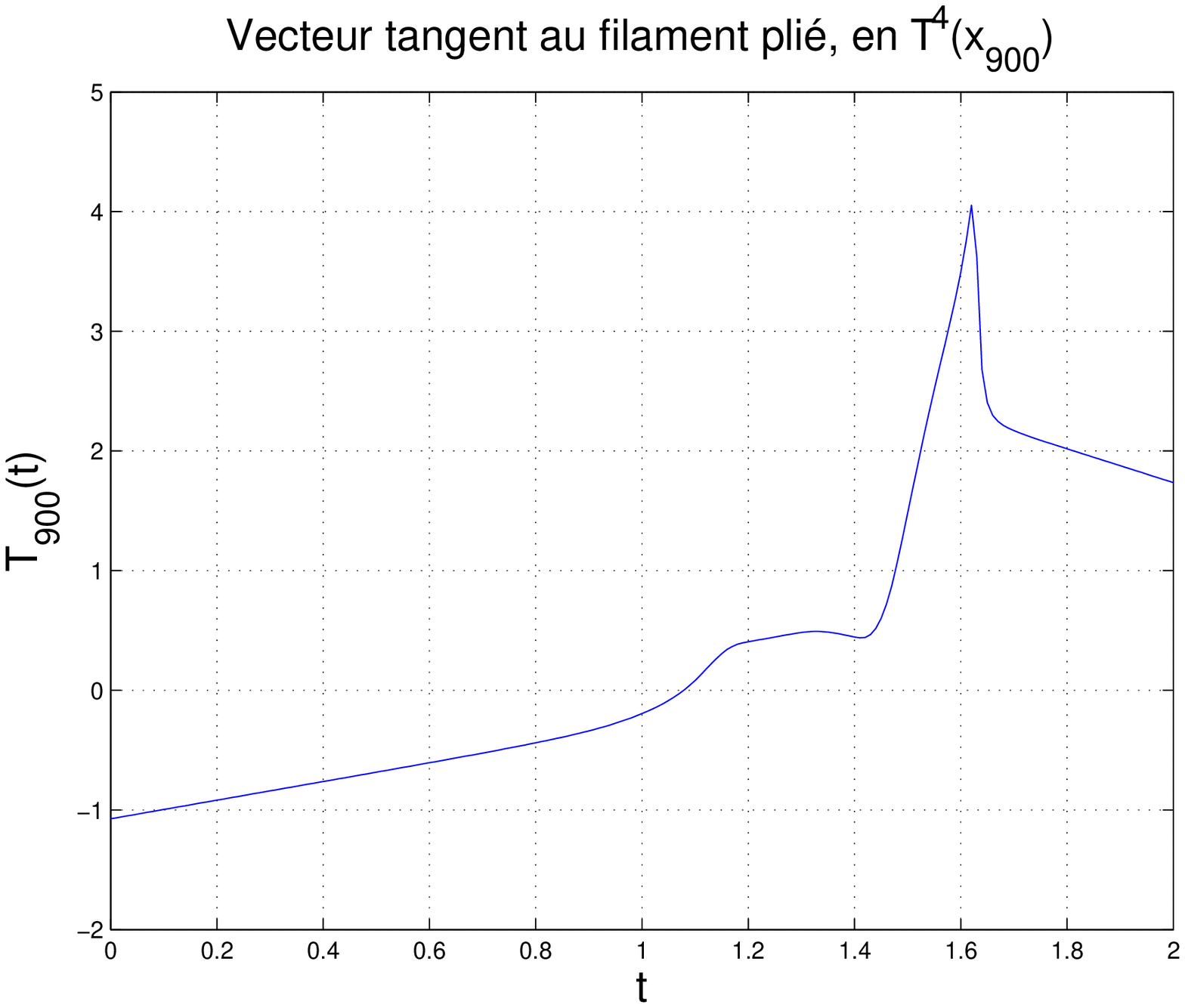}
          \\
     (a) en $T^4(x_{690})$ & (b) en $T^4(x_{900})$
\end{tabular}
\caption{\label{fig:pli_tangent_apres} Vecteurs tangents au
filament apr\`es pliage ($t=4$).}
\end{center}
\end{figure}

Compte-tenu du vecteur tangent au temps $t=4$
(figure~\ref{fig:pli_tangent_apres}a), le filament pli\'e est
faiblement transverse \`a la vari\'et\'e stable en $T^2(x_{690})$.
Apr\`es le pliage, il semble donc que la principale action de
$T^2$ soit de renforcer le pli, par une forte contraction,
l'expansion \'etant tr\`es faible dans la direction tangente au
filament en $T^2(x_{690})$. Il y a par ailleurs assez peu de
diff\'erences entre les diff\'erentes directions instables
repr\'esent\'ees figure~\ref{fig:pli_diff_vp1}. Le pli \'etant
d\'ej\`a form\'e, il y a expansion dans la direction tangente au
filament loin de $T^2(x_{690})$. En effet, le vecteur tangent en
$T^2(x_{900})$ (figure~\ref{fig:pli_tangent_apres}b) co\"incide
avec la direction instable au temps $t=2$. C'est \'egalement le
cas du vecteur tangent en $T^2(x_{400})$, qui est pr\'ecis\'ement
oppos\'e au vecteur tangent en $T^2(x_{900})$.

\paragraph{Visualisation anim\'ee de l'attracteur} Une derni\`ere
fa\c{c}on de comprendre la formation du pli (et la dynamique
g\'en\'erale de l'attracteur) est d'utiliser une
animation\footnote{Il s'agit du fichier
\fichier{film\_delta.avi}.} repr\'esentant la projection dans
$\R^3$ de l'attracteur avec origine des temps $\delta_t$ pour des
valeurs successives de $\delta_t$ (voir section
\ref{sec:origine_temps}).

\begin{figure}
\begin{center}
\begin{tabular}{c@{}c}
     \includegraphics[width=8cm]{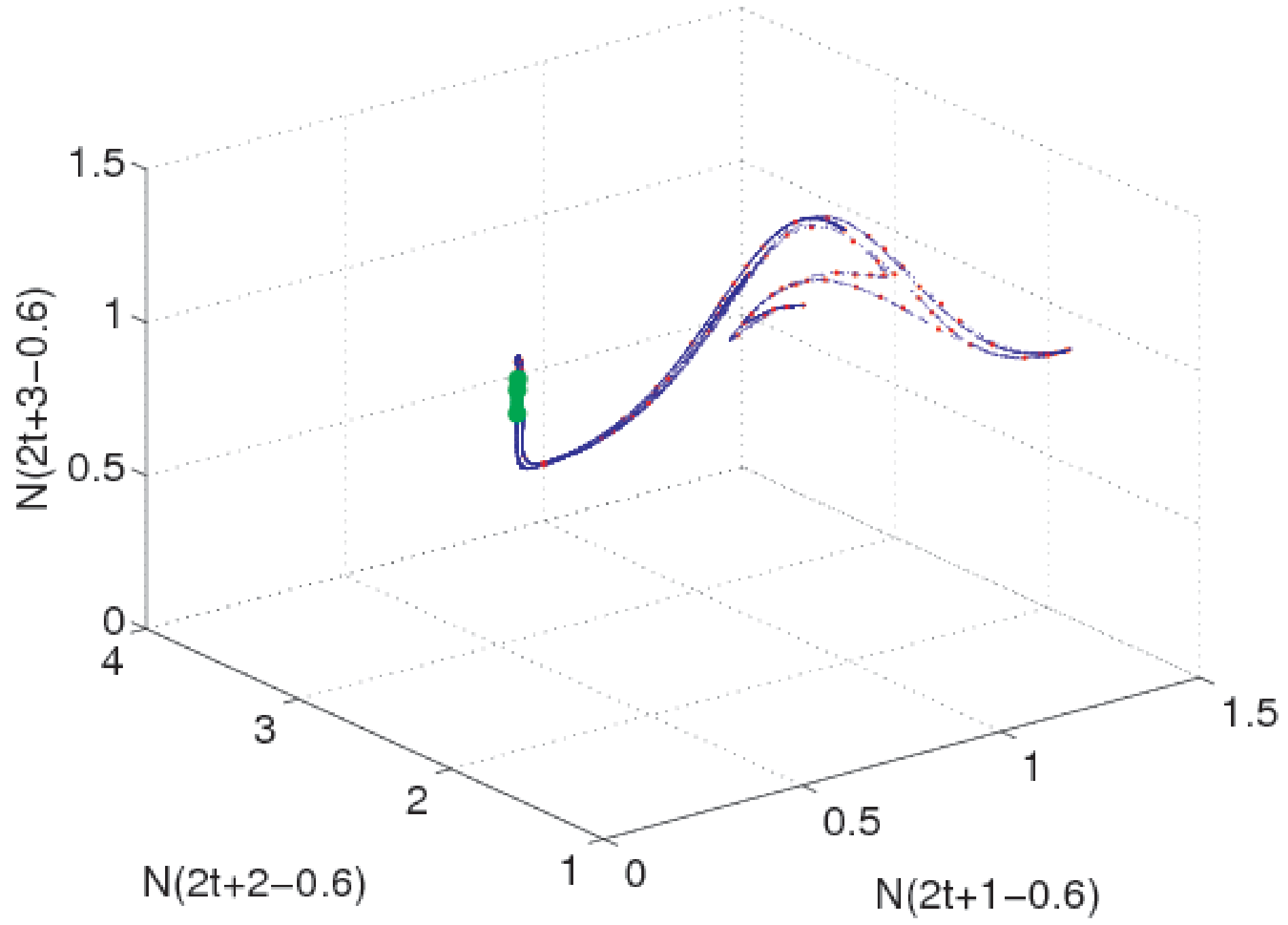}
     &
     \includegraphics[width=8cm]{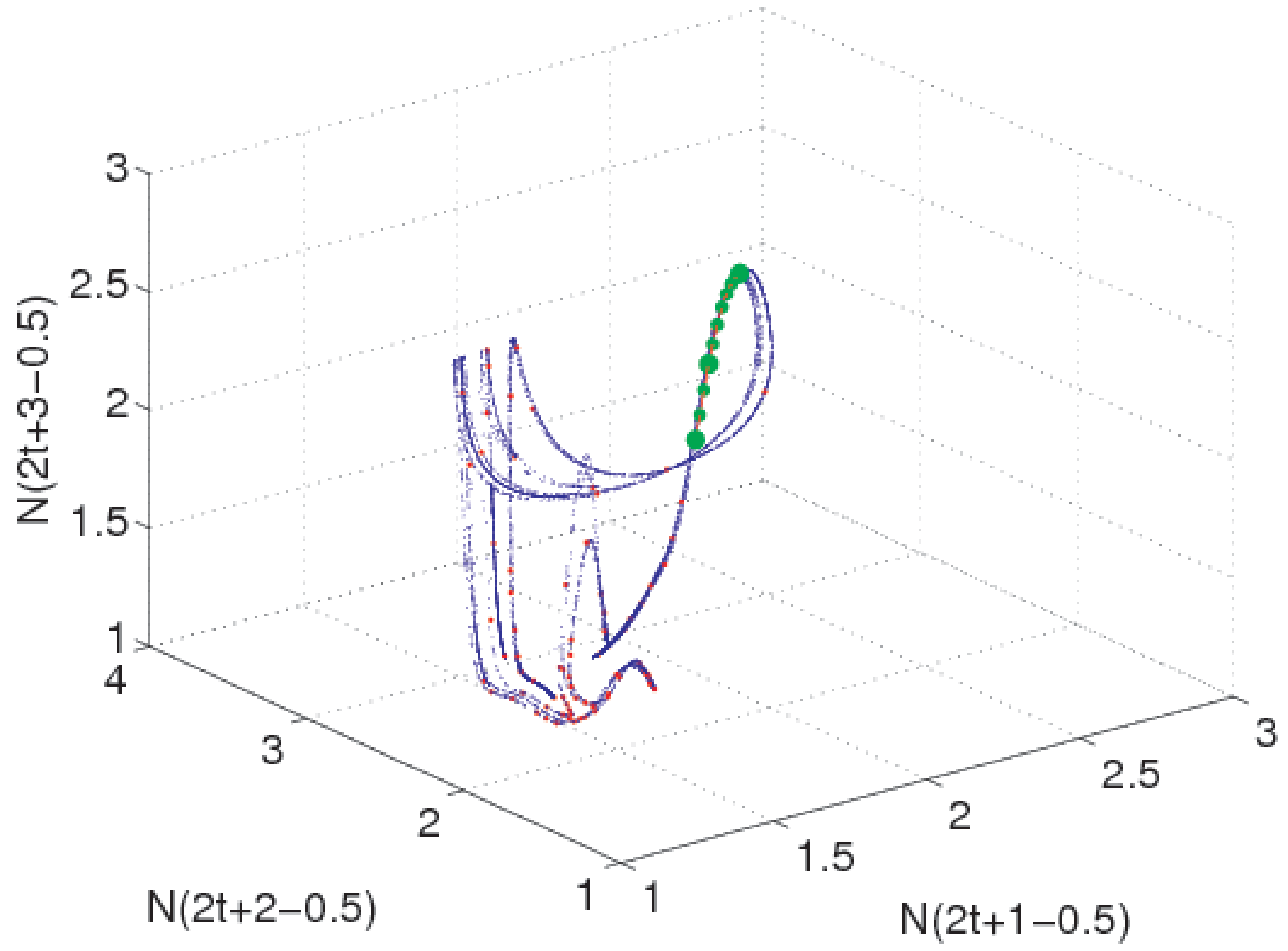}
          \\
     (a) $\delta_t = -0 \virg 6$ & (b) $\delta_t = -0 \virg 5$ \\
     $ $ & $ $ \\
     \includegraphics[width=8cm]{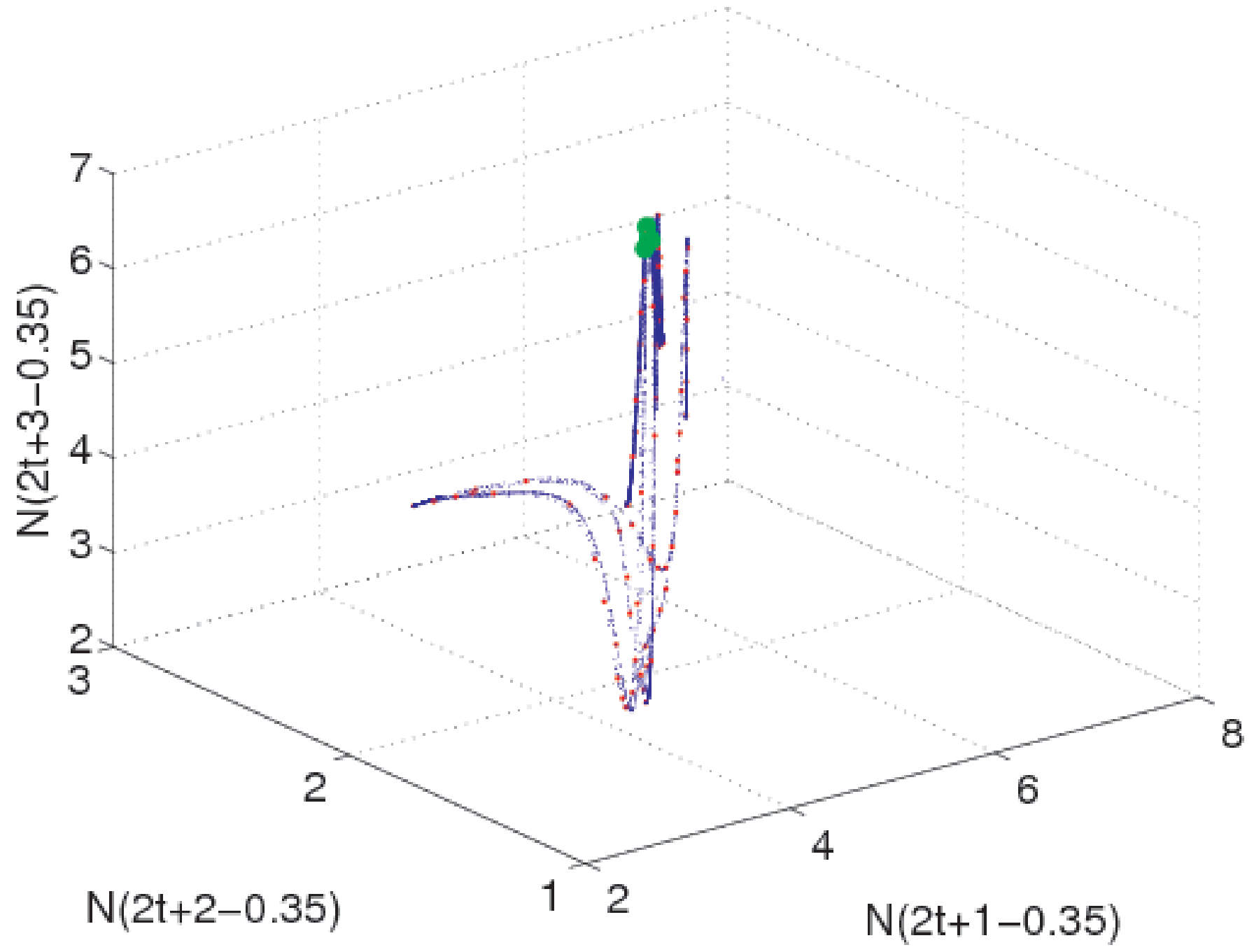}
     &
     \includegraphics[width=8cm]{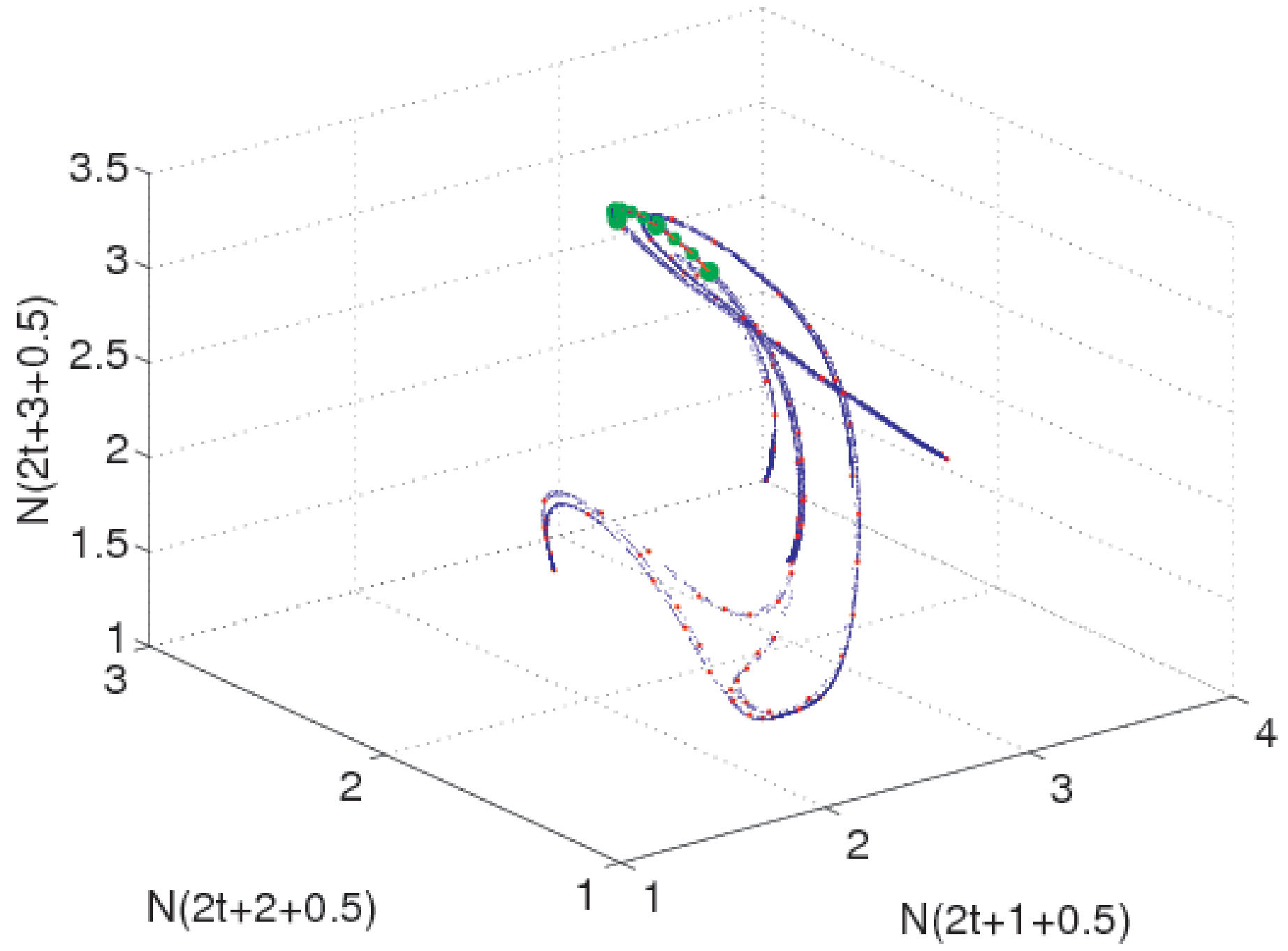}
          \\
     (c) $\delta_t = -0 \virg 35$ & (d) $\delta_t = +0 \virg 5$
\end{tabular}
\caption{\label{fig:pli_delta} L'attracteur et le pli pour
diff\'erentes valeurs de $\delta_t$.}
\end{center}
\end{figure}

On constate pour $-0\virg 6 \leq \delta \leq -0 \virg 35$ qu'une
portion jusque l\`a quasiment rectiligne se plie litt\'eralement
au cours d'un rapide d\'eplacement dans $\R^3$ de la zone
<<inf\'erieure>> (deux coordonn\'ees proches de 1, la troisi\`eme
grande) vers la zone <<sup\'erieure>> (une coordonn\'ee proche de
1, les deux autres grandes). Le pli est inexistant pour $\delta_t
= - 0 \virg 6$ (figure~\ref{fig:pli_delta}a), quasiment form\'e
pour $\delta_t = -0\virg 5$ mais pas encore plac\'e en $x_{690}$
(figure~\ref{fig:pli_delta}b) et l'est totalement \`a $\delta_t =
-0 \virg 35$ (\latin{i.e.} $t=1\virg 65$ avec les conventions de
cette section : voir figure~\ref{fig:pli_delta}c).

La suite de la d\'eformation, pour $\delta_t > -0\virg 4$, ne fait
qu'accentuer ce pli en \'etirant l'attracteur dans une direction
et en le contractant dans les autres (au voisinage du pli). Cette
contraction se ressent plus particuli\`erement autour de $\delta_t
= 0 \virg 5$ (\latin{i.e.} $t \approx 2 \virg 5$ :
figure~\ref{fig:pli_delta}d).

\subsubsection{Non-hyperbolicit\'e de
l'attracteur}\label{sec:hyperbol} 
\begin{figure}
\begin{center}
\includegraphics[width=\textwidth]{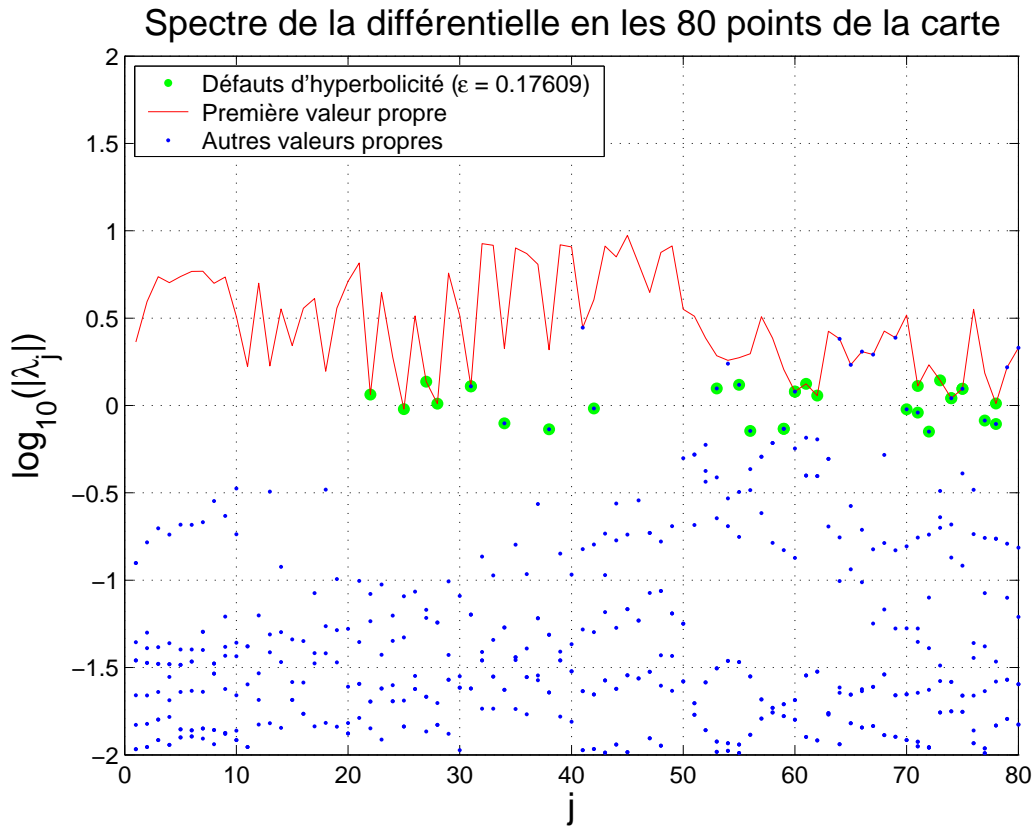}
\caption{\label{fig:hyperbol} Spectre de la diff\'erentielle de
$T^2$ en les 80 points de la carte.}
\end{center}
\end{figure}

La figure~\ref{fig:hyperbol}
repr\'esente le spectre de la diff\'erentielle de $T^2$ en chacun
des 80 points de la carte (en mettant la structure $L^2$ canonique
sur chaque espace tangent). Cela nous permet de caract\'eriser
l'action de $T^2$ sur les diff\'erentes zones de l'attracteur.

\begin{figure}
\begin{center}
    \includegraphics[width=\textwidth]{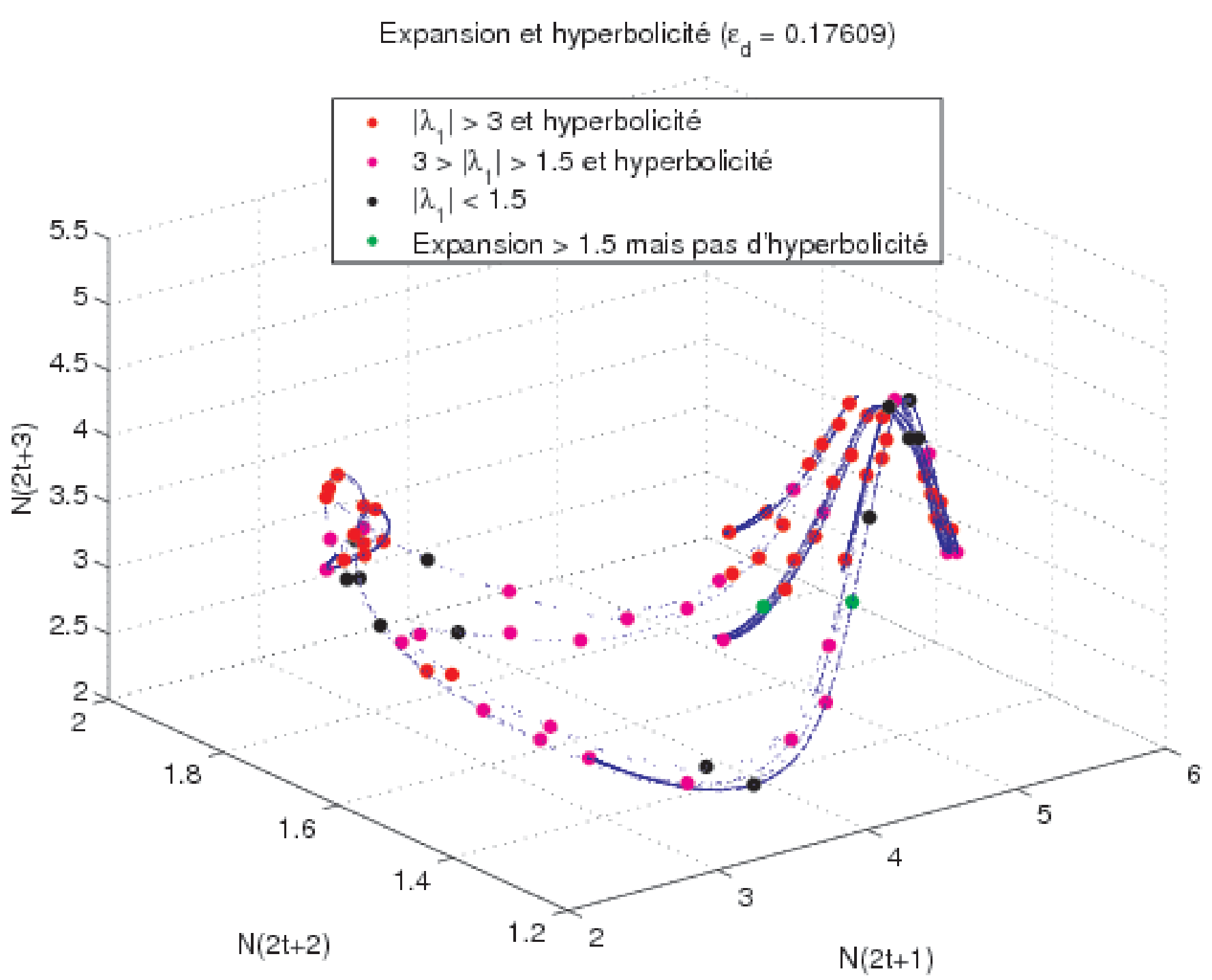}
     \caption{\label{fig:hyperbol_3d}Spectre de la diff\'erentielle selon les zones de l'attracteur.}
\end{center}
\end{figure}

La figure~\ref{fig:hyperbol_3d} permet de situer dans $\R^3$ les
zones d'expansion (dans la ou les direction(s) instable(s)) et de
<<non-hyperbolicit\'e>>\footnote{non-hyperbolicit\'e avec la
structure $L^2$ sur chaque espace tangent. Il est possible qu'une
autre structure rende l'attracteur hyperbolique.}.

L'expansion est assez forte dans l'ensemble, \`a l'exception de
trois r\'egions : le pli (25, 28, 22, 27), le creux (61, 62) et le
triangle (78 ; 65 ; 75, 72, 73). Ces trois zones semblent ainsi
jouer un r\^ole particulier dans la dynamique. Nous avons d\'ej\`a
commenc\'e l'\'etude du pli, il faudra \'egalement s'int\'eresser
aux deux autres zones.

Le <<triangle>> poss\`ede sans aucun doute une orbite de p\'eriode
3 pour $T^2$. Il faudrait \'evaluer plus pr\'ecis\'ement
l'hyperbolicit\'e de cette orbite, qui joue un grand r\^ole dans
la dynamique globale.

Le <<creux>> semble quant \`a lui \^etre une autre zone de pli
dans l'attracteur. Il faudrait d\'eterminer si c'est effectivement
le cas, et si ce pli est ou non distinct du pli que nous avons
d\'ej\`a mis au jour.

Enfin, deux points pr\'esentent \`a la fois une forte expansion et
un d\'efaut d'hyperbolicit\'e (c'est-\`a-dire une seconde valeur
propre proche de 1) : 42 et 70. Leurs images par $T^2$ \'etant
proche du point 62, il est possible que leur particularit\'e soit
simplement li\'ee au <<creux>>.

\subsubsection{Stabilit\'e de la structure dans l'espace des
param\`etres} Apr\`es avoir d\'etaill\'e certains \'el\'ements de
la structure d'un seul attracteur, celui que nous avons observ\'e
pour $(A_0,\rho,\gamma) = (0\virg 15; 0\virg 30; 8\virg 25)$, on
peut s'interroger sur la persistance de cette structure quand on
fait varier l\'eg\`erement les param\`etres. Nous pouvons
d'ores-et-d\'ej\`a esquisser une r\'eponse, visuelle, \`a l'aide
d'une animation\footnote{Voir les fichiers
\fichier{film\_A0\_100\_\_200\_300\_10\_50\_08250\_1\_2coul\_0.1\_0.3.avi}
pour une vue d'ensemble, et
\fichier{film\_zoomA0\_100\_\_200\_300\_10\_50\_08250\_1.avi} pour
un zoom sur la r\'egion pli\'ee (zone 3).} repr\'esentant l'attracteur
$(A_0; 0\virg 30; 8\virg 25)$ dans $\R^3$ quand $A_0$ varie (voir
aussi la section \ref{sec:explor_A0_0.30_8.25}, consacr\'ee \`a
cette exploration, notamment le diagramme de bifurcation
\ref{diag:A0_0.30_8.25}).

Pour $0\virg 135 \leq A_0 \leq 0 \virg 160$, l'attracteur grandit
contin\^ument en deux parties, \`a partir d'une forme tr\`es
simple (uniquement la r\'egion de l'\'equilibre et du pli,
\latin{i.e.} zones 2 et 3). Pour $A_0 \approx 0 \virg 139$, les
autres r\'egions apparaissent brusquement, dans une configuration
assez similaire \`a $A_0 = 0 \virg 15$. La croissance de
l'attracteur se fait par l'allongement de filaments s'achevant
pour une pointe, un peu de la m\^eme fa\c{c}on que lorsqu'on a
observ\'e comment la vari\'et\'e instable s'\'etend \`a
l'int\'erieur de l'attracteur (voir section
\ref{sec:variete_instable}). L'essentiel de la structure est
conserv\'e, m\^eme si la complexit\'e de la dynamique va
apparemment en s'accroissant.

En $A_0 = 0 \virg 16$, les deux composantes connexes se rejoignent
pour n'en former qu'une, qui continue \`a grandir de la m\^eme
fa\c{c}on. Mises \`a part quelques interruptions (une ou deux
fen\^etres de p\'eriodicit\'e), l'attracteur reste intact (tout en
se complexifiant au fur-et-\`a-mesure de sa croissance) jusqu'\`a
$A_0 \approx 0 \virg 19$. Il dispara\^it alors par brusques
paliers pour laisser place \`a une grande fen\^etre de
p\'eriodicit\'e, puis un cycle se refermant sur un point fixe. On
ne le retrouve plus trace ensuite de cet attracteur.

Un zoom sur la r\'egion du pli (zone 3) nous a permis de
d\'eterminer si la pointe o\`u se situent les points 29 et 30
continue \`a grandir en se pliant, parall\`element aux filaments
des points 21 \`a 28. Contrairement \`a ce que l'on aurait pu
penser, ce n'est pas le cas : la double pointe 29--30 reste d'un
seul c\^ot\'e du pli sans le franchir. Il n'y a donc pas de
complexification de la structure de cette mani\`ere-l\`a.
\section{Perspectives}
Il reste beaucoup de travail \`a faire pour comprendre ce
mod\`ele. \`A la lumi\`ere des simulations dont nous venons
d'analyser les r\'esultats, nous pouvons cependant d\'ej\`a
esquisser quelques pistes afin de poursuivre ce travail.

\subsection{Simulations num\'eriques}
Tout d'abord, plusieurs aspects restent \`a \'etudier
num\'eriquement afin de comprendre au mieux le comportement du
mod\`ele. Certaines questions \'etaient d\'ej\`a \'evoqu\'es dans
\cite{Yoccoz:toymodel}, d'autres sont venues \`a la lueur des
nouvelles simulations num\'eriques.

\paragraph{Explorations}
\begin{itemize}
\item Lors de l'augmentation de $A_0$ dans le diagramme $(A_0,
0.30, 8.25)$, les deux filaments se rejoignent-ils
pr\'ecis\'ements ou bien se rapprochent-ils suffisamments pour
induire le m\'elange ?
\item Dans le cas d'une exploration particuli\`ere, chercher
le plus grand nombre possible d'attracteurs, notamment en
<<suivant>> automatiquement chacun jusqu'\`a la perte de
stabilit\'e, et en utilisant plusieurs conditions initiales. Ce
travail pourrait d\'ej\`a \^etre fait dans le cas des explorations
d\'ej\`a faites, notamment pour \'eclaircir la question des
discontinuit\'es du diagramme.
\item Exploration plus exhaustive de l'espace des param\`etres.
\item Tester d'autres transitions de l'\'equilibre stable vers la
dynamique chaotique. \item D\'eterminer des r\'egions de l'espace
des param\`etres (en 2 ou 3 dimensions) o\`u se produisent les
divers comportements observ\'es (orbite p\'eriodique attractive,
cycle attractif, attracteur de type H\'enon et autres attracteurs
\'etranges), et les bifurcations qui se produisent \`a l'interface
entre ces diff\'erentes zones. Pour chaque comportement observ\'e,
on pourrait chercher \`a caract\'eriser un peu plus
pr\'ecis\'ement la dynamique dans un cas particulier au moins
(notamment pour les cycles et les attracteurs de type H\'enon).
\item D\'eterminer par un calcul num\'erique la nature des bifurcations
observ\'ees sur le diagramme, en calculant les valeurs propres de
la diff\'erentielle (soit en lin\'earisant l'\'equation, soit
exp\'erimentalement). Ceci serait particuli\`erement utile au
niveau des discontinuit\'es apparentes du diagramme.
\end{itemize}

\paragraph{\'Etude de l'attracteur \'etrange}
\begin{itemize}
\item G\'eom\'etrie de l'attracteur :
\begin{itemize}
\item D\'ecrire g\'eom\'etriquement l'attracteur
avec plus de pr\'ecision, et notamment comment se s\'eparent les
diff\'erents filaments (en se limitant aux plus gros d'entre eux).
Ceci devrait permettre d'\'elaborer un mod\`ele simple pouvant
g\'en\'erer un tel attracteur. Il faudrait en particulier
caract\'eriser les <<embranchements>>.
\item Y a-t-il des zones plus <<denses>> que
d'autres (au sens de la mesure de Lebesgue et non de la mesure
physique) ?
\end{itemize}
\item Dynamique sur l'attracteur : il reste beaucoup de travail
pour comprendre pr\'ecis\'ement l'int\'eraction entre cette
dynamique en trois dimensions et la dynamique en temps continu.
\item Dimension fractale : est-elle la m\^eme dans toutes les r\'egions
de l'attracteur ?
\item Sensibilit\'e aux conditions initiales : calculer les exposants de
Lyapunov.
\item \'Equilibre et orbites p\'eriodiques : \begin{itemize}
\item S'assurer de la correction de l'estimation de la seconde plus
grande valeur propre de la diff\'erentielle \`a l'\'equilibre,
\'eventuellement en lin\'earisant directement l'\'equation.
\item Localiser quelques orbites p\'eriodiques de faible
p\'eriode (en particulier l'orbite de p\'eriode 3 pour $T^2$,
situ\'ee dans la zone 8 et l'orbite de p\'eriode 2 dans les zones
4-5) et y effectuer la m\^eme \'etude que pour l'\'equilibre,
notamment en \'evaluant les valeurs propres de la diff\'erentielle
de $T^m$ et en tra\c{c}ant la vari\'et\'e instable. \item Y a-t-il
d'autres \'equilibres, associ\'es ou non \`a ces orbites, situ\'es
en-dehors de l'attracteur ?
\item Le trac\'e que nous avons fait repr\'esente l'ensemble
$\omega$-limite : o\`u sont situ\'es les points p\'eriodiques ?
sont-ils denses ?
\end{itemize}
\item Hyperbolicit\'e : utiliser les orbites p\'eriodiques
pr\'esentes sur l'attracteur pour la tester avec plus de
g\'en\'eralit\'e.
\item Pli(s) : \begin{itemize} \item Peut-on \'eviter les
discontinuit\'es de la courbure en augmentant encore la
pr\'ecision de localisation du pli ?
\item Y a-t-il un lien entre pics de courbure et d\'efaut
d'injectivit\'e de la projection ? \item Localiser les autres plis
s'il y en a (notamment en \'etudiant le <<creux>> et sa
pr\'eimage, qui se sont distingu\'es dans l'\'etude de
l'hyperbolicit\'e).
\end{itemize}
\item Mesure physique : \'Evaluer la mesure physique sur l'attracteur :
quelles zones sont plus charg\'ees, quelles zones le sont moins ?
\item Essayer de <<suivre>> pr\'ecis\'ement la structure mise en \'evidence sur
cet attracteur lorsque l'on fait varier l\'eg\`erement les
param\`etres (l'\'equilibre, le pli, les valeurs propres de la
diff\'erentielle dans ces r\'egions, etc.). Jusqu'o\`u peut-on la
suivre ?
\end{itemize}

\paragraph{D'autres attracteurs \'etranges ?} Le cas $(0\virg 18; 0\virg
30; 8\virg 25)$ (figure~\ref{fig:0.18_0.30_8.25}) semble plus
complexe mais peut-\^etre encore plus int\'eressant que
l'attracteur que nous avons \'etudi\'e. Parmi les explorations
d\'ej\`a effectu\'ees (ou \`a venir), on pourrait chercher \`a
approfondir l'\'etude de la dynamique de certains objets
remarquables, \`a la lumi\`ere du travail d\'ej\`a effectu\'e.

\subsection{Conjectures}
Au vu des r\'esultats num\'eriques, on peut \'enoncer quelques
conjectures, en vue de rendre rigoureuses les observations
qualitatives que nous venons de faire, et dont nous sommes \`a peu
pr\`es s\^urs.

La dynamique est chaotique pour certaines valeurs des
param\`etres. L'attracteur est \'etrange. La dynamique chaotique
est persistante, au voisinage de cet attracteur \'etrange.

La dimension fractale de l'attracteur \'etrange est comprise
strictement entre 1 et $1\virg 5$.

L'attracteur est <<quasi-hyperbolique>> (cf. propri\'et\'es d'un
syst\`eme hyperbolique et de l'attracteur de H\'enon en
annexe~\ref{annexe:hyperbolique}).

On peut le d\'ecomposer en un nombre fini de parties dans
lesquelles $T^1$ est transitive (<<d\'ecomposition spectrale>>).

Il y a une cascade sous-harmonique lorsqu'on se d\'eplace dans
l'espace des param\`etres \`a partir de $\gamma$ petit, $\rho$
proche de 0 ou 1, $A_0$ petit, pour se diriger vers des valeurs
plus grandes de $\gamma$ ou de $\rho$, ou bien vers des valeurs
interm\'ediaires de $\rho$.

\subsection{Questions biologiques}
Il y a de nombreuses mani\`eres de complexifier le mod\`ele pour
le rendre plus r\'ealiste, par exemple en ne supposant plus que la
survie est ind\'ependante de l'\^age. On peut \'egalement essayer
d'introduire des effets maternels.
\section{Conclusion}
Revenons tout d'abord au probl\`eme purement biologique que nous
nous sommes pos\'e initialement. Il s'est av\'er\'e qu'un mod\`ele
simple combinant maturation et des saisons r\'eguli\`eres peut
engendrer des comportements chaotiques extr\^emement complexes,
pour des valeurs assez raisonnables des param\`etres. La
principale condition, et qui se trouve remplie pour ce qui
concerne \emph{Microtus epiroticus}, est une tr\`es forte
f\'econdit\'e.

Il est clair que ce mod\`ele pourrait difficilement \^etre
pr\'edictif, tant il a \'et\'e simplifi\'e sans se soucier
finement du cycle de vie des populations qui nous int\'eressent.
En revanche, nous pouvons d'ores-et-d\'ej\`a tirer des conclusions
qualitatives, la plus importante \'etant que l'impr\'evisibilit\'e
des effectifs futurs \`a long terme peut avoir lieu dans un
environnement stable, o\`u tous les hivers sont strictement
identiques.

De plus, de nombreuses difficult\'es techniques soulev\'ees lors
de la mise en \oe uvre de simulations sur ce premier mod\`ele
seront tr\`es utiles pour des calculs num\'eriques sur des
mod\`eles plus complexes d\'eriv\'es ou non de celui-ci.

\medskip

D'un point de vue math\'ematique, nous avons eu un aper\c{c}u de
la richesse des comportements que peut g\'en\'erer un mod\`ele
somme toute assez \'el\'ementaire. En nous attardant sur l'un de
ces <<attracteurs \'etranges>>, nous avons pu \'evaluer la
complexit\'e d'un seul de ces comportements, en mettant en
\'evidence des ph\'enom\`enes tr\`es mal compris, voire jamais
abord\'es encore. Nous avons finalement pos\'e beaucoup plus de
questions que nous n'avons apport\'e de r\'eponses. De nombreuses
simulations restent ainsi encore \`a faire, soit pour confirmer
une hypoth\`ese expliquant les r\'esultats obtenus, soit pour nous
aider \`a en formuler au sujet des questions pour lesquelles nous
n'arrivons m\^eme pas \`a esquisser une r\'eponse.

Le travail th\'eorique restant \`a faire est lui aussi immense.
Avec quelques outils \'el\'ementaires, nous avons pu d\'efinir
l'attracteur global du syst\`eme, mais nous n'avons aucune
information \`a son sujet. Les exp\'erimentations num\'eriques
nous ont permis de formuler quelques conjectures, mais ont surtout
pour objectif de nous indiquer par quelles voies il serait
possible de s'attaquer \`a la r\'esolution de celles-ci. Au vu des
quelques figures que nous avons pu tracer, il semble assez
probable que les <<attracteurs \'etranges>> --- s'il s'av\`erent
en \^etre effectivement --- que nous avons observ\'es ont
certainement un grand int\'er\^et dans l'\'etude th\'eorique des
syst\`emes dynamiques non-uniform\'emement hyperboliques.

\clearpage

\pagenumbering{Roman}
\appendix
\section{Simulations : travail pr\'eliminaire} \label{annexe:sim_preliminaire}

\subsection{Discr\'etisation du mod\`ele}
\label{par:modele_discretisation} On veut passer du mod\`ele
continu d\'ecrit par l'\'equation \eqref{eq:modele} \`a un mod\`ele
discret, qui pourra \^etre simul\'e num\'eriquement.

Fixons un entier $p>0$, ce sera le nombre de classes d'\^ages
consid\'er\'ees par ann\'ee. On note $n_i$ le nombre de naissances qui
ont lieu dans l'intervalle de temps
$\left[\frac{i-1}{p};\frac{i}{p}\right[$, $N_i$ l'effectif mature
moyen au cours de ce m\^eme intervalle, $e_i$ le facteur
saisonnier \'egal \`a la moyenne\footnote{En pratique, on a pris la
moyenne des valeurs aux extr\'emit\'es de l'intervalle.} dans cet
intervalle de $m_{\rho,\epsilon}$ (d\'efinie par
\eqref{eq:saison}, $i$ \'etant consid\'er\'e mod. $p$). Soit $s_i$ la
proportion des individus matures et vivants parmi ceux qui sont
n\'es $i$ pas de temps auparavant, et $m$ la fonction de f\'econdit\'e
d\'efinie par les relations \eqref{eq:fecondite}. On a alors les
relations suivantes :

\begin{equation}
\left\{
\begin{aligned}
n_i &= \frac{m(N_i)\times N_i \times e_i}{p} \\
N_i &= \sum_{k=1}^{2p} s_k \times n_{i-k}
\end{aligned} \right. \end{equation}

La condition initiale $\left(n_i\right)_{1 \leq i \leq 2p}$ \'etant
donn\'ee, ceci permet de calculer les $n_i$ pour tout $i > 2p$.

Remarquons que le calcul de $s_i$ n'est pas toujours \'evident.
Pour les calculs, dans le cas o\`u $A_0$ est un multiple entier de
$1/p$ ($A_0 = i_0 / p$), on a pris $s_i = \left( 1-\frac{i}{p
\times A_1} \right) \1_{i \geq i_0}$. Lorsque $A_0$ n'est pas
adapt\'e au pas de discr\'etisation, on a pris pour $s_i$ la
moyenne des $s_{j,q}$ pour $(i-1)\times k +1 \leq j \leq i\times
k$, o\`u $s_{j,q}$ d\'esigne le coefficient de survie calcul\'e
avec un pas $q=p \times k$ (le plus souvent, $k=100$). C'est le
cas du diagramme \ref{diag:A0_0.30_8.25}.

Pour des raisons pratiques de calcul, on a souvent utilis\'e le
vecteur des naissances $n_i$ au lieu des effectifs matures $N_i$,
car cela \'evite de calculer deux fois $n_i$ au cours de la
simulation. Il est bien s\^ur tr\`es simple de passer des
naissances aux effectifs matures, mais en perdant les deux
premi\`eres ann\'ees. Ce choix explique l'apparition des
naissances dans les r\'esultats, alors que celles-ci ne sont pas
explicitement utilis\'ees dans le mod\`ele.

\subsection{Choix des conditions initiales}

Pour choisir une condition initiale <<au hasard>>, on a choisi de
d\'eterminer un vecteur de naissances al\'eatoire. Les naissances
successives sont tir\'ees suivant des lois uniformes
ind\'ependantes. Deux m\'ethodes ont \'et\'e employ\'ees.

\begin{figure}
\begin{center}
\begin{tabular}{c@{}c}
     \includegraphics[width=7cm]{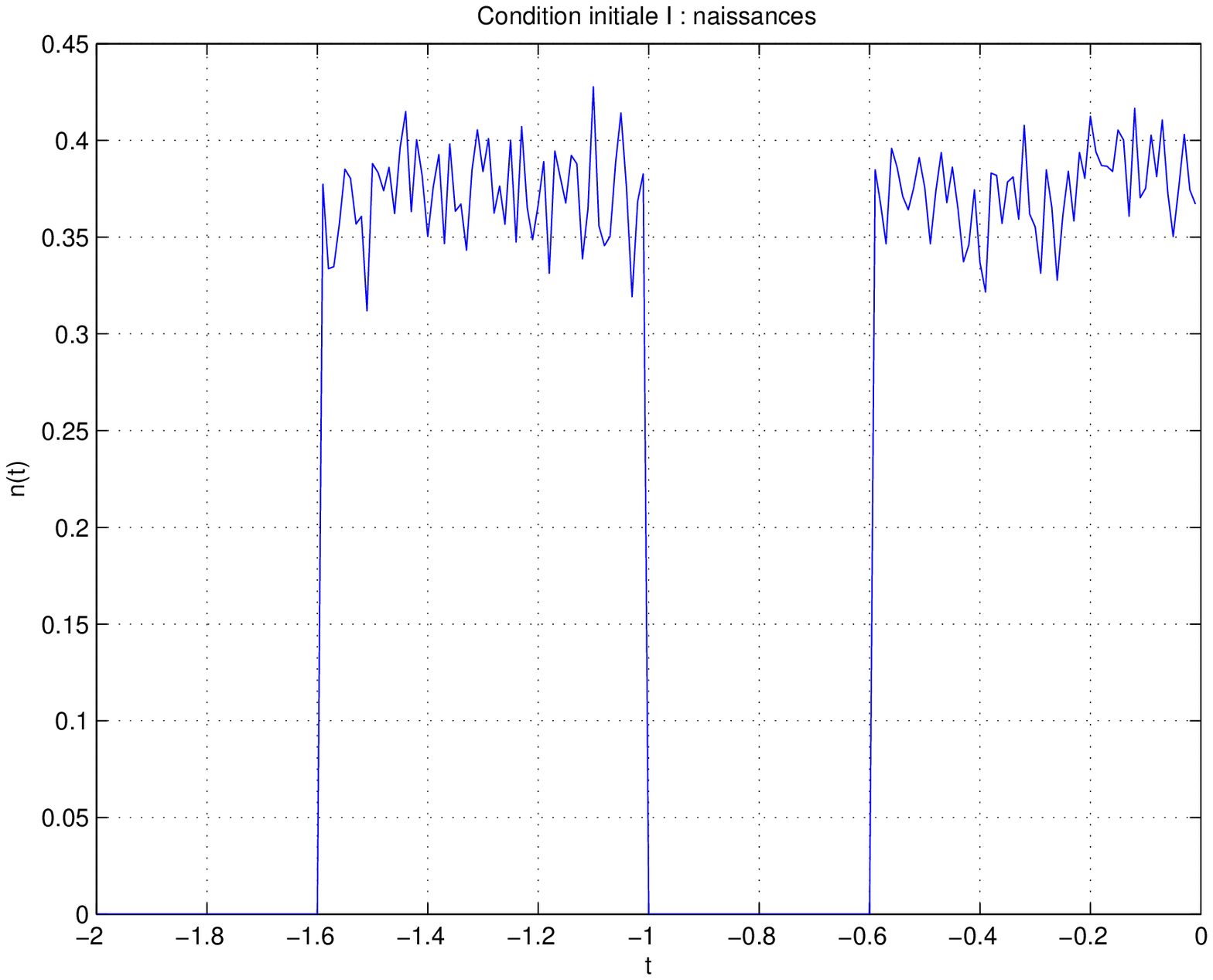}
     &
     \includegraphics[width=7cm]{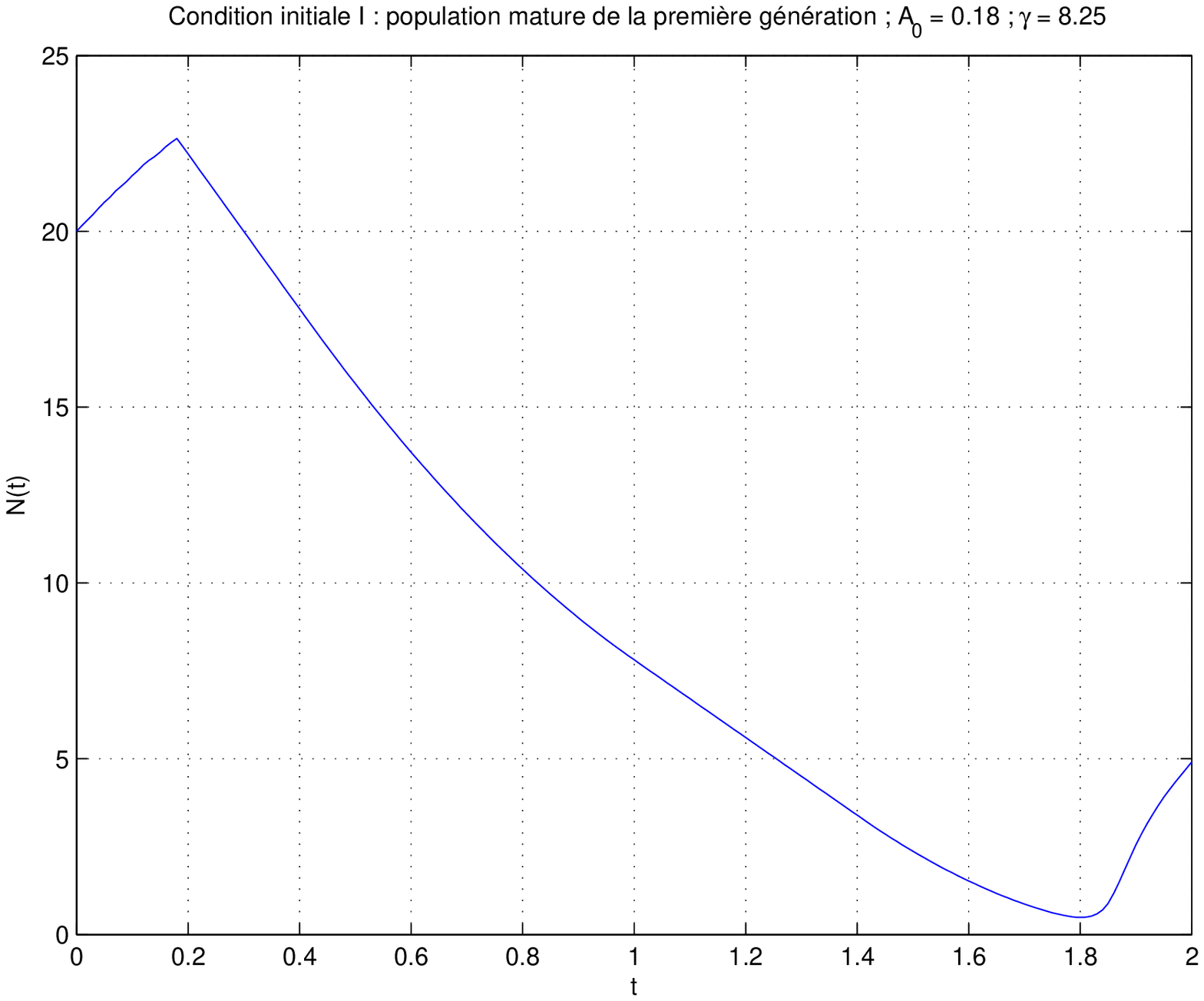}
     \\
     Naissances & Premi\`ere g\'en\'eration
\end{tabular}
\caption{\label{fig:cond_I} Condition initiale (I).}
\end{center}
\end{figure}

Pour la condition initiale (I), on a impos\'e qu'il n'y ait pas de
naissances en hiver (avec $\rho=0\virg 41$), et que l'effectif
mature soit \'egal \`a 20 \`a l'instant $t=0$ (avec $A_0 = 0 \virg
18$, $A_1 = 2$). La figure~\ref{fig:cond_I} repr\'esente les
naissances ($-2 \geq t < 0$) et la premi\`ere g\'en\'eration de
populations matures ($0 \leq t \leq 2$) qui en d\'ecoule (avec
$A_0 = 0 \virg 18$, $A_1 = 2$, $\gamma = 8 \virg 25$).

\begin{figure}
\begin{center}
\begin{tabular}{c@{}c}
     \includegraphics[width=7cm]{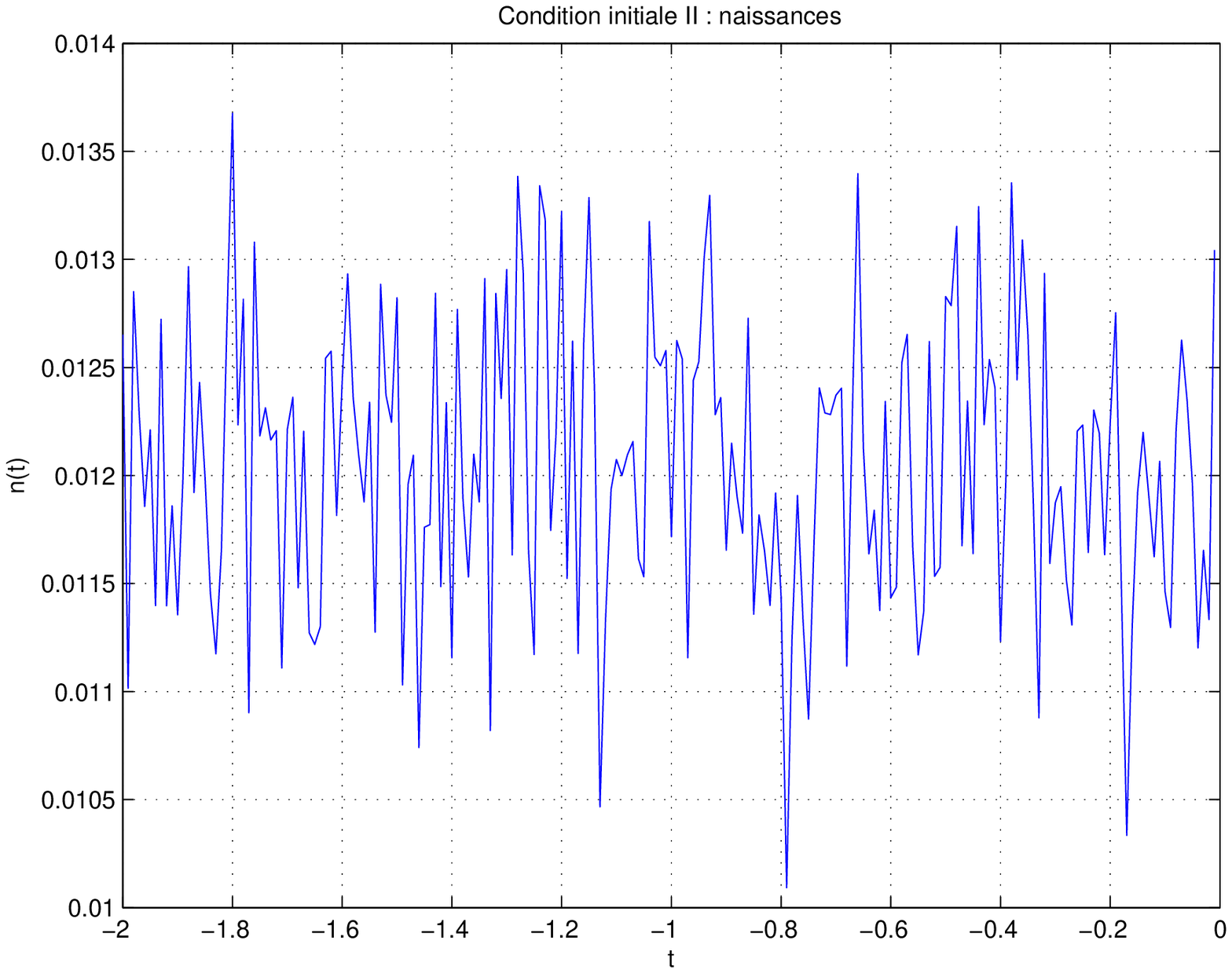}
     &
     \includegraphics[width=7cm]{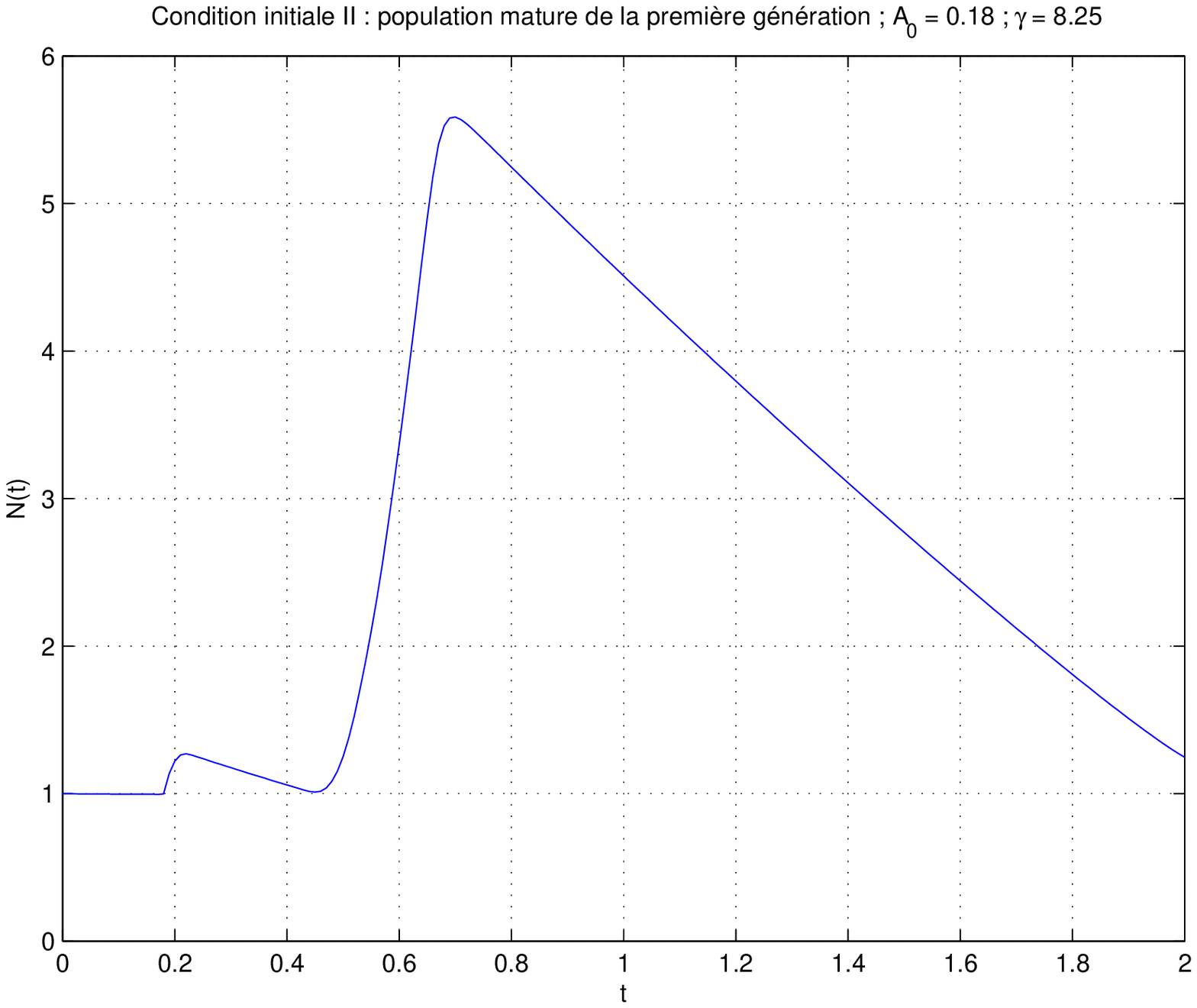}
     \\
     Naissances & Premi\`ere g\'en\'eration
\end{tabular}
\caption{\label{fig:cond_II}Condition initiale (II).}
\end{center}
\end{figure}

Pour la condition initiale (II), qui a \'et\'e utilis\'ee le plus
souvent, on n'a pas tenu compte de l'hiver (le choix d'une valeur
particuli\`ere de $\rho$ n'\'etant ps justifi\'e), et on a
impos\'e un effectif mature \'egal \`a 1 \`a l'instant $t=0$
(cette valeur \'etant plus raisonnable au vu de la dynamique
stationnaire du syst\`eme). La figure~\ref{fig:cond_II}
repr\'esente dans les m\^emes conditions que pr\'ec\'edemment les
naissances et la premi\`ere g\'en\'eration correspondante.

On pourrait bien s\^ur concevoir d'autres m\'ethodes de choix
al\'eatoire d'une condition initiale, \'evitant mieux les biais
possibles, mais ce n'est pas tr\`es important pour l'usage que
nous en avons eu. Cela ne serait utile que dans le cadre d'une
\'etude de la <<taille>> de bassins d'attractions de plusieurs
attracteurs.

\subsection{Choix du pas de discr\'etisation}

Le param\`etre $p$ (nombre de classes d'\^age par ann\'ee,
appel\'e $an$ pour plus de clart\'e) est d\'ecisif pour les
simulations. La complexit\'e de l'algorithme est en effet
proportionnelle au carr\'e de ce pas. Le choix $an=100$ est le
r\'esultat d'un compromis entre rapidit\'e du calcul et
pr\'ecision, apr\`es quelques simulations test. Cette valeur est
de plus raisonnable pour ce mod\`ele : au vu de l'ordre de
grandeur des param\`etres $A_0$ et $\rho$, il ne semble pas utile
d'\^etre plus pr\'ecis. De plus, la situation biologique de
d\'epart \'etant discr\`ete, il est inutile de consid\'erer une
\'echelle de temps inf\'erieure \`a trois jours.

\subsection{Lissage des fonctions}
Le choix des fonctions $m(N)$ et $m_{\rho}(t)$ \'etant assez
arbitraire et peu r\'ealiste, le lissage des fonctions a surtout
\'et\'e fait \`a titre pr\'eventif, pour que d'\'eventuelles
discontinuit\'es ou irr\'egularit\'es dans la dynamique ne soient
pas d\^ues au manque de r\'egularit\'e des fonctions utilis\'ees
dans le mod\`ele. Il a \'et\'e motiv\'e par l'observation d'angles
sur certains attracteurs qui semblaient lisses par ailleurs, au
cours de simulations pr\'eliminaires.

Il semble en r\'ealit\'e que ces changements n'ont pas modifi\'e
qualitativement les propri\'et\'es globales du syst\`eme dans
l'espace des param\`etres. En revanche, il est certain que pour
des valeurs fix\'ees des param\`etres, la plus infime modification
des fonctions utilis\'ees dans le mod\`ele peut modifier
enti\`erement le comportement observ\'e.
\clearpage
\section{Simulations : traitement des donn\'ees} \label{annexe:traitement}

\subsection{Visualisation de l'attracteur en 3 dimensions}
On repr\'esente $(N(t),N(t+1),N(t+2)$ pour les valeurs enti\`eres
de $t \in [t_{\min};t_{\max}]$. Ces valeurs enti\`eres
correspondent pr\'ecis\'ement \`a la fin de l'\'et\'e dans le cas
$\epsilon_{ete}=0$, au milieu de l'automne dans le cas
g\'en\'eral.

\paragraph{Intervalle de temps choisi}
Le plus souvent, on a choisi de se limiter \`a $\nombre{19001}
\leq t \leq \nombre{20000}$. La valeur maximale est choisie pour
que le r\'egime transitoire soit largement d\'epass\'e, et elle fixe
la dur\'ee du calcul : avec 100 pas par an, un calcul sur
\nombre{20000} ans se fait en \`a peu pr\`es 2 minutes avec un
ordinateur r\'ecent\footnote{PC avec un processeur Athlon 2\virg 4
GHz et 768 Mo de m\'emoire vive.}. Avec $t\geq \nombre{19001}$, on
dispose de suffisamment de points pour distinguer clairement le
type d'attracteur (orbite p\'eriodique, cycle, type H\'enon,
\'etrange, etc.) tout en \'evitant le r\'egime transitoire.

Ces deux valeurs ont \'et\'e test\'ees sur quelques exemples de
valeurs des param\`etres, en s'assurant que le r\'egime
transitoire est tr\`es largement pass\'e. Lors des simulations
suivantes, on a v\'erifi\'e\footnote{Cela se d\'etecte facilement
sur la repr\'esentation en trois dimensions.} que c'\'etait
toujours le cas.

\paragraph{Visualisation de la mesure physique}
Pour limiter la taille des figures repr\'esentant l'attracteur
$(0.15,0.30,8.25)$, nous avons repr\'esent\'e une partie seulement
des points que nous avons calcul\'es, de telle sorte qu'il y ait
de l'ordre de %
\nombre{20000}
points sur chaque graphique. Certaines zones \'etant beaucoup plus
denses que d'autres, cette limitation a \'et\'e faite dans des
proportions diff\'erentes suivant les huit zones de l'attracteur
que nous avons d\'efinies. Il devrait donc y avoir environ deux
fois plus de points dans les zones 2 \`a 5, sur tous les
graphiques o\`u figurent les 80 points de la carte. Cette
transformation n'a ainsi pas \'et\'e effectu\'ee sur la
figure~\ref{fig:0.15_0.30_8.25b}, qui permet donc de se faire une
id\'ee plus exacte de la mesure physique sur l'attracteur.

\subsection{Diagrammes de bifurcation}
Afin de comprendre le r\^ole que jouent les param\`etres du
mod\`ele dans la dynamique, nous avons r\'ealis\'e des animations
repr\'esentant les graphiques en 3 dimensions, l'un des
param\`etres variant au cours du temps. Pour retranscrire une
partie de ces informations sur un graphique bidimensionnel, nous
avons d\^u tracer des diagrammes de bifurcation. Le principe est
le suivant : on porte en abscisse l'un des param\`etres du
syst\`eme et $N(t)$ en ordonn\'ee, pour $t$ entier,
$\nombre{19001} \leq t \leq \nombre{20000}$. On s'est content\'e
de \nombre{1000} valeurs de $t$, au lieu de \nombre{10000}, car la
projection unidimensionnelle de l'attracteur ne permet pas de bien
saisir sa g\'eom\'etrie. On distingue tout au plus les ensembles
finis, les ensembles fractaux et les ensembles continus. Une
pr\'ecision accrue n'aurait rien apport\'e.

\subsection{Injectivit\'e de la projection}

Pour \'evaluer l'injectivit\'e de la projection $\pi$ : $\R^{201}$
$\rightarrow$ $\R^3$, nous avons cherch\'e \`a \'evaluer $\sup_{t
\neq t^{\prime} \in \N}
\frac{\norm{x_{201}(t)-x_{201}(t^{\prime})}_{\R^{201}}}
{\norm{x_{3}(t)-x_{3}(t^{\prime})}_{\R^3}}$ au voisinage des 80
points de la <<carte>> de l'attracteur (section
\ref{sec:injectivite}). En chaque point, nous avons d\'etermin\'e
les \'el\'ements de la boule $B$ de rayon $r = 0.1$ dans
$L^p(\R^{201})$ centr\'ee en ce point $x_{201}(t_0)$. Nous avons
alors calcul\'e $\sup_{t \neq t_0 \in B}
\frac{\norm{x_{201}(t)-x_{201}(t^{\prime})}_{\R^{201}}}
{\norm{x_{3}(t)-x_{3}(t^{\prime})}_{\R^3}}$.

Nous avons indiqu\'e une deuxi\`eme information sur les graphiques
ainsi obtenus : le nombre de points dans chaque boule. En effet,
certains points sont situ\'es dans des zones beaucoup plus denses
que d'autres (au sens de la mesure physique), et cela peut fausser
les r\'esultats obtenus puisque nous ne conservons que le maximum
sur les points de la boule. Il est normal d'obtenir un r\'esultat
plus \'elev\'e si la boule contient plus de points, puisque l'on
risque d'avoir des points tr\`es proches dans $R^3$ par
<<accident>>, au vu des approximations que nous sommes oblig\'es
de faire.

Nous avons choisi les normalisations suivantes pour les normes
$L^p$ dans $\R^N$, $p < \infty$ : $\norm{(x_1, \ldots ,
x_N)}_{L^p}^p = \frac{1}{N} \sum_{i=1}^N x_i^p$. La norme
$L^{\infty}$ est simplement le $\sup$ des coordonn\'ees.

Notons que lorsque $\delta_t$ est pris non-nul, le vecteur
$x_{201}(t)$ est d\'ecal\'e d'autant. Les boules consid\'er\'ees
ne sont donc pas le m\^emes lorsque $\delta_t$ est diff\'erent.
Ceci permet de mieux comparer les diff\'erentes valeurs de
$\delta_t$, et de ne pas privil\'egier $\delta_t = 0$.

Pour r\'esumer les nombreux graphiques obtenus, nous avons choisi
d'extraire des 80 valeurs de $\sup$ deux donn\'ees : le maximum et
la valeur m\'ediane. Elles nous permettent d'avoir une assez bonne
id\'ee de la qualit\'e de la projection.

Enfin, lorsque nous avons consid\'er\'e d'autres visualisations,
par exemple avec une \'echelle logarithmique, nous avons
proc\'ed\'e aux m\^emes op\'erations, en rempla\c{c}ant $x_3$ par
$g(x_3)$, $x_{201}$ \'etant inchang\'e. Il y a alors un
in\'evitable changement d'\'echelle homoth\'etique : par exemple,
avec $g(x) = \lambda x$ et $\lambda$ tr\`es grand, le r\'esultat
serait artificiellement bon. Il faut donc rapporter les
quantit\'es calcul\'ees au diam\`etre de la projection de
l'attracteur. Dans le cas de l'\'echelle logarithmique, cela ne
fait que confirmer notre conclusion en sa d\'efaveur.

\subsection{G\'eom\'etrie}
La repr\'esentation simplifi\'ee de la figure~\ref{fig:geometrie}
est simplement un extrait de la vari\'et\'e instable, $n=18$ (voir
annexe~\ref{annexe:var_u}). Le d\'ecoupage a \'et\'e fait morceau
par morceau, et le choix des r\'egions est purement visuel.

\subsection{D\'ecomposition spectrale, m\'elange}
Pour d\'eterminer s'il y a ou non m\'elange topologique, nous
avons utilis\'e plusieurs couleurs suivant la valeur de $t$ modulo
$N$, pour quelques valeurs de $N$ entre 2 et 10. Lorsque des zones
bien s\'epar\'ees se sont distingu\'ees, nous avons conclu \`a la
$N$-p\'eriodicit\'e de la dynamique. Dans le cas contraire, si
nous avions l'impression d'un m\'elange des couleurs pour toutes
les valeurs de $N$ (comme \`a la figure~\ref{fig:0.18_0.30_8.25}),
nous avons conclu au m\'elange topologique. Il ne s'agit donc que
d'impressions visuelles, et non de v\'erifications rigoureuses.

\subsection{Calcul de la dimension fractale}
La dimension fractale est d\'efinie en
annexe~\ref{annexe:theorie_dimfract}. \label{annexe:dimfract}

Calculer la dimension fractale d'un ensemble dont on ne conna\^it
qu'un nombre fini de points, avec une pr\'ecision limit\'ee, est
loin d'\^etre un probl\`eme facile. Nous avons donc d\^u faire de
nombreuses approximations pour tenter d'obtenir une valeur
approch\'ee raisonnable.

\paragraph{Attracteur $(0\virg 15; 0\virg 30; 8\virg 25)$}
Nous avons consid\'er\'e l'attracteur dans $\R^3$, et non dans
$\R^{201}$, avec la projection naturelle, l'origine des temps
\'etant prise $\delta_t=0$. Il s'agit donc d'un ensemble $K$
d'environ $N_0=%
\nombre{100000}
$ points. Pour diff\'erentes valeurs
de $\epsilon$, nous avons calcul\'e le nombre
$\widetilde{N}_{\epsilon}(K)$ de cubes $C_{i,j,k} = [i \epsilon;
(i+1) \epsilon] \times [j \epsilon; (j+1) \epsilon] \times [k
\epsilon; (k+1) \epsilon]$ qui contiennent au moins un point de
$K$. La figure~\ref{fig:dim_fract_0.15_0.30_8.25} repr\'esente
ainsi $\log_{10} \widetilde{N}_{\epsilon}(K)$ en fonction de
$\log_{10} \epsilon$.

En th\'eorie, la dimension fractale est l'oppos\'e de la pente
limite en $-\infty$ de cette courbe (rigoureusement, on sait que
$\widetilde{N}_{\epsilon}(K) \geq N_{\epsilon}(K) \geq
\frac{\widetilde{N}_{\epsilon}(K)}{8}$, la pente limite doit donc
\^etre la m\^eme). Mais d\`es que $\epsilon$ est assez petit,
$N_{\epsilon}$ est constant \'egal \`a $N_0=%
\nombre{100000}
$, car l'ensemble $K$ est fini. Nous avons donc consid\'er\'e la
pente (obtenue par r\'egression lin\'eaire) en nous restreignant
\`a $N_{\epsilon} < \frac{N_0}{10} = N_{\epsilon_0}$ et $\epsilon
> \sqrt{\epsilon_0}$. Le choix de ce domaine a \'et\'e fait au vu
des donn\'ees et correspond \`a la zone o\`u les points sont bien
align\'es.

En raison de ces nombreuses approximations, il ne faut pas
accorder trop d'importance \`a la valeur pr\'ecise que nous avons
obtenue, seul l'ordre de grandeur (entre 1 et $1\virg 5$) a de
bonnes raisons d'\^etre exact.

\paragraph{Diagramme $(0.18,0.41,\gamma)$}
Pour r\'ealiser la figure~\ref{dim_f:0.18_0.41_gamma}, nous avons
d\^u calculer la dimension fractale d'un grand nombre d'ensemble,
via un processus automatis\'e. Au lieu de %
\nombre{100000}
points, nous avons d\^u nous contenter de $N_0=%
\nombre{10000}
$ points. La pr\'ecision de ces calculs est donc encore
inf\'erieure. Nous disposons probablement d'une l\'eg\`ere
sous-estimation des dimensions fractales de ces attracteurs.

\subsection{Sensibilit\'e aux conditions initiales}
On choisit un point sur l'attracteur, consid\'er\'e dans
$\R^{201}$, et l'on d\'etermine l'ensemble des points de
l'attracteur situ\'es \`a une distance inf\'erieure \`a $r$ pour
la norme $L^p$ (sur la figure~\ref{fig:passe_futur_eq_continu},
$r=0 \virg 04 $ et $p=\infty$). Par point <<de l'attracteur>>, on
entend qu'il s'agisse d'un des %
\nombre{20000}
points d'une orbite calcul\'ee pr\'ec\'edemment et dont on a les
coordonn\'ees dans $\R^{201}$. Ces points sont donc r\'epartis
selon la mesure physique.

Pour chacun de ces points, on trace $N(t)-N^{\star}(t)$ pour $t
\in [-20;20]$, o\`u $N^{\star}(t)$ correspond au point
initialement choisi. L'intervalle de temps de r\'ef\'erence est
$[0;2]$, le graphique de la
figure~\ref{fig:passe_futur_eq_continu} n'est donc pas exactement
centr\'e sur le pr\'esent.

\subsection{Point fixe}
\begin{figure}
\begin{center}
\includegraphics[height=7cm]{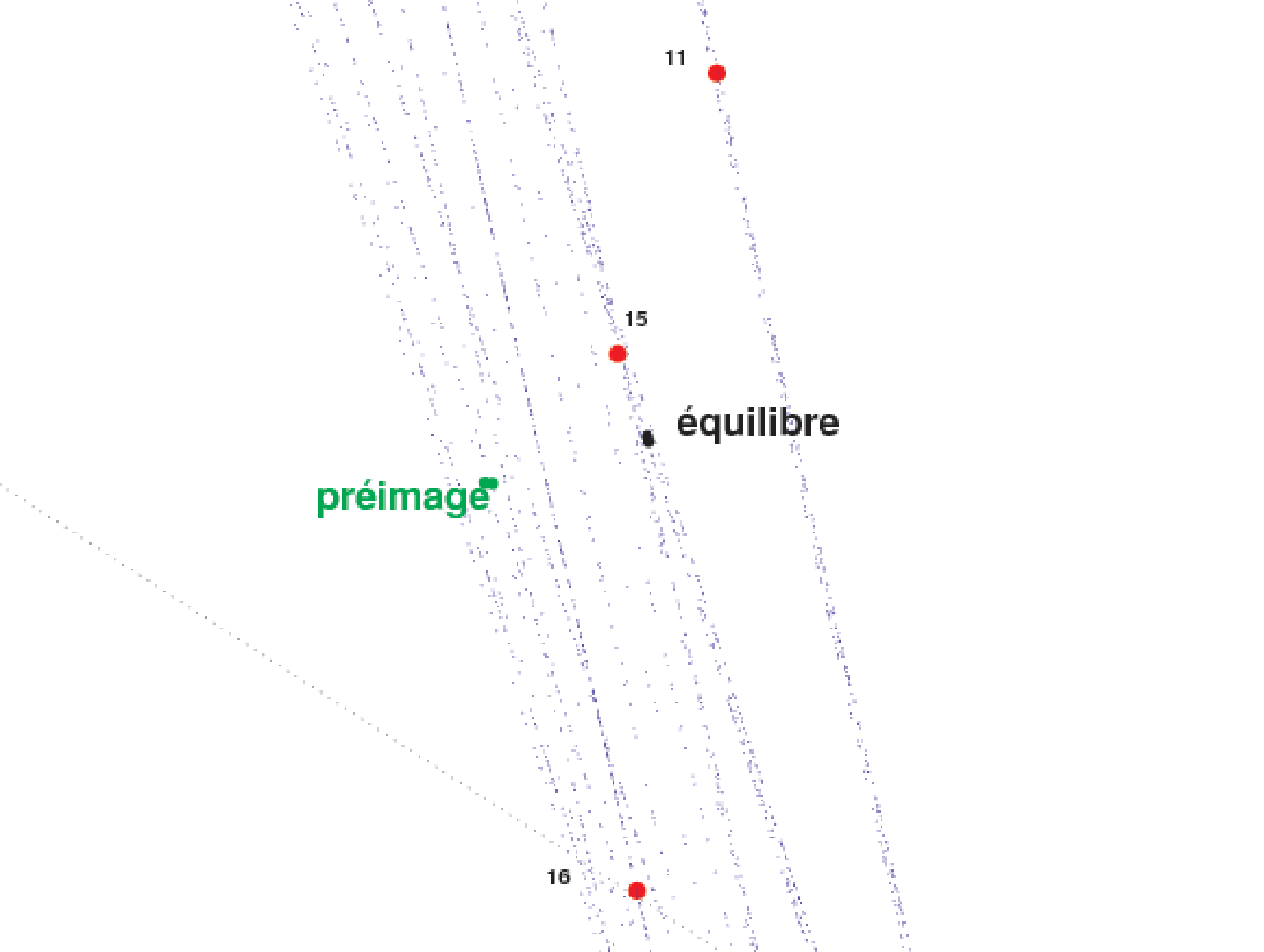}
\caption{\label{fig:preimage_equilibre_zoom} Les <<points fixes>>
et leurs pr\'eimages par $T^2$.}
\end{center}
\end{figure}

\paragraph{Premi\`ere localisation dans $\R^3$}
\`A la recherche d'\'eventuels points fixes, nous avons cherch\'e
\`a minimiser la distance $L^1$ (non-normalis\'ee), dans $\R^3$,
entre $x_3(t)$ et $x_3(t+2)=T^2(x_3(t))$. Sur les %
\nombre{50000}
points consid\'er\'es (nous nous sommes limit\'es \`a la
<<grande>> composante), deux sont \'ecart\'es de leurs images de
moins $4.10^{-3}$. Ces points \'etant tr\`es proches l'un de
l'autre, nous avons consid\'er\'e qu'ils sont \`a proximit\'e d'un
unique \'equilibre. La visualisation de ce point de son image en 3
dimensions ont confirm\'e cette impression dans la mesure o\`u
l'image du point <<fixe>> en est tr\`es proche et sur le m\^eme
filament (figure~\ref{fig:preimage_equilibre_zoom}).

En effet, d'autres points sur l'attracteur sont \'egalement \`a
une petite distance (inf\'erieure \`a $2.10^{-2}$) de leur image
par $T^2$ (\`a proximit\'e du point 43), mais celle-ci est sur un
filament bien distinct. Il n'y aurait donc qu'un seul point fixe
sur l'attracteur.

\paragraph{Localisation plus pr\'ecise, dans $\R^{201}$}
Afin de disposer d'une meilleure approximation de ce point fixe,
aussi bien dans $\R^3$ que dans $\R^{201}$, nous avons eu recours
\`a une m\'ethode un peu plus sophistiqu\'ee.

Nous avons tout d'abord choisi deux points de part et d'autre de
l'\'equilibre pr\'esum\'e, avec une marge assez importante, et
nous avons trac\'e le segment rejoignant ces deux points. Nous
avons calcul\'e les images successives de ce segment par $T^2$, en
les tronquant de telle sorte que l'on reste dans un m\^eme
voisinage de l'\'equilibre pr\'esum\'e. Apr\`es un nombre
suffisant d'it\'erations de ce processus, nous disposons d'une
bonne approximation de la vari\'et\'e instable locale \`a
l'\'equilibre. C'est en effet ce que montre le $\lambda$-lemme
sous certaines conditions\footnote{Le $\lambda$-lemme ne
s'applique sans doute pas directement dans notre cadre non
hyperbolique, en dimension infinie, mais nous donne de bonnes
raisons de penser que sa conclusion est au moins approximativement
vraie.}. Le segment initial ayant \'et\'e choisi transverse \`a la
vari\'et\'e stable, et la troncature nous permettant d\'eviter la
r\'egion de pli, cela devrait \^etre vrai dans notre cas pratique.

Disposant de la vari\'et\'e instable locale, son image par
$T^{-2n}$ doit se contracter autour de l'\'equilibre pour $n$
assez grand. Nous avons ainsi commenc\'e par calculer $T^{4n}$ de
la vari\'et\'e instable, suffisamment tronqu\'ee pour \'eviter le
pli, tout en gardant la m\'emoire du pass\'e (les troncatures
successives nous obligeant \`a rajouter des points sur la ligne
bris\'ee). Pour une valeur suffisamment grande de $n$ (de l'ordre
de 10), nous avons obtenu, en appliquant $T^{-2n}$ \`a la
vari\'et\'e instable locale obtenue, un morceau de vari\'et\'e
instable tr\`es proche de l\'equilibre.

En minimisant la distance $L^1$ entre les points de cette
vari\'et\'e instable locale et leur image (dans $\R^{201}$), on
obtient ainsi une tr\`es bonne approximation de l'\'equilibre.
L'erreur $L^1$ que nous avons obtenue est en effet environ \'egale
\`a $10^{-4}$, ce qui est presque %
\nombre{3000}
 fois mieux que notre
premi\`ere approximation\footnote{dans la mesure o\`u la norme
$L^1$ utilis\'ee n'est pas normalis\'ee}.

\subsection{Vari\'et\'e instable}\label{annexe:var_u}

On part de la vari\'et\'e instable locale \`a l'\'equilibre. La
vari\'et\'e instable globale est alors donn\'ee par la
d\'efinition \ref{def:var_stable_globale}. Comme il s'agit d'un
\'equilibre, il s'agit de calculer les images successives par
$T^2$ de la vari\'et\'e instable locale.

Afin de garder une pr\'ecision finale correcte, on a augment\'e
progressivement le nombre de points d\'efinissant la vari\'et\'e
instable. On a rajout\'e des points interm\'ediaires entre deux
sommets cons\'ecutifs de la ligne bris\'ee d\`es que leur distance
$L^2$ dans $\R^3$ est inf\'erieure \`a $\eta = 10^{-2}$.

On a alors deux mani\`eres de visualiser la croissance de la
vari\'et\'e instable dans l'attracteur. D'une part, les images
successives $T^{2n}(W^u_{loc}(x_{eq}))$ pour $1 \leq n \leq 18$.
D'autre part, en ne consid\'erant que $n = 18$, on peut faire
grandir la vari\'et\'e instable en partant du voisinage de
l'\'equilibre. L'ordre d'apparition des diff\'erentes zones est en
principe le m\^eme. La premi\`ere m\'ethode est plus naturelle et
plus simplement interpr\'etable, mais la seconde donne plus de
d\'etails sur l'ordre dans lequel la vari\'et\'e instable se
d\'eploie dans l'attracteur.

\subsection{Calcul de la diff\'erentielle et de ses valeurs
propres}\label{annexe:differentielle} \'Etant donn\'e un point $x$
de $\R^{201}$, il est ais\'e d'estimer la diff\'erentielle de
$T^2$ en ce point : on fixe $\epsilon=10^{-6}$ et on calcule
$(T^2(x+\epsilon e_i)-T^2(x))/\epsilon$. Cela nous donne la
d\'eriv\'ee partielle dans la direction $e_i$ (de la base
canonique de $\R^{201}$. La matrice des d\'eriv\'ees partielles
nous donne une expression de la diff\'erentielle.

Le calcul des valeurs propres et des vecteurs propres est alors
r\'ealis\'e \`a l'aide des fonctions int\'egr\'ees de
Matlab$^\circledR$. Cela signifie que les espaces tangents en $x$
et $T^2(x)$ sont \'egaux \`a l'espace euclidien $\R^{201}$. La
norme consid\'er\'ee est donc la norme $L^2(\R^{201})$. Pour
diff\'erentes valeurs de $\epsilon$ (allant de $10^{-3}$ \`a
$10^{-6}$), les valeurs des cinq premi\`eres valeurs propres
changent tr\`es peu : $\lambda_1 = -2\virg 2942$, $\lambda_2
\approx 0\virg 0433$, $\lambda_3 \approx -0\virg 0283$,
$\lambda_{4,5} \approx 0\virg 0214 \pm 0\virg 0028 i$. Une
incertitude persiste cependant sur les valeurs propres
$\lambda_j$, $j \geq 2$, dans la mesure o\`u les algorithmes de
calculs sont assez instables. Il serait utile de v\'erifier les
r\'esultats ci-dessus par un autre calcul de la diff\'erentielle,
par exemple en lin\'earisant directement l'\'equation.

Les vecteurs propres repr\'esent\'es sont normalis\'es pour que
leur moyenne ($L^1$) sur $[0;2]$ soit \'egale \`a 1.

\subsection{Hyperbolicit\'e}
Pour le calcul de la diff\'erentielle et de son spectre, on a
utilis\'e la m\'ethode pr\'ec\'edente avec $\epsilon=10^{-3}$. On
a consid\'er\'e qu'une valeur propre met en d\'efaut
l'hyperbolicit\'e de l'attracteur lorsque
($\absj{\log_{10}(\lambda)} < \epsilon_d=\log_{10}(3/2)$). Notons
que les valeurs propres d\'ependent fortement de la structure des
espaces tangents d'arriv\'ee et de d\'epart. Ce n'est pas parce
que la structure canonique du fibr\'e tangent ne rend pas
l'attracteur hyperbolique que celui-ci ne l'est pas.

Pour certains points, les deux premi\`eres valeurs propres sont
\'egales (en module) : il s'agit des points o\`u l'on a deux
valeurs propres complexes conjugu\'ees de module maximal. Dans les
autres cas, il s'agit toujours d'une valeur propre r\'eelle.

\subsection{Pli : Localisation}
La premi\`ere \'etape dans la localisation du pli est arbitraire
et un peu impr\'ecise : on choisit deux points sur l'attracteur
$x_1$ et $x_{1000}$ (figure~\ref{fig:pli}a), \`a peu pr\`es sur le
m\^eme filament et tels que $T^2(x_1)$ et $T^2(x_{1000})$ sont
situ\'es de part et d'autre du pli (figure~\ref{fig:pli}b).

On d\'etermine ensuite un segment dans $\R^{201}$ joignant
$T^{-10}(x_1)$ et $T^{-10}(x_{1000})$, sur lequel on place %
\nombre{1000}
points (y compris les extr\'emit\'es). En r\'eappliquant $T^{10}$
\`a ces %
\nombre{1000}
points, on d\'efinit une ligne bris\'ee $(x_j)_{j=1 \ldots 1000}$
: c'est le filament non-pli\'e de la figure~\ref{fig:pli}a, qui
correspond au temps $t=0$. Le temps $t$ correspond ainsi \`a la
ligne bris\'ee $(T^t(x_j))_{j=1 \ldots 1000}$.

Les d\'etails concernant les calculs de courbure sont donn\'es \`a
la section suivante.

Les calculs des diff\'erentielles ont \'et\'e effectu\'es comme
indiqu\'e section \ref{annexe:differentielle}. La diff\'erentielle
a \'et\'e calcul\'ee \`a deux instants distincts : avant pliage
($t=0$) et \`a l'emplacement du pli ($t=2$).

Pr\'ecisons enfin que lorsqu'on parle du filament \`a l'instant
$t=0$ (ou de sa courbure, etc.), il s'agit des $x_j$. Dans le
cadre d'une visualisation d'un $x_j$ particulier, dans $\R^{201}$,
on note l'intervalle de temps $[-2;0]$ : c'est l'instant final qui
est pris en compte. De m\^eme, le filament \`a l'instant $s$ est
compos\'e des $T^s(x_j)$, et correspond \`a l'intervalle $[s-2;s]$
dans $\R^{201}$.

\paragraph{Deuxi\`eme tentative}
On a r\'ealis\'e une deuxi\`eme tentative de localisation du pli,
plus pr\'ecise, devant en principe supprimer les irr\'egularit\'es
de courbure que nous avons constat\'ees. Pour cela, on utilise la
vari\'et\'e instable globale au point fixe calcul\'ee
pr\'ec\'edemment. Pour $n=13$, une partie de $f^n(W^u(x_{eq}))$
est situ\'ee au niveau du pli. On a alors d\'etermin\'e sa
pr\'eimage par $f^{13}$ et augment\'e le nombre sommets de la
ligne bris\'ee pour disposer d'un filament de plus de %
\nombre{1000}
points dans la zone de pliage (en recalculant son image par
$f^{13}$).

\subsection{Vecteur tangent, courbure d'une ligne bris\'ee}
Soient $x_1 , \ldots , x_N$ les sommets successifs de cette ligne
bris\'ee.

\paragraph{Vecteur tangent}
Le vecteur tangent $T_j$ en un point $x_j$ a \'et\'e calcul\'e
tout simplement \`a l'aide de la formule suivante : \[ T_j =
\frac{x_{j+1}-x_j}{\norm{x_{j+1}-x_j}_{L^1}} \times 201. \] $T_j$
est donc normalis\'e pour avoir une moyenne ($L^1$) \'egale \`a 1
(c'est un \'el\'ement de $\R^{201}$).

\paragraph{Courbure}
En supposant les sommets de la ligne bris\'ee assez proches, on
peut \'evaluer la courbure de la courbe qu'elle approche \`a
l'aide d'une version discr\`ete des formules continues
d\'efinissant la courbure :

\begin{align*}
ds_j &= \norm{x_{j+1}-x_j}_{L^2} \\
T_j &= \frac{x_{j+1}-x_j}{ds_j} \\
\kappa_j &= \absj{\frac{T_{j+1}-T_j}{ds_j}}\end{align*}

Les points $x_j$ sont consid\'er\'es dans $\R^3$ et non dans
$R^{201}$. La norme $L^2$ utilis\'ee pour estimer l'\'el\'ement de
longueur $ds_j$ n'est pas normalis\'e, mais cela ne change rien
pour le calcul de $\kappa_j$. On estime de cette fa\c{c}on la
valeur absolue de la courbure, et non la courbure proprement dite.
\clearpage
\section{Rappels de syst\`emes dynamiques} \label{annexe:rappels}
L'objectif de cette section est de rappeler les notions les plus
importantes de th\'eorie des syst\`emes dynamiques et que nous
avons \'evoqu\'ees pr\'ec\'edemment. Pour plus de d\'etails, on se
reportera par exemple \`a
\cite{KatokHasselblatt:IntroDynamicalSystems} ou
\cite{Benoistpaulin:systemesdynamiques}. Pour les aspects plus
sp\'ecifiquement hyperboliques, on se reportera \`a
\cite{PalisTakens:HyperbolicityChaotic} et
\cite{Yoccoz:HyperbolicDynamics}.

\subsection{Un peu de vocabulaire}
Un \emph{syst\`eme dynamique continu} est la donn\'ee d'un espace
$X$ et d'un groupe \`a un param\`etre de transformation $(f^t)_{t
\in \R}$ (c'est-\`a-dire une famille d'applications continues $X
\rightarrow X$ telles que $f^{t+t^{\prime}} = f^t \circ
f^{t^{\prime}}$). Un \emph{syst\`eme dynamique discret} est la
donn\'ee d'un espace topologique $X$ (l'\emph{espace des phases})
et d'une application $f$ continue $X \rightarrow X$. On peut se
ramener du premier cas au second par le biais de l'application de
retour de Poincar\'e $f^1$. Dans la suite, sauf indication
contraire, on se placera toujours dans le cas discret.

\begin{Def}[Orbite] Si $x$ est un point de $X$, l'\emph{orbite
(positive) de $x$} est l'ensemble $\{ f^n(x) \telque n \geq 0 \}$.
Si $f$ est bijective, l'\emph{orbite de $x$} est $\{ f^n(x)
\telque n \in \Z \}$.  \end{Def}

La th\'eorie des syst\`emes dynamiques s'int\'eresse
particuli\`erement au comportement des orbites. Il est souvent
utile de consid\'erer aussi des <<pseudo-orbites>> au sens de la
d\'efinition suivante.

\begin{Def}[$\delta$-pseudo-orbite] Soit $\delta >0$. Une suite
$(x_n)_{n \in \N}$ (ou $(x_n)_{n \in \Z}$) est une
$\delta$-pseudo-orbite si $\forall i \in \N$ (ou $\Z$),
$d(f(x_i),x_{i+1}) < \delta$.
\end{Def}

Par exemple, une orbite simul\'ee num\'eriquement est une
$\delta$-pseudo-orbite, puisque les calculs sont effectu\'es avec
une pr\'ecision limit\'ee.

\begin{Def}[Partie invariante] \label{def:invariante} Une partie $A$ de $X$ est dite
\emph{invariante\footnote{attention, cette d\'efinition est
parfois remplac\'ee par $f(A)\subset A$.}} par $f$ si $f(A) = A$.
\end{Def}

Dans la suite, on supposera $X$ m\'etrique compact.

\medskip

La notion de conjugaison topologique est extr\^emement
importantes. Elle traduit l'id\'ee que deux syst\`emes dynamiques
sont topologiquement \'equivalents.
\begin{Def}[Conjugaison topologique] Soit $r \geq 0$. Deux applications $C^r$ $f: X
\rightarrow X$ et $g: Y \rightarrow Y$ sont \emph{topologiquement
conjugu\'ees} lorsqu'il existe un hom\'eomorphisme $h : X
\rightarrow Y$ tel que $f = h^{-1} \circ g \circ h$. \end{Def} 

Lorsque $h$ est un $C^m$ diff\'eomorphisme ($m \leq r$), on parle
de \emph{conjugaison lisse}. Parfois, on peut seulement trouver $h
: X \rightarrow Y$ continue surjective telle que $h \circ f = g
\circ h$. On parle alors de \emph{semi-conjugaison}.

Nous pouvons d\'esormais d\'efinir la stabilit\'e structurelle
d'un syst\`eme dynamique.
\begin{Def}[Stabilit\'e structurelle]
\label{def:struc_stable} Une application $f$ $C^r$ est \emph{$C^m$
structurellement stable} ($1 \leq m \leq r \leq \infty$) s'il
existe un voisinage $U$ de $f$ pour la $C^m$ topologie telle que
toute application $g \in U$ est topologiquement conjugu\'ee \`a
$f$.

Si de plus on peut choisir $h = h_g$ dans la conjugaison de $f$ et
$g$ tel que $h_g$ et $h_g^{-1}$ convergent uniform\'ement vers
l'identit\'e lorsque $g$ converge vers $f$ pour la topologie
$C^m$, alors on dit que $f$ est \emph{$C^m$ fortement
structurellement stable}.
\end{Def}

\subsection{R\'ecurrence}
\begin{Def}[Point p\'eriodique] Un point $x$ de $X$ est dit
\emph{p\'eriodique} lorsqu'il existe $n \geq 1$ tel que $f^n (x) =
x$. On note $Per(f)$ l'ensemble des points p\'eriodiques.
\end{Def}
\[ f(Per(f)) = Per(f) \text{ et } f(\adherence{Per(f)}) =
\adherence{Per(f)}. \]

\begin{Def}[Point r\'ecurrent] Un point $x$ de $X$ est dit
\emph{positivement r\'ecurrent} (resp. \emph{n\'egativement
r\'ecurrent}) si $x$ est un point d'accumulation de la suite
$(f^n(x))_{n \geq 0}$ (resp. $(f^{-n}(x))_{n \geq 0}$). On note
$R^+(f)$ l'ensemble des points positivement r\'ecurrents, $R^-(f)$
l'ensemble des points n\'egativement r\'ecurrents et $R(f) =
R^+(f) \cup R^-(f)$ l'ensemble des points \emph{r\'ecurrents}.
\end{Def}

Autrement dit, partant d'un point r\'ecurrent, on revient une
infinit\'e de fois dans son voisinage. On d\'emontre :
\begin{gather*}
Per(f) \subset R^+(f) \cap R^-(f) \neq \emptyset \\
f(R^+(f))=R^+(f) \text{ et } f(R^-(f))=R^-(f) \end{gather*}

\begin{Def}[Point limite] Pour tout point $x$, on note $\omega (x)$
(\emph{ensemble $\omega$-limite de $x$}) l'ensemble des points
d'accumulation de $(f^n(x))_{n \geq 0}$ et $\alpha (x)$
(\emph{ensemble $\alpha$-limite de $x$})l'ensemble des points
d'accumulation de $(f^n(x))_{n \leq 0}$. On d\'efinit alors
l'\emph{ensemble $\omega$-limite} $L^+ (f) = \overline{\bigcup_{x
\in X} \omega(x)}$, l'\emph{ensemble $\alpha$-limite} $L^- (f) =
\overline{\bigcup_{x \in X} \alpha(x)}$ et l'\emph{ensemble
limite} $L(f) = L^+(f) \cup L^-(f)$.
\end{Def}
Un point de $\omega(x)$ est un point dont l'orbite issue de $x$
visite le voisinage une infinit\'e de fois. On a les
propri\'et\'es suivantes :
\begin{gather*} R^+(f) \subset L^+(f) \text{ et } R^-(f) \subset
L^-(f) \\ f(L^+(f))=L^+(f) \text{ et } f(L^-(f))=L^-(f)
\end{gather*}

\begin{Def}[Point errant\footnote{en anglais : wandering}]
Un point $x \in X$ est dit \emph{errant} s'il poss\`ede un
voisinage $V$ tel que $f^n(V) \cap V = \emptyset$ pour tout $n
\geq 1$. Sinon, on dit que $x$ est \emph{non-errant}, et on note
$\Omega(f)$ l'ensemble des points non-errants.
\end{Def}
Un point est donc non-errant lorsque tout voisinage se recoupe au
moins une fois dans le futur. On peut d\'emontrer les r\'esultats
suivants :
\begin{gather*} \adherence{L^+(f)} \subset \Omega(f) \text{ et }
\adherence{L^-(f)} \subset \Omega(f)
\\ f(\Omega(f)) = \Omega(f) \end{gather*}

\begin{Def}[Point r\'ecurrent par cha\^ine\footnote{en anglais :
chain-recurrent}] Un point est r\'ecurrent par cha\^ine si pour
tout $\delta >0$ il existe une $\delta$-pseudo-orbite p\'eriodique
issue de $x$. On note $C(f)$ l'ensemble des points r\'ecurrents
par cha\^ine. \end{Def} Un point r\'ecurrent par cha\^ine est un
point qui peut revenir exactement en lui-m\^eme en autorisant des
<<erreurs>> d'amplitude aussi petites que l'on veut.
\begin{gather*}
f(C(f)) = C(f) \text{ et } \Omega(f) \subset C(f) \\
C(f_{|\Lambda})= \Lambda \text{ si } \Lambda \in \{
\adherence{Per(f)} , \adherence{R(f)} , \adherence{L(f)} , C(f) \}
\end{gather*}

\begin{Pro} En r\'esum\'e, on a : \[ Per(f) \subset
R^{\pm}(f) \subset L^{\pm}(f) \subset L(f) \subset \Omega(f)
\subset C(f) \] et chacune de ces inclusions peut \^etre stricte.
De plus, chacun de ces ensembles (ainsi que leurs adh\'erences)
est une partie invariante par $f$, au sens de la d\'efinition
\ref{def:invariante}.
\end{Pro}


\begin{Def}[Transitif] Un hom\'eomorphisme $f$ est
\emph{transitif} si pour tout ouvert non-vide $U$, $\bigcup_{n \in
\N} f^n(U)$ est dense dans $X$. \end{Def} Ceci \'equivaut \`a dire
qu'il existe $x \in X$ dont l'orbite est dense (\latin{i.e.}
$\omega(x) = X$).

\begin{Def}[M\'elange topologique] \label{def:melange}
Un hom\'eomorphisme $f$ d'un espace m\'etrique compact $X$ est
\emph{topologiquement m\'elangeant} si pour tous $U,V$ ouverts
non-vides, il existe $n_0 \in \Z$ tel que $\forall n \geq n_0$,
$f^n(U) \cap V \neq \emptyset$.
\end{Def}

Si $f$ est topologiquement m\'elangeant, alors $f^k$ est transitif
pour tout entier $k \neq 0$. La r\'eciproque est fausse (voir
l'exemple de la rotation du cercle : \ref{ex:rotation_cercle}). Il
n'y a en effet pas n\'ecessairement de <<m\'elange>> de l'espace
des phases sous l'action de la dynamique dans le cas d'un
syst\`eme transitif. Il y a ainsi transitivit\'e de tout syst\`eme
restreint \`a une orbite p\'eriodique, mais jamais de m\'elange.

\begin{Def}[Minimalit\'e] Une partie ferm\'ee $A \subset X$ est
\emph{minimale} pour $f$ si elle est non-vide, invariante par $f$
et si $A$ ne contient pas de ferm\'e non-vide invariant par $f$
autre que $A$.\end{Def} Ceci \'equivaut \`a dire que l'orbite
(positive) de tout point $x \in A$ est dense dans $X$. En
particulier, tout point de $A$ est r\'ecurrent.

\begin{ex} \label{ex:rotation_cercle} Le cercle $S^1$ est minimal pour la rotation $R_{\alpha}$ : $\theta$
$\rightarrow$ $\theta + \alpha$ (mod. 1) si $\alpha$ est
irrationel. Ce syst\`eme est transitif, mais pas topologiquement
m\'elangeant.\end{ex}
\begin{ex}
Le doublement de l'angle $\theta$ $\rightarrow$ $2 \theta $ (mod.
1) sur $S^1$ est topologiquement m\'elangeant, donc
transitif.\end{ex}

\begin{Pro}
Si $X$ est m\'etrique compact non-vide, $f$ continue $X
\rightarrow X$, alors $X$ contient un ferm\'e minimal pour $f$. En
particulier, $R(f) \neq \emptyset$.
\end{Pro}

Les d\'efinitions suivantes pr\'ecisent les notions intuitives
d'attracteur et de bassin d'attraction.
\begin{Def}[Attracteur] \label{def:attracteur} Une partie compacte $A \subset X$ est un
\emph{attracteur} pour $f$ s'il existe un voisinage $V$ de $A$ et
un entier $N \in \N$ tel que $f^N (V) \subset V$ et $A =
\bigcap_{n \in N} f^n(V)$. \end{Def}

\begin{Def}[Bassin d'attraction] \label{def:bassin} Soit $A$ un attracteur. Le
\emph{bassin d'attraction} de $A$, not\'e $B(A)$, est l'ensemble
des points $x \in X$ tels que $\omega(x) \subset A$. \end{Def}

\subsection{Dynamique hyperbolique} \label{annexe:hyperbolique}

Le cas particulier de la dynamique hyperbolique est extr\^ement
important, notamment parce qu'elle est pr\'esente dans la plupart
des syst\`emes dynamiques. 
Le comportement hyperbolique est, comme nous allons le voir, le
comportement <<typique>> d'un syst\`eme dynamique. Nous parlerons
ici de syst\`emes dynamiques uniform\'ement hyperboliques.

\subsubsection{D\'efinitions}

Commen\c{c}ons par traiter le cas d'une application lin\'eaire.

\begin{Def}[Application lin\'eaire hyperbolique] Une
bijection lin\'eaire $T : E=\R^n \rightarrow \R^n$ est dite
\emph{hyperbolique} s'il existe une d\'ecomposition $E = E_s
\oplus E_u$ en somme directe de sous-espaces $T$-invariants
(\latin{i.e.} $T E_s = E_s$ et $T E_u = E_u$) tels que, en notant
$S = T |_{E_s}$ et $U = T|_{E_u}$, il existe $n \geq 1$ tel que
$\norm{S^n} <1$ et $\norm{U^{-n}}<1$.
\end{Def} Cette d\'efinition est ind\'ependante du choix de la
norme, et est \'equivalente \`a dire que $T$ n'a pas de valeur
propre de module 1.

On dit que $\norm{\cdot}$ est \emph{adapt\'ee} \`a $T$ lorsque $n=1$
convient dans la d\'efinition pr\'ec\'edente et si $\forall x_s
\in E_s$ et $x_u \in E_u$, on a $\norm{x_s + x_u} = \max \{
\norm{x_s} , \norm{x_u} \} $. On appelle alors \emph{constante
d'hyperbolicit\'e} la constante \[ ch(T) = \max \{ \norm{S} ,
\norm{U^{-1}} \} <1. \]

Dans le reste de cette section, on consid\'erera $X=M$ une
vari\'et\'e lisse, munie d'une distance $d$, $U$ un ouvert de $M$
et $f : U \rightarrow M$ un $C^1$ diff\'eomorphisme sur son image.

\begin{Def}[Point p\'eriodique hyperbolique] Un point p\'eriodique
$p$ de $f$, de p\'eriode $n$, est \emph{hyperbolique} si $Df^n_p:
T_p M \rightarrow T_p M$ est une application lin\'eaire
hyperbolique. Son orbite est appell\'ee \emph{orbite p\'eriodique
hyperbolique}. \end{Def}

Sans perte de g\'en\'eralit\'e, on peut se limiter au cas d'un
point fixe. Le th\'eor\`eme de Grobman-Hartman affirme alors qu'au
voisinage d'un point fixe hyperbolique, un diff\'eomorphisme est
topologiquement conjugu\'e \`a sa diff\'erentielle.

\begin{The}[Grobman-Hartman]
Soit $\Omega$ un ouvert de $\R^N$, $f: \Omega \rightarrow \R^N$ un
$C^1$-diff\'eomorphisme local, $x_0$ un point fixe hyperbolique de
$f$ et $T=Df(x_0)$ la diff\'erentielle de $f$ en $x_0$.  Alors il
existe des voisinages ouverts $U$ de 0 dans $\R^N$ et $V$ de $x_0$
dans $\Omega$, et un hom\'eomorphisme $H : U \rightarrow V$ tel
que, pour tout $x$ dans $U$ avec $T(x)$ dans $U$, on a \[ H \circ
T(x) = f \circ H(x) .\]
\end{The}

G\'en\'eralisons cette notion au cas d'un ensemble invariant
quelconque.

\begin{Def}[Ensemble hyperbolique] Une partie
$f$-invariante $\Lambda \subset U$ est \emph{hyperbolique} si pour
tout $x \in \Lambda$ il existe une d\'ecomposition $TM_x = E^s_x
\oplus E^u_x$ et des constantes $\lambda < 1 < \mu$ et une
m\'etrique Riemannienne sur $M$ v\'erifiant les propri\'et\'es
suivantes : \begin{itemize}
\item $\forall x \in \Lambda$, $Tf_x(E^s_x) = E^s_{f(x)}$ et
$Tf_x(E^u_x) = E^u_{f(x)}$. \item $\forall x \in \Lambda$,
$\norm{Tf_{x \mid E^s_x}} \leq \lambda$ et $\norm{Tf^{-1}_{x \mid
E^u_x}} \leq \mu^{-1}$ (les normes \'etant induites par la
m\'etrique de $M$).
\end{itemize}
\end{Def}

On peut alors montrer que les sous-espaces $E^s_x$ et $E^u_x$
d\'ependent contin\^ument de $x$, ont des dimensions localement
constantes, et sont uniform\'ement transverses (il existe
$\alpha_0 >0$ tel que pour tous $x \in \Lambda$, $\xi \in E^s_x$,
$\eta \in E^u_x$, l'angle entre $\xi$ et $\eta$ est au moins
$\alpha_0$).

\begin{Def}[Anosov] Un $C^1$ diff\'eomorphisme $f : M \rightarrow
M$ d'une vari\'et\'e compacte $M$ est appel\'e un
\emph{diff\'eomorphisme Anosov} si $M$ est hyperbolique pour $f$.
\end{Def}

L'ensemble des diff\'eomorphismes Anosov sur $M$ est un ouvert de
$C^1(M,M)$.

Pour d\'eterminer si un ensemble est hyperbolique, en g\'en\'eral, on
regarde s'il v\'erifie la condition de c\^one suivante.

\begin{Pro}[Condition de c\^one] Soit $\Lambda$ une partie
$f$-invariante, $U$ un voisinage de $\Lambda$, $\alpha > 1$. On
suppose qu'il existe en tout point $x \in U$ une d\'ecomposition
$TM_x = E^1_x \oplus E^2_x$.

Supposons que $\forall x\in U$, $\forall v \in T_x M$, en posant
$v=v_1+v_2$ et $w = T_x f(v) = w_1+w_2$ (d\'ecomposition sur $E^1_x$
et $E^2_x$), on a : $$ \norm{v_2} \geq \norm{v_1} \Rightarrow
\norm{w_2} \geq \alpha \norm{w_1} \text{ et } \norm{w_2} \geq
\alpha \norm{v_2}$$ $$ \norm{w_2} \geq \norm{w_1} \Rightarrow
\norm{v_1} \geq \alpha \norm{w_1} \text{ et } \norm{v_1} \geq
\alpha \norm{v_2}.$$

Sous ces conditions, $\Lambda$ est hyperbolique.
\end{Pro}

Bien s\^ur, cette condition est suffisante, mais pas n\'ecessaire.

\subsubsection{Propri\'et\'es fondamentales}
\begin{Def}[Expansivit\'e]
Un hom\'eomorphisme $f : X \rightarrow X$ est \emph{expansif} s'il
existe une constante $\delta_0 >0$ telle que pour tous $x,y \in
X$, $x \neq y$, il existe $n \in \Z$ tel que
$d(f^n(x),f^n(y))>\delta_0$.
\end{Def}

\begin{Pro}[Expansivit\'e] \label{pro:expansivite}
La restriction d'un diff\'eomorphisme \`a un ensemble hyperbolique
est expansive.
\end{Pro}

Le lemme de pistage est fondamental pour justifier la validit\'e
de simulations num\'eriques, o\`u l'on fait des calculs avec une
pr\'ecision limit\'ee sur un syst\`eme dynamique chaotique. En
revanche, il ne garantit pas que les orbites pistant les
pseudo-orbites sont typiques. Ainsi, pour l'application $f: x
\rightarrow 2x$ (mod. 1), une orbite calcul\'ee par ordinateur
s'achevera toujours en 0, car la condition initiale est donn\'ee
par un nombre fini de bits. L'ordinateur calcule ainsi toujours
une vraie orbite, mais syst\'ematiquement attir\'ee par 0, ce qui
n'est pas le comportement typique du syst\`eme.

\begin{Pro}[Lemme de pistage\footnote{shadowing lemma}] Si $\Lambda$
est un compact hyperbolique pour $f$, alors il existe un voisinage
$U(\Lambda)$ de $\Lambda$ tel que pour tout $\delta >0$, il existe
$\epsilon >0$, $\forall (x_n)_{n \in \Z}$ $\epsilon$-pseudo-orbite
contenue dans $U(\Lambda)$, il existe $x \in X$ v\'erifiant
$\forall n \in \Z$, $d(f^n(x),x_n) < \delta$.
\end{Pro}

Les orbites sous l'action d'une petite perturbation de $f$ sont
ainsi proches des <<vraies>> orbites de $f$, ce qui nous am\`ene
\`a la question de la stabilit\'e structurelle. En fait, on peut
m\^eme montrer qu'en un certain sens, la stabilit\'e structurelle
est \'equivalente \`a la notion d'hyperbolicit\'e.

\begin{The}[Stabilit\'e structurelle]
Si $\Lambda$ est hyperbolique pour $f : U \rightarrow M$, alors
pour tout voisinage $V \subset U$ de $\Lambda$ et tout $\delta
>0$, il existe $\epsilon >0$ tel que si $f^{\prime} : U \rightarrow
X$ et $d_{C^1} (f _{| V}, f^{\prime}) < \epsilon$, il existe un
ensemble hyperbolique $\Lambda^{\prime} =
f^{\prime}(\Lambda^{\prime}) \subset V$ pour $f^{\prime}$ et un
hom\'eomorphisme $h : \Lambda^{\prime} \rightarrow \Lambda$ avec
$d_{C^0} (Id, h) + d_{C^0} (Id, h^{-1}) < \delta$ tel que le
diagramme suivant commute : \[\begin{CD}
\Lambda^{\prime} @>f^{\prime}>> \Lambda^{\prime}\\
@V{h}VV @VV{h}V\\
\Lambda @>f>> \Lambda
\end{CD} \] De plus, $h$ est
unique si $\delta$ est assez petit.

En particulier, les diff\'eomorphismes Anosov sont fortement
structurellement stables (voir d\'efinition
\ref{def:struc_stable}). \end{The}

\subsubsection{Vari\'et\'es stables et instables}
La d\'efinition d'un ensemble hyperbolique d\'egage au voisinage
de chaque point deux directions, l'une stable, l'autre instable,
provenant de l'\'etude du syst\`eme lin\'earis\'e. Intuitivement
(le th\'eor\`eme de Grobman-Hartman fait d\'ej\`a un pas dans ce
sens), le syst\`eme non-lin\'earis\'e devrait pr\'esenter le
m\^eme type de d\'ecomposition, au moins localement : cela conduit
\`a d\'efinir les vari\'et\'es stables et instables. Le
th\'eor\`eme suivant justifie leur d\'efinition, dans un cadre
local.

\begin{The}\label{the:var_stab_inst}Soit $\lambda$ hyperbolique pour un $C^1$
diff\'eomorphisme $f : V \rightarrow M$, avec des constantes
$\lambda < 1 < \mu$. Alors, pour tout point $x \in \Lambda$, il
existe $W^s(x)$ et $W^u(x)$, images de disques par des
$C^1$-plongements, appel\'es \emph{vari\'et\'es stables et
instables locales} en $x$, telles que \begin{enumerate} \item $T_x
W^s(x) = E^s_x$ et $T_x W^u(x) = E^u_x$. \item $f(W^s(x)) \subset
W^s(f(x))$ et $f^{-1}(W^u(x)) \subset W^u(f^{-1}(x))$. \item pour
tout $\delta
>0$, il existe $C(\delta)$ telle que pour tout $n \in \N$,
\begin{gather*} \forall y \in W^s(x), \, d(f^n(x),f^n(y)) <
C(\delta) (\lambda + \delta)^n d(x,y) \\ \forall y \in W^u(x), \,
d(f^{-n}(x),f^{-n}(y)) < C(\delta) (\mu - \delta)^{-n} d(x,y).
\end{gather*} \item il existe $\beta >0$ et une famille de
voisinages $O_x$ contenant une boule autour de $x \in \Lambda$ de
rayon $\beta$ telle que \begin{align*} W^s(x) &= \{ y | f^n(y) \in
O_{f^n(x)} , n \in \N \} \\ W^u(x) &= \{ y | f^{-n}(y) \in
O_{f^{-n}(x)} , n \in \N \} \end{align*} \end{enumerate} \end{The}

Les vari\'et\'es locales stables et instables ne sont pas uniques,
mais l'intersection de $W^s_1(x)$ et $W^s_2(x)$ contient toujours
un voisinage de $x$. On peut ainsi d\'efinir ind\'ependamment du
choix des vari\'et\'es locales les \emph{vari\'et\'es stables et
instables globales}.

\begin{Def}[Vari\'et\'es stables et instables globales]
\label{def:var_stable_globale}
\begin{align*} \widetilde{W}^s(x) &= \bigcup_{n=0}^{\infty}
f^{-n}(W^s(f^n(x))) \\ \widetilde{W}^u(x) &=
\bigcup_{n=0}^{\infty} f^n(W^u(f^{-n}(x)))
\end{align*}
\end{Def}

On a alors une caract\'erisation topologique des vari\'et\'es
stables et instables correspondant \`a la notion intuitive.

\begin{Pro}[Vari\'et\'es stables et instables]
\begin{align*} \widetilde{W}^s(x) &= \{ y \in U \telque
d(f^n(x), f^n(y)) \cv{} 0 \} \\ \widetilde{W}^u(x) &= \{ y \in U
\telque d(f^{-n}(x), f^{-n}(y)) \cv{} 0 \}
\end{align*}\end{Pro}

Ainsi, deux vari\'et\'es stables (resp. instables) globales sont
disjointes ou \'egales.

\subsubsection{Produit local, ensemble localement maximal}
Lorsque des vari\'et\'es stables et instables globales en un point
s'intersectent, il est int\'eressant de consid\'erer leurs points
d'intersection. Ceci am\`ene \`a d\'efinir le produit local, \`a
l'aide de la proposition suivante (qui d\'ecoule en partie du
th\'eor\`eme \ref{the:var_stab_inst}).

\begin{Pro} \label{pro:prod_local}Soit $x \in \Lambda$. On note
$W^s_{\epsilon} (x)$ et $W^u_{\epsilon} (x)$ les boules de rayon
$\epsilon$ dans $\widetilde{W}^s(x)$ et $\widetilde{W}^u(x)$.
Alors il existe $\epsilon >0$ tel que pour tous $x \in \Lambda$,
$W^s_{\epsilon} (x) \cap W^u_{\epsilon} (y)$ contient au plus un
point $[x, y]$, et il existe $\delta >0$ tel que si $d(x,y) <
\delta$ avec $x,y \in \Lambda$, alors $W^s_{\epsilon} (x) \cap
W^u_{\epsilon} (y) \neq \emptyset$.
\end{Pro}

\begin{Def}[Produit local]
On dit qu'un ensemble hyperbolique $\Lambda$ a une \emph{structure
de produit local} si pour $\epsilon >0$ assez petit, le point
$[x,y]$ d\'efini par la proposition \ref{pro:prod_local}
appartient \`a $\Lambda$.
\end{Def}

En fait, cette propri\'et\'e est \'equivalente \`a la notion
suivante (que l'on distinguera bien de la notion d'attracteur,
d\'efinition \ref{def:attracteur}).

\begin{Def}[Localement maximal] Soit $\Lambda$ hyperbolique pour
$f: U \rightarrow M$. S'il existe un voisinage ouvert $V$ de
$\Lambda$ tel que $\Lambda = \Lambda^f_V := \bigcap_{n \in \Z}
f^n(\overline{V})$, on dit que $\Lambda$ est \emph{localement
maximal}.
\end{Def}

Dans ce cas, on montre que les points p\'eriodiques de $f$ sont
denses dans l'ensemble des points non-errants de $f_{|\Lambda}$.

\begin{The}
Soit $\Lambda$ un ensemble compact hyperbolique. Il a une
structure de produit local si et seulement si il est localement
maximal.
\end{The}

Une autre notion importante reli\'ee au produit local est celle
d'\emph{intersection homocline}.
\begin{Def}[Intersection homocline] \label{def:homocline}
Soit $p$ un point fixe hyperbolique. Un point $q\neq p$ est
homocline \`a $p$ si $q \in \widetilde{W}^s(p) \cap
\widetilde{W}^u(p)$. Cette intersection est dite \emph{transverse
homocline} si les deux vari\'et\'es stables et instables
s'intersectent transversalement. \end{Def}

Les exemples les plus importants sont le pendule (en temps
continu) et le fer-\`a-cheval (en temps discret). En g\'en\'eral,
un syst\`eme poss\'edant une intersection homocline pr\'esente une
dynamique tr\`es complexe. \`A ce sujet, on pourra consulter
\cite{PalisTakens:HyperbolicityChaotic}.

\subsubsection{D\'ecomposition spectrale}
Il s'agit d'un r\'esultat important qui permet d'\'etudier la
r\'ecurrence des orbites dans le cas d'ensemble localement
maximaux hyperboliques.

\begin{The}[D\'ecomposition spectrale] \label{the:decompo_spectrale}
Soit $\Lambda$ compact hyperbolique localement maximal pour un
diff\'eomorphisme $f : U \rightarrow M$. Alors il existe une
famille finie des compacts invariants disjoints $\Lambda_1,
\ldots, \Lambda_m$ telle que $\Omega(f_{| \Lambda}) =
\bigcup_{i=1}^m \Lambda_i$. Les ensembles $\Lambda_i)=
\Lambda_{i}$ sont <<irr\'eductibles>> au sens o\`u
$f_{|\Lambda_i}$ est transitif. De plus, $\Lambda_i =
\bigcup_{j=1}^{m_i} \Lambda_{i,j}$ avec $f(\Lambda_{i,j}) =
\Lambda_{i,j+1}$, et $f^{M_i}_{| \Lambda_{i,1}}$ est
topologiquement m\'elangeante.
\end{The}

Si de plus $\Lambda = \overline{C(f_{|\Lambda})}$ (\latin{i.e.}
$\Lambda$ est \emph{r\'ecurrent par cha\^ine}) , alors le
th\'eor\`eme suivant montre que $\Lambda =
\overline{\Omega(f_{|\Lambda})}$, et donc la d\'ecomposition
spectrale s'applique directement \`a $\Lambda$.

\begin{The}
Supposons que $\Lambda = \adherence{Per(f)}$, $\adherence{R(f)}$,
$\adherence{L(f)}$ ou $\adherence{C(f)}$ est hyperbolique. Alors
cet ensemble est localement maximal et r\'ecurrent par cha\^ine.
De plus, cet ensemble est \'egal \`a $\adherence{Per(f)}$.
\end{The}

\begin{Def} Si $\Lambda$ est compact invariant par $f$,
les ensembles stables et instables de $K$ sont \begin{align*}
W^s(\Lambda) &= \{ y \in X \telque \omega(y) \subset \Lambda \}
\\ W^u(\Lambda) &= \{ y \in X \telque \alpha(y) \subset \Lambda
\}. \end{align*} \end{Def}

Cette d\'efinition est coh\'erente avec les notions
pr\'ec\'edentes de vari\'et\'es stables et instables globales,
comme le justifie la proposition suivante.

\begin{Pro} Si $\Lambda$ est un compact invariant hyperbolique localement
maximal, alors $W^s(\Lambda) = \bigcup_{x \in \Lambda}
\widetilde{W}^s(x)$ et
$W^u(\Lambda) = \bigcup_{x \in \Lambda} \widetilde{W}^u(x)$. \end{Pro}

\begin{Def}
On appelle \emph{ensemble basique} l'un des $\Lambda_i$ de la
d\'ecomposition spectrale. Trois situations peuvent alors se
pr\'esenter : \begin{itemize} \item si $W^s(\Lambda_i)$ est
ouvert, on dit que $\Lambda_i$ est un \emph{attracteur}. \item si
$W^u(\Lambda_i)$ est ouvert, on dit que $\Lambda_i$ est un
\emph{r\'epulseur}\footnote{repeller en anglais}. \item si aucun
de ces deux ensembles n'est ouvert, alors $\Lambda_i$ est de
\emph{type-selle}. \end{itemize}
\end{Def}

\subsubsection[Exemple du sol\'eno\"ide]{Exemple de dynamique uniform\'ement hyperbolique : le
sol\'eno\"ide}
\label{annexe:solenoide}

\begin{figure}
\begin{center}
\begin{tabular}{c@{}c}
    \includegraphics[width=7cm]{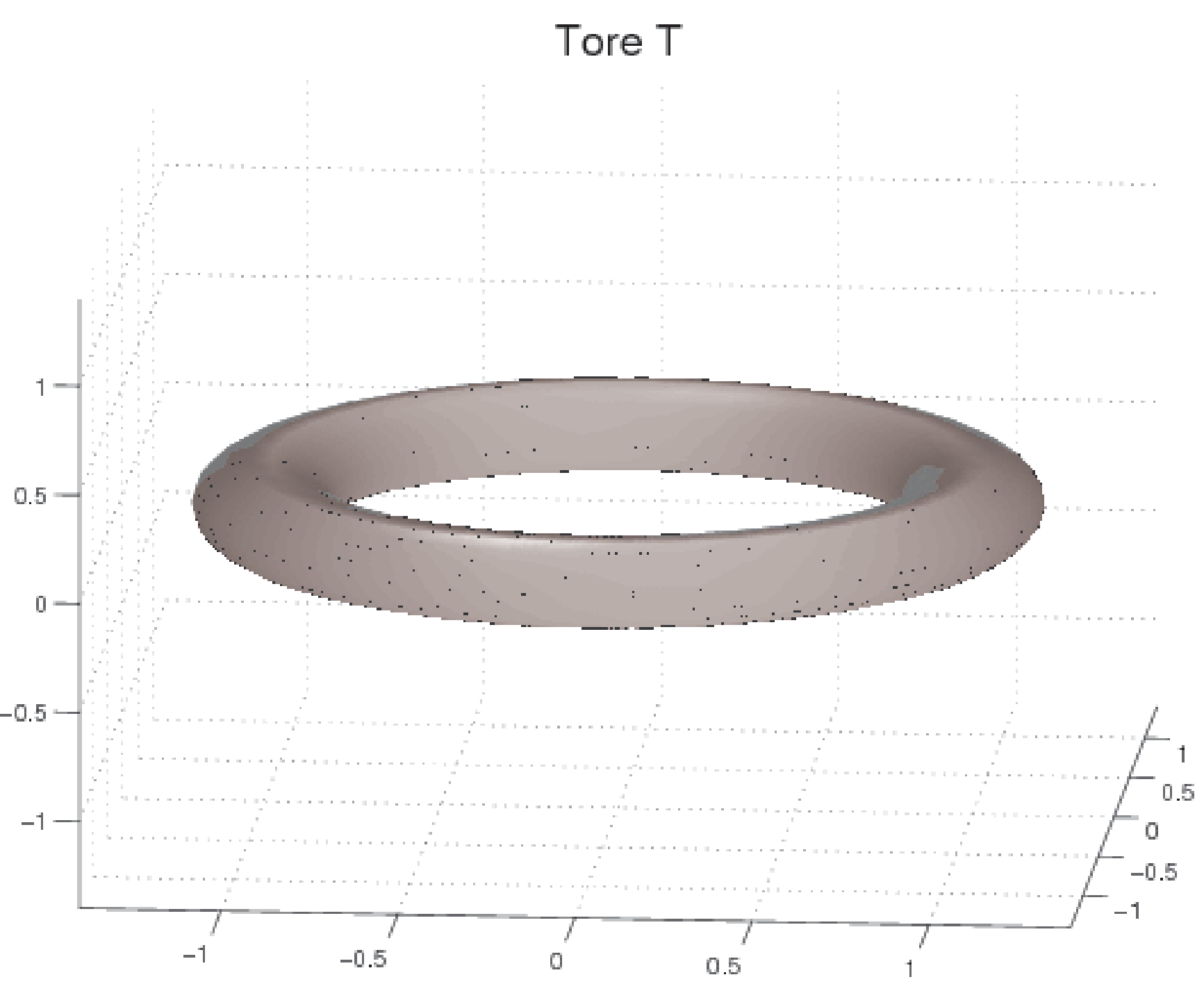}
     &
    \includegraphics[width=7cm]{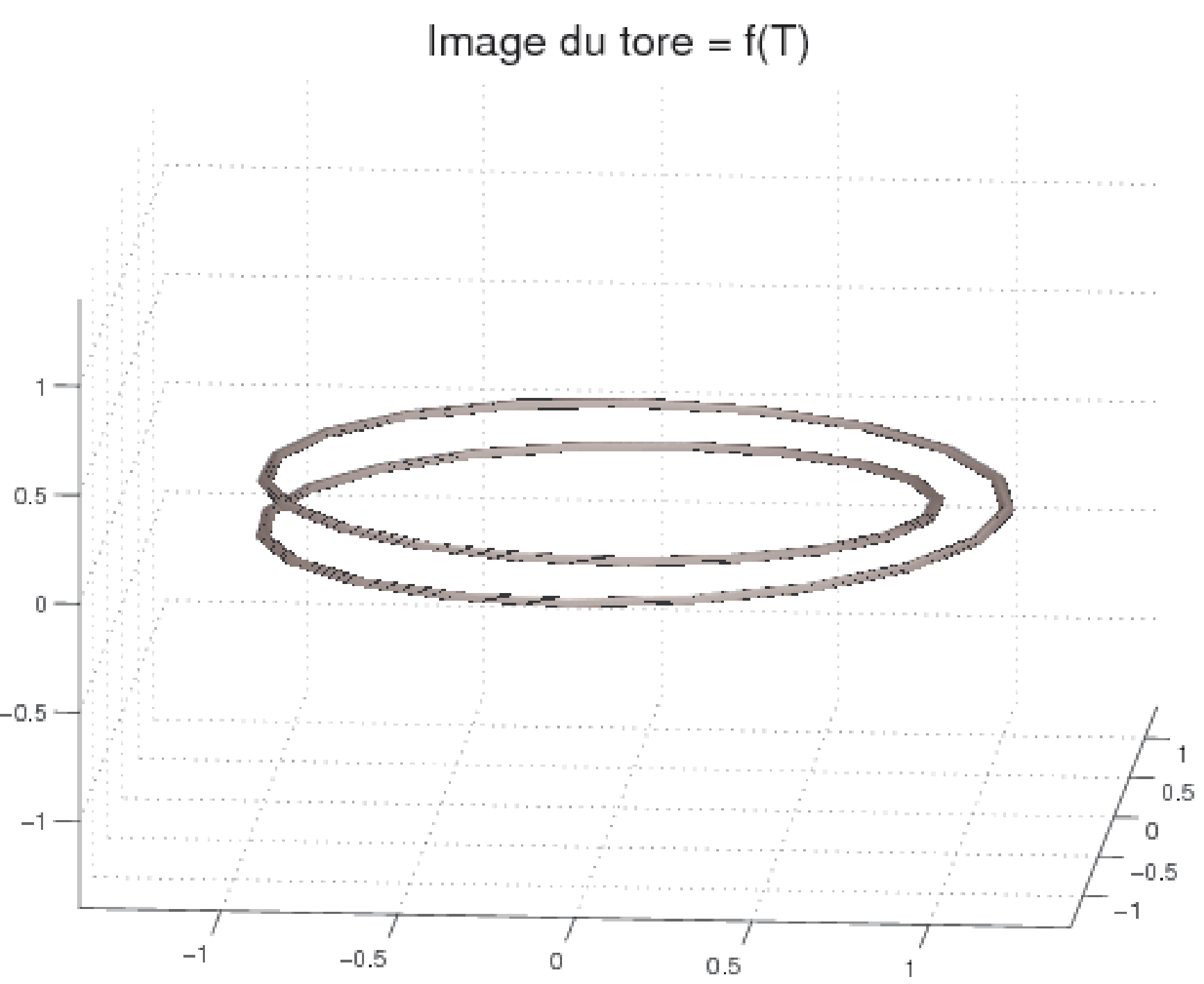}
\\
(a) $T$ & (b) $f(T)$
\end{tabular}
\caption{\label{fig:solenoide_T_fT} Le sol\'eno\"ide : $T$ et
$f(T)$.}
\end{center}
\end{figure}

\begin{figure}
\begin{center}
\begin{tabular}{c@{}c}
    \includegraphics[width=7cm]{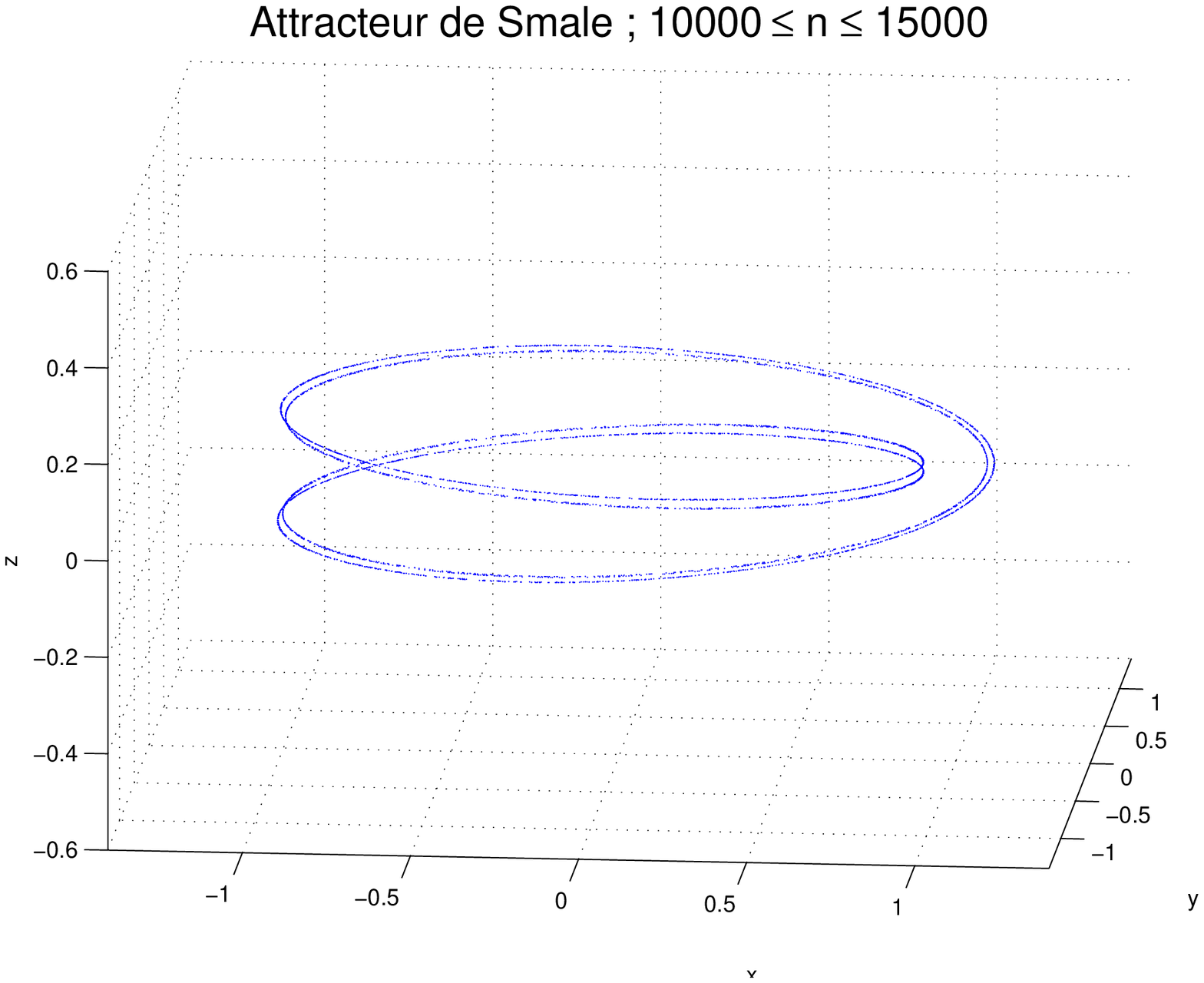}
     &
    \includegraphics[width=7cm]{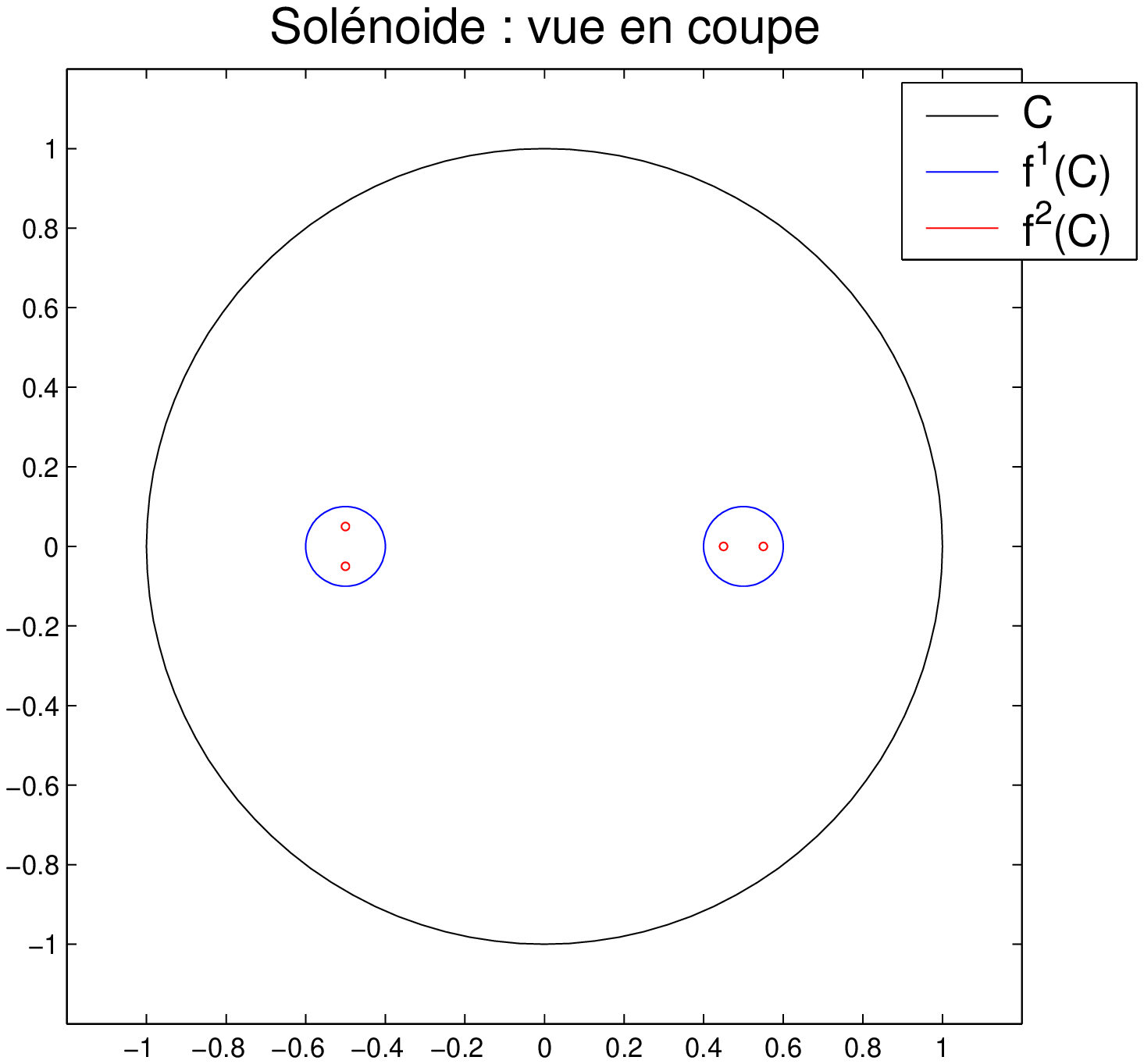}
     \\
     (a) $\Lambda$ & (b) coupe de $\Lambda$
\end{tabular}
\caption{\label{fig:solenoide_attr} Le sol\'eno\"ide.}
\end{center}
\end{figure}

Soit le tore $T = \{ (\theta, z) \telque \theta \in \R/\Z, z \in
\C,
\absj{z} \leq 1 \}$ et l'application \[ f : \begin{CD} T @>>> T \\
(\theta,z) @>>> \left(2\theta \text{ mod. 1} , \frac{1}{2} e^{2
\pi i \theta} + \frac{1}{10} z \right) \end{CD}
\] Le tore et son image par $f$ sont repr\'esent\'es figure~\ref{fig:solenoide_T_fT}. L'ensemble maximal invariant de $f$ est
$\Lambda = \bigcap _{n \in \N} f^n(T)$, et il est
hyperbolique\footnote{Les sous-espaces $E^s_x$ sont les plans
$\theta = cte$, les $E^u_x$ sont de dimension 1, tangents \`a
$\Lambda$.}, on l'appelle \emph{attracteur de Smale}
(figure~\ref{fig:solenoide_attr}a). Localement, $\Lambda$ est le
produit d'un ensemble de Cantor dyadique (inclus dans le disque
unit\'e de $\R^2$) et d'une droite
(figure~\ref{fig:solenoide_attr}b).

Les vari\'et\'es stables sont les sections $C_{\theta_0} =
\{(\theta ,z) \in \Lambda | \theta = \theta_0 \}$. Les
vari\'et\'es instables sont plus difficiles \`a d\'ecrire, et on
peut montrer que chacune est dense dans $\Lambda$.

Dans le cas du sol\'eno\"ide, il est \'egalement possible de
montrer l'existence d'une mesure physique (d\'efinition
\ref{def:mesure_physique}) qui donne la r\'epartition statistique
dans $\Lambda$ des points de presque toutes les orbites.

\subsubsection{Dynamique non-uniform\'ement hyperbolique}
\label{annexe:Henon} Consid\'erons l'exemple de l'attracteur de
H\'enon. Soit
l'application \[ H = H_{b,c} : \begin{CD} \R^2 @>>> \R^2 \\
(x,y) @>>> (x^2+c-by,x) \end{CD}
\] avec $0<b\ll 1$ et $c$ un peu plus grand que -2.

\begin{figure}
\begin{center}
\begin{tabular}{c@{}c@{}c}
    \includegraphics[width=4.5cm]{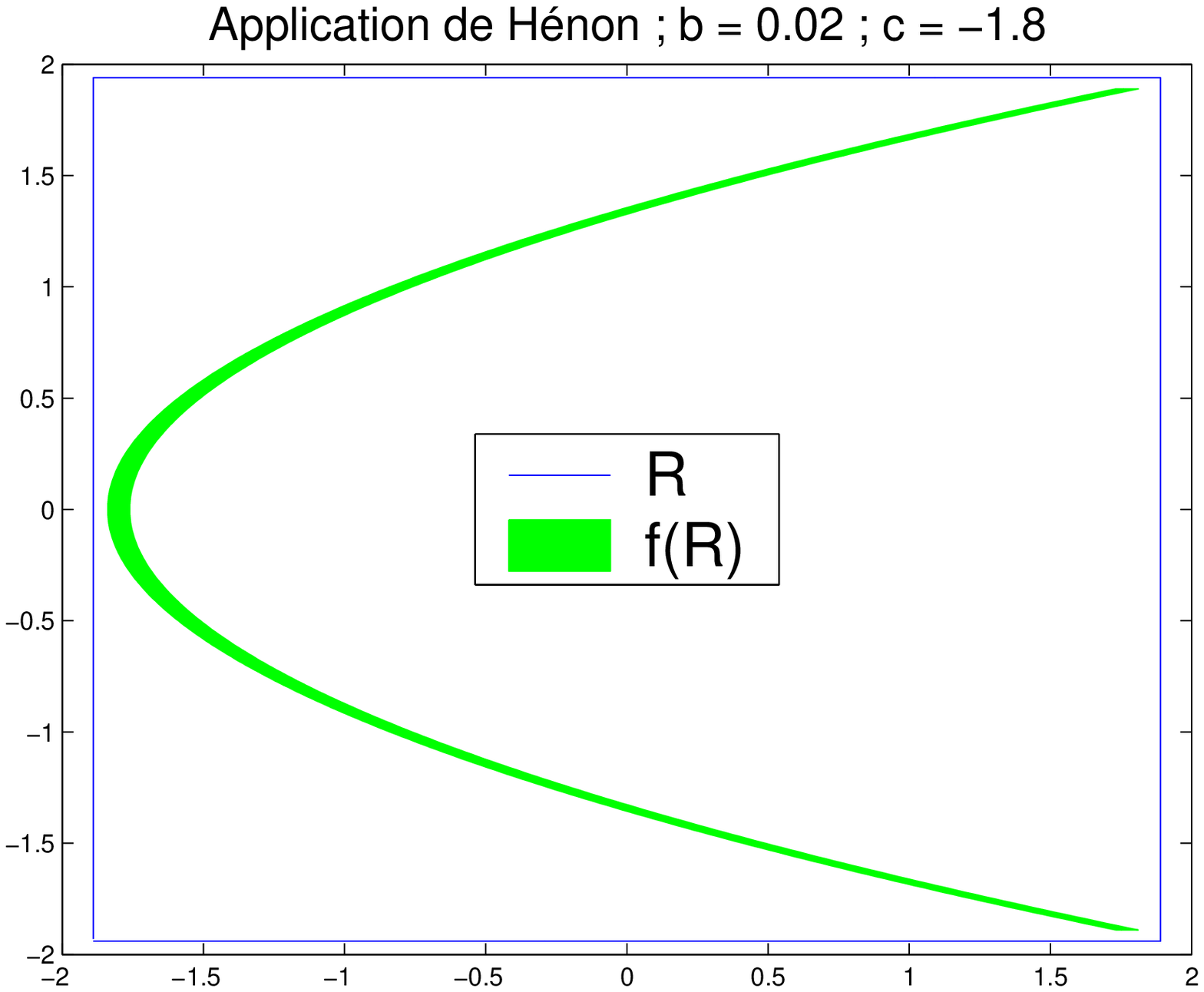} &
    \includegraphics[width=4.5cm]{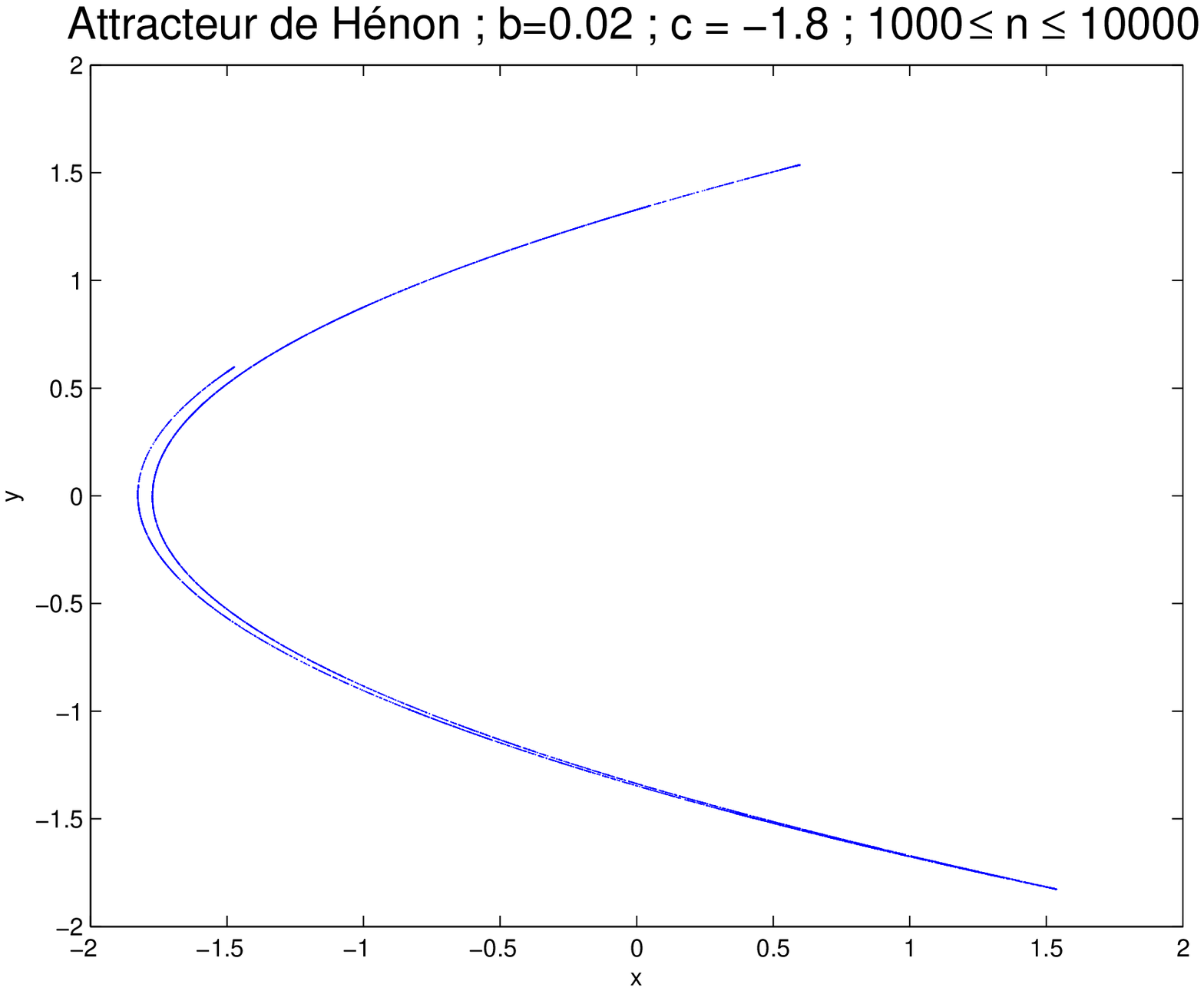} &
    \includegraphics[width=4.5cm]{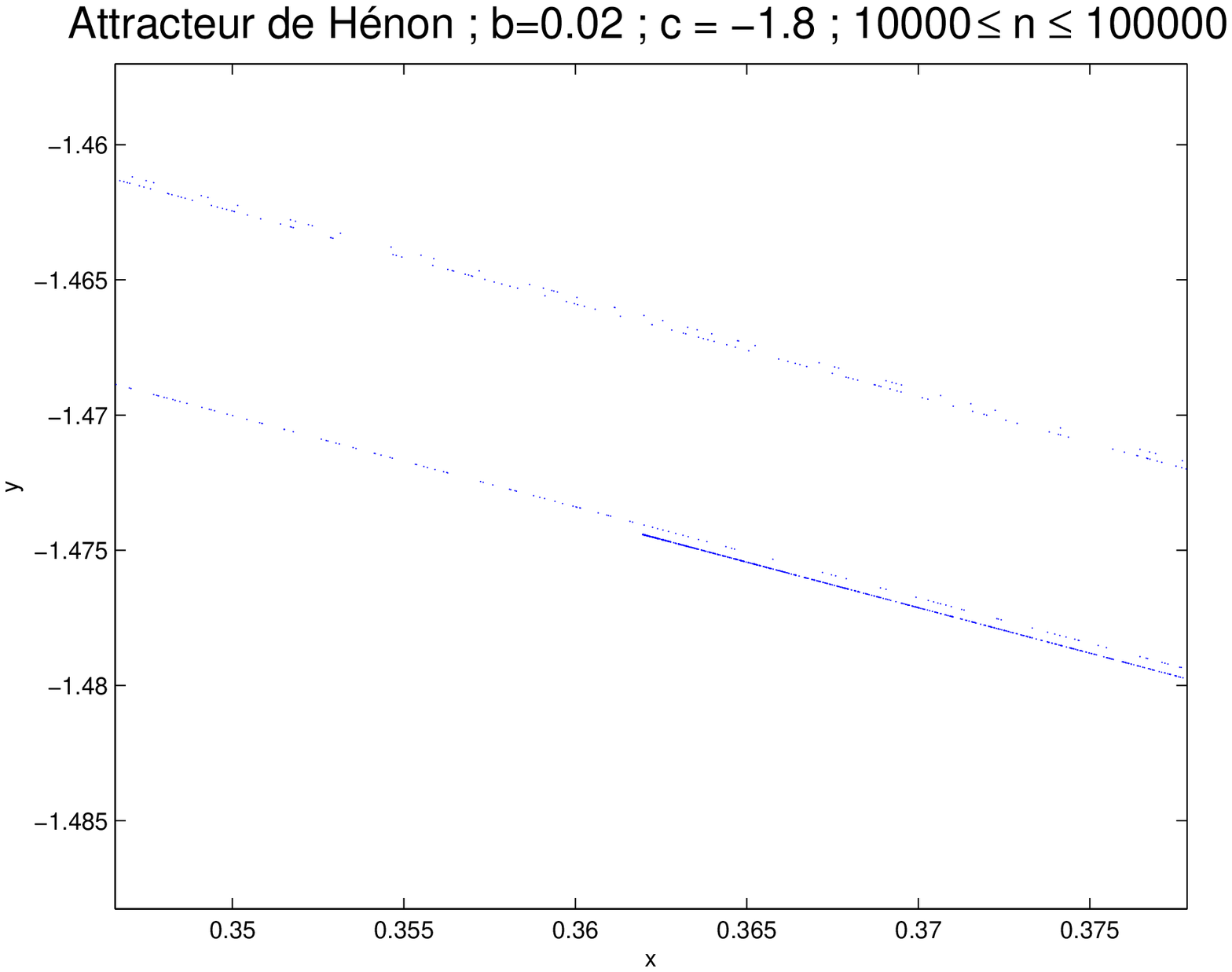}
     \\
     (a) $R$ et $f(R)$ & (b) $\Lambda$ & (c) zoom dans $\Lambda$
\end{tabular}
\caption{\label{fig:henon} L'application de H\'enon et son
attracteur ($b=0.02$, $c=-1.8$).}
\end{center}
\end{figure}

Il existe un rectangle $R$ tel que $H(R) \subset R$ et $H(R)$
ressemble \`a un arc de parabole <<\'epaissi>>
(figure~\ref{fig:henon}a) : le rectangle est fortement pinc\'e,
\'etir\'e, et pli\'e (dans le cas du sol\'eno\"ide, il n'y avait
pas de pli). L'attracteur est $\Lambda = \bigcup_{n \in \N}
H^n(R)$, et poss\`ede en presque tout point une structure de
Cantor $\times$ une droite (comme le sol\'eno\"ide)
(figure~\ref{fig:henon}a). En revanche, il existe un ensemble de
points (dense dans $\Lambda$) o\`u ce n'est pas le cas : ce sont
les <<pointes>>. Ainsi, figure~\ref{fig:henon}c, on visualise une
petite zone de l'attracteur, qui semblait rectiligne sur la vue
d'ensemble, et on distingue une pointe, \latin{i.e.} un filament
qui ne se poursuit pas vers la gauche.

On a alors : \begin{itemize} \item une d\'ecomposition $\R^2 =
E^s_x \oplus E^u_x$ en \emph{presque tous} les points $x \in
\Lambda$, la contraction et l'expansion n'\'etant \emph{pas
uniformes}, et les espaces $E^s_x$ et $E^u_x$ \emph{ne d\'ependent
pas contin\^ument de $x$}.
\item sensibilit\'e aux conditions initiales. \item le lemme
de pistage reste valide pour \emph{presques toutes} les
pseudo-orbites (mais pas toutes). \item \emph{pas de stabilit\'e
structurelle} : avec des perturbations arbitrairement petites de
$c$, on peut obtenir une dynamique compl\`etement diff\'erente.
Cependant, pour presque toute perturbation, on a un attracteur du
m\^eme type.
\item la mesure physique existe.
\end{itemize}

\subsection{Dynamique symbolique}
La dynamique symbolique fournit un exemple de syst\`eme dynamique
tr\`es important car il permet de mod\'eliser la dynamique de
tr\`es larges classes de syst\`emes. Nous la d\'ecrivons ici
tr\`es bri\`evement. On trouvera une introduction plus compl\`ete
dans \cite{Benoistpaulin:systemesdynamiques} et
\cite{KatokHasselblatt:IntroDynamicalSystems}.

Soit $\mathcal{A}$ un alphabet (\latin{i.e.} un ensemble) fini,
muni de la topologie discr\`ete, et $X = \mathcal{A}^{\Z}$ muni de
la topologie produit, l'ensemble des mots sur $\mathcal{A}$. C'est
un espace m\'etrique compact, muni de la distance
\[ d(\omega,\omega^{\prime}) = \sup_{i \in \Z} \left( 2^{-\absj{i}} \1_{\omega_i
\neq \omega^{\prime}_i} \right). \] On note $\sigma : X
\rightarrow X$ et on appelle \emph{d\'ecalage \`a
gauche}\footnote{shift en anglais} l'application d\'efinie par
$\sigma(\omega_i) = \omega_{i+1}$. Le syst\`eme dynamique
$(X,\sigma)$ est appel\'e \emph{syst\`eme de Bernoulli} ou
\emph{syst\`eme symbolique}.

Soit $\mathcal{B} \subset \mathcal{A}^2$ un ensemble de r\`egles
de compatibilit\'e. On appelle \emph{sous-d\'ecalage de type fini}
une partie de $X$ \[ \Sigma = \{ (\omega_i)_{i \in \Z} \telque
(\omega_i,\omega_{i+1}) \in \mathcal{B}, \, \forall i \in \Z \}.
\] $\Sigma$ est invariante par $\sigma$.

Ces notions permettent de coder de fa\c{c}on combinatoire la
dynamique de certains syst\`emes, \`a l'aide de partitions de
Markov. Une \emph{partition de Markov} d'un ensemble invariant
compact hyperbolique localement maximal $\Lambda$ est un
recouvrement fini par des rectangles propres\footnote{Un rectangle
$R$ est une partie de $\Lambda$ de diam\`etre plus petit que
$\epsilon/10$ ($\epsilon$ \'etant donn\'e par la proposition
\ref{pro:prod_local}) et telle que $[x,y] \in R$ si $x,y \in R$.
Il est propre si $R = \adherence{\interieur{R}}$.} $(R_a)_{a \in
\mathcal{A}}$ d'int\'erieurs disjoints et tels que si $x \in
\interieur{R_{a}}$ et $f(x) \in \interieur{R_b}$, alors
$f(W^s_{R_a}(x)) \subset W^s_{R_b}(f(x))$ et $W^u_{R_b}(f(x))
\subset f(W^u_{R_a}(x))$.
On d\'efinit alors un ensemble de transitions \[ \mathcal{B} = \{
(a,b) \in \mathcal{A}^2 \telque f(\interieur{R_a}) \cap
\interieur{R_b} \neq \emptyset \} \] et $\Sigma$ le
sous-d\'ecalage de type fini associ\'e. Le th\'eor\`eme suivant
pr\'ecise dans quelle mesure $(\Sigma,\sigma)$ <<code>> la
dynamique de $f$ sur $\Lambda$.

\begin{The}
Soient $\mathcal{R} = (R_a)_{a \in \mathcal{A}}$ une partition de
Markov de $\Lambda$ pour $f$ et $(\Sigma,\sigma)$ le
sous-d\'ecalage de type fini associ\'e. \begin{enumerate} \item
Pour tout $\omega \in \Sigma$, l'intersection $\bigcap_{n \in \Z}
f^{-n}(R_{\omega_n})$ est r\'eduite \`a un point $\pi(\omega)$.
\item L'application $\pi : \Sigma \rightarrow \Lambda$ est
continue, surjective, et le diagramme suivant commute : \[
\begin{CD} \Sigma @>\sigma>> \Sigma \\@V{\pi}VV @VV{\pi}V \\ \Lambda @>f>> \Lambda \end{CD} \]
\item Pour toute mesure de probabilit\'e $\sigma$-invariante et
ergodique $\mu$, de support $\Sigma$, on a \[ \mu \left( \{ \omega
\in \Sigma \telque \card \pi^{-1}(\pi(\omega))>1 \} \right) = 0
.\]
\end{enumerate}
\end{The}

Un ensemble basique $\Lambda_i$ de la d\'ecomposition spectrale
(Th\'eor\`eme \ref{the:decompo_spectrale}) poss\`ede une partition
de Markov de diam\`etre arbitrairement petit (voir
\cite{KatokHasselblatt:IntroDynamicalSystems}).

\subsection{Chaos} \label{annexe:chaos}
La notion de <<chaos>> en syst\`emes dynamiques, contrairement \`a
sa signification usuelle de d\'esordre total, se r\'ef\`ere \`a
une situation o\`u les orbites ne convergent pas vers une orbite
p\'eriodique ou quasi-p\'eriodique, et o\`u l'\'evolution des
orbites est impr\'evisible \`a un certain point, ou leur
comportement est sensible aux conditions initiales. Les premiers
exemples \'etudi\'es furent --- entre autres --- l'attracteur de
Lorenz, l'application logistique et l'application de H\'enon.

\begin{Def}[Orbite chaotique] \label{def:orbite_chaotique}
L'orbite de $x$, $\{ f^n(x) / n \geq 0 \}$, est \emph{sensible}
(ou \emph{chaotique}) s'il existe une constante $C>0$ telle que
\begin{equation} \begin{split} \forall q \in \omega(x), \, \forall
\epsilon >0 , \, \exists n_1, n_2, n > 0 / \, & d(f^{n_1}(x),q)<
\epsilon , \\ d(f^{n_2}(x),q)< \epsilon &\text{ et }
d(f^{n_1+n}(x),f^{n_2+n}(x))> C.
\end{split} \end{equation}
\end{Def}

Une orbite asymptotique \`a une orbite p\'eriodique ou
quasi-p\'eriodique n'est pas chaotique au sens o\`u si
$f^{n_1}(x)$ et $f^{n_2}(x)$ sont proches, alors $f^{n_1+n}(x)$ et
$f^{n_2+n}(x)$ restent proches pour tout $n \geq 0$.

Une orbite sensible est \'egalement impr\'evisible dans la mesure
o\`u savoir qu'un point $y$ de l'orbite est extr\^emement proche
de $q \in \omega(x)$ n'est pas suffisant pour pr\'edire le futur
de $y$ \`a une distance $C$ pr\`es.

Dans l'ensemble stable d'un attracteur hyperbolique non-trivial,
de m\^eme que l'on a une forte sensibilit\'e aux conditions
initiales\footnote{Comme l'indique la propri\'et\'e
d'expansivit\'e \ref{pro:expansivite}}, on peut montrer que
l'ensemble des points ayant une orbite chaotique a une mesure de
Lebesgue totale.

\begin{Def}[Dynamique chaotique] \label{def:dynamique_chaotique}
Un syst\`eme dynamique $(X,f)$ est \emph{sensible} (ou a une
\emph{dynamique chaotique}) lorsque l'ensemble des points ayant
une orbite chaotique a une mesure de Lebesgue
non-nulle\footnote{Cette d\'efinition n'a de sens que lorsque $X$
est une vari\'et\'e, pour que les ensembles de mesure de Lebesgue
nulle soient d\'efinis.}.
\end{Def}

Cependant, le chaos ainsi d\'efini ne doit pas \^etre
interpr\'et\'e comme une totale impr\'edictibilit\'e. En effet, on
observe num\'eriquement, pour certains syst\`emes chaotiques, que
pour toute condition initiale prise dans un certain ouvert, on
obtient le m\^eme ensemble $\omega$-limite. Ceci conduit \`a la
notion d'attracteur \'etrange.

\begin{Def}[Attracteur \'etrange] Une partie compacte $A$ de $X$
est un \emph{attracteur \'etrange} s'il existe un ouvert $U$ et $N
\subset U$ de mesure de Lebesgue nulle tel que $\forall x \in U
\backslash N$, $\omega(x) = A$ et l'orbite de $x$ est chaotique.
\end{Def}

Un exemple d'attracteur \'etrange est l'\emph{attracteur de
H\'enon} (section \ref{annexe:Henon}). On appelle parfois
\'egalement attracteur \'etrange un attracteur $A$ tel que $f$ a une
d\'ependance sensible aux conditions initiales avec probabilit\'e
totale sur $B(A) \times B(A)$ (o\`u $B(A)$ est le bassin
d'attraction de $A$ : voir d\'efinition \ref{def:bassin}).

Une derni\`ere notion importante est celle de \emph{dynamique
chaotique persistante}, qui traduit que de petites perturbations
de $f$ ont, avec une probabilit\'e positive, une dynamique
chaotique. Cette d\'efinition a un sens lorsque par exemple $f =
f_{\alpha}$ est param\'etr\'ee par $\alpha \in \R^n$, car alors on
dispose de la mesure de Lebesgue sur l'espace des param\`etre
$\alpha$. De fa\c{c}on plus restrictive, on peut demander la
persistance d'une dynamique chaotique dans un voisinage ouvert de
$f$.

Une notion que nous avons d\'ej\`a introduite est \'etroitement
reli\'ee au chaos. Il s'agit de celle d'\emph{intersection
homocline} (d\'efinition \ref{def:homocline}). Il y a
\'equivalence entre l'existence d'une orbite chaotique (voir
section \ref{annexe:chaos}, d\'efinition
\ref{def:orbite_chaotique}) et l'existence d'une orbite homocline.
En revanche, cela n'entra\^ine pas forc\'ement que la dynamique
est chaotique (d\'efinition \ref{def:dynamique_chaotique}).

\paragraph{Chaos et simulations num\'eriques} Il est
probl\'ematique de vouloir observer ou m\^eme caract\'eriser un
comportement chaotique lors d'une simulation num\'erique. Comment
en effet mettre en \'evidence un tel ph\'enom\`ene malgr\'e la
pr\'ecision finie d'un ordinateur ? Celle-ci a plusieurs
cons\'equences majeures.

Tout d'abord, les erreurs d'arrondi font que l'on n'observe que
des pseudo-orbites. Si le syst\`eme \'etudi\'e poss\`ede une
propri\'ete de pistage, comme c'est le cas avec les syst\`emes
uniform\'ement hyperboliques, on a de quoi \^etre partiellement
rassur\'e. Il reste cependant des cas (par exemple le doublement
de l'angle) o\`u les orbites qu'un ordinateur peut pister ne sont
pas des orbites typiques du syst\`eme. De m\^eme, lorsque les
orbites calcul\'ees sont born\'ees, toutes les pseudo-orbites
observ\'ees sont p\'eriodiques (m\^eme si la p\'eriode est tr\`es
longue), en raison du nombre fini de d\'ecimales que l'on peut
calculer. Il faut donc fixer (arbitrairement) un seuil pour
s\'eparer orbites p\'eriodiques et non-p\'eriodiques.

Un deuxi\`eme effet est que l'on ne peut observer que le
comportement en temps fini. Comment alors \^etre s\^urs qu'il
s'agit bien du comportement stationnaire, et non d'un r\'egime
transitoire tr\`es long ? Il nous faut en effet fixer un seuil \`a
partir duquel on observe la dynamique <<\`a l'infini>>. Le choix
de ce seuil est crucial pour \'eviter des erreurs, tout en
limitant la dur\'ee des calculs.

Enfin, lorsque l'on \'etudie un syst\`eme d\'ependant de
param\`etres r\'eels, il faut garder \`a l'esprit que l'on ne peut
observer celui-ci que sur un ensemble de mesure nul, l'ensemble
des rationnels. C'est tout l'int\'er\^et de consid\'erer la
persistance de la dynamique dans un voisinage ouvert, $\Q$ \'etant
dense dans $\R$. Ce probl\`eme peut cependant se ramener \`a celui
du lien entre pseudo-orbites et vraies orbites si la famille
$(f_{\alpha})_{\alpha \in \R}$ d\'epend contin\^ument de $\alpha$
pour la topologie de la convergence uniforme sur $X$, car alors
une orbite sous $f_{\alpha+\epsilon}$ est une pseudo-orbite sous
$f_{\alpha}$ si $\epsilon$ est assez petit.

\subsection{Bifurcations}
Consid\'erons une famille de syst\`emes dynamiques d\'ependant
d'un ou plusieurs param\`etres. M\^eme si pour presque toutes les
valeurs des param\`etres, le syst\`eme a un comportement
transverse (par exemple structurellement stable), il peut y avoir
des valeurs particuli\`eres de ceux-ci o\`u se produit une
transition entre deux diff\'erents types d'orbites. De tels
changements sont appel\'es \emph{bifurcations}. Leur \'etude ---
qui est une branche \`a part enti\`ere de la th\'eorie des
syst\`emes dynamiques --- est fondamentale pour comprendre les
propri\'et\'es d'un syst\`eme typique car les bifurcations
montrent comment diff\'erents comportements transverses peuvent
appara\^itre.

Nous ne parlons ici que de quelques cas simples de bifurcations,
en petite dimension. Il en existe bien s\^ur beaucoup d'autres
types. Nous nous limitons de plus \`a des bifurcations
\emph{locales}, c'est-\`a-dire pouvant \^etre d\'efinies seulement
au voisinage d'un point, par opposition aux bifurcations
\emph{globales}. Nous consid\'erons plus particuli\`erement le cas
des bifurcations structurellement stables, d\'efinies de la
mani\`ere suivante dans le cas de syst\`emes discrets.

\begin{Def}[Bifurcation structurellement stable] Une famille $\{
f_{\tau} \}$ de $C^{\infty}$ diff\'eomorphismes d\'efinis
localement a une \emph{bifurcation structurellement stable} \`a
$\tau = \tau_0$ si $f_{\tau_0}$ n'est pas localement
structurellement stable et si pour toute famille $\{ g_{\tau} \}$
de $C^{\infty}$ diff\'eomorphismes d\'efinis localement
suffisamment $C^2$-proche de $\{ f_{\tau} \}$, il existe une
reparam\'etrisation $\phi(\tau)$ de $\{ g_{\tau} \}$ et une
famille continue $\{ h_{\tau} \}$ d'hom\'eomorphismes d\'efinis
localement telle que \[ g_{\phi(\tau)} = h_{\tau}^{-1} \circ
f_{\tau} \circ h_{\tau} \] partout o\`u cela est d\'efini.
\end{Def}

\subsubsection{Diagramme de bifurcations}
Il existe un moyen simple de visualiser une bifurcation, appel\'e
diagramme de bifurcation. On trace l'ensemble $\omega$-limite $L^+
(f_{\epsilon})$ pour les diff\'erentes valeurs du param\`etre
$\epsilon$, que l'on porte sur l'axe des abscisses. Un tel
diagramme peut ais\'ement \^etre trac\'e num\'eriquement, en
prenant pour ensemble $\omega$-limite les valeurs de
$f^n_{\epsilon}(x)$ pour $n$ <<grand>> et pour un ou plusieurs $x$
choisis al\'eatoirement.

Il y a cependant une diff\'erence entre un diagramme obtenu par
simulations et un diagramme th\'eorique : les objets instables, ou
de <<petit>> bassin d'attraction, n'apparaissent que dans le
second cas. Il n'est ainsi pas forc\'ement simple de d\'eterminer
la nature d'une bifurcation en comparant son diagramme empirique
avec les diagrammes th\'eoriques des bifurcations classiques.

\subsubsection{Cas discret, dimension 1}\label{annexe:bif1}

En dimension 1, on peut classifier les bifurcations
structurellement stables autour d'un point d'\'equilibre $p$. En
effet, dans ce cas, la d\'eriv\'ee de $f_{\tau_0}$ en $p$ doit
valoir $\lambda = \pm 1$.

Commen\c{c}ons par le cas $\lambda = 1$. La famille
$(f_{(+1),\tau})_{\tau \in \R}$, d\'efinie par
\begin{equation} \label{eq:bif_dim1_+1} \forall x \in \R , \,
f_{(+1),\tau} (x) = x + x^2 + \tau \end{equation} a une
bifurcation structurellement stable en $\tau_0 = 0$, avec
d\'eriv\'ee 1, et est caract\'eristique de cette situation.
\begin{Pro}\label{pro:bif_dim1_+1}La bifurcation de la famille \eqref{eq:bif_dim1_+1}
en $\tau_0 = 0$ est structurellement stable, et toute bifurcation
locale structurellement stable en dimension 1 ayant lieu en un
point fixe avec d\'eriv\'ee 1 est (topologiquement) \'equivalente
(apr\`es reparam\'etrisation) \`a cette bifurcation. \end{Pro}

\begin{figure}
\begin{center}
\begin{tabular}{c@{}c@{}c}
     \includegraphics[width=4.5cm]{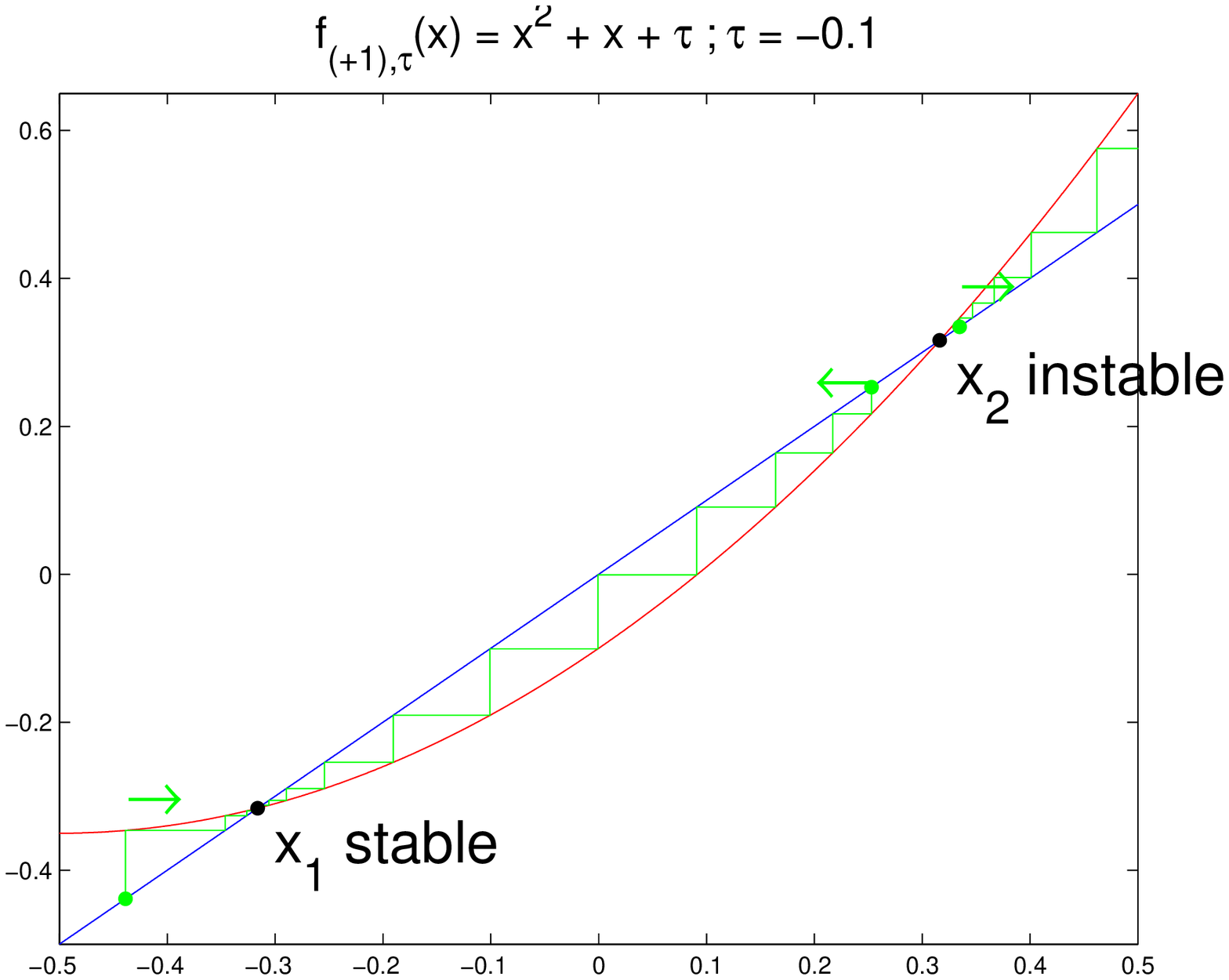} &
     \includegraphics[width=4.5cm]{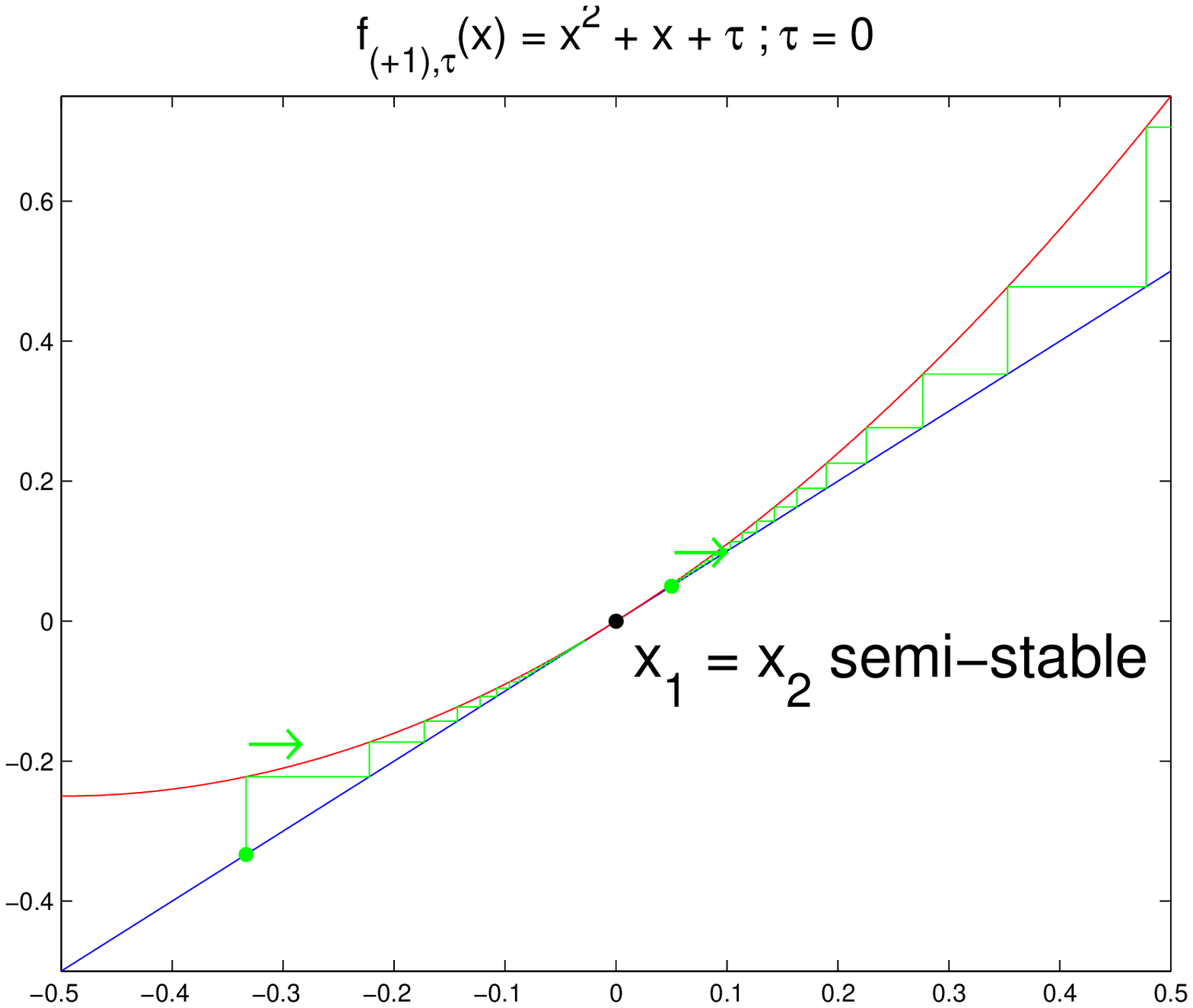} &
     \includegraphics[width=4.5cm]{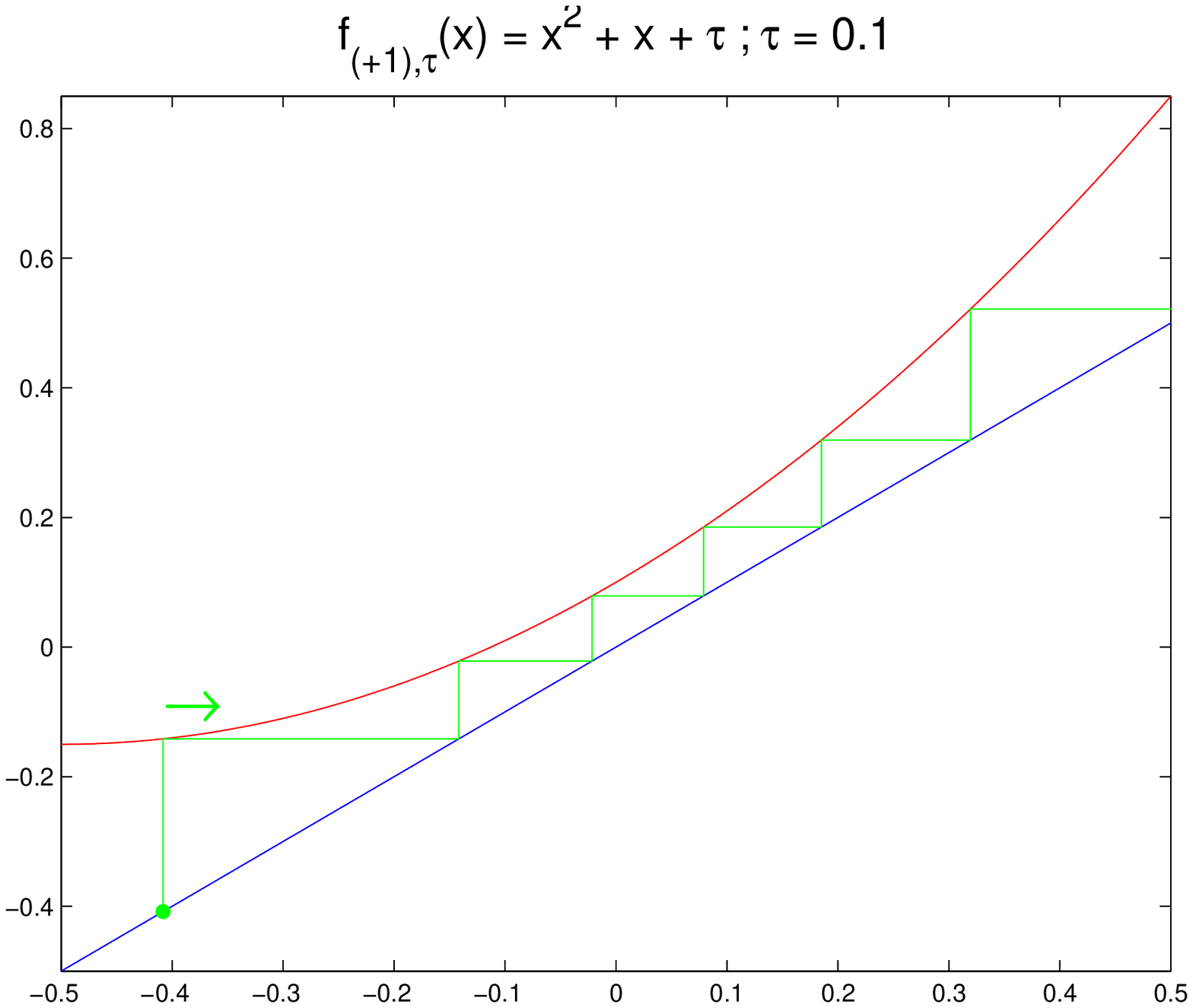}
     \\
     $\tau < \tau_0$ & $\tau = \tau_0$ & $\tau>\tau_0$
\end{tabular}
\caption{\label{fig:bif1_+1_f} Bifurcation de la famille
$f_{(+1),\tau}(x)$.}
\end{center}
\end{figure}

\begin{figure}
\begin{center}
     \includegraphics[width=7cm]{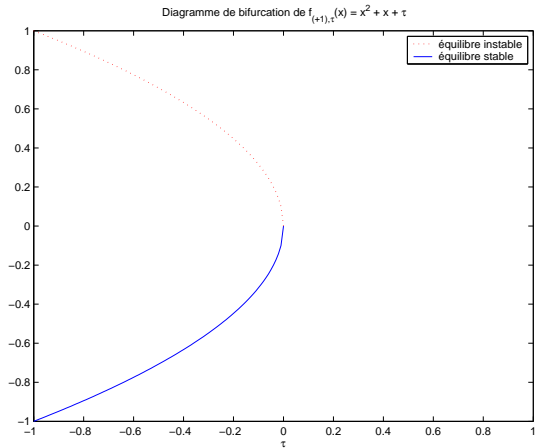}
\caption{\label{fig:bif1_+1_diag} Diagramme de bifurcation de la
famille $f_{(+1),\tau}(x)$, autour de $\tau_0 = 0$.}
\end{center}
\end{figure}

Ainsi, pour $\tau < \tau_0$, $f_{(+1),\tau}$ poss\`ede un point
fixe stable $-\sqrt{-\tau}$ et un point fixe instable
$\sqrt{-\tau}$ ; pour $\tau=\tau_0$, ces deux points fixes sont
confondus, et l'\'equilibre qui en r\'esulte est semi-stable ;
enfin, d\`es que $\tau > \tau_0$, $f_{(+1),\tau}$ n'a plus de
point fixe (figure~\ref{fig:bif1_+1_f}). Le diagramme de
bifurcation correspondant est repr\'esent\'e \`a la
figure~\ref{fig:bif1_+1_diag}.

\begin{figure}
\begin{center}
\begin{tabular}{c@{}c@{}c}
     \includegraphics[width=4.5cm]{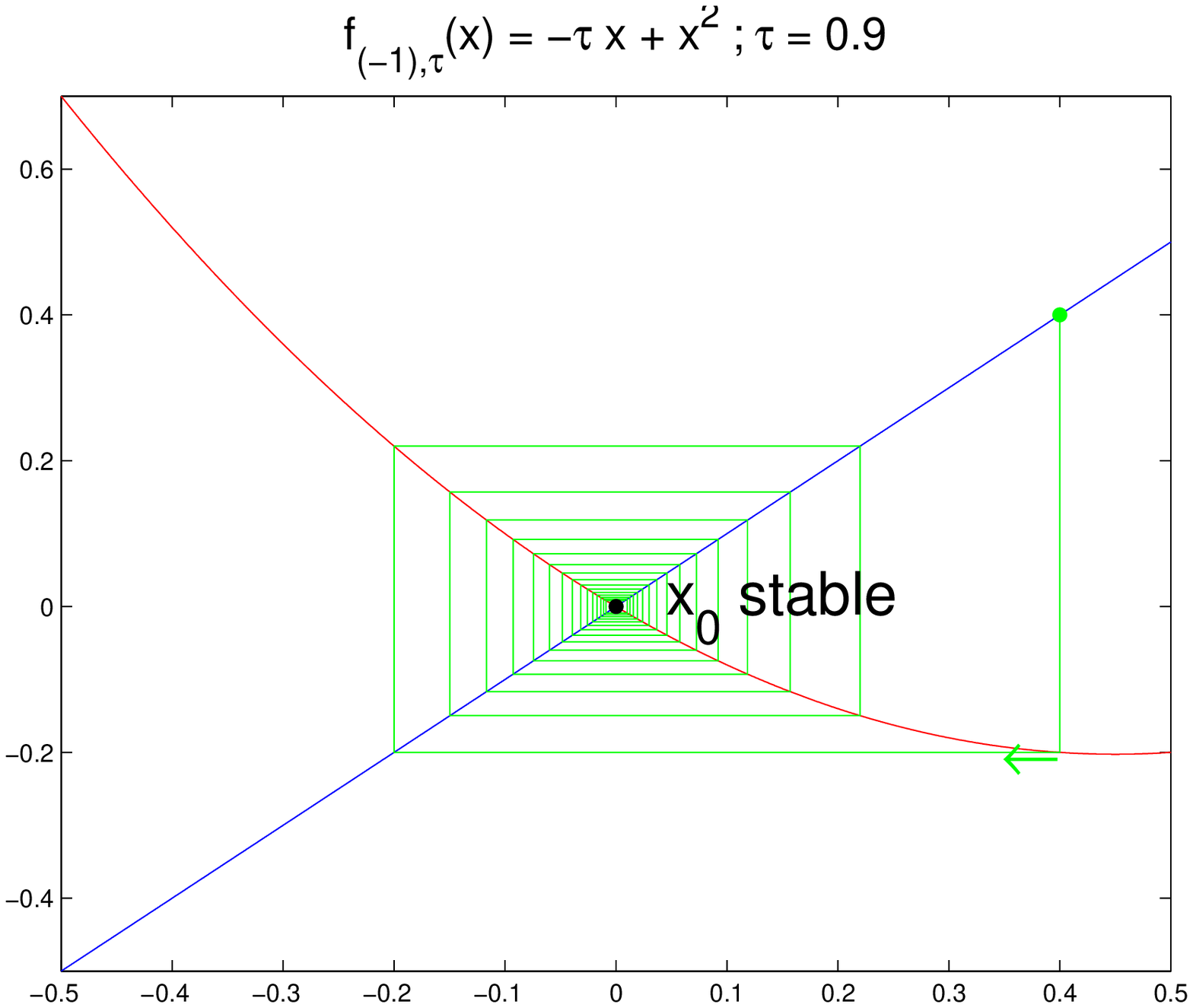} &
     \includegraphics[width=4.5cm]{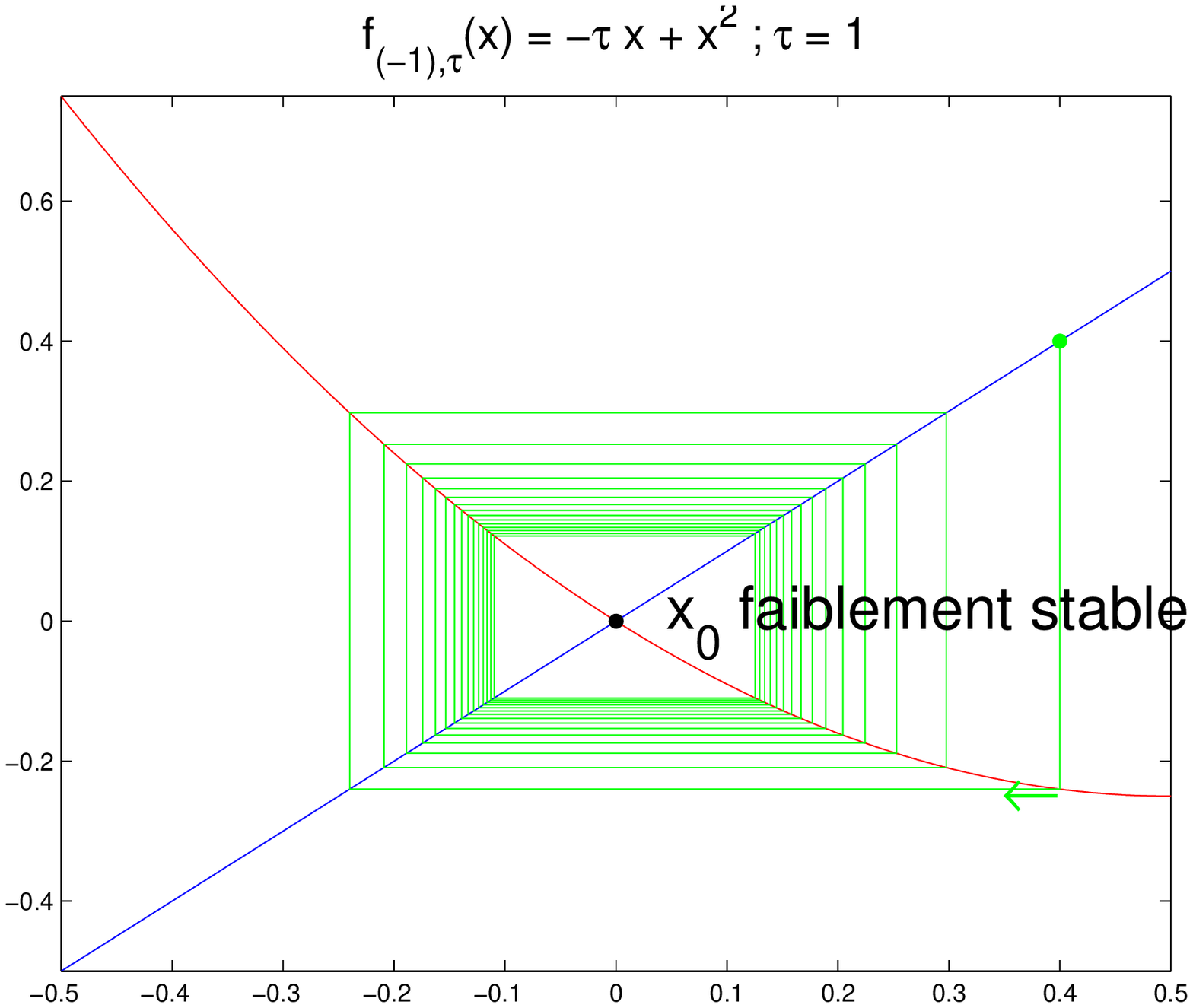} &
     \includegraphics[width=4.5cm]{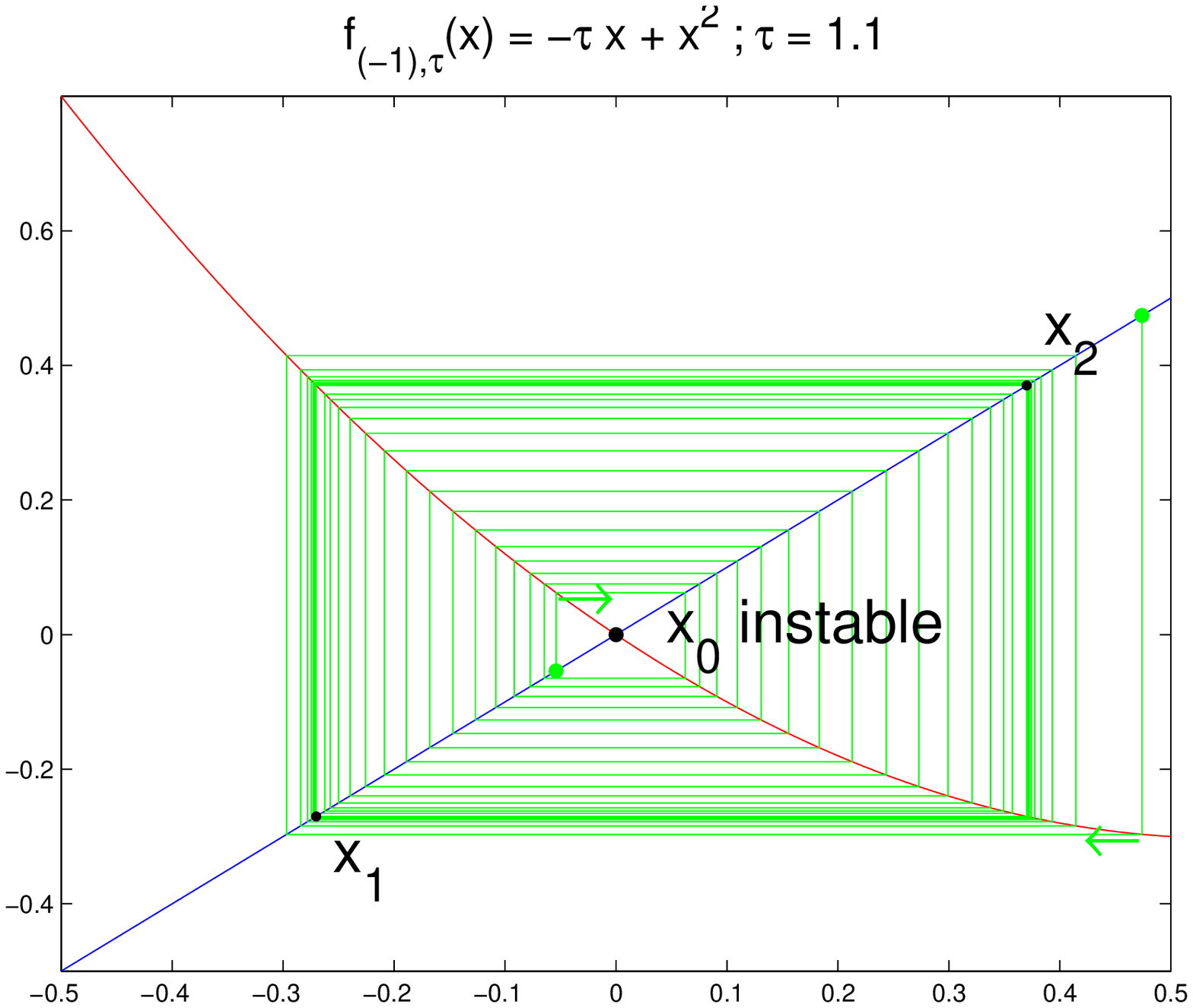}
     \\
     $\tau < \tau_0$ & $\tau = \tau_0$ & $\tau>\tau_0$
\end{tabular}
\caption{\label{fig:bif1_-1_f} Bifurcation subie par
$f_{(-1),\tau}(x)$.}
\end{center}
\end{figure}

\begin{figure}
\begin{center}
\begin{tabular}{c@{}c@{}c}
     \includegraphics[width=4.5cm]{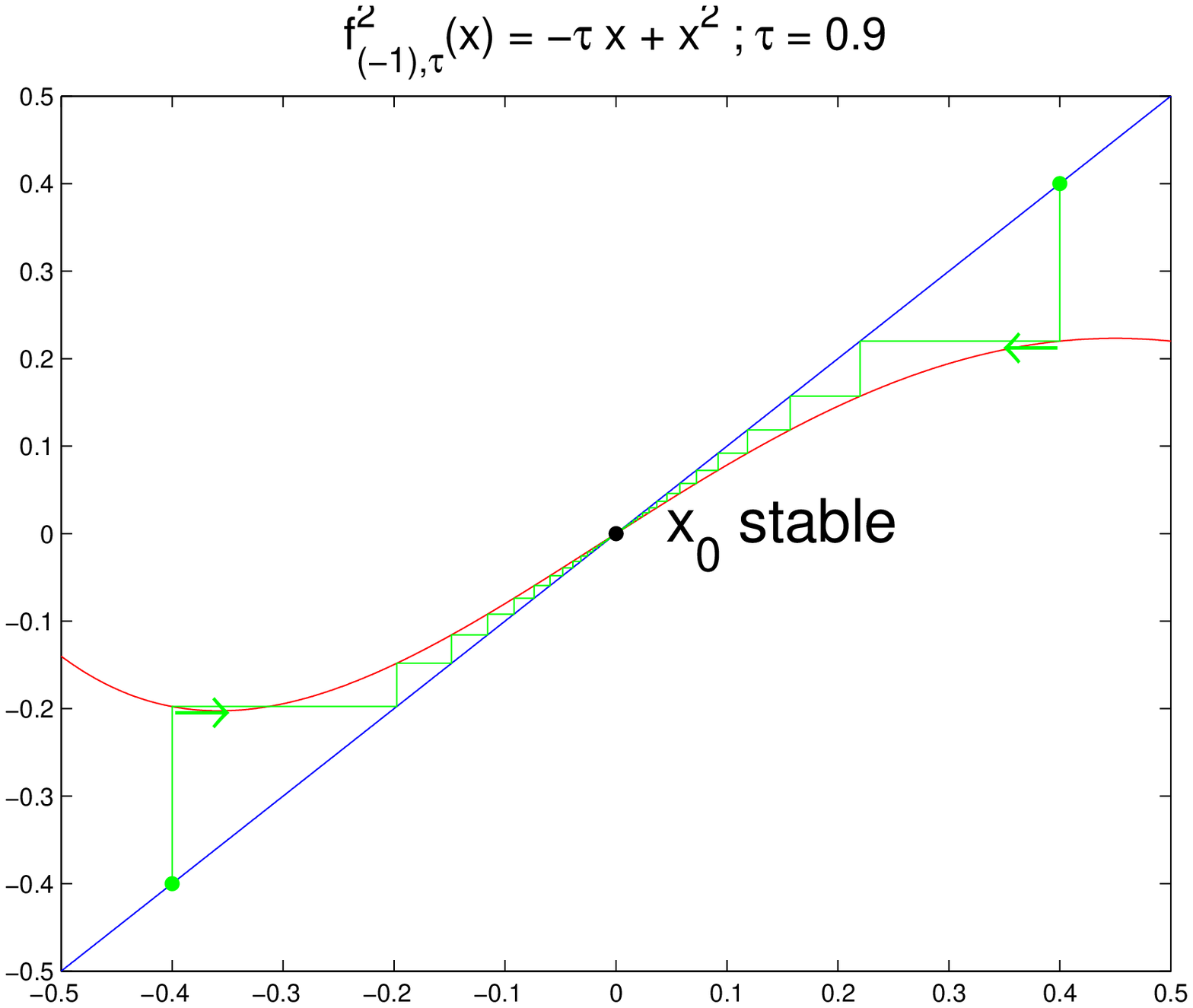} &
     \includegraphics[width=4.5cm]{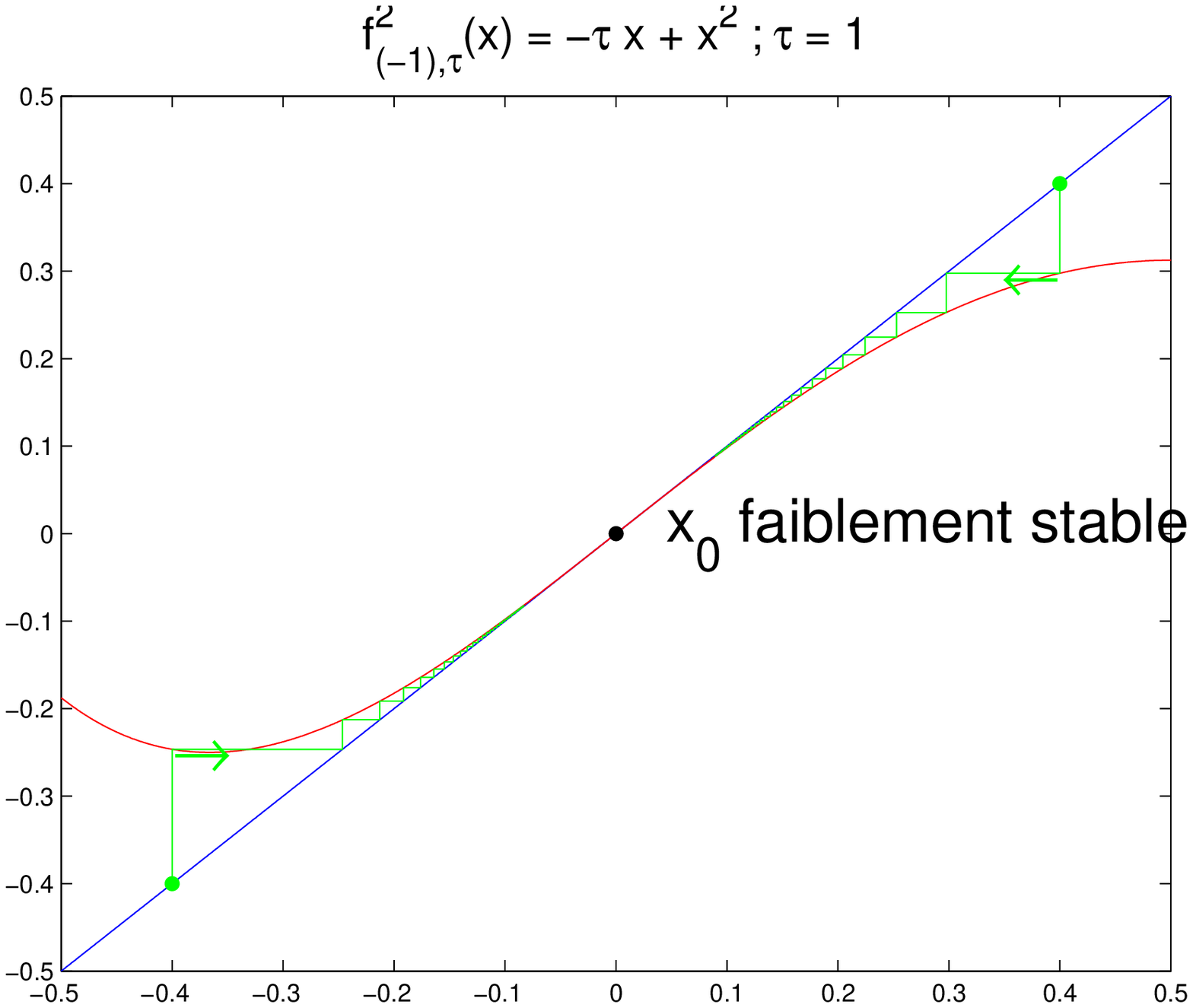} &
     \includegraphics[width=4.5cm]{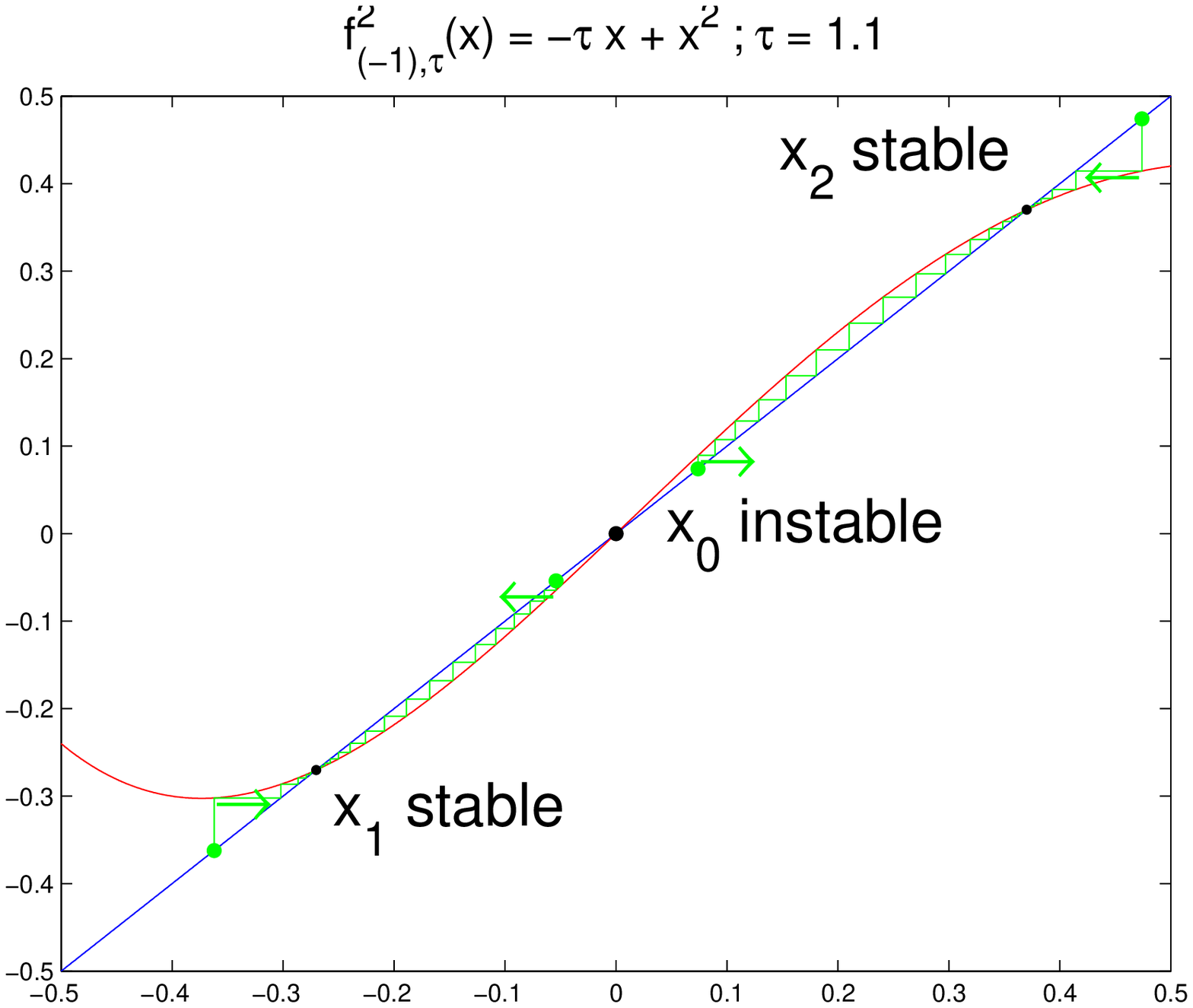}
     \\
     $\tau < \tau_0$ & $\tau = \tau_0$ & $\tau>\tau_0$
\end{tabular}
\caption{\label{fig:bif1_-1_f2} Bifurcation subie par
$f^{2}_{(-1),\tau}(x)$.}
\end{center}
\end{figure}

\begin{figure}
\begin{center}
     \includegraphics[width=7cm]{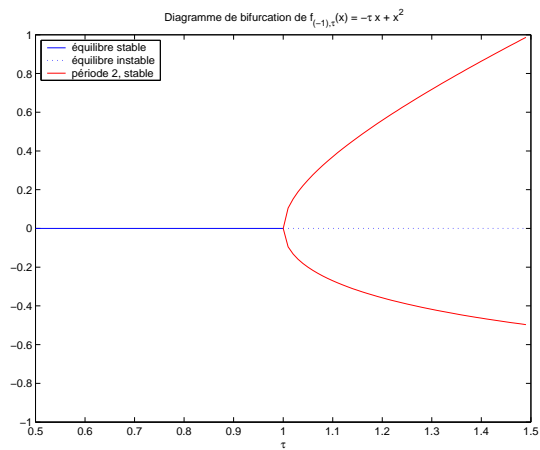}
\caption{\label{fig:bif1_-1_diag} Diagramme de bifurcation de la
famille $f_{(-1),\tau}(x)$, autour de $\tau_0 = 1$.}
\end{center}
\end{figure}

Dans le cas o\`u $\lambda = -1$, le point fixe $p$ est transverse
et donc persistent. La valeur de la d\'eriv\'ee en $p$ en
sup\'erieure \`a $-1$ pour $\tau < \tau_0$ et inf\'erieure \`a
$-1$ pour $\tau > \tau_0$, le point fixe restant isol\'e. Cela
s'accompagne de la cr\'eation d'une orbite stable de p\'eriode 2,
tandis que le point fixe devient instable. On parle de bifurcation
par \emph{doublement de p\'eriode}\label{def:doublement_periode},
dont l'exemple typique est le suivant : \begin{equation}
f_{(-1),\tau}(x) = -\tau x + x^2 \end{equation} au voisinage de
$x_0=0$, $\tau_0 = 1$. On montre alors une proposition similaire
\`a la proposition \ref{pro:bif_dim1_+1}, ce qui ach\`eve la
classification dans le cas de la dimension 1. Pour visualiser
cette bifurcation, on peut tracer $f_{(-1),\tau}$
(figure~\ref{fig:bif1_-1_f}), mais aussi $f^2_{(-1),\tau}$
(figure~\ref{fig:bif1_-1_f2}) pour mieux comprendre les orbites de
p\'eriode 2. Le diagramme de cette bifurcation est repr\'esent\'e
figure~\ref{fig:bif1_-1_diag}.


\subsubsection{Bifurcation selle-n\oe ud}
En dimension sup\'erieure, les bifurcations structurellement
stables se produisent lorsqu'une valeur propre vaut $\pm 1$ et les
autres sont en-dehors du cercle unit\'e.

Un exemple classique, en dimension deux ou plus, est le suivant :
deux point fixes, un noeud (point fixe attractif) et une selle
(attractif dans une direction, r\'epulsif dans une autre) se
rencontrent. Apr\`es bifurcation, il n'y a plus aucun point fixe
(localement). Une telle bifurcation est appel\'ee \emph{selle-n\oe
ud}\footnote{saddle-node en anglais}.

\begin{figure}
\begin{center}
\begin{tabular}{c@{}c@{}c}
     \includegraphics[width=4.5cm]{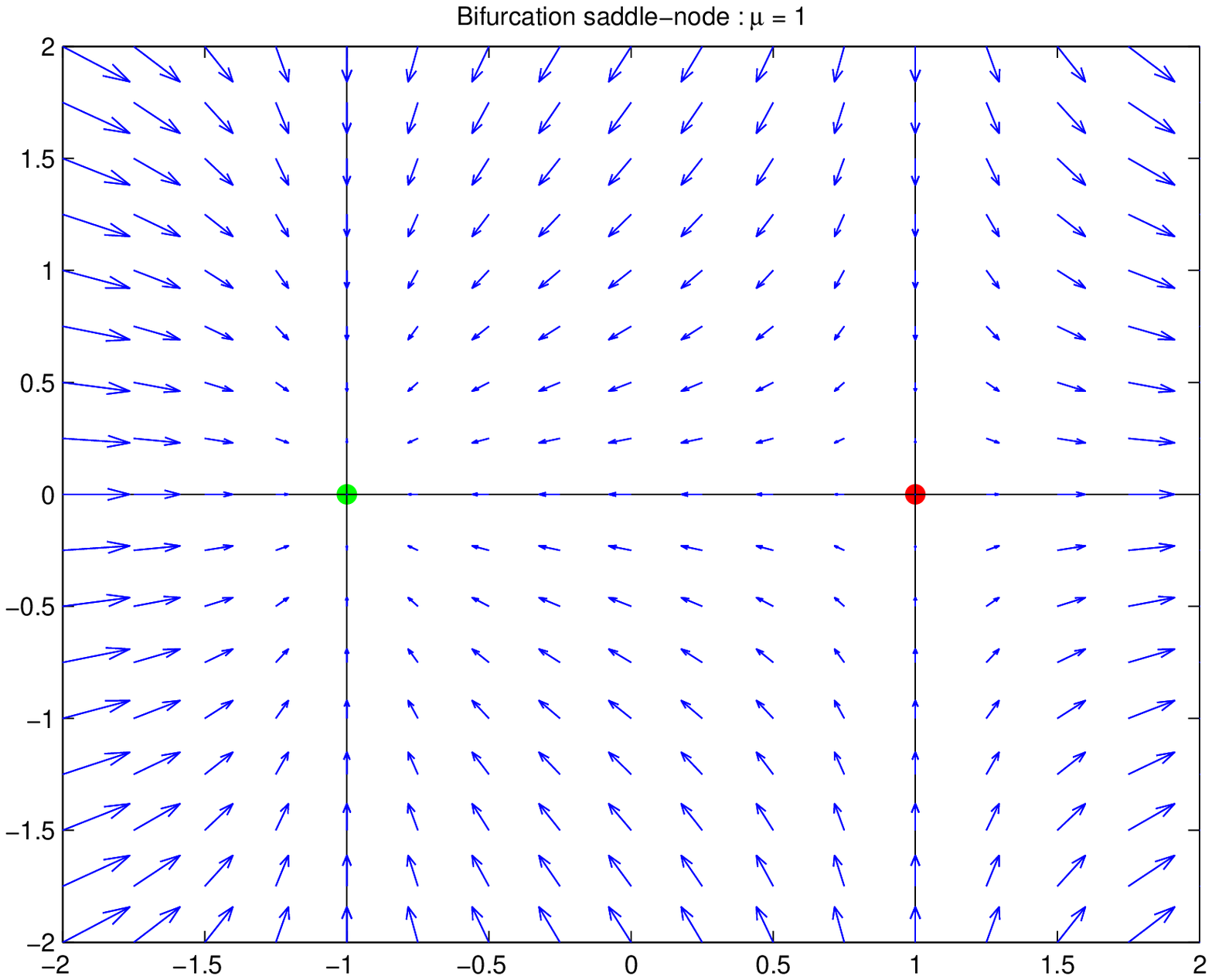}
     &
     \includegraphics[width=4.5cm]{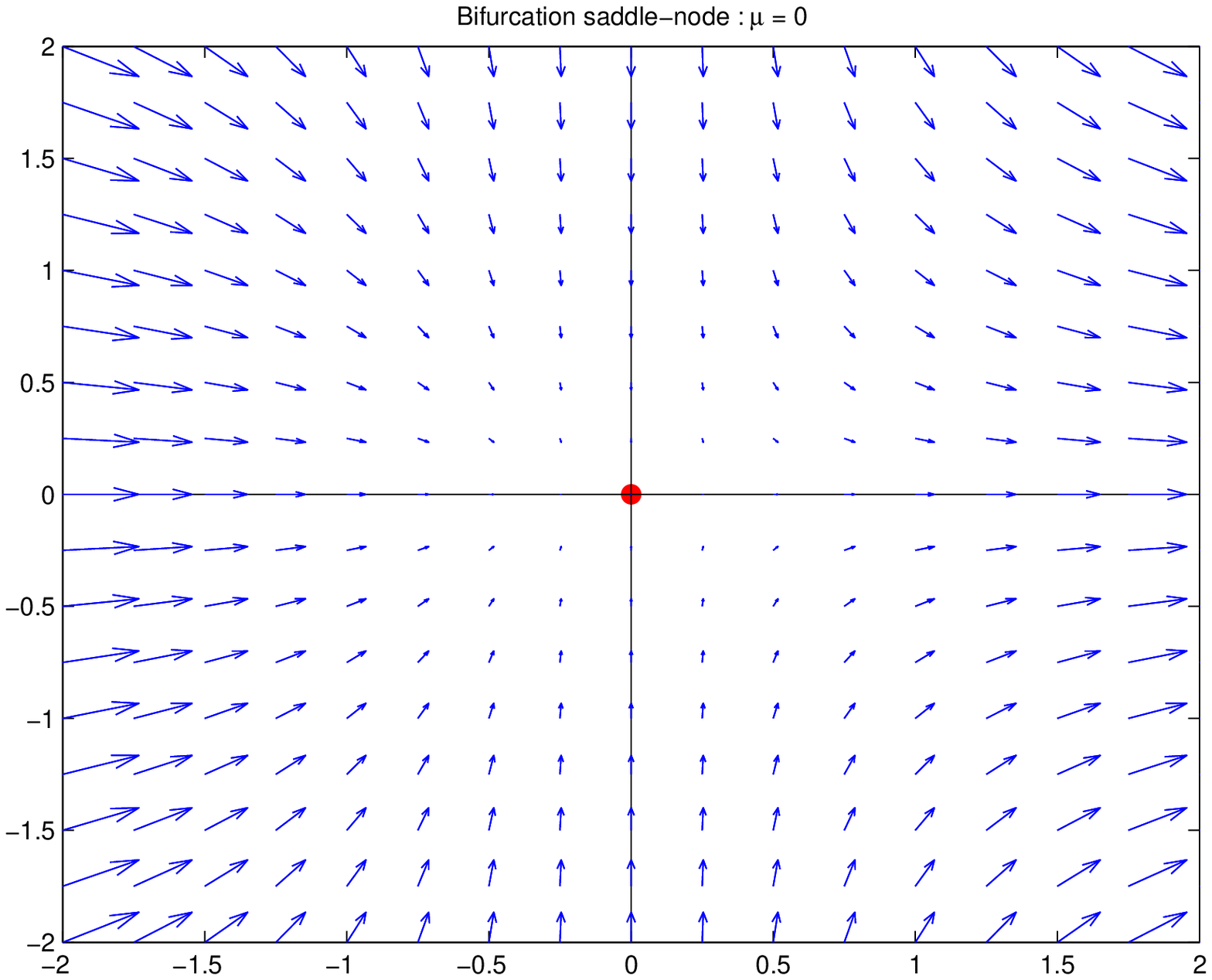}
     &
     \includegraphics[width=4.5cm]{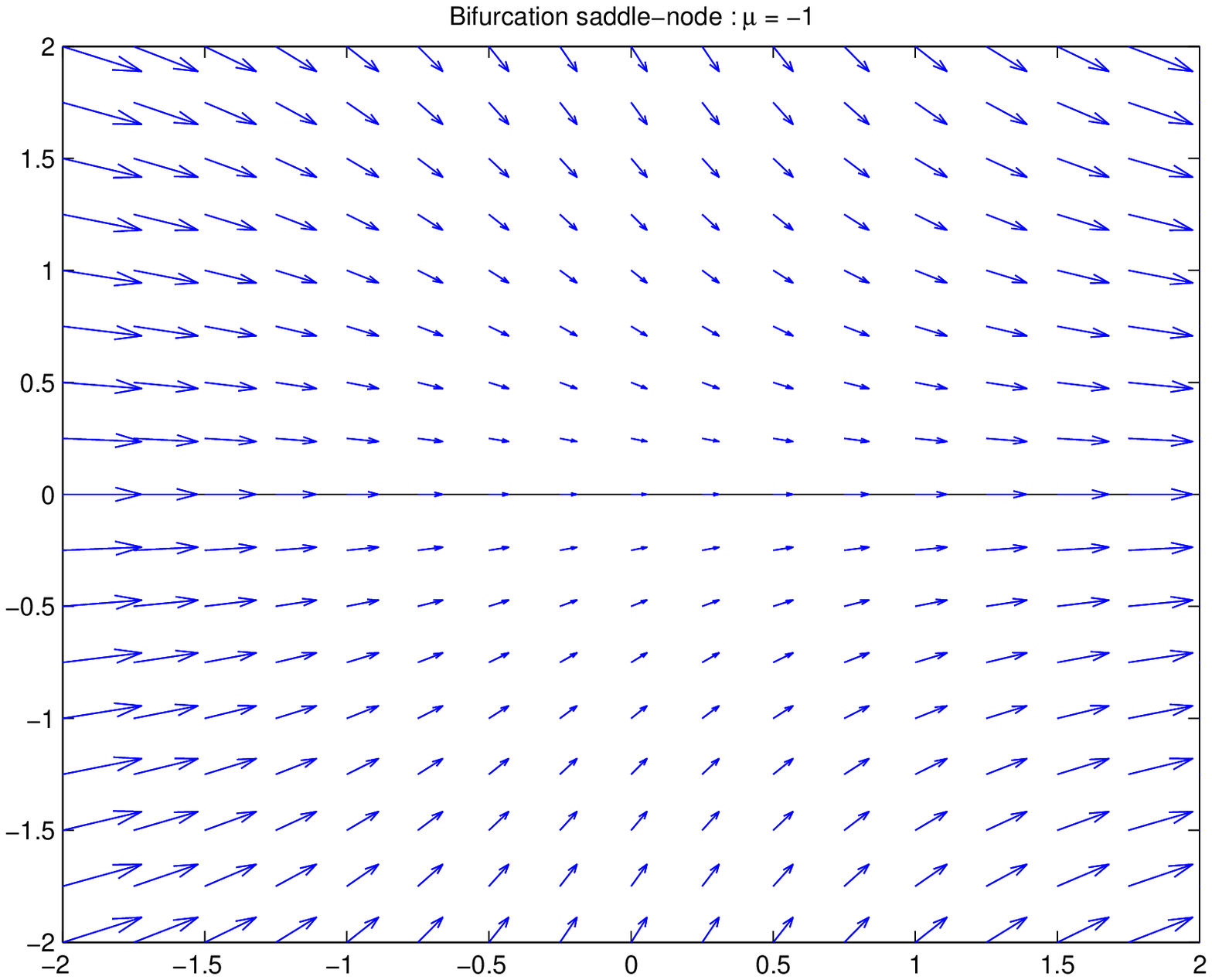}
     \\
     $\mu = 1 > \mu_0$ & $\mu = \mu_0$ & $\mu = -1 < \mu_0$
\end{tabular}
\caption{\label{fig:saddle_node_f} Bifurcation selle-noeud :
espace des phases de part et d'autre de $\mu_0=0$.}
\end{center}
\end{figure}

Le syst\`eme diff\'erentiel suivant donne un exemple de
bifurcation selle-n\oe ud :
\begin{equation}\label{eq:exemple_saddle_node}
\left\{
\begin{aligned}
 \frac{dx}{dt} &= x^2 - \mu \\
\frac{dy}{dt} &= - y
\end{aligned} \right. \end{equation}
Il ne s'agit en fait que d'une l\'eg\`ere modification par rapport
\`a la bifurcation \ref{eq:bif_dim1_+1}, qui se produit sur la
premi\`ere coordonn\'ee de ce syst\`eme. La deuxi\`eme
coordonn\'ee est l\`a pour que le point fixe instable devienne une
selle (il ne peut pas y avoir de selle en dimension 1). Le
diagramme de bifurcation est donc exactement le m\^eme que celui
de la figure~\ref{fig:bif1_+1_diag}. L'espace des phases de part
et d'autre de la bifurcation ($\mu = 0$) est repr\'esent\'e
figure~\ref{fig:saddle_node_f}.

\subsubsection{Bifurcation de Hopf}\label{annexe:Hopf}

\paragraph{\'Etude d'un exemple dans $\R^2$}
Consid\'erons l'exemple de la famille de syst\`emes dynamiques
continus suivante :
\begin{equation}\label{eq:exemple_Hopf}
\left\{
\begin{aligned}
 \frac{dx}{dt} &= - \lambda y + \epsilon x - a x (x^2 + y^2) \\
\frac{dy}{dt} &= \lambda x + \epsilon y - a y (x^2 + y^2)
\end{aligned} \right. \end{equation}
o\`u $\lambda$ et $a$ sont des constantes strictement positives.
Pour tout $\epsilon$, $(0,0)$ est un \'equilibre du syst\`eme, les
valeurs propres de la d\'eriv\'ee en 0 sont $\mu_{\epsilon} = i
\lambda + \epsilon$ et $\overline{\mu_{\epsilon}}$. L'\'equilibre
est donc stable si $\epsilon <0$ et instable si $\epsilon >0$.

En coordonn\'ees <<polaires>> (un peu modifi\'ees), $R = x^2 +
y^2$ et $\theta = \arctan \frac{y}{x}$, \eqref{eq:exemple_Hopf}
devient : \begin{equation} \left\{
\begin{aligned}
 \frac{dR}{dt} &= 2R (\epsilon - a R) \\
\frac{d\theta}{dt} &= \lambda
\end{aligned} \right. \end{equation}

\begin{figure}
\begin{center}
\begin{tabular}{c@{}c@{}c}
     \includegraphics[width=4.5cm]{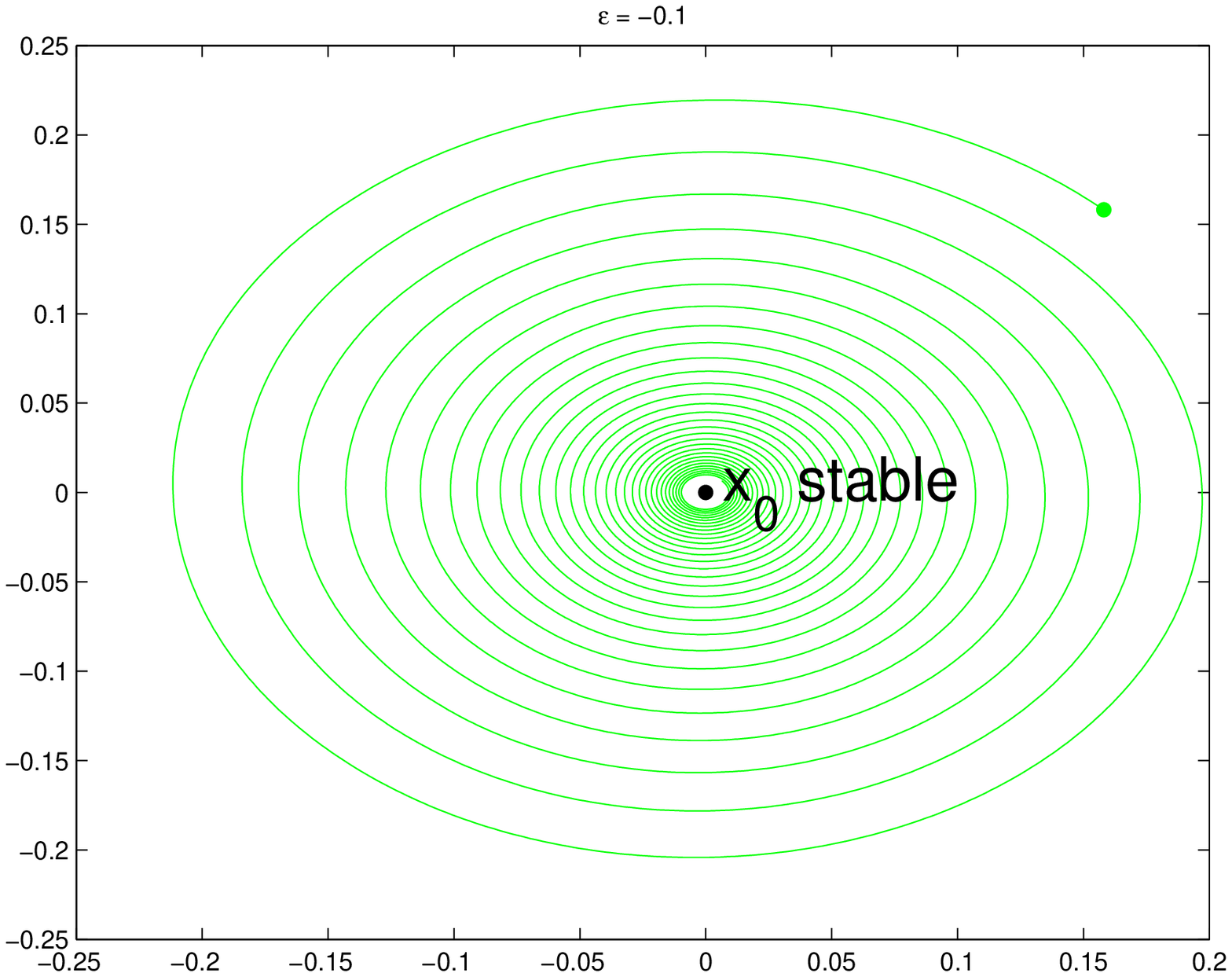}
     &
     \includegraphics[width=4.5cm]{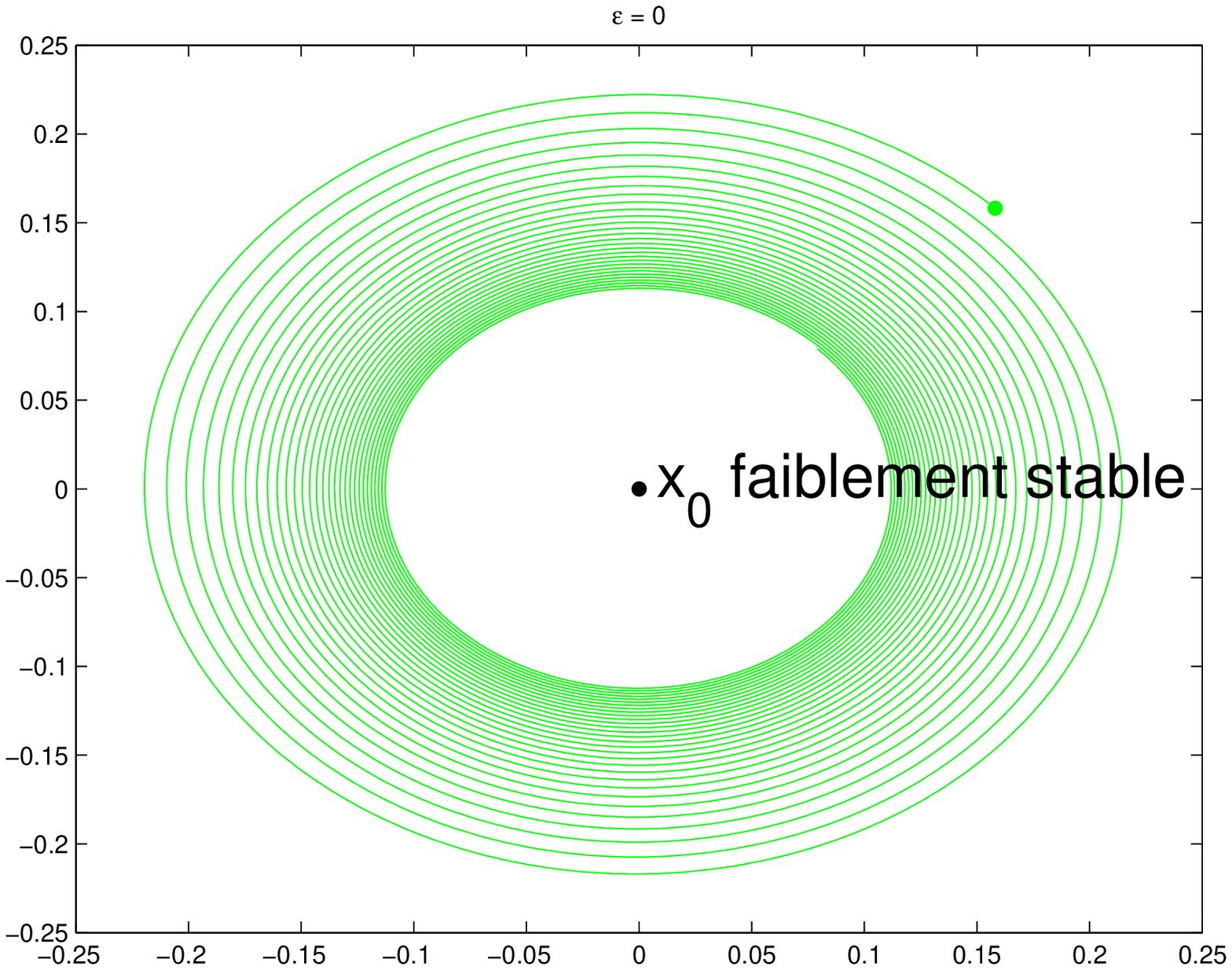}
     &
     \includegraphics[width=4.5cm]{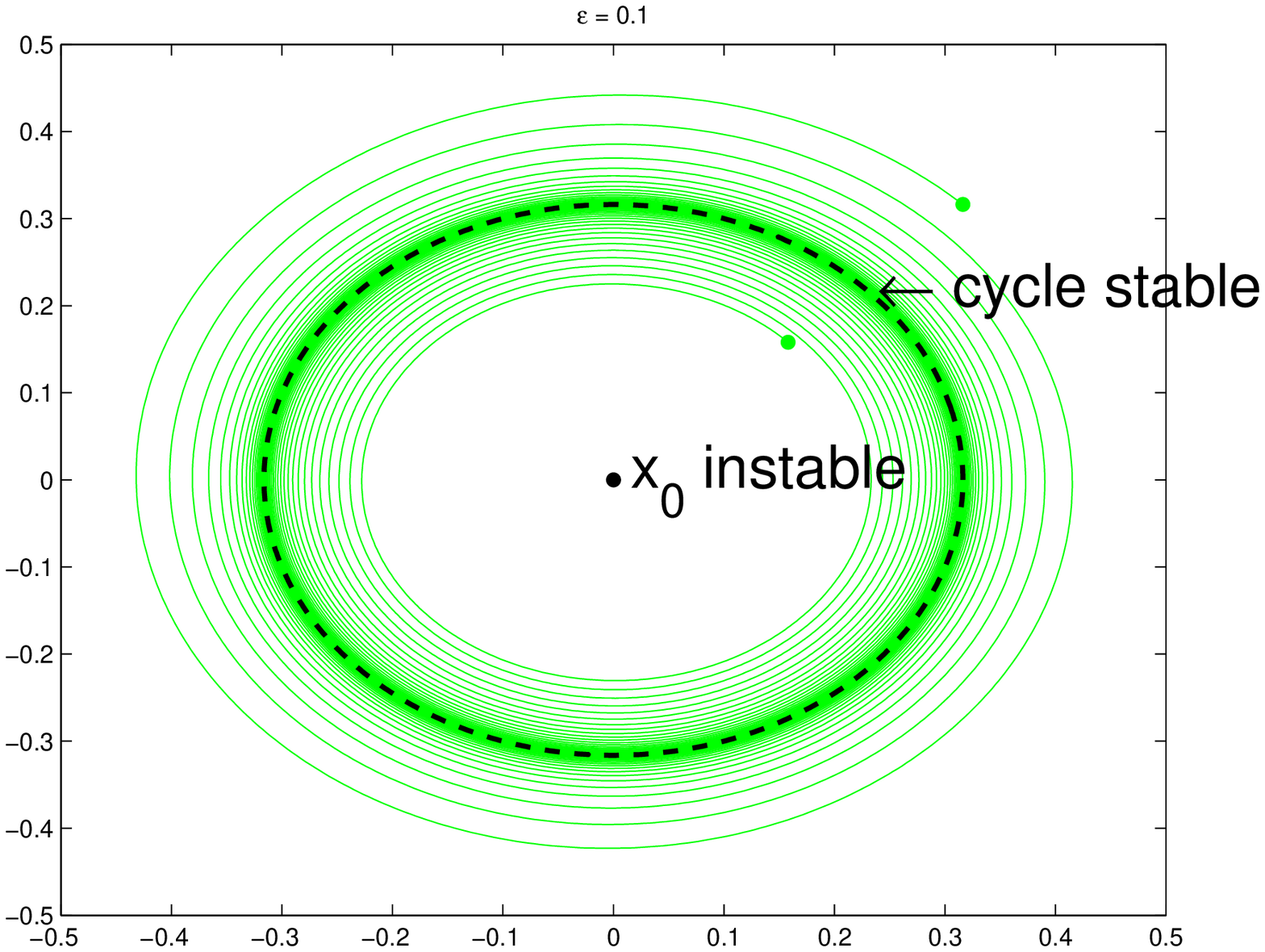}
     \\
     $\epsilon < \epsilon_0$ & $\epsilon = \epsilon_0$ & $\epsilon > \epsilon_0$
\end{tabular}
\caption{\label{fig:Hopf_f} Bifurcation de Hopf : orbites de part
et d'autre de $\epsilon_0=0$.}
\end{center}
\end{figure}

Ce syst\`eme se r\'esout explicitement (voir des exemples
d'orbites figure~\ref{fig:Hopf_f}, dans le cas $a=1$,
$\lambda=2\pi$), d'o\`u
\begin{itemize}
\item si $\epsilon <0$, toutes les solutions convergent vers
l'\'equilibre. \item si $\epsilon >0$, toutes les solutions (sauf
la solution constante nulle) convergent vers l'orbite p\'eriodique
\begin{equation*} \left\{
\begin{aligned}
 R &= \frac{\epsilon}{a} \\
\dot{\theta} &= \lambda
\end{aligned} \right. \end{equation*}
\end{itemize}

\begin{figure}
\begin{center}
     \includegraphics[width=8cm]{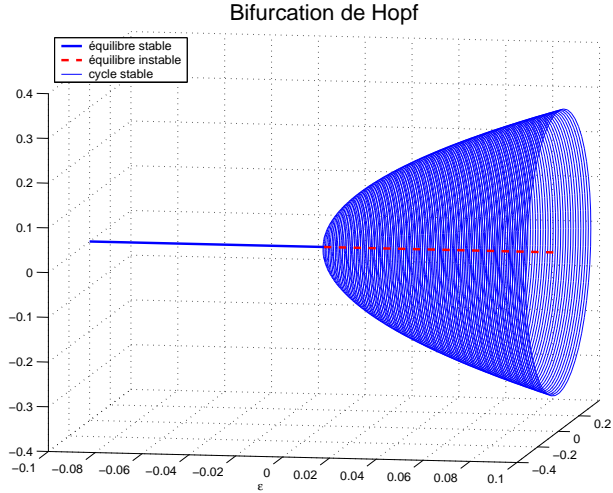}
\caption{\label{fig:Hopf_diag} Diagramme d'une bifurcation de
Hopf.}
\end{center}
\end{figure}

La figure~\ref{fig:Hopf_diag} repr\'esente le diagramme de
bifurcations de cette famille de syst\`emes dynamiques en
$\epsilon_0=0$.

\paragraph{Description du ph\'enom\`ene g\'en\'eral} 
La bifurcation d\'ecrite au paragraphe pr\'ec\'edent est une
\emph{bifurcation de Hopf}. Plus g\'en\'eralement, consid\'erons
la famille \`a un param\`etre d'\'equations diff\'erentielles dans
$\R^N$ \begin{equation}\label{eq:Hopf_general} \frac{dx}{dt} =
F_{\epsilon} (x) .\end{equation} Nous faisons l'hypoth\`ese
$\mathbf{(H_0)}$ : $F_0(0) = 0$ et $D_0(F_0)$ n'a que des valeurs
propres de partie r\'eelle strictement n\'egative, sauf deux qui
sont imaginaires pures et non-nulles : $\mu_0 = i \lambda$ et
$\overline{\mu_0} = - i \lambda$, avec $\lambda >0$. Dans un
voisinage de l'origine, le syst\`eme peut se r\'eecrire (apr\`es
changement de variable), \`a des termes n\'egligeables pr\`es :
\begin{equation}
\left\{
\begin{aligned}
 \frac{dx_0}{dt} &= - \lambda x_1 - a x_0 (x_0^2 + x_1^2) \\
\frac{dx_1}{dt} &= \lambda x_0 - a x_1 (x_0^2 + x_1^2) \\
\frac{dx^{\prime}}{dt} &= A x^{\prime}
\end{aligned} \right. \end{equation}

Nous supposerons d\'esormais $\mathbf{(H_1)}$ : $a>0$.

Nous avons enfin besoin d'une derni\`ere hypoth\`ese relative \`a
la d\'ependence en $\epsilon$. Dans un voisinage de 0, on peut
suivre l'\'equilibre et les valeurs propres $\mu_{\epsilon}$,
$\overline{\mu_{\epsilon}}$ proches de l'axe imaginaire. On
suppose $\mathbf{(H_2)}$ : $\frac{\partial}{\partial \epsilon} \Re
\mu_{\epsilon} > 0$ en $\epsilon = 0$.

Sous les hypoth\`eses $\mathbf{(H_0)}$, $\mathbf{(H_1)}$ et
$\mathbf{(H_2)}$, la dynamique de l'\'equation
\eqref{eq:Hopf_general} pr\'esente une \emph{bifurcation de Hopf}
au voisinage de l'origine en $\epsilon = 0$ : \begin{itemize}
\item pour $\epsilon <0$ petit, il y a un \'equilibre stable.
\item pour $\epsilon=0$, l'\'equilibre reste stable mais plus
faiblement.
\item pour $\epsilon >0$ petit, l'\'equilibre est instable, mais une
orbite p\'eriodique quasi-circulaire de rayon $\simeq
\sqrt{\epsilon/a}$ est stable. \end{itemize}

\paragraph{Cas des diff\'eomorphismes}
Un ph\'enom\`ene semblable peut se produire pour des syst\`emes
dynamiques discrets $x \mapsto f_{\epsilon}(x)$, $x \in \R^N$. On
fait les hypoth\`eses suivantes : \begin{enumerate} \item
$f_0(0)=0$, et les valeurs propres de $D_0 f_0$ ont toutes un
module strictement inf\'erieur \`a 1 sauf deux, $\mu_0$ et
$\overline{\mu_0}$ pour lesquelles $\absj{\mu_0}=1$. \item Pour
$k=1,2,3,4$, $\mu_0^k \neq 1$, \latin{i.e.} $\mu_0 \notin \{ \pm
1, \pm i , \pm j \}$. \item $\mathbf{(H_1^{\prime})}$ et
$\mathbf{(H_2^{\prime})}$ comme dans le paragraphe pr\'ec\'edent.
\end{enumerate}

La dynamique pour $\epsilon$ proche de 0 est alors la m\^eme que
dans le cas pr\'ec\'edent. Un exemple de tel diff\'eomorphisme est
donn\'e par \begin{equation} f_{\epsilon} (z) = \lambda
(1+\epsilon) z - a z \absj{z}^2 ,\, z \in \C \end{equation} avec
$\absj{\lambda}=1$, $\lambda \neq \pm 1$, $a>0$. L'\'equilibre 0 est
stable pour $\epsilon<0$, faiblement stable pour
$\epsilon=\epsilon_0=0$, instable pour $\epsilon > 0$ et alors le
cercle $\absj{z} = \left(\frac{\epsilon}{a}\right)^{1/2}$ est
invariant et attire toutes les orbites proches de 0 sauf
l'\'equilibre lui-m\^eme.

Remarquons \'egalement que si la dynamique sur la courbe
invariante est proche d'une rotation, elle ne se comporte pas
toujours comme une rotation. C'est le cas pour presque tous les
param\`etres, mais pas n\'ecessairement pour tous.

\paragraph{Cas des orbites p\'eriodiques}
On se ram\`ene en fait aux diff\'eomorphismes via
l'\emph{application de retour de Poincar\'e}. En effet, soit
l'\'equation diff\'erentielle $\frac{dx}{dt} = F_0(x)$ dans $\R^N$
poss\'edant une solution p\'eriodique $x_0$. Consid\'erons une
section $\Sigma$ transverse \`a l'orbite $x_0$ en $x_0(t_0)$. Une
condition initiale suffisamment proche de $x_0(t_0)$ retourne sur
$\Sigma$ en un temps fini, ce qui d\'efinit (dans un voisinage de
$x_0(t_0)$) un diff\'eomorphisme $f_0$ de $\Sigma$. La m\^eme
op\'eration pouvant \^etre faite pour une petite perturbation
$F_{\epsilon}$ de $F_0$, cela d\'efinit une famille $f_{\epsilon}$
de diff\'eomorphismes, comme dans le paragraphe pr\'ec\'edent.

\subsubsection{Autres bifurcations}
Nous n'avons bien s\^ur pas abord\'e ici toutes les bifurcations
possibles, m\^eme en nous limitant \latin{a priori} \`a un cadre
restreint. Un exemple particuli\`erement int\'eressant est celui
de la \emph{bifurcation homocline}, reli\'ee \`a celle
d'intersection homocline (d\'efinition \ref{def:homocline}) : deux
intersections transverses homoclines se rencontrent, forment une
tangence \`a cet instant, puis disparaissent. Une r\'ef\'erence
\`a ce sujet est \cite{PalisTakens:HyperbolicityChaotic}.

\subsection{Dynamique des polyn\^omes quadratiques} \label{annexe:poly_quadratiques}
La r\'ef\'erence pour cette section est
\cite{Yoccoz:PolynomesQuadratiques}.

On consid\`ere la famille d'applications\footnote{Tout polyn\^ome
complexe de degr\'e deux est conjugu\'e par une application affine
\`a une application de cette forme. C'est en particulier le cas de
la famille logistique $x \mapsto rx(1-x)$, bien connue en
dynamique des populations.} $P_c : z \mapsto z^2 +c$ pour $z \in
\C$ et $c \in \C$.

Cette famille de syst\`emes dynamiques est l'une des plus simples
qui, en dimension 1, peut g\'en\'erer un comportement chaotique.
Son \'etude est de plus particuli\`erement int\'eressante car on y
observe des ph\'enom\`enes que l'on retrouve dans de nombreux
autres cas.

\subsubsection{Ensembles de Julia et de Mandelbrot}
Il est int\'eressant de se placer dans $\C$ au lieu de $\R$ car on
peut alors utiliser de nombreux r\'esultats d'analyse complexe.
Nous reviendrons ensuite au cas r\'eel.

\begin{figure}
\begin{center}
\begin{tabular}{c@{}c@{}c}
     \includegraphics[width=4.5cm]{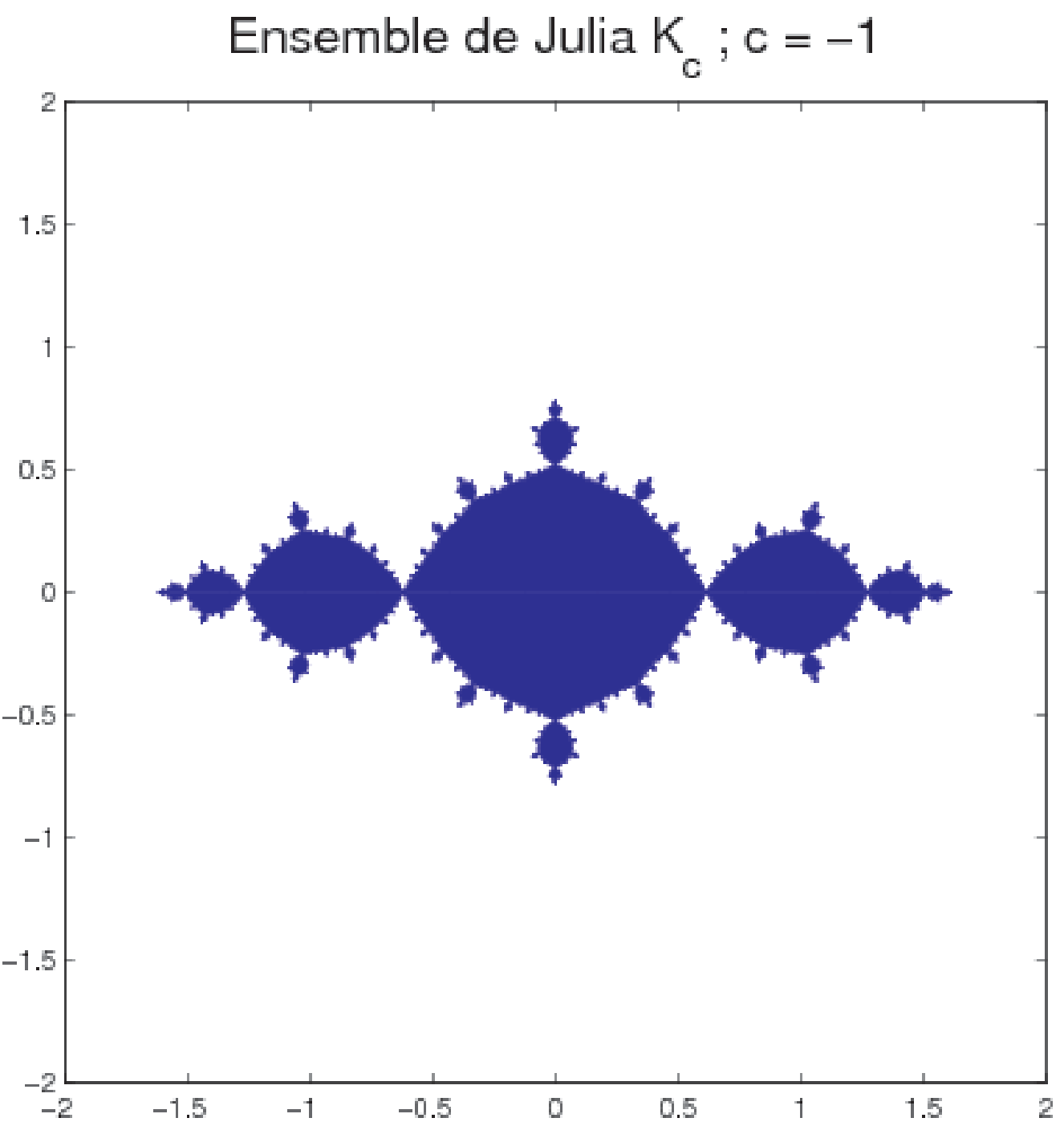}
     &
     \includegraphics[width=4.5cm]{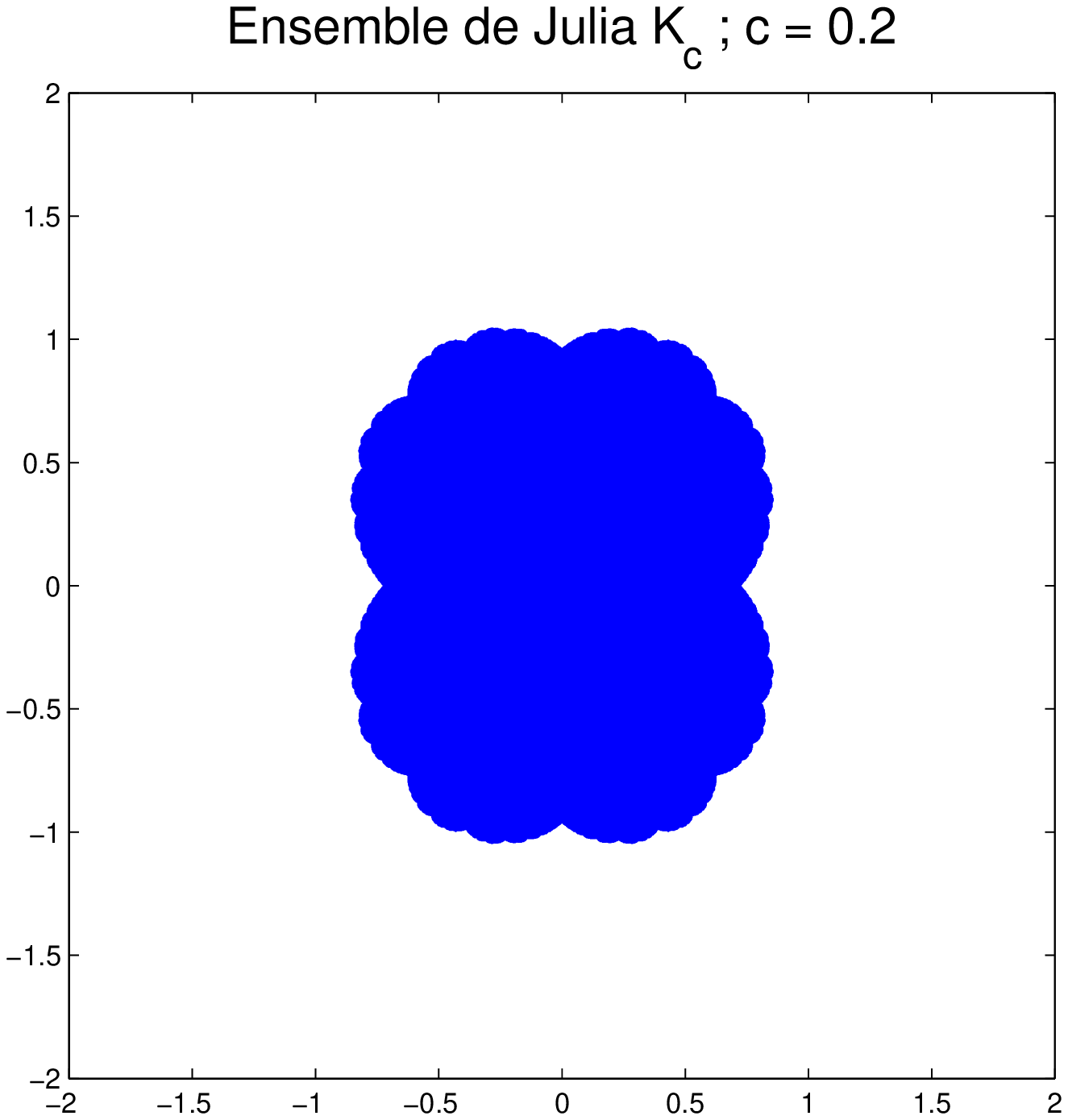}
     &
     \includegraphics[width=4.5cm]{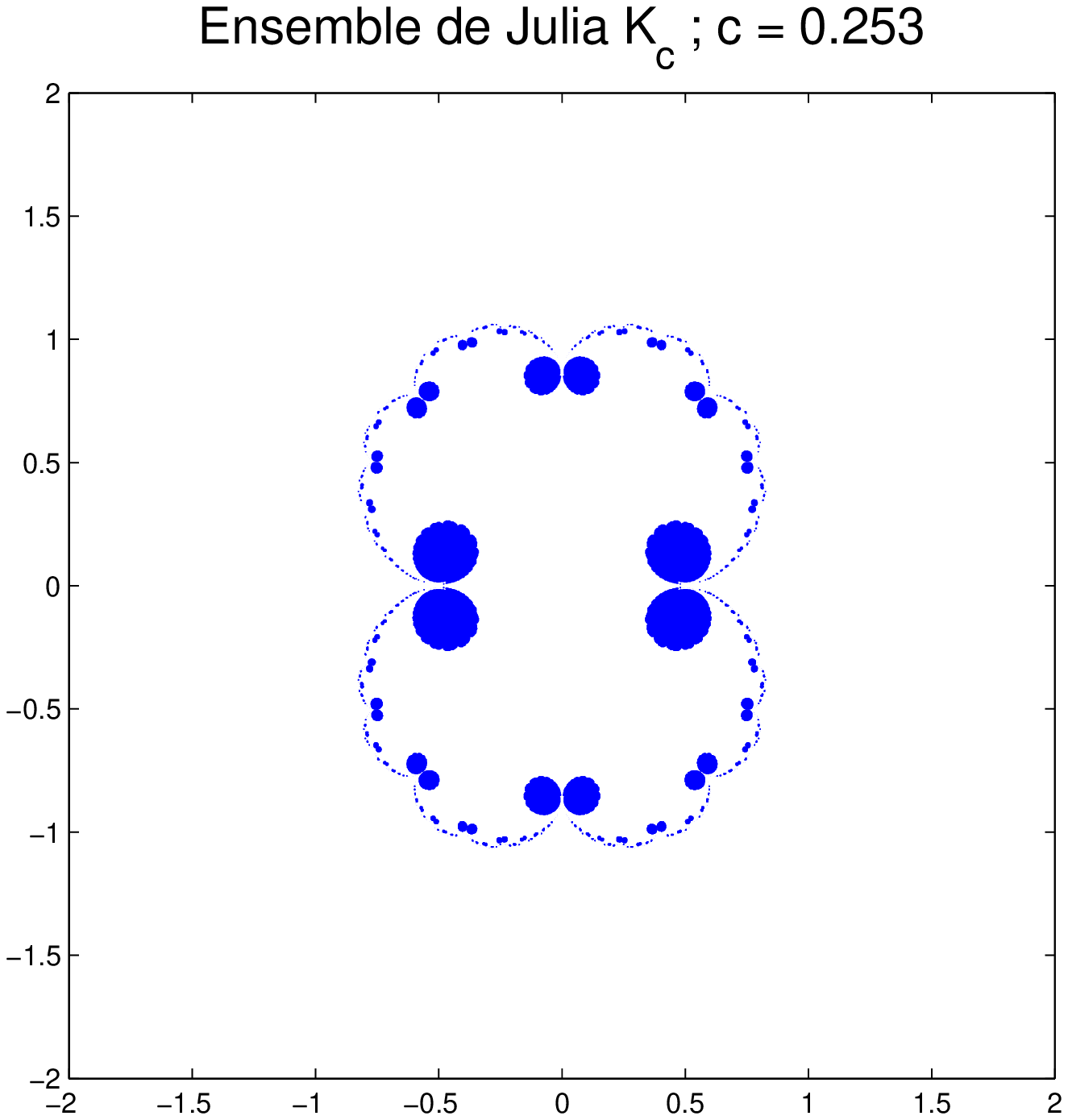}
     \\
     (a) $c = -1 \in M$ & (b) $c = 0.2 \in M$ & (c) $c = 0.253 \notin M$
\end{tabular}
\caption{\label{fig:julia} Ensemble de Julia rempli $K_c$.}
\end{center}
\end{figure}

Pour $c \in \C$, on consid\`ere l'\emph{ensemble de Julia rempli}
(figure~\ref{fig:julia})
\[ K_c = \{ z \in \C \telque P_c^n(z) \text{ est born\'e} \} \] que l'on
peut \'egalement \'ecrire \[ K_c = \bigcap_{n \geq 0}
P_c^{-n}(\adherence{\mathbb{D}(O,R)}) \] en ayant pos\'e $R = (1+
\sqrt{1+4\absj{c}})/2$. Ainsi, $K_c$ est : \begin{itemize} \item
compact, \item non-vide et il contient tous les points
p\'eriodiques de $P_c$, \item totalement invariant, \latin{i.e.}
$P_c(K_c) = K_c = P_c^{-1}(K_c)$, \item plein, \latin{i.e.} $\C
\backslash K_c$ est connexe. \end{itemize}

Le bord $J_c = \partial K_c$ est l'\emph{ensemble de Julia}.
D'apr\`es un th\'eor\`eme montr\'e ind\'ependamment par Julia et
Fatou, c'est aussi l'adh\'erence de l'ensemble des points
p\'eriodiques r\'epulsifs.

\begin{figure}
\begin{center}
    \includegraphics[height=8cm]{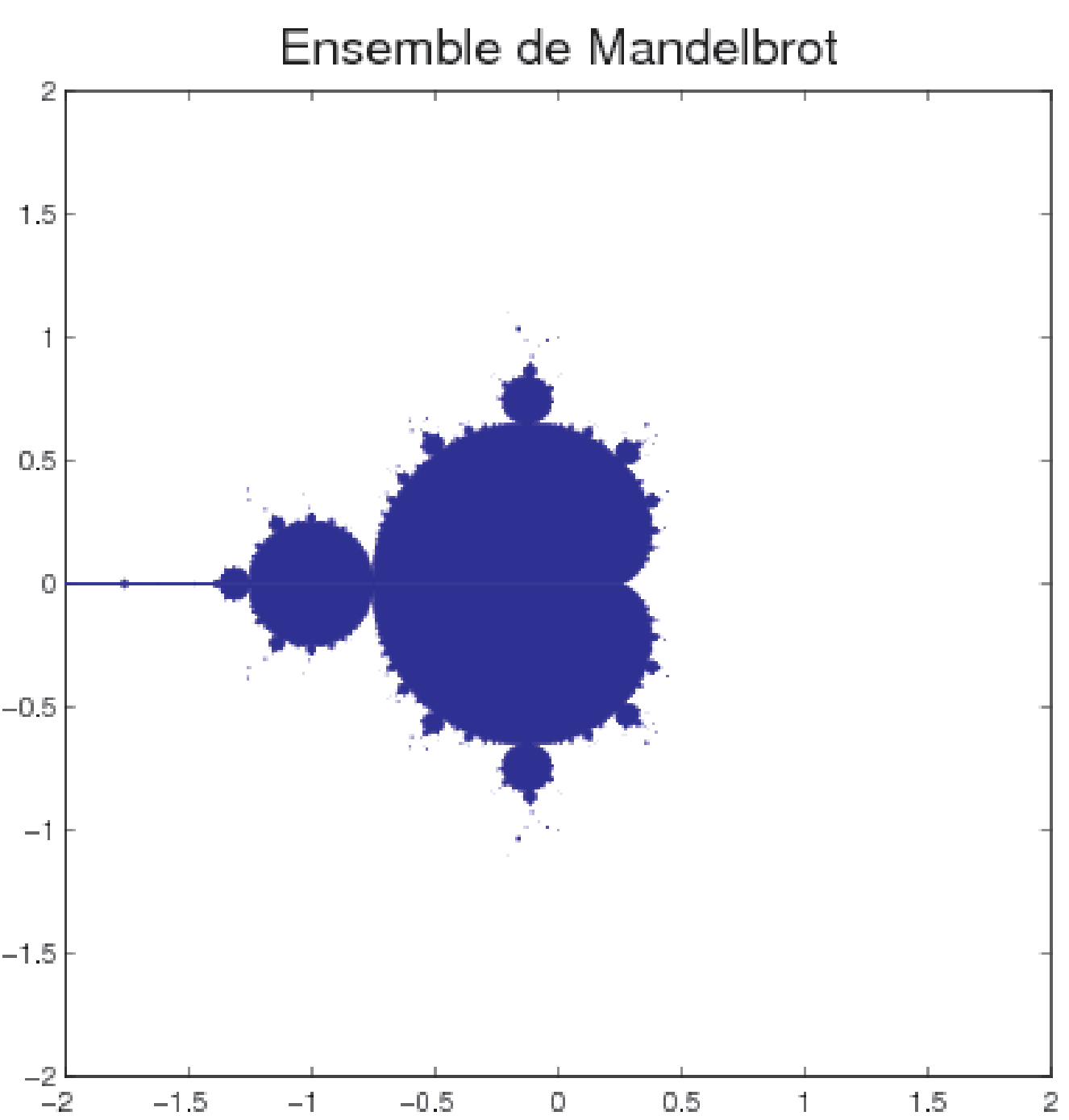}
    \caption{\label{fig:mandelbrot} L'ensemble de Mandelbrot :
    approximation num\'erique.}
\end{center}
\end{figure}

Un th\'eor\`eme de Fatou (1919) montre que $0 \in K_c$ si et
seulement si $K_c$ est connexe. Dans l'espace des param\`etres, on
consid\`ere l'\emph{ensemble de Mandelbrot}
(figure~\ref{fig:mandelbrot})
\[ M = \{ c \in \C \telque K_c \text{ est connexe} \}. \] On
montre que $M = \{ c \telque \absj{P_c^n(0)}\leq 2 , \, \forall n
>0 \}$ et donc $M$ est compact. De plus, $M$ est plein,
sym\'etrique par rapport \`a l'axe r\'eel qu'il coupe suivant
l'intervalle $[-2,1/4]$. Sur la figure~\ref{fig:mandelbrot}, on
distingue des \^ilots disjoints de la grande composante de $M$. Un
calcul plus pouss\'e montrerait qu'ils lui sont en r\'ealit\'e
reli\'es par des filaments extr\^emement fins.

Si $c \notin M$, $K_c = J_c$ est un ensemble de Cantor
(figure~\ref{fig:julia}.c). La dynamique est de type d\'ecalage,
et $P_c$ est hyperbolique.

Si $c \in M$, la pr\'esence d'orbites p\'eriodiques attractives
(\latin{i.e.} $z_0$ tel que $P_c^m(z_0)=z_0$ et
$\absj{(P_c^m)^{\prime}(z_0)} <1$) est d\'eterminante pour la
structure de $K_c$. Douady a montr\'e (1982) que $P_c$ a au plus
une orbite p\'eriodique attractive. Lorsque c'est le cas, le
bassin d'attraction $W = \{ z \telque \lim_{n \rightarrow \infty}
d(P_c^n(z),O(z_0))=0 \}$ de l'orbite est l'int\'erieur de $K_c$,
et $P_{c|J_c}$ est hyperbolique. Un tel $c$ est alors dans
l'int\'erieur de $M$, et la composante connexe de l'int\'erieur de
$M$ contenant $c$ est appel\'ee \emph{composante hyperbolique} de
$M$.

Par exemple, l'ensemble des $c$ tels que $P_c$ poss\`ede un point
fixe attractif est l'int\'erieur d'une cardio\"ide dite principale
contenant 0. Pour un tel $c$, l'int\'erieur de $K_c$ a une seule
composante et $J_c$ est un quasi-cercle
(figure~\ref{fig:julia}.b). La figure~\ref{fig:julia}.a donne un
exemple d'ensemble de Julia rempli lorsque $c$ est dans une autre
composante hyperbolique.

L'int\'erieur de $M$ est dense dans $M$ et contient toutes les
composantes hyperboliques. La conjecture d'hyperbolicit\'e dit que
l'union des composantes hyperboliques est en fait exactement
l'int\'erieur de $M$. Parmi les r\'esultats partiels obtenus dans
cette direction, on a montr\'e que les composantes hyperboliques
de $M$ rencontrent $M \cap \R = [-2,1/4]$ suivant un ensemble
dense.

En revanche, cette intersection n'est pas de mesure totale, comme
le montre le th\'eor\`eme de Jakobson \cite{Yoccoz:Jakobson}.

\subsubsection{Dynamique sur la droite r\'eelle}\label{annexe:quadr_reelle}
Supposons $c \in m \cap \R = [-2,1/4]$. Nous venons de voir que
pour un ensemble dense (mais pas de mesure totale) de valeurs de
$c$, $c$ est dans une composante hyperbolique et donc il existe
une unique orbite p\'eriodique attractive, et son bassin
d'attraction est l'int\'erieur de $K_c$.

\begin{figure}
\begin{center}
    \includegraphics[height=8cm]{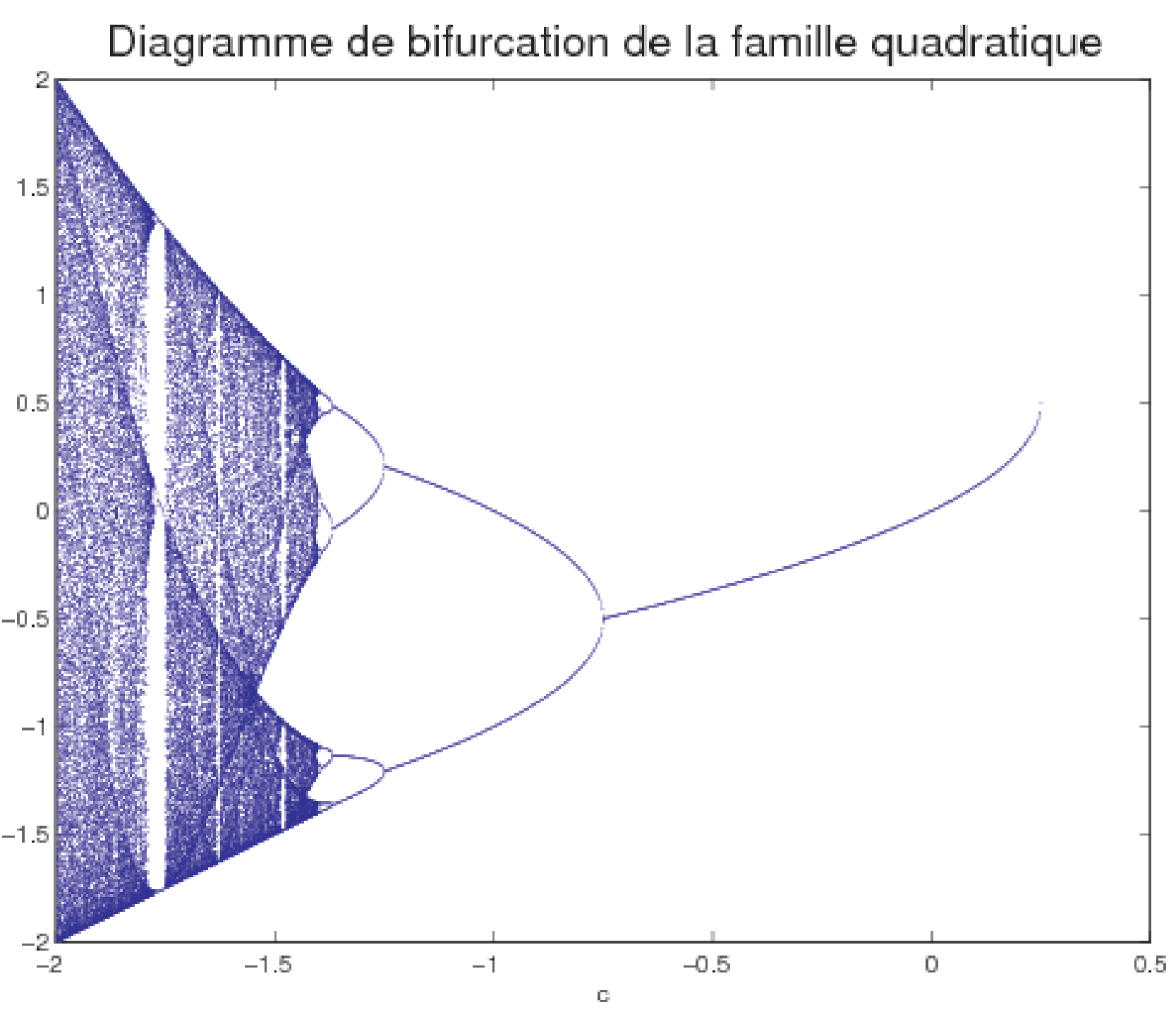}
    \caption{\label{fig:quadr_bif} Dynamique de $P_c$ pour $c\in [-2;\frac{1}{4}]$.}
\end{center}
\end{figure}

Partons de $c=1/4$ et faisons diminuer $c$ (voir le diagramme de
bifurcations, figure~\ref{fig:quadr_bif}). On a tout d'abord un
point fixe attractif, puis une orbite attractive de p\'eriode 2
(apr\`es une bifurcation doublement de p\'eriode en $c=c^{(1)}$).
Les bifurcations doublement de p\'eriode se succ\`edent ainsi
jusqu'\`a atteindre $c = c^{(\infty)}$ o\`u il n'y a plus d'orbite
p\'eriodique attractive. Cette succession de bifurcations est
appel\'ee \emph{cascade sous-harmonique directe}.

Le point $c = c^{(\infty)}$ est appel\'e le \emph{point de
Feigenbaum}, o\`u la dynamique peut encore \^etre d\'ecrite assez
simplement. Un intervalle $I$ est stable, et poss\`ede deux
sous-intervalles $I_0$ et $I_1$ disjoints tels que $P_c(I_0)
\subset I_1$ et $P_c(I_1) \subset I_0$. Dans chacun de ces
intervalles $I_{\alpha}$, on trouve deux sous-intervalles
$I_{\alpha,0}$ et $I_{\alpha,1}$ disjoints tels que
$P_c^2(I_{\alpha,0}) \subset I_{\alpha,1}$ et $P_c^2(I_{\alpha,1})
\subset I_{\alpha,0}$. On retrouve cette dynamique en faisant un
double changement d\'echelle, spatial et temporel. L'ensemble
limite a ainsi une structure d'ensemble de Cantor. Cependant, la
dynamique de $P_c$ sur cet ensemble n'est pas chaotique.

Au del\`a du point de Feigenbaum, la dynamique devient chaotique.
Il se produit alors un processus, miroir de la cascade directe, de
regroupement par bandes : les composantes connexes de l'attracteur
fusionnent successivement. On parle alors de \emph{cascade
inverse}.

La complexit\'e du comportement du syst\`eme lorsque $c$ varie
provient de l'alternance de r\'egimes p\'eriodiques et chaotiques.
En effet, la situation pr\'esent\'ee avec une p\'eriode 1 initiale
se reproduit pour toutes les valeurs de p\'eriodes impaires.
Ainsi, des fen\^etres de p\'eriodicit\'e s'installent brutalement
\`a la suite de r\'egimes chaotiques, s'ach\`event par une cascade
sous-harmonique suivie d'une cascade inverse et d'un r\'egime
chaotique. Ces fen\^etres de p\'eriodicit\'e sont denses dans
$[-2; 1/4]$, ce qui montre bien toute la complexit\'e de la
dynamique dans cette r\'egion de l'espace des param\`etres. De
plus, le compl\'ementaire de cet ensemble ayant une mesure
non-nulle, il reste possible d'observer un comportement chaotique
en choisissant le param\`etre $c$ al\'eatoirement suivant la
mesure de Lebesgue.

Une \'etude plus d\'etaill\'ee (et plus exp\'erimentale) est faite
dans \cite{Manneville:chaos}.

\subsection{Comportement statistique des orbites} 
L'\'etude des syst\`emes dynamiques mesurables est l'objet de la
th\'eorie ergodique, qui est notamment introduite dans
\cite{Benoistpaulin:systemesdynamiques}. Nous ne donnons ici que
quelques d\'efinitions utiles pour notre \'etude.

\begin{Def}[Mesure invariante] Une mesure $\mu$ est invariante par
l'application (mesurable) $f$ si pour toute partie mesurable $A$
\[ \mu(f^{-1}(A)) = \mu(A). \]
\end{Def}

\begin{Def}[Ergodicit\'e] Une application $f : (X,\mathcal{B},\mu)
\rightarrow (X,\mathcal{B},\mu)$ qui pr\'eserve $\mu$ est
\emph{ergodique} si \[ \forall A \in \mathcal{B}, \, f^{-1}(A) = A
\Rightarrow \mu(A) \in \{ 0 , 1 \} .\] \end{Def}

\begin{Pro} Une application $f$ est ergodique si et seulement si
toute application $\phi : X \rightarrow \C$ mesurable, telle que
$\phi \circ f = \phi$ presque partout, est presque partout
constante. \end{Pro} On peut remplacer dans cette proposition
<<mesurable>> par $L^1(X,\mu)$ ou $L^2(X,\mu)$.

\begin{ex} Les rotations d'angle $\alpha$ irrationnel et le doublement de
l'angle sont ergodiques sur le cercle $S^1$, pour la mesure de
Lebesgue sur le cercle. \end{ex}

La proposition suivante fait le lien avec la notion physique
d'ergodicit\'e.
\begin{Pro} Si $X$ est un espace m\'etrique s\'eparable, $\mu$ une
probabilit\'e bor\'elienne sur $X$, $f : X \rightarrow X$ continue
pr\'eservant $X$. Si $f$ est ergodique, alors $\mu$-presque toute
orbite est dense dans $X$. \end{Pro}

On a alors une estimation quantitative de la <<densit\'e>> des
orbites : pour toute partie mesurable $A$, la proportion de temps
pass\'ee dans $A$ par presque toutes les orbites est \'egale \`a
$\mu(A)$.
\begin{The}[Th\'eor\`eme ergodique de Birkhoff] Soit $(X,\mathcal{B},\mu)$
un espace mesur\'e, $f : X \rightarrow X$ mesurable pr\'eservant
$\mu$. Pour tout $\phi$ dans $L^1(X,\mu)$, on note \[ S_n \phi (x)
= \frac{1}{n} \sum_{k=0}^{n-1} \phi(f^k(x)) \text{ (somme de
Birkhoff de $\phi$).} \] La limite $\widetilde{\phi}(x) = \lim_{n
\rightarrow \infty} S_n \phi (x)$ existe pour $\mu$-presque tout
$x$, $\widetilde{\phi} \circ f = \widetilde{\phi}$ presque
partout. Pour toute partie $f$-invariante $A$ mesurable, de mesure
finie, on a \[ \int_A \phi d\mu = \int_A \widetilde{\phi} d\mu. \]

En particulier, si $\mu$ est une mesure de probabilit\'e
ergodique, alors \[ \widetilde{\phi}(x) = \int_X \phi d\mu \] pour
$\mu$-presque tout $x$.
\end{The}

\begin{Def}[Mesure physique, Ruelle--Bowen] \label{def:mesure_physique}
C'est une mesure de probabilit\'e $\mu$ invariante par $f$, telle
que pour toute application $\phi$ continue sur $X$, pour
$\lambda$-presque tout $x \in X$, $\widetilde{\phi}(x) = \int_X
\phi d\mu$.
\end{Def} 
Cette condition est bien plus forte que l'ergodicit\'e, puisque
contrairement au r\'esultat du th\'eor\`eme de Birkhoff, le
r\'esultat de convergence est vrai $\lambda$-p.p. ($\lambda$ est
la mesure de Lebesgue sur $X$, dont le support est $X$ tout
entier), et non $\mu$-p.p., $\mu$ pouvant avoir un support bien
moins grand que $X$ tout entier. En particulier, si le support de
$\mu$ a une mesure de Lebesgue nulle, le th\'eor\`eme ergodique de
Birkhoff \'enonce un r\'esultat que l'on n'observera jamais (p.s.)
si l'on choisit une condition initiale $x$ al\'eatoirement suivant
$\lambda$. De plus, lorsque la mesure physique existe (cela a
\'et\'e prouv\'e dans le cas du sol\'eno\"ide), elle est unique
(ce n'est pas toujours le cas pour les mesures ergodiques).

La mesure physique (lorsqu'elle existe) donnant la densit\'e de
$\lambda$-presque toute orbite, c'est elle que l'on observe
empiriquement au cours des simulations num\'eriques.

\subsection{Dimension fractale} \label{annexe:theorie_dimfract}

Certains des attracteurs que nous avons \'evoqu\'es ont --- au
moins partiellement --- une structure d'ensemble de
Cantor\footnote{Notamment le sol\'eno\"ide, section
\ref{annexe:solenoide}, et l'attracteur de H\'enon, section
\ref{annexe:Henon}.}, de dimension non-enti\`ere. Nous allons
donner un sens \`a cette affirmation, en d\'efinissant la
dimension fractale d'un compact $K$. Il existe plusieurs autres
notions de dimension non-enti\`ere (reli\'ees les unes aux
autres), par exemple la dimension de Hausdorff\footnote{on montre
en g\'en\'eral que la dimension de Hausdorff $HD(K)$ est
inf\'erieure o\`u \'egale \`a $D_f(K)$. Il y a \'egalit\'e pour
des classes assez g\'en\'erales d'ensemble, par exemple pour
l'exemple d'ensemble de Cantor d\'ecrit dans ce paragraphe. On
pourra se r\'ef\'erer \`a \cite{PalisTakens:HyperbolicityChaotic}
pour le cas des ensembles de Cantor d\'efinis dynamiquement.} ; la
dimension fractale poss\`ede l'avantage d'\^etre la plus simple
\`a \'evaluer num\'eriquement.

\begin{Def}[Dimension fractale]
Soit $K$ un compact d'un espace m\'etrique $(X,d)$. Pour tout
$\epsilon >0$, on note $N_{\epsilon}(K)$ le nombre minimal de
boules de rayon $\epsilon$ n\'ecessaires pour recouvrir $K$. La
\emph{capacit\'e limite} ou \emph{dimension fractale} de $K$ est
d\'efinie par \[ D_f (K) = \limsup_{\epsilon \rightarrow 0}
\frac{\log N_{\epsilon}(K)} {- \log \epsilon} \]
\end{Def}

Lorsque $K$ est une sous-vari\'et\'e de dimension finie, la
dimension fractale est \'egale \`a la dimension topologique.

Un autre cas classique est celui des ensembles de Cantor.
Consid\'erons un exemple o\`u la dimension se calcule facilement :
pour $I=[a,b]$ intervalle, on note $f(I)=[a,a+\frac{b-a}{3}] \cup
[a+\frac{2(b-a)}{3},b]$. Le compact $K = \bigcap_{n \in \N} f^n
([0,1])$ est un ensemble de Cantor. On a $N_{3^n}(K) = 2^n$, et
$\epsilon \mapsto N_{\epsilon}(K)$ est croissante, donc si $3^n
\leq \epsilon \leq 3^{n+1}$, \[ \frac{n \log 2}{(n+1) \log 3} \leq
\frac{\log N_{\epsilon}(K)} {- \log \epsilon} \leq \frac{(n+1)
\log 2}{n \log 3}. \] On en d\'eduit que \[D_f(K) = \frac{\log
2}{\log 3}.\]

\subsection{Th\'eor\`eme de Whitney} \label{annexe:whitney}
Pour visualiser les r\'esultats des simulations num\'eriques, nous
avons projet\'e en dimension 3 les points de $\R^{N}$ ($N$ grand)
que nous avions calcul\'es. L'une des justifications \latin{a
posteriori} de la validit\'e de la m\'ethode est th\'eorique et
passe par le th\'eor\`eme de Whitney. En effet, l'attracteur
semblant avoir une dimension (fractale) strictement inf\'erieure
\`a $1\virg 5$, il est possible de le plonger dans $\R^3$.

\begin{The}[Whitney] \label{the:whitney} Toute vari\'et\'e compacte lisse de dimension
$n \in \N$ se plonge dans $\R^{2n+1}$.
\end{The}

Ce r\'esultat est d\'emontr\'e dans \cite{Lafontaine:geodiff}. Il
se g\'en\'eralise au cas d'un compact de dimension fractale $d$,
qui se plonge dans $\R^N$ d\`es que $N > 2d$.
\section{R\'esultats d\'etaill\'es}\label{annexe:details}

\carteB{1}{0} \carteB{2}{1} \carteB{3}{2} \carteB{4}{3}
\carteB{5}{4} \carteB{6}{5} \carteB{7}{6} \carteB{8}{7} \clearpage

\begin{figure}
\begin{center}
\begin{tabular}{c}
     \includegraphics[height=8.5cm]{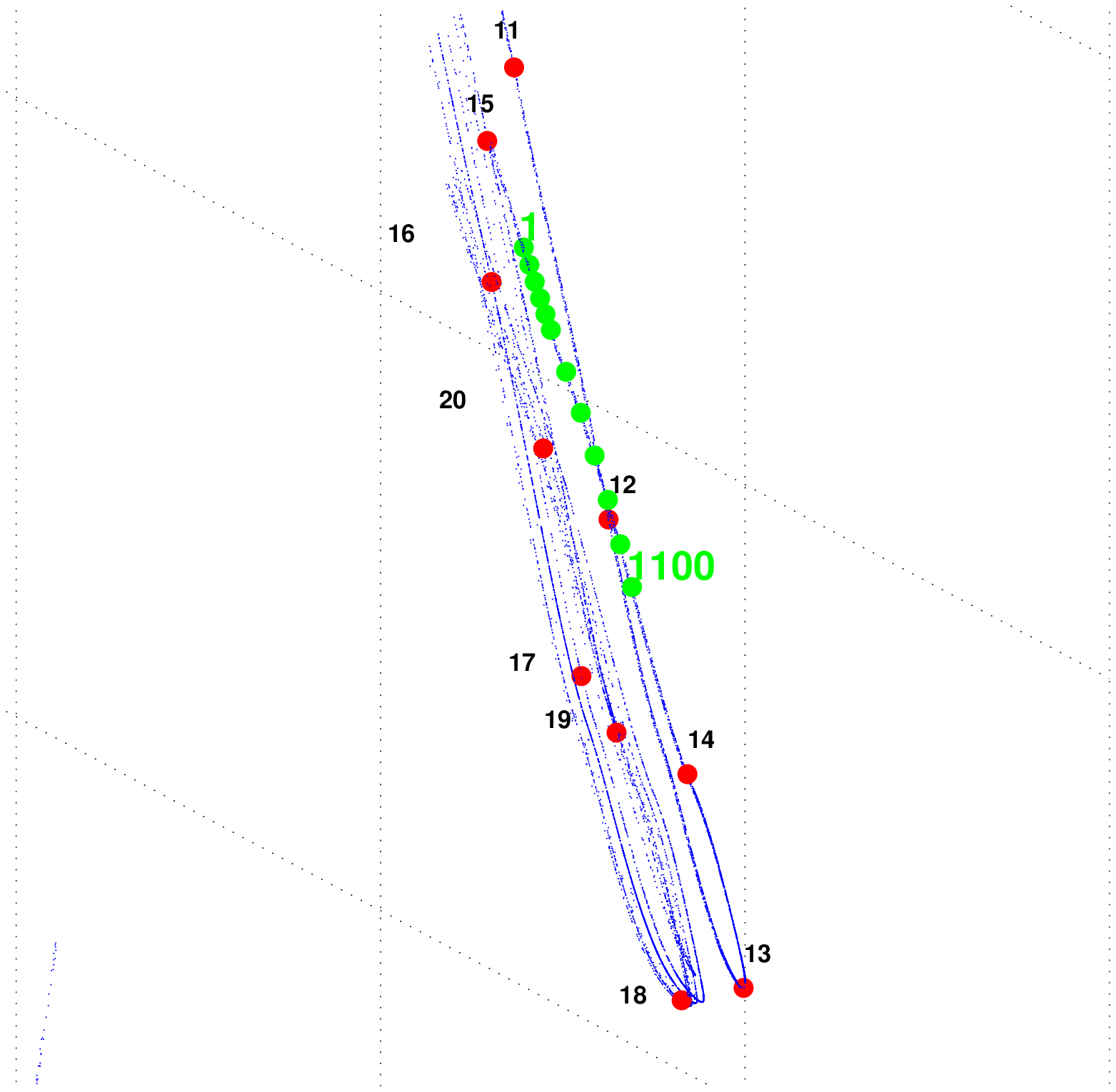} \\
    (a) $x_1 , \ldots , x_{1100}$ \\
     \includegraphics[height=8.5cm]{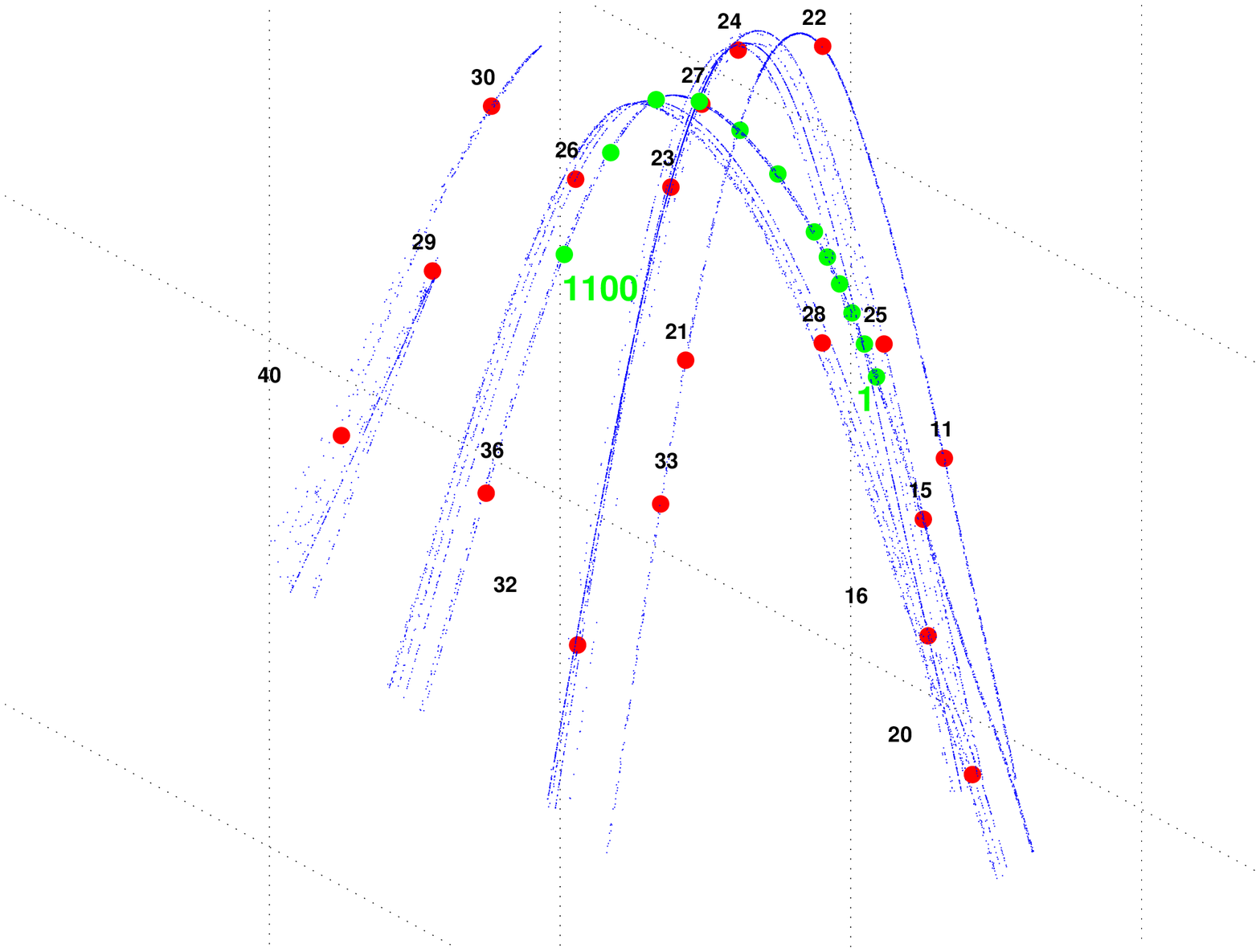} \\
    (b) $T^2(x_1), \ldots, T^2(x_{1100})$
\end{tabular}
\caption{\label{fig:pli2} Localisation du pli et de sa pr\'eimage
 (deuxi\`eme m\'ethode).}
\end{center}
\end{figure}

\begin{figure}
\begin{center}
\begin{tabular}{c@{}c}
     \includegraphics[width=7.5cm]{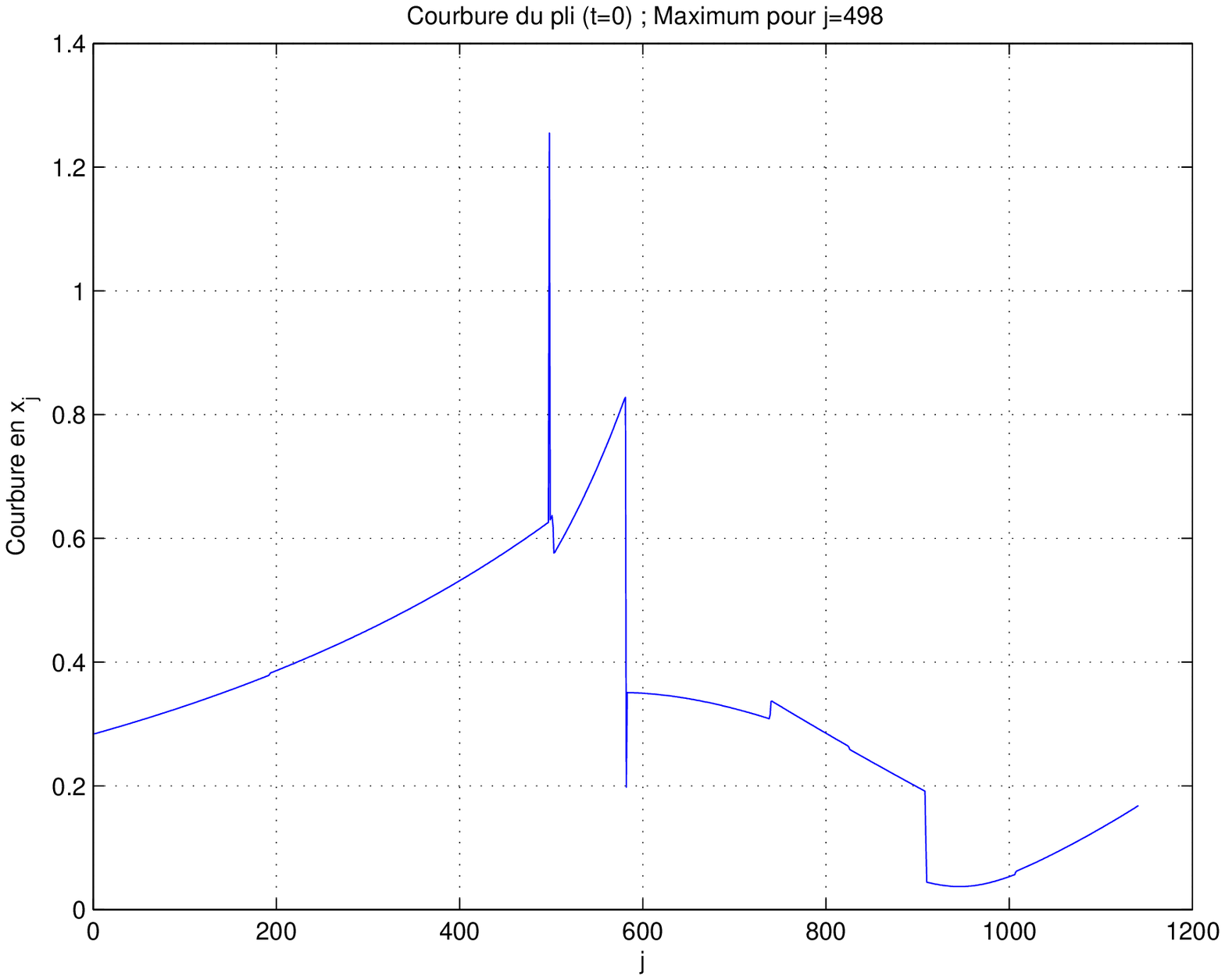}
     &
     \includegraphics[width=7.5cm]{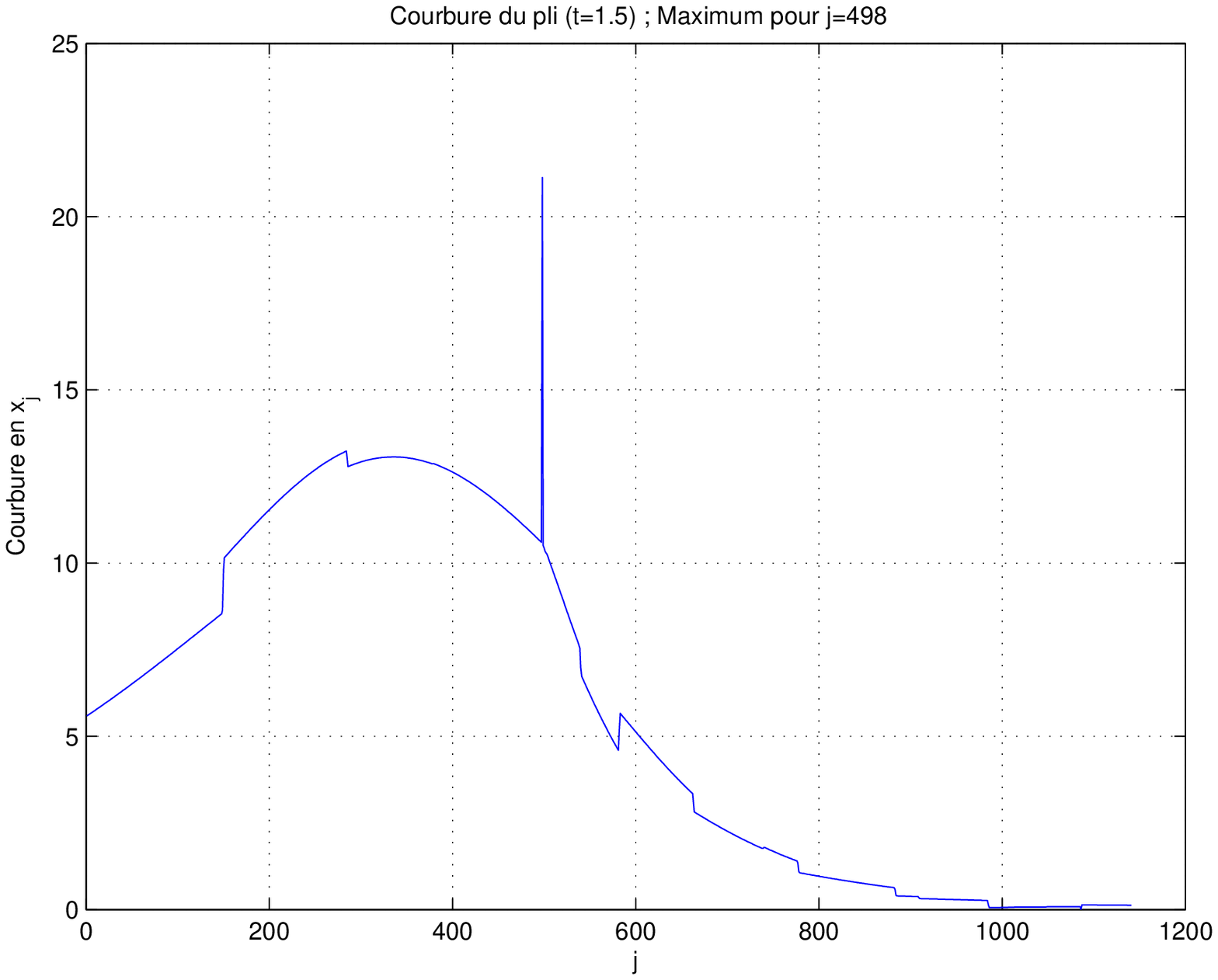}
     \\
     (a) $t=0$ & (b) $t=1\virg 5$
     \\
     \includegraphics[width=7.5cm]{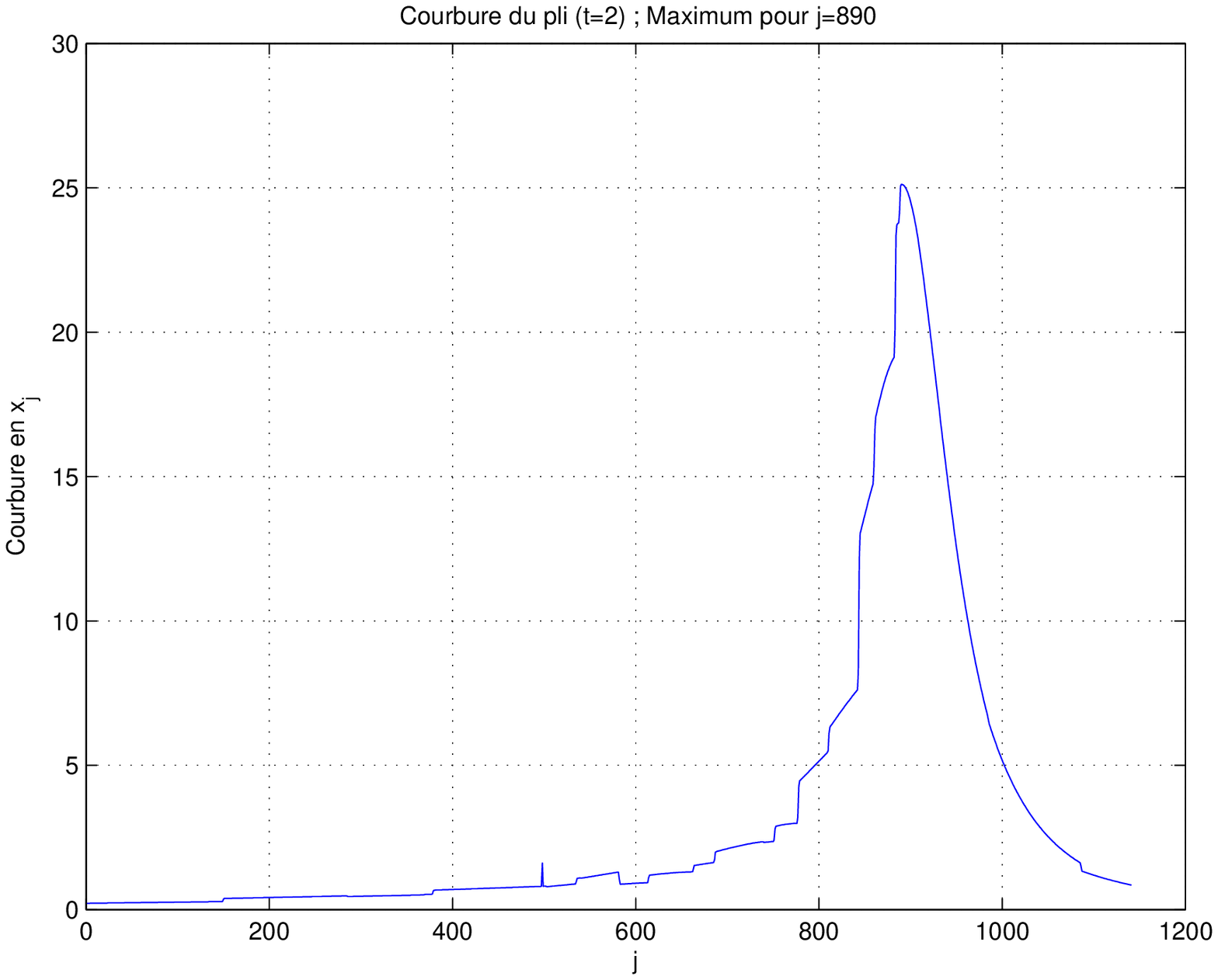}
     &
     \includegraphics[width=7.5cm]{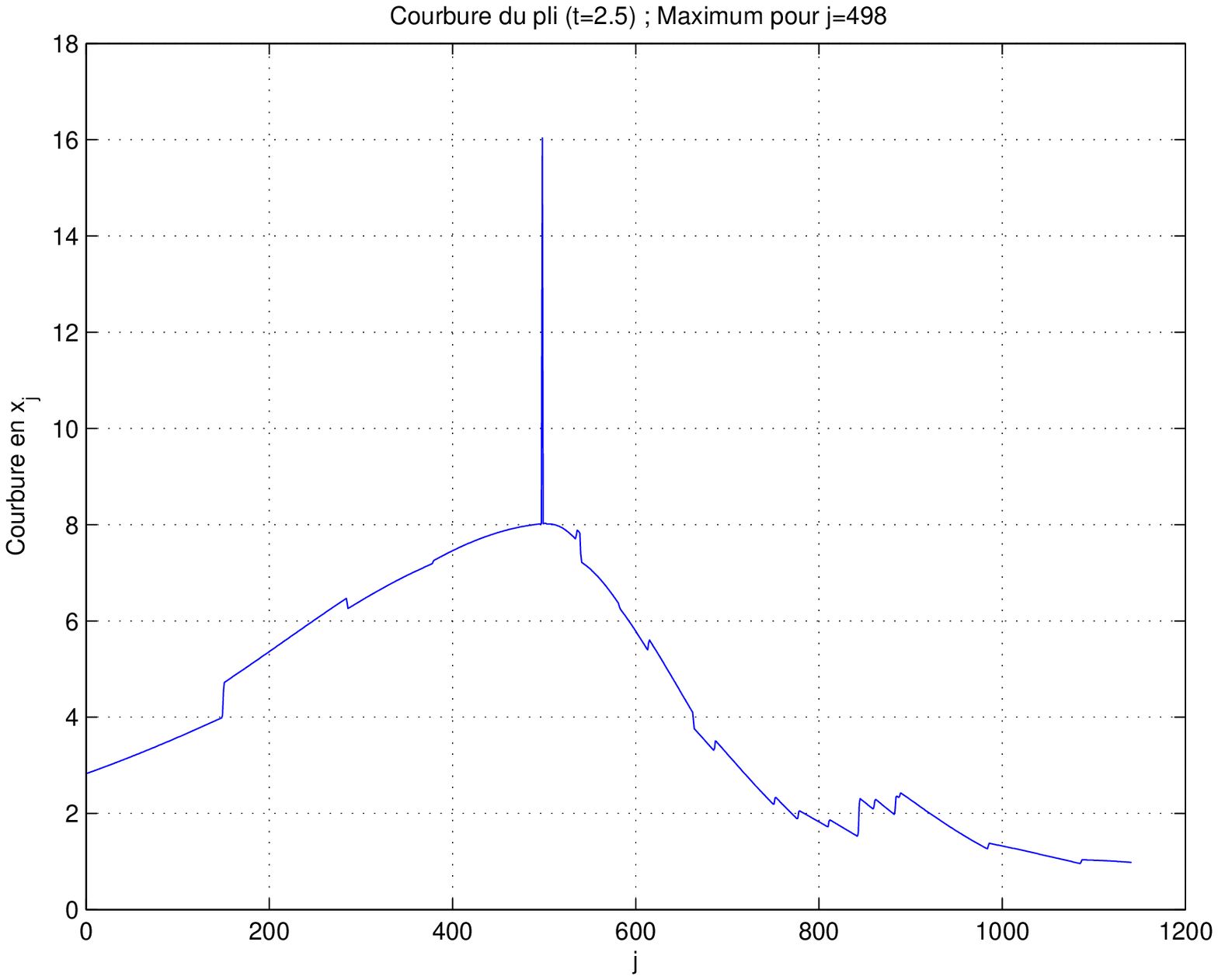}
     \\
     (c) $t=2$ & (d) $t=2\virg 5$
     \\
     \includegraphics[width=7.5cm]{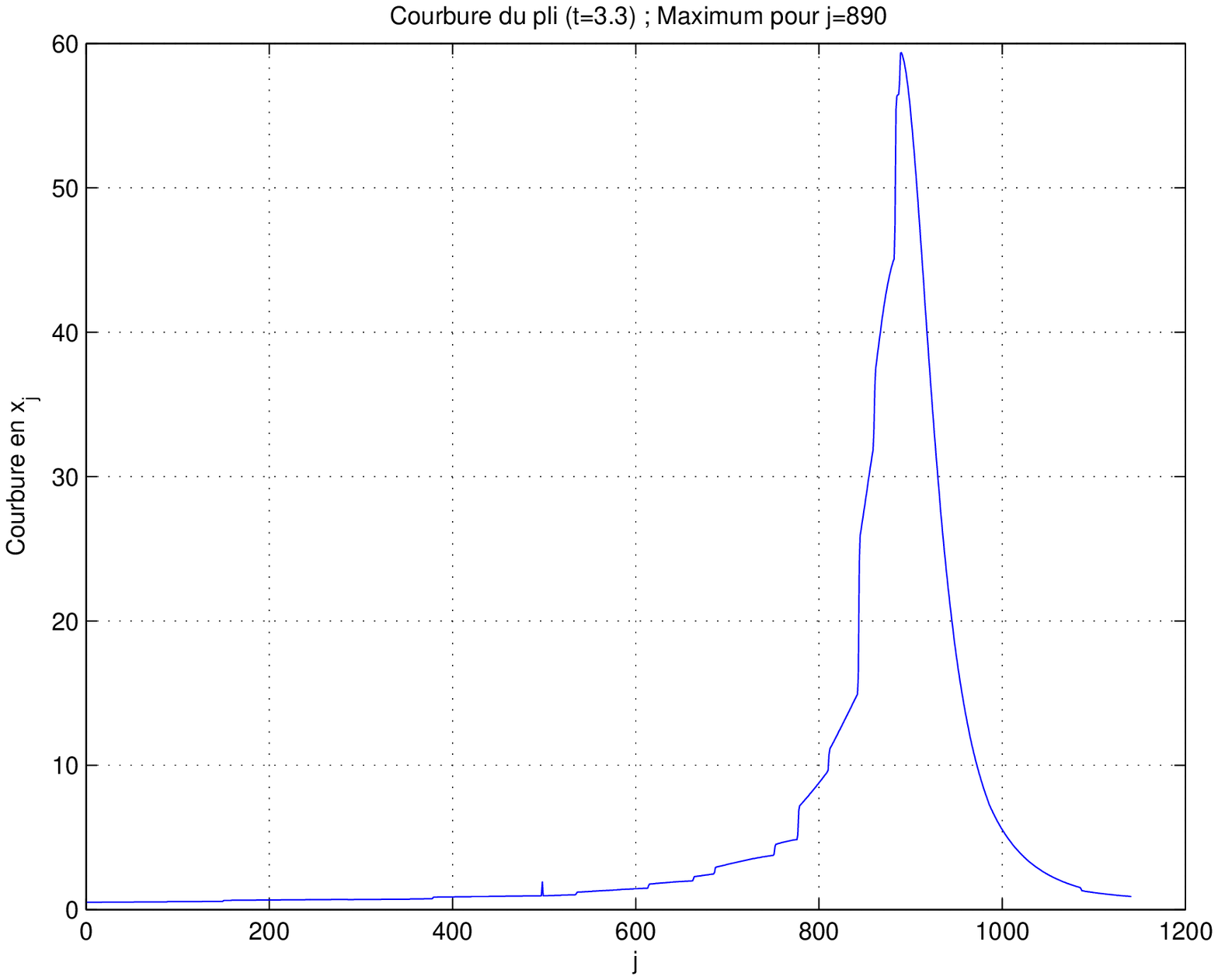}
     &
     \includegraphics[width=7.5cm]{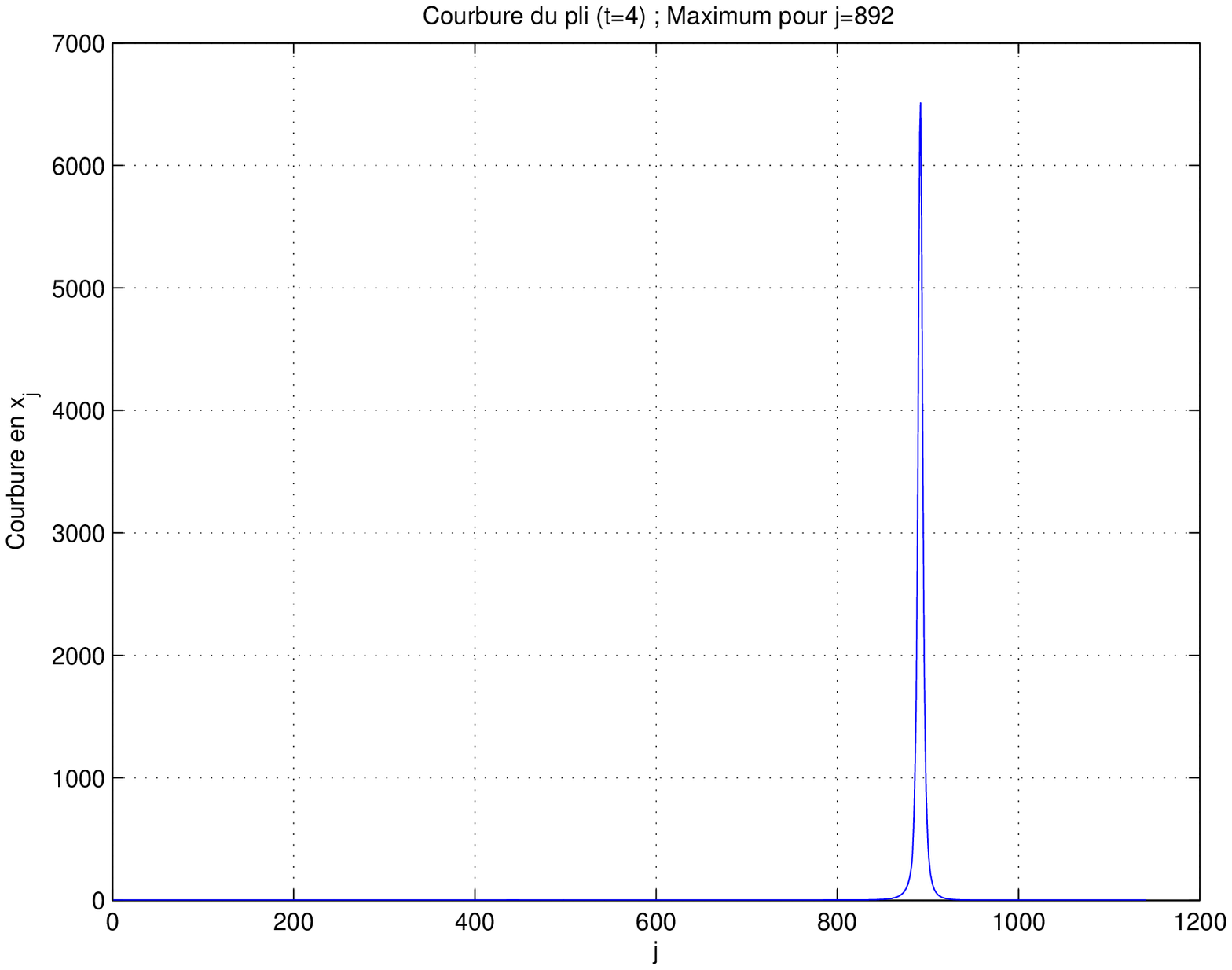}
     \\
     (e) $t=3\virg 3$ & (f) $t=4$
\end{tabular}
\caption{\label{fig:courbure_pli2} Formation du pli : courbure
(deuxi\`eme m\'ethode).}
\end{center}
\end{figure}

\begin{figure}
\begin{center}
\begin{tabular}{c@{}c}
     \includegraphics[height=7cm]{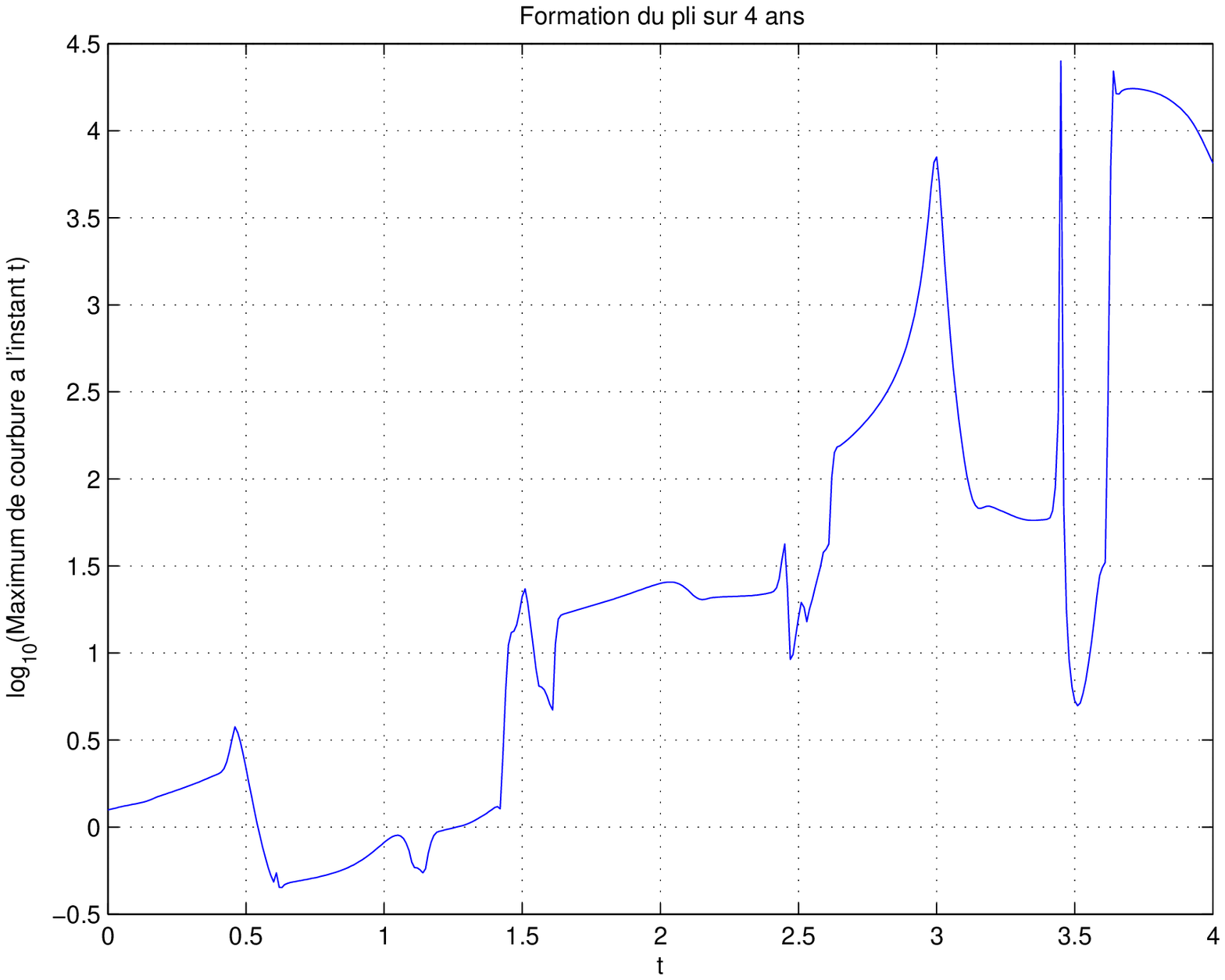}
     \\
     (a) Valeur maximale
     \\
     \includegraphics[height=7cm]{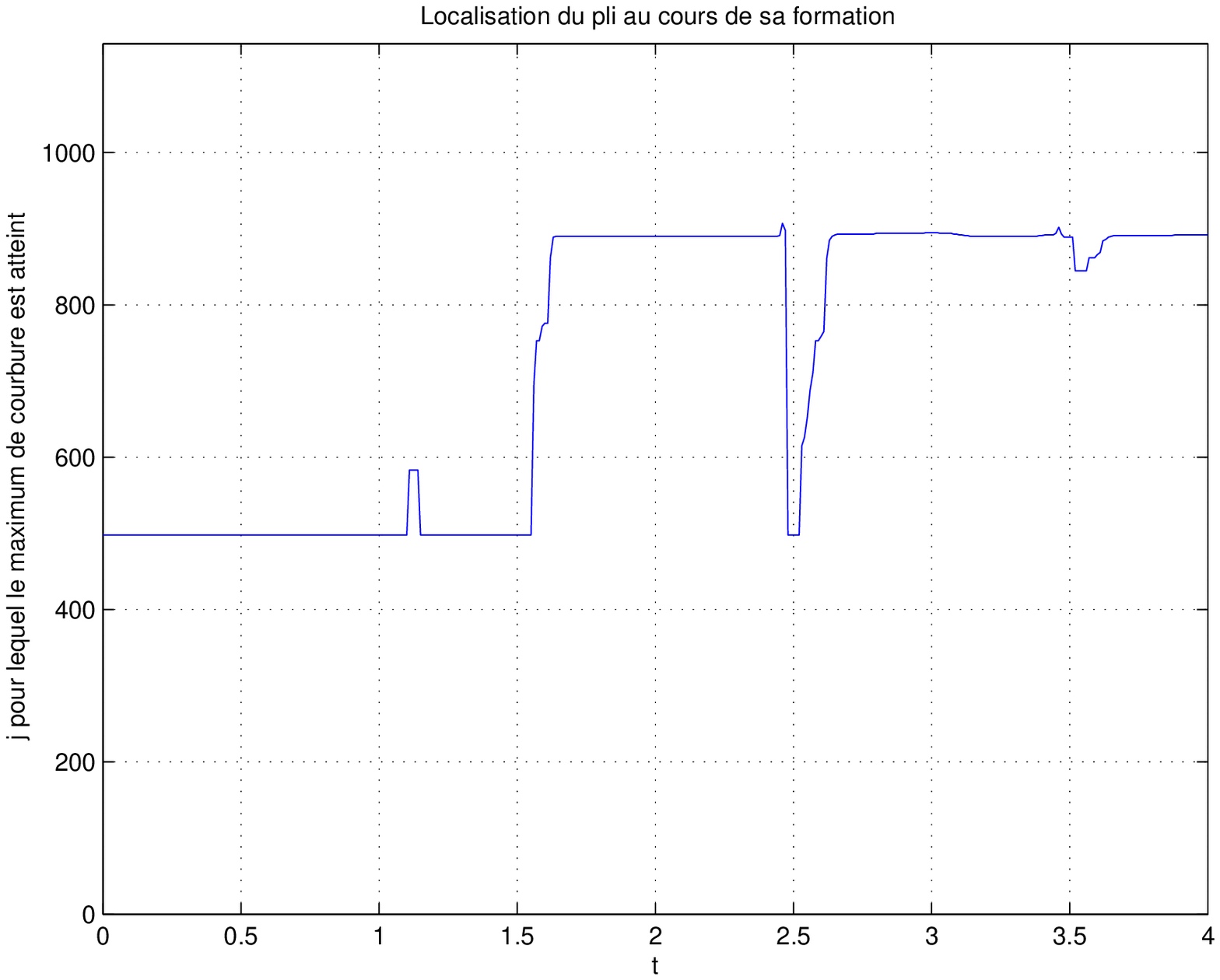}
     \\
     (b) Localisation du maximum
\end{tabular}
     \caption{\label{fig:courbure_maxi2}Formation du pli : maximum de courbure (deuxi\`eme m\'ethode).}
\end{center}
\end{figure}
\clearpage
\bibliographystyle{alpha}

\bibliography{maths,biologie}

\end{document}